\documentclass[AHL,XML,longabstracts,Unicode]{cedram}
\usepackage{tikz}
\usetikzlibrary{arrows,calc,positioning,decorations.pathreplacing,cd}
\usepackage{paralist}
\usepackage[mathscr]{euscript}
\usepackage{accents}
\usepackage{xparse}
\usepackage{array}
\usepackage{multirow}
\usepackage{enumitem}

\datereceived{2021-11-17}
\daterevised{2023-12-01}
\dateaccepted{2024-01-15}
\editor{G. Powell}
\DOI{10.NNNN/ahl.XXXX}

\AtBeginDocument{
\newtheorem{athm}{Theorem}

\theoremstyle{remark}
\newtheorem{rmk}[cdrthm]{Remark}
\newtheorem{eg}[cdrthm]{Example}
\newtheorem{convention}[cdrthm]{Convention}
\newtheorem{construction}[cdrthm]{Construction}
\newtheorem{condition}[cdrthm]{Condition}
\newtheorem{hypothesis}[cdrthm]{Hypothesis}
}

\newcommand{\Emb}{\mathrm{Emb}}
\newcommand{\Diff}{\mathrm{Diff}}
\newcommand{\Aut}{\mathrm{Aut}}
\newcommand{\Map}{\mathrm{Map}}
\newcommand{\mbar}{\ensuremath{{\,\,\overline{\!\! M\!}\,}}}
\newcommand{\nbar}{\ensuremath{{\,\,\overline{\!\! N}}}}
\newcommand{\B}{\ensuremath{\mathbf{B}}}
\newcommand{\N}{\mathscr{N}}
\newcommand{\longtwoheadrightarrow}{\longrightarrow\mathrel{\mkern-14mu}\rightarrow}
\newcommand{\longhookrightarrow}{\lhook\joinrel\longrightarrow}
\newcommand{\too}{\ensuremath{\longrightarrow}}
\newcommand{\cat}{\ensuremath{\mathrm{Cat}}}
\newcommand{\topo}{\ensuremath{\mathrm{Top}}}
\newcommand{\obj}{\ensuremath{\mathrm{ob}}}
\newcommand{\pr}{\ensuremath{\mathrm{pr}}}
\newcommand{\opp}{\ensuremath{\mathrm{op}}}
\newcommand{\fib}{\ensuremath{\mathsf{fib}}}
\newcommand{\id}{\mathrm{id}}
\newcommand{\ab}{\mathrm{ab}}
\newcommand{\BM}{\mathrm{BM}}
\newcommand{\MCGo}{\mathbf{\Gamma}}
\newcommand{\MCGno}{\boldsymbol{\mathcal{N}}}
\newcommand{\MCG}{\mathrm{MCG}}
\newcommand{\LCS}{\varGamma}
\newcommand{\Sym}{\mathfrak{S}}
\newcommand{\bn}{\mathbf{n}}
\newcommand{\LB}{\mathfrak{LB}}
\newcommand{\Int}{\mathring}
\newcommand{\lB}{\mathbf{LB}}
\newcommand{\tw}{\mathrm{tw}}
\newcommand{\colim}{\mathrm{colim}}
\newcommand{\Mot}{\mathrm{Mot}}
\newcommand{\mt}{\mathrm{t}}
\newcommand{\Aalg}{\ensuremath{\mathbb{A}\text{-}\mathrm{Alg}}}
\newcommand{\Rings}{\ensuremath{\mathrm{Rings}}}
\newcommand{\groups}{\ensuremath{\mathrm{Grp}}}
\newcommand{\fU}{\ensuremath{\mathfrak{U}}}
\newcommand{\diffdec}{\ensuremath{\mathrm{Diff}_{\mathrm{dec}}}}
\newcommand{\diffdecbr}{\ensuremath{\mathrm{Diff}_{\mathrm{dec}}^{\mathrm{br}}}}

\newcommand{\diffdecbrplus}{\ensuremath{\mathrm{Diff}_{\mathrm{dec}}^{\mathrm{br},+}}}
\newcommand{\diffdece}{\ensuremath{\mathrm{Diff}_{\mathrm{dec},\epsilon}}}
\newcommand{\diffdeceprime}{\ensuremath{\mathrm{Diff}_{\mathrm{dec},\epsilon'}}}
\newcommand{\diffdecet}{\ensuremath{\mathrm{Diff}_{\mathrm{dec},\epsilon,t}}}
\newcommand{\diffdecetprime}{\ensuremath{\mathrm{Diff}_{\mathrm{dec},\epsilon',t'}}}
\newcommand{\embdec}{\ensuremath{\mathrm{Emb}_{\mathrm{dec}}}}
\newcommand{\embdece}{\ensuremath{\mathrm{Emb}_{\mathrm{dec},\epsilon}}}
\newcommand{\cdec}{\ensuremath{C^\infty_{\mathrm{dec}}}}
\newcommand{\embgm}{\ensuremath{\mathrm{Emb}_{\langle \cG,\cM \rangle}}}
\newcommand{\Beta}{\boldsymbol{\beta}}
\newcommand{\triv}{\mathrm{tr}}
\newcommand{\col}{\mathrm{col}}
\newcommand{\unt}{\mathrm{u}}
\newcommand{\sR}{\mathscr{R}}
\newcommand{\sL}{\mathscr{L}}
\newcommand{\fF}{\mathfrak{F}}
\newcommand{\im}{\mathrm{Im}}
\newcommand{\zero}{\mathbf{0}}
\newcommand{\sF}{\mathscr{F}}
\newcommand{\tF}{\mathtt{F}}
\newcommand{\eF}{\mathsf{F}}
\newcommand{\fG}{\mathbf{G}}
\newcommand{\sG}{\mathsf{G}}

\NewDocumentCommand{\covlr}{O{\bullet} O{\bullet}}{%
\ensuremath{{}_{#1}\mathrm{Cov}_{#2}}%
}
\NewDocumentCommand{\covl}{O{\bullet}}{%
\ensuremath{{}_{#1}\mathrm{Cov}}%
}
\NewDocumentCommand{\covr}{O{\bullet}}{%
\ensuremath{\mathrm{Cov}_{#1}}%
}
\NewDocumentCommand{\toplr}{O{\bullet} O{\bullet}}{%
\ensuremath{{}_{#1}\mathrm{Top}_{#2}}%
}
\NewDocumentCommand{\topl}{O{\bullet}}{%
\ensuremath{{}_{#1}\mathrm{Top}}%
}
\NewDocumentCommand{\topr}{O{\bullet}}{%
\ensuremath{\mathrm{Top}_{#1}}%
}
\NewDocumentCommand{\modlr}{O{\bullet} O{\bullet}}{%
\ensuremath{{}_{#1}\mathrm{Mod}_{#2}}%
}
\NewDocumentCommand{\modl}{O{\bullet}}{%
\ensuremath{{}_{#1}\mathrm{Mod}}%
}
\NewDocumentCommand{\modr}{O{\bullet}}{%
\ensuremath{\mathrm{Mod}_{#1}}%
}
\NewDocumentCommand{\topmodlr}{O{\bullet} O{\bullet} O{\bullet}}{%
\ensuremath{{}_{#1}\mathrm{Top}_{#2}\mathrm{Mod}_{#3}}%
}
\NewDocumentCommand{\topmodlrpr}{O{\bullet} O{\bullet} O{\bullet}}{%
\ensuremath{{}_{#1}\mathrm{Top}^{\mathrm{pr}}_{\,\, #2}\mathrm{Mod}_{#3}}%
}
\NewDocumentCommand{\topmodl}{O{\bullet} O{\bullet}}{%
\ensuremath{{}_{#1}\mathrm{Top}_{#2}\mathrm{Mod}}%
}
\NewDocumentCommand{\topmodr}{O{\bullet} O{\bullet}}{%
\ensuremath{\mathrm{Top}_{#1}\mathrm{Mod}_{#2}}%
}
\NewDocumentCommand{\setlr}{O{\bullet} O{\bullet}}{%
\ensuremath{{}_{#1}\mathrm{Set}_{#2}}%
}
\NewDocumentCommand{\setl}{O{\bullet}}{%
\ensuremath{{}_{#1}\mathrm{Set}}%
}

\newcommand{\cB}{\mathcal{B}}
\newcommand{\cC}{\mathcal{C}}
\newcommand{\cD}{\mathcal{D}}
\newcommand{\cE}{\mathcal{E}}
\newcommand{\cF}{\mathcal{F}}
\newcommand{\cG}{\mathcal{G}}
\newcommand{\cH}{\mathcal{H}}

\newcommand{\cL}{\mathcal{L}}
\newcommand{\cM}{\mathcal{M}}

\newcommand{\cO}{\mathcal{O}}
\newcommand{\cP}{\mathcal{P}}
\newcommand{\cQ}{\mathcal{Q}}

\newcommand{\cT}{\mathcal{T}}
\newcommand{\cU}{\mathcal{U}}
\newcommand{\cV}{\mathcal{V}}

\newcommand{\bA}{\mathbb{A}}
\newcommand{\bB}{\mathbb{B}}

\newcommand{\bD}{\mathbb{D}}

\newcommand{\bN}{\mathbb{N}}

\newcommand{\bR}{\mathbb{R}}
\newcommand{\bS}{\mathbb{S}}

\newcommand{\bZ}{\mathbb{Z}}

\renewcommand{\geq}{\geqslant}
\renewcommand{\leq}{\leqslant}

\newcommand{\incl}[3][right]%
{%
\draw[<-,>=#1 hook] #2 to ($ #2!0.5!#3 $);
\draw[->] ($ #2!0.5!#3 $) to #3;%
}

\newenvironment{itemizeb}%
{\begin{compactitem}

}%
{\end{compactitem}}

\author[\initial{M.} \lastname{Palmer}]{\firstname{Martin} \lastname{Palmer}}
\address{Institutul de Matematică Simion Stoilow al Academiei Române, 21 Calea Griviței, 010702 Bucharest, Romania}
\email{mpanghel@imar.ro}

\author[\initial{A.} \lastname{Soulié}]{\firstname{Arthur} \lastname{Soulié}}
\address{Normandie Univ., UNICAEN, CNRS, LMNO, 14000 Caen, France}
\email{artsou@hotmail.fr, arthur.soulie@unicaen.fr}

\title[Unified topological representations]{Topological representations of motion groups and mapping class groups -- a unified functorial construction}
\alttitle{Représentations topologiques des groupes de mouvement et des groupes de difféotopie -- Une construction fonctorielle unifiée}

\begin{abstract}
For groups of a topological origin, such as braid groups and mapping class groups, an important source of interesting and highly non-trivial representations is given by their actions on the twisted homology of associated spaces; these are known as \emph{homological representations}. Representations of this kind have proved themselves especially important for the question of \emph{linearity}, a key example being the family of topologically-defined representations introduced by Lawrence and Bigelow, and used by Bigelow and Krammer to prove that braid groups are linear.
In this paper, we give a unified foundation for the construction of homological representations using a functorial approach. Namely, we introduce \emph{homological representation functors} encoding a large class of homological representations, defined on categories containing all \emph{mapping class groups} and \emph{motion groups} in a fixed dimension. These source categories are defined using a topological enrichment of the Quillen bracket construction applied to categories of decorated manifolds. This approach unifies many previously-known constructions, including those of Lawrence-Bigelow, and yields many new representations.
\end{abstract}

\begin{altabstract}
Une source majeure de représentations intéressantes et grandement non-triviales pour les groupes ayant une origine topologique, tels que les groupes de tresses et les groupes de difféotopies, est donnée par leurs actions sur de l'homologie tordue d'espaces associés ; celles-ci sont connues sous le nom de \emph{représentations homologiques}. Les représentations de ce type se sont montrées particulièrement importantes pour les questions de \emph{linéarité}, un exemple clef étant celui de la famille de représentations définies topologiquement par Lawrence et Bigelow, et utilisée par Bigelow et Krammer pour démontrer que les groupes de tresses sont linéaires.
Dans cet article, nous établissons des fondations unifiées pour la construction de représentations homologiques en utilisant des méthodes fonctorielles. Plus précisément, nous introduisons des \emph{foncteurs de représentations homologiques}, qui encodent de larges classes de représentations homologiques, et qui sont définis sur des catégories contenant tous les \emph{groupes de difféotopie} et tous les \emph{groupes de mouvement} pour une dimension fixée. Ces catégories sources sont définies à partir d'un enrichissement topologique de la construction de support due à Quillen, que l'on applique à des catégories de variétés décorées. Cette approche unifie de nombreuses constructions déjà connues, y compris celles de Lawrence et Bigelow, et produit beaucoup de nouvelles représentations.
\end{altabstract}

\thanks{The first author was partially supported by a grant of the Romanian Ministry of Education and Research, CNCS - UEFISCDI, project number PN-III-P4-ID-PCE-2020-2798, within PNCDI III. The second author was partially supported by a Rankin-Sneddon Research Fellowship of the University of Glasgow, by the Institute for Basic Science IBS-R003-D1 and by the ANR Projects ChroK ANR-16-CE40-0003 and AlMaRe ANR-19-CE40-0001-01.}

\begin{document}
\shorthandoff{!}
\maketitle

\section*{Introduction}

For a smooth manifold $M$ with non-empty boundary and a closed submanifold $Y$ of its interior, the mapping class group $\MCG(M,Y)$ is the group of isotopy classes of diffeomorphisms of $M$ fixing its boundary pointwise and fixing $Y$ setwise.
The corresponding \emph{motion group} $\Mot_{Y}(M)$ is the fundamental group of the space of embeddings of $Y$ into the interior of $M$ up to diffeomorphisms of $Y$.
(When $Y$ is orientable, there is also an oriented variant $\Mot_{Y}^{+}(M)$ of the motion group, for which we require diffeomorphisms of $Y$ to be orientation-preserving.)
For example, if $\bD^{2}$ denotes the closed unit $2$-disc and $\underline{n}$ a set of $n \geq 1$ distinct points in its interior, the classical braid group on $n$ strands $\B_{n}$ is isomorphic both to the mapping class group $\MCG(\bD^{2},\underline{n})$ and to the motion group $\Mot_{\underline{n}}(\bD^{2})$.

The representation theory of \emph{mapping class groups} and of \emph{motion groups} is very rich, and the subject of much active research -- see Birman and Brendle's survey \cite[\S 4]{BirmanBrendlesurvey} for braid groups or Margalit's expository paper \cite{Margalit} for more general mapping class groups of surfaces. A reason for this is that representation theory transforms abstract group theory problems into linear algebra questions, which are generally more accessible.
In particular, a group is said to be \emph{linear} if it acts faithfully on a finite-dimensional vector space. Whether or not a group satisfies this property is a fundamental question: if it does, the group is isomorphic to a subgroup of a general linear group over a field and falls into the class of \emph{matrix groups}, and thus automatically satisfies many remarkable properties for free; see for instance \cite{Suprunenko}.
Classical braid groups provide a famous positive answer to this question by Bigelow and Krammer~\cite{bigelow2001braid,KrammerLK}, as well as the mapping class group of a torus which is classically known to be isomorphic to the symplectic group $\mathrm{Sp}_{2}(\bZ)$ (see for example \cite[Th.~2.5]{farbmargalit}), or the mapping class groups of a genus two surface thanks to the work of Bigelow and Budney \cite{BigelowBudney}. 
In contrast, there are also families with a negative answer, such as automorphism groups of free groups --- which may be viewed as mapping class groups of certain $3$-manifolds --- by Formanek and Procesi~\cite{FormanekProcesi}. It is known by \cite{Button2016} (see also \cite[Cor.~1.4(b)]{KoberdaLuoSun2019}) that any faithful, finite-dimensional representation of mapping class groups of genus at least $3$ must be in characteristic zero. However, in general, the question of linearity remains a wide open fundamental problem for the vast majority of mapping class groups and motion groups; see \cite{BrendleHamidi-Tehrani} or \cite[\S 1]{Margalit}.

\paragraph*{Homological representations.}

In order to understand the representation theory of mapping class groups and motion groups, we view them not just as abstract groups, but use their geometric and topological structures. Combining this philosophy (of studying these groups via their \emph{topology}) with the philosophy of representation theory itself (of studying groups via \emph{linear algebra}), one is led naturally to homology. More precisely, the idea consists in constructing representations of these groups via their actions on the homology of some topological spaces naturally associated to them. Such representations are said to be \emph{homological}.

As an archetypal example, for each braid group $\B_{n}$, Lawrence \cite{Lawrence1} and Bigelow \cite{bigelow2001braid} constructed a by now well-known representation $\LB_{k}(n)$ for each $k\geq 1$, called the $k$-th \emph{Lawrence-Bigelow representation}, defined as follows. Let $\bD_{n}$ be the $n$-punctured disc and $C_{k}(\bD_{n})$ the configuration space of $k$ unordered points in $\bD_{n}$. Lawrence \cite[\S 2]{Lawrence1} and Bigelow \cite[\S 1.2]{bigelow2001braid} define via a geometrical method regular covering spaces $C_{k}(\bD_{n})^{\phi_{k}}$ of $C_{k}(\bD_{n})$, with deck transformation group $Q_{k}=\bZ$ if $k=1$ and $\bZ^{2}$ if $k\geq2$; see \S\ref{ss:applications_motion_groups} for more details. The representation $\LB_{k}(n)$ is then defined by the natural action of $\B_{n}$ on the homology group $H_{k}(C_{k}(\bD_{n})^{\phi_{k}};\bZ)$.
In particular, $\LB_{1}(n)$ is the reduced \emph{Burau} representation originally introduced by Burau \cite{burau}, while $\LB_{2}(n)$ corresponds to the \emph{Lawrence-Krammer-Bigelow} representation, which Bigelow \cite{bigelow2001braid} and Krammer \textup{\cite{KrammerLK}} independently proved to be faithful. Moreover, the representations used to prove the linearity of the mapping class groups of a torus and of a genus two surface are also homological; see \cite[Th.~2.5]{farbmargalit} and \cite[\S 3]{BigelowBudney}.
In particular, the representations that have thus far provided positive answers for the linearity question have all been homological representations. This suggests that a systematic treatment of constructing homological representations that works for all motion groups and all mapping class groups is an important and natural avenue of investigation related to the linearity question for these families of groups.

\paragraph*{A functorial approach.}

In parallel, another point of view adopted in this paper is to treat mapping class groups and motion groups simultaneously \emph{in families}. A \emph{family of groups} means a collection of groups $G_{n}$ indexed by the natural numbers $\bN$ equipped with morphisms $G_{n} \to G_{n+1}$. Typically, we consider a set of geometrically-consistent manifolds $\{M_{n}\}_{n\in\mathbb{N}}$ in order to construct a family of mapping class groups $\MCG(M_{n})\to \MCG(M_{n+1})$, or a set of geometrically-consistent submanifolds $\{Y_{n}\}_{n\in\mathbb{N}}$ of a fixed manifold $M$ to construct a family of motion groups $\Mot_{Y_{n}}(M)\to \Mot_{Y_{n+1}}(M)$. For example, we consider braid groups as a family thanks to the inclusions $\B_{n}\hookrightarrow \B_{n+1}$ induced by adding a strand on the left.
The purpose of considering groups in families is to require representations to respect the coherences that naturally arise between the different groups in the family. More precisely, for a family of groups $\{g_{n} \colon G_{n} \to G_{n+1}\}_{n\in\mathbb{N}}$ and an associative, unital ring $R$, a collection of representations $\{\rho_{n} \colon G_{n}\rightarrow \Aut_{R}(V_{n})\}_{n\in\mathbb{N}}$ is \emph{coherent} if it comes equipped with module homomorphisms $v_{n} \colon V_{n} \to V_{n+1}$ such that $v_{n}$ is \emph{equivariant} with respect to $g_{n}$, meaning that for all $g \in G_{n}$ and $x \in V_{n}$ we have $v_{n}(g.x) = g_{n}(g).v_{n}(x)$.

Furthermore, we consider a slightly larger notion of representation than the classical one of a group homomorphism $\rho_{n} \colon G_{n}\rightarrow \Aut_{R}(V_{n})$: we generally allow $\rho_{n}$ to take values in the larger group $\Aut_{\bZ}(V_{n}) \supseteq \Aut_{R}(V_{n})$ and assume that it comes equipped with an action $r_{n}\colon G_{n}\rightarrow \Aut_{\Rings}(R)$ encoding its failure to respect the $R$-module structure of $V_n$; see Definition~\ref{def:bimodule_category}. The representation $\rho_{n}$ is then said to be \emph{twisted}; if $r_{n}$ is trivial, then we recover the classical notion, and $\rho_{n}$ is said to be \emph{genuine} or \emph{untwisted}.

Now, the notion of coherent (twisted or genuine) representations of a family of groups can be encoded in a functorial way as follows. We denote by $\modr[R]$ the category of right $R$-modules. We introduce the larger category $\modr[R]\subseteq \modr[R]^{\tw}$ in which morphisms are permitted to act on the underlying ring $R$; we call this the category of \emph{twisted} right $R$-modules; see Definition~\ref{def:bimodule_category}.
Let $\fG$ be the groupoid with objects indexed by non-negative integers, with the groups $\{G_{n}\}_{n\in\mathbb{N}}$ as automorphism groups and with no morphisms between distinct objects. Let us tentatively assume that there exists a category $\langle \cG_{\circ}, \cM_{\circ} \rangle$ containing $\fG$ as its underlying groupoid and with a preferred morphism $\iota_{n} \colon n\rightarrow n+1$ for each object $n$, satisfying $\iota_{n} \circ g = g_{n}(g) \circ \iota_{n}$ for each $g \in G_{n}$.
Then, functors
\begin{equation}
\label{eq:coherent_rep_functors}
\langle \cG_{\circ}, \cM_{\circ} \rangle \too \modr[R]^{\tw} \quad \textrm{and} \quad \langle \cG_{\circ}, \cM_{\circ} \rangle \too \modr[R]
\end{equation}
give us coherent representations of $\{G_{n}\}_{n\in\mathbb{N}}$, which are respectively \emph{twisted} and \emph{untwisted} over $R$.
It is an important question whether or not our representations take values in the untwisted subcategory $\modr[R]\subseteq \modr[R]^{\tw}$. Indeed, if this is not the case, the encoded representations are not \emph{genuine} representations of the family of groups, since the groups also act on the ground ring $R$. Although one may always find a subring $R'\subset R$ over which the representations become genuine (i.e.~unwisted), the representations typically also become infinite-dimensional as $R'$-modules. See \cite{BPS,PSIIp} for examples of such situations, where the representations are finitely generated as $R$-modules with $R=\bZ[Q]$ for $Q$ an infinite group, while $R'=\bZ$.
We introduce a general criterion for the representations that we construct to be untwisted in \S\ref{sss:col_coeff_unwtisted_rep}.

The category $\langle \cG_{\circ}, \cM_{\circ} \rangle$ is defined by the \emph{Quillen bracket construction} $\langle -,- \rangle$, introduced in \cite[p.219]{graysonQuillen}, applied to a monoidal groupoid $\cG_{\circ}$ and a $\cG_{\circ}$-module groupoid $\cM_{\circ}$. That there exist such groupoids $\cG_{\circ}$ and $\cM_{\circ}$ such that the category $\langle \cG_{\circ}, \cM_{\circ} \rangle$ contains the groupoid $\fG$ as its underlying groupoid is a mild assumption, naturally satisfied in all of our examples. For instance, we consider the braid groupoid $\Beta$ (see \cite[Chap.~XI, \S 4]{MacLane1}) and take $\cG_{\circ} = \cM_{\circ} = \fG = \Beta$ to deal with the braid groups; see \S\ref{sss:category_surface_braid_groups}.
We provisionally defer the treatment of this question until after Theorem~\ref{athm:thm_B}. For now, we assume that this category exists and mention that, for the purpose of constructing homological representations, we first introduce a \emph{topologically-enriched} semicategory $\langle \cG_{\circ}, \cM_{\circ} \rangle^{\mt}$, and define $\langle \cG_{\circ}, \cM_{\circ} \rangle := \pi_{0}(\langle \cG_{\circ}, \cM_{\circ} \rangle^{\mt})$ by taking path-components of all morphism spaces (this is \emph{a priori} just a semicategory, but turns out to be a category). \emph{From now on, we consider a family $\{g_{n} \colon G_{n} \to G_{n+1}\}_{n\in\mathbb{N}}$ of mapping class groups or motion groups in a fixed dimension $d$, contained in a category $\langle \cG_{\circ}, \cM_{\circ} \rangle=\pi_{0}(\langle \cG_{\circ}, \cM_{\circ} \rangle^{\mt})$ as above.}

In addition, thanks to our functorial point of view, notions of \emph{polynomiality} may be introduced for functors of the form \eqref{eq:coherent_rep_functors}.
We refer the reader to \cite[\S 4.1]{RWW} or \cite[\S 4.1]{PSIIp} for a detailed introduction to this notion.
In particular, a significant application of polynomiality is homological stability with twisted coefficients for several families of mapping class groups and motion groups: see the work of Randal-Williams and Wahl \cite[Theorems D, 5.26 and I]{RWW}; this was also proven earlier by Ivanov \cite{Ivanov} and Boldsen \cite{Boldsen} for mapping class groups of orientable surfaces.
Also, certain notions of polynomiality offer a classifying tool for functors (see \cite[\S 1.3]{DV3} for instance), while the representation theory of mapping class groups and motion groups is known to be \emph{wild} (for example, this is proven by Erdmann and Nakano~\cite{ErdmannNakano} for braid groups on $n\geq 6$ strands). This is therefore another argument supporting our functorial approach to representations.
The study of polynomiality for homological representation functors of certain mapping class groups and motion groups is the central topic of the sequel paper \cite{PSIIp}.

In summary, our goal is to develop functorial constructions of homological representations that:
\begin{itemizeb}
\item apply naturally to all mapping class groups and motion groups in a fixed dimension (\emph{globality}),
\item respect the natural coherences between these groups (\emph{functoriality}),
\item produce a wide range of new, interesting and highly non-trivial representations (\emph{richness}), for the purpose of studying questions of linearity and of polynomiality.
\end{itemizeb}

\paragraph*{Functorial constructions of homological representations.}

Our first main result is Theorem~\ref{athm:thm_A} below, which introduces functorial constructions of coherent representations for the family of groups $\{g_{n} \colon G_{n} \to G_{n+1}\}_{n\in\mathbb{N}}$. They depend on the following parameters:
\begin{itemizeb}
\item the dimension $d$ of the family of mapping class groups or motion groups, which we assume from now on to be different from $4$. The reason for this restriction is related to establishing the existence of a certain decomposition of embedding spaces; see Lemma~\ref{lem:decomposition-disjoint-union} and Remark~\ref{rmk:d_neq_4_explanation};
\item a closed submanifold $Z \subset \bR^d$;
\item an open subgroup $\sG$ of the group $\Diff(Z)$ of diffeomorphisms of $Z$;
\item an integer $\ell \geq 1$, indexing the $\ell$-th lower central series functor $\LCS_{\ell}\colon \groups \to \groups$ defined by sending a group $G$ to $\LCS_{\ell}(G)$, see Example~\ref{eg:examples_FQG};
\item a non-negative integer $i \geq 0$, indexing the degree of the twisted homology functor $H_{i}$, see Proposition~\ref{prop:twisted-homology}.
\end{itemizeb}

\begin{athm}[{Theorems~\ref{thm:construction}, \ref{thm:global_functor_motion_groups} and \ref{thm:global_functor_mcg}}]
\label{athm:thm_A}
For a fixed set of parameters $\{i,Z,\sG,\ell\}$ as described above, the action of the family of groups $\{g_{n} \colon G_{n} \to G_{n+1}\}_{n\in\mathbb{N}}$ on the twisted homology of certain embedding spaces determines functors
\begin{equation}
\label{eq:ThmA_hom_rep_functors_twisted}
L_{i}(\cF_{(Z,\sG,\ell)})\colon \langle \cG_{\circ}, \cM_{\circ} \rangle \too \modr[\bZ[\cQ]]^{\tw} \quad \textrm{and} \quad L_{i}(\sF_{(Z,\sG,\ell)})\colon \langle \cG_{\circ}, \cM_{\circ} \rangle \too \modr[\bZ[\cQ]]^{\tw},
\end{equation}
where $\cQ$ denotes a group built out of the deck transformation groups of the regular covering spaces corresponding to the coefficients in the twisted homology.
In particular, the functors \eqref{eq:ThmA_hom_rep_functors_twisted} encode coherent homological representations of the family of groups $\{g_{n} \colon G_{n} \to G_{n+1}\}_{n\in\mathbb{N}}$. They are thus called \textbf{homological representation functors}.

Additionally, there exists a universal quotient $\cQ^{\unt}$ of the group $\cQ$, together with variants of the homological representation functors \eqref{eq:ThmA_hom_rep_functors_twisted}, such that the encoded representations of the family of groups $\{g_{n} \colon G_{n} \to G_{n+1}\}_{n\in\mathbb{N}}$ are \emph{untwisted} over the ground ring $\bZ[\cQ^{\unt}]$:
\begin{equation}
\label{eq:ThmA_hom_rep_functors_untwisted}
L_{i}(\cF^{\unt}_{(Z,\sG,\ell)})\colon \langle \cG_{\circ}, \cM_{\circ} \rangle \too \modr[\bZ[\cQ^{\unt}]] \quad \textrm{and} \quad L_{i}(\sF^{\unt}_{(Z,\sG,\ell)})\colon \langle \cG_{\circ}, \cM_{\circ} \rangle \too \modr[\bZ[\cQ^{\unt}]].
\end{equation}
These alternative homological representation functors are therefore called \textbf{untwisted}.
\end{athm}

We note that one could simply replace $\cQ$ with the trivial group, in which case the functors \eqref{eq:ThmA_hom_rep_functors_twisted} would already be untwisted (since $\bZ$ has no non-trivial ring automorphisms). However, this would multiply the dimensions of all of the underlying modules of \eqref{eq:ThmA_hom_rep_functors_twisted} by $\lvert \cQ \rvert$, which is typically infinite, and for applications -- especially to questions of linearity -- one would like to have \emph{finite-dimensional} representations.

The main geometric input in the construction of the homological representation functors \eqref{eq:ThmA_hom_rep_functors_twisted} and \eqref{eq:ThmA_hom_rep_functors_untwisted} is contained in the construction of two continuous semifunctors
\begin{equation}
\label{eq:ThmA_semifunctors}
\cF_{(Z,\sG,\ell)}\colon \langle \cG_{\circ}, \cM_{\circ} \rangle^{\mt} \too \covr \quad \textrm{and} \quad \sF_{(Z,\sG,\ell)}\colon \langle \cG_{\circ}, \cM_{\circ} \rangle^{\mt} \too \covr
\end{equation}
as well as their \emph{untwisted} variants $\cF^{\unt}_{(Z,\sG,\ell)}$ and $\sF^{\unt}_{(Z,\sG,\ell)}$. We emphasise that these are \emph{semifunctors}, not genuine \emph{functors}: this is because their source is the topological lift $\langle \cG_{\circ}, \cM_{\circ} \rangle^{\mt}$ of the category $\langle \cG_{\circ}, \cM_{\circ} \rangle$, which is a topologically-enriched \emph{semicategory}, not a genuine \emph{category}; see Theorem~\ref{athm:thm_B}. We also use a topologically-enriched category $\covr$ of topological spaces equipped with regular coverings (introduced in Definition~\ref{def:bicoverings}). The idea to construct the semifunctors \eqref{eq:ThmA_semifunctors} consists in defining regular coverings of certain embedding spaces depending on the manifold $Z$ and the group $\sG$, induced by the lower central series $\LCS_{\ell}$ via the commutative diagrams \eqref{eq:split-short-exact-sequence-quotient} and \eqref{eq:split-short-exact-sequence-quotient-mcg} respectively. As one may see from these diagrams, the fundamental difference between the two constructions lies in the fact that the regular coverings for $\cF_{(Z,\sG,\ell)}$ are induced by the actions of some \emph{motion groups}, whereas they are induced by the actions of some \emph{mapping class groups} for $\sF_{(Z,\sG,\ell)}$. See \S\ref{ss:homological_representation_functor_motion_groups} for the detailed construction of $\cF_{(Z,\sG,\ell)}$ and $\cF^{\unt}_{(Z,\sG,\ell)}$, and \S\ref{ss:homological_representation_functor_mcg} for that of $\sF_{(Z,\sG,\ell)}$ and $\sF^{\unt}_{(Z,\sG,\ell)}$.

Having constructed the continuous semifunctors \eqref{eq:ThmA_semifunctors}, the overall procedure of Theorem~\ref{athm:thm_A} to define the homological representation functors \eqref{eq:ThmA_hom_rep_functors_twisted} and \eqref{eq:ThmA_hom_rep_functors_untwisted} is summarised in the commutative diagram \eqref{eq:diagram-homological-representations} below. We consider any one of the semifunctors \eqref{eq:ThmA_semifunctors} or their untwisted variants, which we denote by $\eF_{(Z,\sG,\ell)}$.
\begin{equation}
\label{eq:diagram-homological-representations}
\centering
\begin{split}
\begin{tikzpicture}
[x=1mm,y=1mm]
\node (c) at (0,0) {$\langle \cG_{\circ}, \cM_{\circ} \rangle^{\mt}$};
\node (c0) at (0,-12) {$\langle \cG_{\circ}, \cM_{\circ} \rangle$};
\node (top) at (45,0) {$\topr[{\bZ[\cQ]}]^{\tw}$};
\node (mod) at (90,0) {$\modr[{\bZ[\cQ]}]^{\tw}$};
\draw[->] (c) to node[above,font=\small]{$\tilde{\eF}_{(Z,\sG,\ell)}$} (top);
\draw[->] (top) to node[above,font=\small]{$H_{i}$} (mod);
\draw[->>] (c) to node[left,font=\small]{$\pi_{0}$} (c0);
\draw[->,dashed] (c0) to[out=0,in=195] node[below,font=\small]{$L_{i}(\eF_{(Z,\sG,\ell)})$} (mod);
\end{tikzpicture}
\end{split}
\end{equation}
We first promote the semifunctor $\eF_{(Z,\sG,\ell)}$ to a semifunctor $\tilde{\eF}_{(Z,\sG,\ell)}\colon \langle \cG_{\circ}, \cM_{\circ} \rangle^{\mt} \to \topr[\bZ[\cQ]]^{\tw}$. Here, $\topr[\bZ[\cQ]]^{\tw}$ denotes the topologically-enriched category of topological spaces equipped with \emph{twisted} bundles of $\bZ[\cQ]$-bimodules (i.e.~the morphisms are permitted to act on the ground ring $\bZ[\cQ]$); see Definition~\ref{def:bundles_of_bimodules}. More precisely, the semifunctor $\tilde{\eF}_{(Z,\sG,\ell)}$ is defined from $\eF_{(Z,\sG,\ell)}$ via a procedure that involves a \emph{linearisation} of $\eF_{(Z,\sG,\ell)}$ (see \S\ref{sss:lift}) followed by its \emph{fibrewise tensor product} (see \S\ref{sss:fibrewise_tensor}) with an associated \emph{colimit coefficient system} (see \S\ref{sss:col_coeff_unwtisted_rep}). In particular, the group $\cQ$ corresponds to a colimit of the deck transformation groups of the covering spaces defining the semifunctor $\eF_{(Z,\sG,\ell)}$; see Notation~\ref{notation-Qcol}.
The second step simply consists in applying the twisted homology functor $H_{i}$ (see \S\ref{sss:twisted-homology}). Finally, the composite $H_{i}\circ \tilde{\eF}_{(Z,\sG,\ell)}$ factors uniquely through the projection functor $\pi_{0}\colon \langle \cG_{\circ}, \cM_{\circ} \rangle^{\mt}\twoheadrightarrow \langle \cG_{\circ}, \cM_{\circ} \rangle$ given by taking path-components of morphism spaces; see Lemma~\ref{lem:factorisation_pi_{0}}.
The desired output is the functor $L_{i}(\eF_{(Z,\sG,\ell)}) \colon \langle \cG_{\circ}, \cM_{\circ} \rangle \to \modr[{\bZ[\cQ]}]^{\tw}$; this is a coherent family of (possibly twisted) representations of $\{G_{n}\}_{n\in\bN}$.
A more elaborate version of diagram \eqref{eq:diagram-homological-representations}, allowing one to take fibrewise tensor products with different coefficient systems $V$, is described in diagram \eqref{eq:construction} and Definition~\ref{def:construction}, but the underlying ideas are the same. We refer the reader to \S\ref{sss:construction} for the detailed procedure defining diagram \eqref{eq:diagram-homological-representations} and its elaboration \eqref{eq:construction}.

\paragraph*{Categorical framework.}

Let us now tackle the question of defining the topologically-enriched semicategory $\langle \cG_{\circ}, \cM_{\circ} \rangle^{\mt}$, which is a topological lift of the Quillen bracket category $\langle \cG_{\circ}, \cM_{\circ} \rangle$. This is essentially solved by introducing a certain topologically-enriched semicategory $\langle \cD_{d},\cD_{d} \rangle$ as follows.

For each fixed dimension $d\geq 2$, we construct a semi-monoidal topologically-enriched groupoid $(\cD_{d},\natural)$, called the \emph{decorated $d$-manifolds groupoid}; see Definition~\ref{def:Dd}. The objects of $\cD_{d}$ are pairs $(M,A)$ with $M$ a smooth $d$-manifold and $A$ a closed submanifold embedded in the interior of $M$ (together with some auxiliary data needed to construct boundary connected sums); see Definition~\ref{def:decorated-manifolds}. Its morphism spaces consist of \emph{decorated diffeomorphisms}, namely diffeomorphisms of pairs of manifolds compatible with the auxiliary data; see Definitions~\ref{def:decorated-manifolds-morphisms} and \ref{def:decorated-manifolds-morphism-spaces}. The semi-monoidal structure $\natural$ is given by the boundary connected sum of manifolds; see Proposition~\ref{prop:Dd-semi-monoidal}. See also Remark~\ref{r:semi-monoidal} for an explanation of why $\natural$ is only a \emph{semi}-monoidal structure; this is also why $\langle \cD_{d},\cD_{d} \rangle$ is only a \emph{semi}category.
Furthermore, we introduce in \S\ref{sss:Topological-Quillen} a topological enrichment of the Quillen bracket construction $\langle -,- \rangle$ \cite[p.219]{graysonQuillen}; see Proposition~\ref{prop:topological-Quillen}. The topologically-enriched semicategory $\langle \cD_{d},\cD_{d} \rangle$ is then constructed by applying the topologically-enriched Quillen bracket construction to the semi-monoidal topologically-enriched groupoid $(\cD_{d},\natural)$.

Our second main result establishes crucial properties of the semicategory $\langle \cD_{d},\cD_{d} \rangle$, by giving a geometric description of its morphism spaces and by describing its interaction with the path-component functor $\pi_{0}$:

\begin{athm}[{Corollary~\ref{cor:description_UD_{d}}}]
\label{athm:thm_B}
For $d \neq 4$, the morphism spaces $\langle \cD_{d} , \cD_{d} \rangle (M,N)$ of the topologically-enriched semicategory $\langle \cD_{d} , \cD_{d} \rangle$ may be identified with the spaces $\embdec(M,N)$ of decorated embeddings introduced in Definition~\ref{def:embdec}.

Applying the path-component functor $\pi_{0}$ to all morphism spaces, there is an isomorphism of categories:
\begin{equation}
\label{eq:ThmB_iso_pi0}
\pi_{0}(\langle \cD_{d},\cD_{d} \rangle) \cong \langle\pi_{0}(\cD_{d}),\pi_{0}(\cD_{d})\rangle.
\end{equation}
Moreover, $\pi_{0}(\cD_{d})$ is the underlying groupoid of $\pi_{0}(\langle \cD_{d},\cD_{d} \rangle)$. 
\end{athm}

The above description of the morphism spaces of $\langle \cD_{d},\cD_{d} \rangle$ as decorated embedding spaces (see Proposition~\ref{prop:morphism-spaces-bracket}), as well as the isomorphism \eqref{eq:ThmB_iso_pi0} (see Proposition~\ref{prop:topological-Quillen_morphisms}), crucially require a technical Serre fibration result for decorated diffeomorphism groups given in Theorem~\ref{thm:fibre-bundle}.

The topologically-enriched groupoid $\cD_{d}$ contains, by construction, all diffeomorphism groups of $d$-manifolds equipped with configurations of submanifolds. It therefore follows that its groupoid of path-components $\pi_{0}(\cD_{d})$, and hence also the category $\pi_{0}(\langle \cD_{d},\cD_{d} \rangle)$ of which $\pi_{0}(\cD_{d})$ is the underlying groupoid, contains all mapping class groups of $d$-manifolds, as well as all $d$-dimensional motion groups as normal subgroups; see Remark~\ref{rmk:globality}. We may thus realise the aforementioned topologically-enriched semicategory $\langle \cG_{\circ}, \cM_{\circ} \rangle^{\mt}$ as a certain subsemicategory $\langle \cG , \cM \rangle$ of $\langle \cD_{d},\cD_{d} \rangle$; see Construction~\ref{const:category_families_of_groups}. In brief, we consider subgroupoids $\cG$ and $\cM$ of $\cD_{d}$, to which we apply the topologically-enriched Quillen bracket construction of \S\ref{sss:Topological-Quillen}, and Theorem~\ref{athm:thm_B} repeats verbatim for $\langle\cG,\cM\rangle \subseteq \langle \cD_{d},\cD_{d} \rangle$; see Corollary~\ref{cor:description_UD_{d}}.
The category $\langle \cG_{\circ}, \cM_{\circ} \rangle$ is then defined to be $\pi_{0}(\langle \cG , \cM \rangle)\cong \langle\pi_{0}(\cG),\pi_{0}(\cM)\rangle$. (More precisely, it is a certain skeleton of this category; see \S\ref{sss:skeleta}.)
In particular, our mild assumption on each family of groups $\{g_{n} \colon G_{n} \to G_{n+1}\}_{n\in\mathbb{N}}$ is to require that the groups $\{G_{n}\}_{n\in\bN}$ are the automorphism groups of a category of the form $\langle\pi_{0}(\cG),\pi_{0}(\cM)\rangle$, which also contains the maps $\{g_{n}\}_{n\in\bN}$. See Construction~\ref{const:category_families_of_groups} for more details and \S\ref{ss:categories_for_families_of_groups} for concrete examples.

\paragraph*{Applications.}

Finally, we survey the representations generated by the constructions of Theorem~\ref{athm:thm_A} for some illustrative examples introduced in \S\ref{sss:category_mcg}--\S\ref{sss:category_loop_braid_groups}, which represent the archetypal examples of mapping class groups and motion groups.
First, let us indicate which parameters $(Z,\sG,\ell)$ we use to apply Theorem~\ref{athm:thm_A} in these examples. From now on, we fix an integer $k\geq 1$ and a partition $\lambda \vdash k$. We generically denote by $\underline{k}\subset \bR^{2}$ a closed submanifold consisting of $k$ distinct points, and by $\underline{k}\bS^{1} \subset \bR^{3}$ a closed submanifold consisting of a collection of disjoint circles forming a trivial link of $k$ components. We delineate our results for each family of groups in the paragraphs below; they are also summarised in Table~\ref{table:recovering-representations}.

\underline{\emph{Classical braid groups:}} the relevant functors of Theorem~\ref{athm:thm_A} to consider for the groups $\B_{n}=\MCG(\bD^{2},\underline{n})=\Mot_{\underline{n}}(\bD^{2})$ are those of the form $L_{i}(\cF_{(Z,\sG,\ell)})$ with $d=2$, $Z:=\underline{k}$, $\sG:=\Sym_{\lambda}$, $\ell\geq 1$ and $i\geq 0$; see \eqref{eq:output_of_general_construction_classical_braids}.

\begin{athm}[{Theorem~\ref{thm:LawrenceBigelowFunctors}, \S\ref{ss:applications_motion_groups}}]
\label{athm:thm_C_classical_braid_groups}
For $\lambda = \boldsymbol{\{} k \boldsymbol{\}}$ the trivial partition, $\ell = 2$ and $i=k$, the functor $L_{k}(\cF_{(\underline{k},\Sym_{k},2)})$ recovers the aforementioned family $\LB_{k}(n)$ of Lawrence-Bigelow representations \textup{\cite{Lawrence1,bigelow2001braid}} of the braid groups.

Otherwise, the $\B_{n}$-representations encoded by the functors $L_{i}(\cF_{(\underline{k},\Sym_{\lambda},\ell)})$ appear to be new.
\end{athm}

\underline{\emph{Surface braid groups:}} we denote by $\Sigma_{g,1}$ (respectively $\N_{h,1}$) a compact, connected, smooth, orientable (resp.~non-orientable) surface of genus $g$ (resp. $h$) and with one boundary component. The relevant functors of Theorem~\ref{athm:thm_A} to consider for the surface braid groups $\B_{n}(\Sigma_{g,1})=\Mot_{\underline{n}}(\Sigma_{g,1})$ and $\B_{n}(\N_{h,1})=\Mot_{\underline{n}}(\N_{h,1})$ are those of the form $L_{i}(\cF_{(Z,\sG,\ell)})$ with $d=2$, $Z:=\underline{k}$, $\sG:=\Sym_{\lambda}$, $\ell\geq 1$ and $i\geq 0$; see \eqref{eq:output_of_general_construction_surface_braids}.

\begin{athm}[{Example~\ref{eg:An-Ko}, \S\ref{ss:applications_motion_groups}}]
\label{athm:thm_C_surface_braid_groups}
For $\lambda = \boldsymbol{\{} k \boldsymbol{\}}$ the trivial partition, $\ell = 3$ and $i=k$, the variant using Borel-Moore homology of the functor $L_{k}(\cF_{(\underline{k},\Sym_{k},3)})$ recovers the family of An-Ko representations \textup{\cite{AnKo}} of the surface braid groups $\B_{n}(\Sigma_{g,1})$.

Otherwise, the representations of the groups $\B_{n}(\Sigma_{g,1})$ and $\B_{n}(\N_{h,1})$ encoded by the functors $L_{i}(\cF_{(\underline{k},\Sym_{\lambda},\ell)})$ appear to be new.
\end{athm}

\underline{\emph{Loop braid groups:}} for the extended and non-extended loop braid groups $\lB'_{n}$ and $\lB_{n}$ (see \S\ref{sss:category_loop_braid_groups} for their definitions), the relevant functors of Theorem~\ref{athm:thm_A} to consider are those of the form $L_{i}(\cF_{(Z,\sG,\ell)})$ with $d=3$, $\ell\geq 1$ and $i\geq 0$, and:
\begin{itemizeb}
    \item either $Z:=\underline{k}$ together with $\sG:=\Sym_{\lambda}$, see \eqref{eq:output_of_general_construction_loop_braids_points};
    \item or else $Z:=\underline{k}\bS^{1}$ together with $\sG:=\Diff^{+}(\underline{\lambda}\bS^{1})$ or $\Diff(\underline{\lambda}\bS^{1})$, see \eqref{eq:output_of_general_construction_loop_braids_unlinks}, \eqref{eq:output_of_general_construction_loop_braids_unlinks_not_preserved}.
\end{itemizeb}

\begin{athm}[{Example~\ref{eg:loop_Burau}, \S\ref{ss:applications_motion_groups}}]
\label{athm:thm_C_loop_braid_groups}
For $k=1$ and $\lambda = \boldsymbol{\{} 1 \boldsymbol{\}}$ its trivial partition, $\ell = 2$ and $i=1$, the functor $L_{1}(\cF_{(\underline{1},0,2)})$ recovers the family of {loop Burau} representations \textup{\cite{PS0}} of the of extended and non-extended loop braid groups.

Otherwise, the representations of the groups $\lB'_{n}$ and $\lB_{n}$  defined by the functors $L_{i}(\cF_{(\underline{k},\Sym_{\lambda},\ell)})$, $L_{i}(\cF_{(\underline{k}\bS^{1},\Diff^{+}(\underline{\lambda}\bS^{1}),\ell)})$ and $L_{i}(\cF_{(\underline{k}\bS^{1},\Diff(\underline{\lambda}\bS^{1}),\ell)})$ appear to be new.
\end{athm}

\underline{\emph{Mapping class groups of surfaces:}} for the groups $\MCGo_{g,1}=\MCG(\Sigma_{g,1})$ and $\MCGno_{h,1}=\MCG(\N_{h,1})$, we consider the homological representation functors $L_{i}(\cF_{(Z,\sG,\ell)})$ and $L_{i}(\sF_{(Z,\sG,\ell)})$ of Theorem~\ref{athm:thm_A} with $d=2$, $Z:=\underline{k}$, $\sG:=\Sym_{\lambda}$, $\ell\geq 1$ and $i\geq 0$; see \eqref{eq:output_of_general_construction_mapping_class_motion} and \eqref{eq:output_of_general_construction_mcg_mcg}.

\begin{athm}[{Proposition~\ref{prop:Moriyama_recovering}, \S\ref{ss:mcg_construction}}]
\label{athm:thm_C_MCG}
For $\lambda = \boldsymbol{\{} 1^{k}\boldsymbol{\}}$ the discrete partition of $k$, $\ell = 1$ and $i=k$, the variant using Borel-Moore homology of the functor $L_{k}(\cF_{(\underline{k},0,1)})$ encodes the family of Moriyama representations \textup{\cite{Moriyama}} of the mapping class groups $\MCGo_{g,1}$.

Otherwise, the representations of the groups $\MCGo_{g,1}$ and $\MCGno_{h,1}$ encoded by the functors $L_{i}(\cF_{(\underline{k},\Sym_{\lambda},\ell)})$ and $L_{i}(\sF_{(\underline{k},\Sym_{\lambda},\ell)})$ appear to be new.
\end{athm}

\begin{table}[t]
\begin{center}
\begingroup
\small
\def\arraystretch{1.4}
\begin{tabular}{!{\vrule width 1pt}c|c|c|c|c|c|c|c|c!{\vrule width 1pt}}
\noalign{\hrule height 1pt}

\multirow{2}{*}{\normalsize Thm.} & \multirow{2}{*}{\normalsize Family of groups} & \multicolumn{5}{c|}{\normalsize Parameters} & {\normalsize Family of} & \multirow{2}{*}{\normalsize Reference} \\
\cline{3-7}

&& $d$ & $Z$ & $\sG$ & $\ell$ & $i$ & {\normalsize representations} & \\
\hline

\ref{athm:thm_C_classical_braid_groups} & $\B_n$ & $2$ & $\underline{k}$ & $\Sym_k$ & $2$ & $k$ & Lawrence-Bigelow & \cite{Lawrence1,bigelow2001braid} \\
\hline

\ref{athm:thm_C_surface_braid_groups} & $\B_n(\Sigma_{g,1})$ ($g$ fixed) & $2$ & $\underline{k}$ & $\Sym_k$ & $3$ & $k$ & An-Ko & \cite{AnKo} \\
\hline

\ref{athm:thm_C_loop_braid_groups} & $\lB'_n$ and $\lB_n$ & $3$ & $\underline{1}$ & $\{\id\}$ & $2$ & $1$ & Loop Burau & \cite{PS0} \\
\hline

\ref{athm:thm_C_MCG} & $\MCGo_{g,1}$ & $2$ & $\underline{k}$ & $\{\id\}$ & $1$ & $k$ & Moriyama & \cite{Moriyama} \\

\noalign{\hrule height 1pt}
\end{tabular}
\endgroup
\end{center}
\caption{A summary of Theorems~\ref{athm:thm_C_classical_braid_groups}--\ref{athm:thm_C_MCG}, describing the special cases of our construction that recover previously-known homological representations.}
\label{table:recovering-representations}
\end{table}

\paragraph*{Unified foundations.}

As detailed above, several examples of the representations arising from the construction of Theorem~\ref{athm:thm_A} (see Theorems~\ref{athm:thm_C_classical_braid_groups}--\ref{athm:thm_C_MCG} and Table~\ref{table:recovering-representations}) have previously been defined and studied in the literature, at least at the level of individual groups, i.e.~when restricted to the automorphism groups of $\langle \cG_{\circ}, \cM_{\circ} \rangle$. One purpose of setting up the general procedure for constructing homological representations in this paper is to give a \emph{unified foundation} that encompasses all known homological representations, amongst many others, and that inspires the discovery of new representations by applying this unified context to novel settings. Regarding this last point, we mention that the examples described in Theorems~\ref{athm:thm_C_classical_braid_groups}--\ref{athm:thm_C_MCG} represent only a small fragment of the potential of the general construction of Theorem~\ref{athm:thm_A}.

\paragraph*{Perspectives on the construction.}

Although we do not discuss it in more detail in this paper, one may equally well consider examples analogous to those of Theorems~\ref{athm:thm_C_classical_braid_groups}--\ref{athm:thm_C_MCG}, given by the homological representation functors of Theorem~\ref{athm:thm_A} on subcategories of $\fU\cD_{2}$ and $\fU\cD_3$ (or $\fU\cD_{d}$ for $d\geq5$) relevant for automorphism groups of free groups, Torelli groups, handlebody mapping class groups, pure braid groups, as well as higher-dimensional mapping class groups and motion groups.
In addition, there are several natural variations of our constructions of Theorem~\ref{athm:thm_A}. First, one may change the ``flavour'' of (twisted) homology that we use: notably \emph{Borel-Moore} homology rather than ordinary homology (see Theorem~\ref{thm:construction}).
Second, we may naturally package together the homological representation functors \eqref{eq:ThmA_hom_rep_functors_untwisted} by fixing the parameters $(Z,\sG)$, and by considering simultaneously all the lower central series parameters $\{\ell \geq 1\}$, which control the ground ring over which the representation is defined (for each $\ell$, it is the group-ring of a group of nilpotency class at most $\ell - 1$). More precisely, these functors fit together into a tower that may be thought of as a single \emph{pro-nilpotent} (functorial) representation, which may be truncated to recover the representation corresponding to each level $\ell$; see \cite[\S 1]{PSIN}.
In particular, there are many interesting cases where we obtain in this way an infinite tower of distinct representations as $\ell \to \infty$; see \cite[\S 4--\S 6]{PSIN}.

\paragraph*{Outline.}

The aim of \S\ref{s:categorical_framework} is to set up the source categories for homological representation functors. Namely, we define the topologically-enriched groupoids $\cD_{d}$ of decorated manifolds in \S\ref{ss:topological_groupoids_of_decorated_manifolds}. Then, we prove some fundamental properties of decorated diffeomorphisms and embeddings of manifolds: in \S\ref{ss:split-ses} we exhibit several split short exact sequences associated to these topological groups, which are crucial for the construction of homological representation functors in \S\ref{s:general_construction}, and then we prove in \S\ref{ss:fibration-condition} that certain quotients of decorated diffeomorphism groups are Serre fibrations (see Theorem~\ref{thm:fibre-bundle}); these are decisive in the proof of Theorem~\ref{athm:thm_B}. Finally, we introduce the topological enrichment of the Quillen bracket construction in \S\ref{sss:Topological-Quillen} and then in \S\ref{sss:quillen-bracket-categories} we apply it to categories of decorated manifolds and study its main properties, in particular proving Theorem~\ref{athm:thm_B}.

In \S\ref{s:general_construction}, we introduce the various tools and topological constructions to define homological representation functors, which between them imply Theorem~\ref{athm:thm_A}. In more detail, we present in \S\ref{ss:ingredients_construction_representations} the different ingredients and steps of the construction illustrated in diagram \eqref{eq:diagram-homological-representations}, considering an abstract continuous semifunctor of the form \eqref{eq:ThmA_semifunctors}; see Theorem~\ref{thm:construction}. We then construct the two continuous semifunctors $\cF_{(Z,\sG,\ell)}$ and $\sF_{(Z,\sG,\ell)}$ of \eqref{eq:ThmA_semifunctors} in \S\ref{ss:homological_representation_functor_motion_groups} and \S\ref{ss:homological_representation_functor_mcg} respectively (see Theorems~\ref{thm:global_functor_motion_groups} and \ref{thm:global_functor_mcg}), thus completing the proof of Theorem~\ref{athm:thm_A}.

We finally apply our constructions of homological representation functors to the four different families of groups of Theorems~\ref{athm:thm_C_classical_braid_groups}--\ref{athm:thm_C_MCG} in \S\ref{s:applications}. Namely, we first describe (specialising \S\ref{s:categorical_framework} to each case) the appropriate categorical framework encoding these families of groups in \S\ref{ss:categories_for_families_of_groups}. We deal with the applications for the families of motion groups in \S\ref{ss:applications_motion_groups}, proving Theorems~\ref{athm:thm_C_classical_braid_groups}, \ref{athm:thm_C_surface_braid_groups} and \ref{athm:thm_C_loop_braid_groups}, and we treat the case of the mapping class groups in \S\ref{ss:mcg_construction}, proving Theorem~\ref{athm:thm_C_MCG}.

\paragraph*{General notation.}

For a category $\cC$, we use the abbreviation $\mathrm{ob}(\cC)$ to denote its class of objects.
We recall that \emph{semi}categories and \emph{semi}functors between them are defined in the exact same way as (genuine) categories and functors by relaxing all structures or conditions that involve identity morphisms of objects.
Throughout the paper, the word \emph{ring} always refers to an \emph{associative, unital ring} and we assume that ring homomorphisms preserve units. For a ring $R$, we denote by $\modr[R]$ the category of (right) $R$-modules. For an $R$-module $M$, we denote by $\Aut_{R}(M)$ the group of $R$-module automorphisms of $M$. When $R=\bZ$, we omit it from the notation as long as there is no ambiguity.

For a manifold $X$, we denote by $\Diff(X)$ its group of diffeomorphisms. If $X$ is orientable, we denote by $\Diff^{+}(X)$ the subgroup of orientation-preserving diffeomorphisms. If $X$ has non-empty boundary, $\Int{X}$ denotes its interior. For each $d\geq 1$, we denote by $\bS^{d}$ the $d$-sphere.

We denote by $\Sym_{n}$ the symmetric group on a set of $n$ elements. For an integer $n \geq 1$, an \emph{ordered partition of $n$} means an ordered $r$-tuple $\bn = \boldsymbol{\{} n_{1}, \ldots, n_r\boldsymbol{\}}$ of integers $n_{i} \geq 1$ (for some $r \geq 1$ called the \emph{length} of $\bn$) such that $n = \sum_{1\leq i \leq r} n_{i}$ (and without the condition $n_{i}\geq n_{i+1}$).
We recall that the \emph{lower central series} of a group $G$ is the descending chain of subgroups $\{ \LCS_{\ell}(G)\} _{\ell\geq 1}$ defined by $\LCS_{1}(G):=G$ and $\LCS_{\ell + 1}(G):=[G,\LCS_{\ell}(G)]$, the subgroup of $G$ generated by the commutators $[g,h]$ for $g \in G$ and $h \in \LCS_{\ell}(G)$. For the sake of simplicity, each quotient $G/\LCS_{\ell}(G)$ is denoted by $G/\LCS_{\ell}$.

\paragraph*{Acknowledgements.}

The authors would like to thank Paolo Bellingeri, Tara Brendle, Oscar Randal-Williams, Antoine Touz{\'e}, Christine Vespa and Emmanuel Wagner for illuminating discussions and questions. In particular, they thank Antoine Touz{\'e} for pointing out the reference \cite{ErdmannNakano}. They would also like to thank Oscar Randal-Williams for inviting the first author to the University of Cambridge in November 2019, where the authors were able to make significant progress on the present article.
The authors would also like to thank the anonymous referee for very helpful comments and suggestions on earlier versions of this work.

\tableofcontents

\section{Categories of decorated manifolds and their embeddings}\label{s:categorical_framework}

In this section, we define the categories that will serve as the domain of the homological representation functors that we will construct in \S\ref{s:general_construction}. They are obtained from certain monoidal groupoids by the \emph{Quillen bracket construction}, an operation that enlarges a given monoidal groupoid to a category having the original monoidal groupoid as its underlying groupoid.

More precisely, we will start with certain \emph{topologically-enriched} monoidal groupoids, so we describe, in \S\ref{sss:Topological-Quillen}, a topological enrichment of the Quillen bracket construction and show -- subject to a \emph{Serre fibration condition} -- that it behaves well with respect to the functor $\pi_{0}$ that replaces all morphism spaces with their sets of path-components. In \S\ref{ss:topological_groupoids_of_decorated_manifolds}, we define the topologically-enriched monoidal groupoids that we wish to consider, and prove in \S\ref{ss:fibration-condition} that they satisfy this Serre fibration condition.

Informally, the idea is that the domain category, for a given dimension $d\geq 2$, will be a topologically-enriched category $\fU\cD_{d}$ having the property that the automorphism groups of $\pi_{0}(\fU\cD_{d})$ contain all mapping class groups and motion groups in dimension $d$. To construct this, we define in \S\ref{ss:topological_groupoids_of_decorated_manifolds} a topologically-enriched groupoid $\cD_{d}$ whose automorphism groups are the diffeomorphism groups of all $d$-dimensional \emph{decorated manifolds}. The topologically-enriched Quillen bracket construction $\fU$ of \S\ref{sss:Topological-Quillen} then gives us a topologically-enriched category $\fU\cD_{d}$ such that $\pi_{0}(\fU\cD_{d}) \cong \fU(\pi_{0}(\cD_{d}))$, where the latter category is defined by appyling the construction $\fU$ to the path-component category $\pi_{0}(\cD_{d})$ of $\cD_{d}$. The underlying groupoid of $\pi_{0}(\fU\cD_{d})$ is therefore $\pi_{0}(\cD_{d})$, consisting of all mapping class groups of $d$-dimensional \emph{decorated manifolds}, which contain all $d$-dimensional motion groups as normal subgroups.

In \S\ref{sss:quillen-bracket-categories}, we then give an explicit description of morphism spaces of \emph{Quillen bracket categories of manifolds} in terms of embedding spaces, which is a crucial ingredient in the construction of \S\ref{s:general_construction}.

\paragraph*{Outline.}

To summarise, the main points of this section are the following.
\begin{itemizeb}
    \item We construct the topologically-enriched groupoid $\cD_{d}$ of decorated manifolds in \S\ref{ss:topological_groupoids_of_decorated_manifolds}.
    \item We obtain the topologically-enriched category $\fU\cD_{d}$ using the topologically-enriched Quillen bracket construction introduced in \S\ref{sss:Topological-Quillen}.
    \item We identify the morphism spaces of $\fU\cD_{d}$ with certain embedding spaces in Proposition~\ref{prop:morphism-spaces-bracket} (in \S\ref{sss:quillen-bracket-categories}). This crucially uses a Serre fibration result for decorated diffeomorphism groups proved in Theorem~\ref{thm:fibre-bundle} (in \S\ref{ss:fibration-condition}).
    \item The categories that we shall use, in \S\ref{s:applications}, to encode families of mapping class groups and motion groups, are defined by certain topologically-enriched Quillen bracket constructions $\langle \cG,\cM\rangle$, which come with natural inclusions $\langle \cG,\cM\rangle \hookrightarrow \fU\cD_{d}$. The homological representation functors (defined on $\fU\cD_{d}$) constructed in \S\ref{s:general_construction} therefore restrict naturally along these inclusions. In \S\ref{ss:categories_for_families_of_groups}, we describe in detail the relevant categories $\langle \cG,\cM\rangle$ for mapping class groups and motion groups.
\end{itemizeb}

\paragraph*{Preliminaries on categorical tools.}

We refer to \cite[Chap.~VII]{MacLane1} for a complete introduction to the notions of monoidal categories and modules over them. We generically denote a monoidal category by $(\cC,\natural,\zero)$, where $\cC$ is a category, $\natural$ is the monoidal product and $\zero$ is the monoidal unit.
A left-module $(\cM,\sharp)$ over $(\cC,\natural,\zero)$ is a category $\cM$ with a functor $\sharp\colon\cC\times\cM\to\cM$ that is unital and associative. Similarly, \emph{semi-}monoidal categories and left-modules over them are defined in the exact same way as in the genuine monoidal case by relaxing all the structures or conditions involving (left or right) units; in other words they only involve a binary operation admitting an associator that satisfies the pentagon condition.
Note that a (semi-)monoidal category $(\cC,\natural)$ is equipped with a left-module structure over itself, induced by its monoidal product. Each left-module structure $\sharp$ in this paper is defined from some underlying (semi-)monoidal structure $\natural$ (see \S\ref{ss:topological_groupoids_of_decorated_manifolds} and \S\ref{ss:categories_for_families_of_groups}), so we abuse notation by using the same symbol $\natural$ for $\sharp$.
The reason for introducing this \emph{semi-monoidal} setting is that the boundary connected sum of manifolds naturally equips the topologically-enriched groupoid $\cD_{d}$ with a non-unital such structure; see Remark~\ref{r:semi-monoidal} for more details.

We recall that, by \emph{topologically-enriched category}, we mean a category enriched over the symmetric monoidal category of topological spaces with its Cartesian product. The functor $\pi_{0}$ replaces each morphism space with its set of path-components, and so defines a functor from (small) topologically-enriched categories to (small) categories. For any topologically-enriched category $\cC$, there is a natural functor $\cC \to \pi_{0}(\cC)$, also denoted by $\pi_{0}$, sending each point of a morphism space to the path-component that it lies in.

\begin{rmk}\label{rmk:braidings_symmetries}
Although most of the (semi-)monoidal categories and their modules that we consider in this paper are \emph{symmetric} or \emph{braided}, we will not dwell on these considerations, as these properties are unnecessary for our work, in particular to define the topologically-enriched bracket construction of \S\ref{ss:Quillen_bracket_construction} or the homological representation functors in \S\ref{s:general_construction}.
\end{rmk}

\subsection{Decorated manifolds and their diffeomorphisms}\label{ss:topological_groupoids_of_decorated_manifolds}

In this section, we introduce the topologically-enriched groupoids $\cD_{d}$ of \emph{decorated manifolds} and their semi-monoidal structures induced by boundary connected sum. We fix an integer $d \geq 2$.

\begin{convention}\label{convention:small_categories}
All of the (topological) categories that we shall consider are \emph{essentially small}, i.e.~they are each equivalent to a \emph{small} (topological) category. One standard way to see this is as follows; see \cite[Rem.~1.2.7]{Galatius}. In each case, the objects of the categories that we shall consider consist of manifolds, submanifolds and some additional data, such as collar neighbourhoods, etc. Fixing a set $\Omega$, we may consider the full subcategory whose objects are only those where the underlying set of the manifold is a subset of $\Omega$. This is a small category, and if the cardinality of $\Omega$ is at least $\lvert \bR \rvert$, its inclusion into the whole (large) category is essentially surjective on objects; thus it is an equivalence. (We assume that manifolds are second-countable, which implies that their cardinality is no larger than $\lvert \bR \rvert$.) In the following, we implicitly fix a sufficiently large set $\Omega$ and we shall consider the corresponding small subcategories, whenever we need our categories to be small.
\end{convention}

\paragraph*{Open path-components.}

For a topological space $X$, if we give its set of path-components $\pi_0(X)$ the discrete topology, the natural map $X \to \pi_0(X)$ sending each point to the path-component that it lies in is continuous if and only if all path-components of $X$ are open. Similarly, for a topologically-enriched category $\cC$, the natural functor $\cC \to \pi_0(\cC)$ mentioned above is continuous if and only if all path-components of all morphisms spaces of $\cC$ are open. It will be important for us to ensure that this property holds, so we collect here some facts about this property that we shall make use of. We first record some basic facts:

\begin{lemm}
\label{lem:open-path-components}
If $X$ has the property that its path-components are open, then the quotient space $X/{\sim}$ also has this property, for any equivalence relation $\sim$ on $X$. If a collection $X_i$ of spaces all have this property, then so does their disjoint union $\bigsqcup_i X_i$.
\end{lemm}

There is a canonical way to refine the topology on an arbitrary space in order to force this property to hold.

\begin{construction}
\label{construction:open-path-components}
For a topological space $X$, denote by $o(X)$ the topological space with the same underlying set as $X$, equipped with the topology generated by the base consisting of $C \cap U$ for $U$ an open subset of $X$ and $C$ a path-component of $X$.
\end{construction}

One may easily verify the following properties of this construction.

\begin{lemm}
\label{lem:properties-of-o}
The path-components of $o(X)$ are open, and they are the same as the path-components of $X$. The canonical map $o(X) \to X$ defined by the identity map of the underlying sets is continuous, and becomes a homeomorphism when restricted to any path-component. Consequently, it is also a weak equivalence. In addition, the functor $o(-)$ preserves products.
\end{lemm}

\begin{rmk}
\label{rmk:coreflective-subcategory}
In fact, the operation $o(-)$ is a right adjoint to the inclusion of the full subcategory of topological spaces having the property that their path-components are open; this subcategory is therefore \emph{coreflective}.
\end{rmk}

Finally, we will need the following lemma, which explains how the operation $o(-)$ interacts with quotient maps that are Serre fibrations.

\begin{lemm}
\label{lem:open-path-components-Serre-fibrations}
Let $X$ be a topological space and $\sim$ an equivalence relation on $X$. If the quotient map $X \to X/{\sim}$ is a Serre fibration, then so is the quotient map $o(X) \to o(X)/{\sim}$.
\end{lemm}
\begin{proof}
The continuous bijection $o(X) \to X$ induces a continuous bijection $o(X)/{\sim} \to X/{\sim}$, fitting into a commutative square with the two quotient maps. Let us consider the lifting problem:
\begin{equation}
\label{eq:lifting-problem-o}
\begin{tikzcd}
{[0,1]^{d-1}} \ar[rr] \ar[d] && o(X) \ar[d] \ar[r] & X \ar[d] \\
{[0,1]^d} \ar[rr] \ar[densely dashed,rru,"?"] && o(X)/{\sim} \ar[r] & X/{\sim}.
\end{tikzcd}
\end{equation}
Since $X \to X/{\sim}$ is a Serre fibration by hypothesis, there is a diagonal map $[0,1]^d \to X$ making the two outer triangles commute. Since the cube $[0,1]^d$ is path-connected, its image in $X$ lies in a single path-component. Restricted to this path-component, the continuous bijection $o(X) \to X$ is a homeomorphism by Lemma~\ref{lem:properties-of-o}, so we may lift the diagonal map $[0,1]^d \to X$ to a diagonal map $[0,1]^d \to o(X)$ for the left-hand square, as desired. The two triangles into which this splits the left-hand square of \eqref{eq:lifting-problem-o} commute, since the triangles that the outer rectangle was split into commute and the two right-hand horizontal maps in \eqref{eq:lifting-problem-o} are bijective.
\end{proof}

\begin{rmk}
Lemma~\ref{lem:open-path-components-Serre-fibrations} holds more generally for any class of fibrations defined by the right lifting property with respect to a collection of maps whose codomains are path-connected. Moreover, it also holds more generally for any two topologies on the same underlying set that have the same path-components and that agree on each path-component.
\end{rmk}

\paragraph*{Categories of decorated manifolds.}

First, we define the notion of \emph{decorated manifolds} of dimension $d$, their morphisms and associated categories. A motivation to work with this type of manifolds is that the groups $\pi_{0}(\Diff(-))$ of decorated manifolds will contain as normal subgroups all of the mapping class groups and motion groups that we are interested in; see Remark~\ref{rmk:globality}. Also, the notion of decorated manifolds is formed so that their associated categories have a well-defined semi-monoidal structure induced by the boundary connected sum, which would not be the case if we just considered manifolds without decorations; see Remark~\ref{rmk:boundary_connected_sum_comparison}. The terminology ``decorated'' originates from \cite{BodigheimerTillmann}; see also \cite{BasualdoBonatto}.

\begin{notation}[Solid cylinders.]
\label{notation:cylinder}
We denote by $\bD^{d-1}$ the closed unit $(d-1)$-dimensional disc in $\bR^{d-1}$ in the usual $L^2$ metric. For a real number $t > 0$, we write $\bB^d_t = \bD^{d-1} \times [0,t]$
for the \emph{solid $d$-dimensional cylinder of height $t$}. We also write
\[
\partial_{l} \bB^d_t = (\partial \bD^{d-1} \times [0,t]) \cup (\bD^{d-1} \times \{0\})
\]
and call this the \emph{lower boundary} of $\bB^d_t$, as well as $b\bB^d_t = \bD^{d-1} \times \{0\}$
and call this the \emph{base} of $\bB^d_t$. This is illustrated in Figure~\ref{fig:boundary-cylinder}.
\end{notation}

\begin{figure}[t]
\centering
\begin{tikzpicture}
[x=1mm,y=1mm,scale=1]
\fill[green!40] (0,0) ellipse (3 and 10);
\fill[yellow!10] (20,0) ellipse (3 and 10);
\fill[yellow!80,opacity=0.5] (0,10) arc (90:270:3 and 10) -- (20,-10) arc (270:90:3 and 10) -- cycle;
\draw[thick] (0,10) -- (20,10);
\draw[thick] (0,-10) -- (20,-10);
\draw[very thick,blue!70] (0,10) arc (90:270:3 and 10);
\draw[very thick,blue!70,densely dashed] (0,-10) arc (-90:90:3 and 10);
\draw[thick] (20,0) ellipse (3 and 10);
\node at (0,-10) [anchor=north,font=\footnotesize] {$0$};
\node at (20,-10) [anchor=north,font=\footnotesize] {$t$};
\end{tikzpicture}
\caption{An illustration of the notation for the solid cylinder $\bB^d_t$ from Notation~\ref{notation:cylinder} for $d=3$. Its lower boundary $\partial_{l} \bB^3_t$ is coloured yellow, its base $b\bB^3_t$ is yellow-green and the codimension-$2$ stratum $\partial (b\bB^3_t) = \partial \bD^{2} \times \{0\}$ is blue.}
\label{fig:boundary-cylinder}
\end{figure}

\begin{defi}[Boundary-cylinders.]
\label{def:boundary-cylinder}
Let $M$ be a smooth $d$-manifold with non-empty boundary. A \emph{boundary-cylinder} for $M$ is a topological embedding $e \colon \bB^d_{1} \hookrightarrow M$ such that $e^{-1}(\partial M) = \partial_{l} \bB^d_{1}$ and $e$ is a smooth embedding on $\bB^d_{1} \smallsetminus \partial( b\bB^d_{1})$. Two boundary-cylinders $e,e'$ are \emph{equivalent} if they are equal when restricted to $\bB^d_\epsilon \subseteq \bB^d_{1}$ for some $\epsilon > 0$. An equivalence class of boundary-cylinders is called a \emph{boundary-cylinder germ}.
\end{defi}

\begin{defi}[Decorated manifolds.]
\label{def:decorated-manifolds}
A \emph{decorated manifold} is a smooth $d$-manifold $M$ with non-empty boundary, equipped with a closed submanifold $A \subset \Int{M}$ and a pair $(e_{1},e_{2})$ of boundary-cylinder germs for $M \smallsetminus A$ such that $e_{1}(b\bB^d_{1})$ and $e_{2}(b\bB^d_{1})$ are disjoint. See Figure~\ref{fig:decorated-manifolds} for a schematic picture.
\end{defi}

\begin{figure}[t]
\centering
\begin{tikzpicture}
[x=1mm,y=1mm,scale=1]
\begin{scope}
\fill[blue!20] (0,0) -- (10,0) .. controls (15,0) and (15,-5) .. (20,-5) -- (21,-5) -- (29,-5) -- (30,-5) .. controls (35,-5) and (35,0) .. (40,0) -- (50,0) -- (50,10) -- (40,10) .. controls (35,10) and (35,15) .. (30,15) -- (29,15) -- (21,15) -- (20,15) .. controls (15,15) and (15,10) .. (10,10) -- (0,10) -- (0,0);
\draw[thick] (0,0) -- (10,0) .. controls (15,0) and (15,-5) .. (20,-5) -- (21,-5);
\draw[thick,densely dashed] (21,-5) -- (29,-5);
\draw[thick] (29,-5) -- (30,-5) .. controls (35,-5) and (35,0) .. (40,0) -- (50,0) -- (50,10) -- (40,10) .. controls (35,10) and (35,15) .. (30,15) -- (29,15);
\draw[thick,densely dashed] (29,15) -- (21,15);
\draw[thick] (21,15) -- (20,15) .. controls (15,15) and (15,10) .. (10,10) -- (0,10) -- (0,0);
\draw[blue!50,thick] (18,2) circle (2);
\draw[blue!50,thick] (28,7) circle (2);
\fill[blue!50,thick] (22,9) circle (0.5);
\node[blue!70] at (30,0) {$A$};
\draw (10,0) -- (10,10);
\draw (40,0) -- (40,10);
\node at (0,0) [anchor=north,font=\footnotesize] {$0$};
\node at (10,0) [anchor=north,font=\footnotesize] {$\vphantom{0}\epsilon$};
\node at (40,0) [anchor=north,font=\footnotesize] {$\vphantom{0}\epsilon$};
\node at (50,0) [anchor=north,font=\footnotesize] {$0$};
\node at (5,5) {$e_{1}$};
\node at (45,5) {$e_{2}$};
\node at (0,5) [anchor=east] {$M=$};
\end{scope}
\begin{scope}[xshift=60mm]
\fill[green!20] (0,0) -- (10,0) .. controls (15,0) and (15,-5) .. (20,-5) -- (21,-5) -- (29,-5) -- (30,-5) .. controls (35,-5) and (35,0) .. (40,0) -- (50,0) -- (50,10) -- (40,10) .. controls (35,10) and (35,15) .. (30,15) -- (29,15) -- (21,15) -- (20,15) .. controls (15,15) and (15,10) .. (10,10) -- (0,10) -- (0,0);
\draw[thick] (0,0) -- (10,0) .. controls (15,0) and (15,-5) .. (20,-5) -- (21,-5);
\draw[thick,densely dashed] (21,-5) -- (29,-5);
\draw[thick] (29,-5) -- (30,-5) .. controls (35,-5) and (35,0) .. (40,0) -- (50,0) -- (50,10) -- (40,10) .. controls (35,10) and (35,15) .. (30,15) -- (29,15);
\draw[thick,densely dashed] (29,15) -- (21,15);
\draw[thick] (21,15) -- (20,15) .. controls (15,15) and (15,10) .. (10,10) -- (0,10) -- (0,0);
\draw[green!80!black,thick] (17,7) circle (2);
\draw[green!80!black,thick] (23,2) circle (2);
\draw[green!80!black,thick] (28,8) circle (1.5);
\fill[green!80!black,thick] (29,5) circle (0.5);
\node[green!80!black] at (30,0) {$A'$};
\draw (10,0) -- (10,10);
\draw (40,0) -- (40,10);
\node at (0,0) [anchor=north,font=\footnotesize] {$0$};
\node at (10,0) [anchor=north,font=\footnotesize] {$\vphantom{0}\epsilon$};
\node at (40,0) [anchor=north,font=\footnotesize] {$\vphantom{0}\epsilon$};
\node at (50,0) [anchor=north,font=\footnotesize] {$0$};
\node at (5,5) {$e'_{1}$};
\node at (45,5) {$e'_{2}$};
\node at (50,5) [anchor=west] {$=M'$};
\end{scope}
\begin{scope}[xshift=5mm,yshift=-30mm]
\fill[blue!20] (0,0) -- (10,0) .. controls (15,0) and (15,-5) .. (20,-5) -- (21,-5) -- (29,-5) -- (30,-5) .. controls (35,-5) and (35,0) .. (40,0) -- (50,0) -- (50,10) -- (40,10) .. controls (35,10) and (35,15) .. (30,15) -- (29,15) -- (21,15) -- (20,15) .. controls (15,15) and (15,10) .. (10,10) -- (0,10) -- cycle;
\draw[thick] (0,0) -- (10,0) .. controls (15,0) and (15,-5) .. (20,-5) -- (21,-5);
\draw[thick,densely dashed] (21,-5) -- (29,-5);
\draw[thick] (29,-5) -- (30,-5) .. controls (35,-5) and (35,0) .. (40,0) -- (50,0);
\draw[thick] (50,10) -- (40,10) .. controls (35,10) and (35,15) .. (30,15) -- (29,15);
\draw[thick,densely dashed] (29,15) -- (21,15);
\draw[thick] (21,15) -- (20,15) .. controls (15,15) and (15,10) .. (10,10) -- (0,10) -- (0,0);
\draw[blue!50,thick] (18,2) circle (2);
\draw[blue!50,thick] (28,7) circle (2);
\fill[blue!50,thick] (22,9) circle (0.5);
\node[blue!70] at (30,0) {$A$};
\draw (10,0) -- (10,10);
\node at (0,0) [anchor=north,font=\footnotesize] {$0$};
\node at (10,0) [anchor=north,font=\footnotesize] {$\vphantom{0}\epsilon$};
\node at (5,5) {$e_{1}$};
\node at (0,5) [anchor=east] {$M \natural M'=$};
\end{scope}
\begin{scope}[xshift=55mm,yshift=-30mm]
\fill[green!20] (0,0) -- (10,0) .. controls (15,0) and (15,-5) .. (20,-5) -- (21,-5) -- (29,-5) -- (30,-5) .. controls (35,-5) and (35,0) .. (40,0) -- (50,0) -- (50,10) -- (40,10) .. controls (35,10) and (35,15) .. (30,15) -- (29,15) -- (21,15) -- (20,15) .. controls (15,15) and (15,10) .. (10,10) -- (0,10) -- cycle;
\draw[thick] (0,0) -- (10,0) .. controls (15,0) and (15,-5) .. (20,-5) -- (21,-5);
\draw[thick,densely dashed] (21,-5) -- (29,-5);
\draw[thick] (29,-5) -- (30,-5) .. controls (35,-5) and (35,0) .. (40,0) -- (50,0) -- (50,10) -- (40,10) .. controls (35,10) and (35,15) .. (30,15) -- (29,15);
\draw[thick,densely dashed] (29,15) -- (21,15);
\draw[thick] (21,15) -- (20,15) .. controls (15,15) and (15,10) .. (10,10) -- (0,10);
\draw[green!80!black,thick] (17,7) circle (2);
\draw[green!80!black,thick] (23,2) circle (2);
\draw[green!80!black,thick] (28,8) circle (1.5);
\fill[green!80!black,thick] (29,5) circle (0.5);
\node[green!80!black] at (30,0) {$A'$};
\draw (40,0) -- (40,10);
\node at (40,0) [anchor=north,font=\footnotesize] {$\vphantom{0}\epsilon$};
\node at (50,0) [anchor=north,font=\footnotesize] {$0$};
\node at (45,5) {$e'_{2}$};
\end{scope}
\end{tikzpicture}
\caption{Two decorated manifolds and their boundary connected sum.}
\label{fig:decorated-manifolds}
\end{figure}

\begin{defi}[Morphisms of decorated manifolds.]
\label{def:decorated-manifolds-morphisms}
A morphism of decorated manifolds from $(M,A,e_{1},e_{2})$ to $(M',A',e'_{1},e'_{2})$ is a smooth, proper (preimages of compact subspaces are compact) map $\varphi \colon M \to M'$ such that $\varphi(A) \subseteq A'$ and such that, for some $\epsilon > 0$ and for each $i \in \{1,2\}$, we have $\varphi(e_{i}(\bB^d_\epsilon)) = e'_{i}(\bB^d_\epsilon)$ and the composition $(e'_{i})^{-1} \circ \varphi \circ e_{i} \colon \bB^d_\epsilon \rightarrow \bB^d_\epsilon$ is the identity. Write $C^\infty_{\mathrm{dec}}(M,M')$ for the set of such morphisms, where by abuse of notation we are abbreviating $(M,A,e_{1},e_{2})$ to $M$ and $(M',A',e'_{1},e'_{2})$ to $M'$.
\end{defi}

For decorated manifolds $M$ and $M'$, we equip the set $C^\infty_{\mathrm{dec}}(M,M')$ with a topology induced by a colimit of Whitney topologies, as follows.
\begin{itemizeb}
    \item Choose representative boundary cylinders for the boundary-cylinder germs $e_{i}$ and $e'_{i}$. This ensures that the condition in Definition~\ref{def:decorated-manifolds-morphisms} makes sense for a fixed $\epsilon \in (0,1)$, not just for an unspecified $\epsilon \in (0,1)$ that is quantified over.
    \item Fix $\epsilon \in (0,1)$ and write $C^\infty_{\mathrm{dec},\epsilon}(M,M')$ for the subset of $C^\infty_{\mathrm{dec}}(M,M')$ where the condition in Definition~\ref{def:decorated-manifolds-morphisms} holds for this fixed $\epsilon$. Equip each subset $C^\infty_{\mathrm{dec},\epsilon}(M,M')$ with the subspace topology induced from the smooth Whitney topology on the set $C^\infty(M,M')$ of \emph{all} smooth maps from $M$ to $M'$; for details of the Whitney topology, see for example \cite[Chap.~2]{Hirsch1976Differentialtopology}. As a set, note that $C^\infty_{\mathrm{dec}}(M,M')$ is the union of $C^\infty_{\mathrm{dec},\epsilon}(M,M')$ over all choices of $\epsilon \in (0,1)$.
\end{itemizeb}

\begin{lemm}
\label{lem:topology_well_defined}
The colimit topology determined by the increasing filtration $\{ C^\infty_{\mathrm{dec},\epsilon}(M,M') \}_{\epsilon \in (0,1)}$ does not depend on the choices of representative boundary-cylinders for the boundary-cylinder germs $e_{i}$ and $e'_{i}$.
\end{lemm}
\begin{proof}
Different choices of representative boundary-cylinders for the boundary-cylinder germs $e_{i}$ and $e'_{i}$ result in different filtrations, which are cofinal in each other, so the induced colimit topologies are equal.
\end{proof}

\begin{defi}[Morphism spaces.]
\label{def:decorated-manifolds-morphism-spaces}
The morphism set $C^\infty_{\mathrm{dec}}(M,M')$ is topologised as follows: we first consider the colimit topology induced by the increasing filtration $\{ C^\infty_{\mathrm{dec},\epsilon}(M,M') \}_{\epsilon \in (0,1)}$ and then we apply the functor $o(-)$ of Construction~\ref{construction:open-path-components} to ensure that all path-components are open. Thus we may write $C^\infty_{\mathrm{dec}}(M,M') := o\bigl(\underset{\epsilon\to 0}{\mathrm{colim}}(C^\infty_{\mathrm{dec},\epsilon}(M,M'))\bigr)$.
\end{defi}

We note that composition is continuous in this topology since the functor $o(-)$ preserves products (by Lemma~\ref{lem:properties-of-o}) and because composition of smooth, proper maps is continuous in the Whitney topology (see \cite[\S 2, Prop.~1]{Mather1969}).

\begin{rmk}
Depending on $M$ and $M'$, the topology of Definition~\ref{def:decorated-manifolds-morphism-spaces} may differ from the subspace topology inherited directly from the Whitney topology on $C^\infty(M,M')$: the colimit topology is in general finer than the directly inherited Whitney topology, and the operation $o(-)$ refines the topology further. However, these three topologies on $\cdec(M,M')$ are all weakly equivalent (see Lemma~\ref{lem:properties-of-o}). In particular, they have the same $\pi_{0}$.
\end{rmk}

We can now introduce the topologically-enriched categories associated to decorated manifolds:

\begin{defi}[Decorated manifold categories.]
\label{def:Decd}\label{def:Dd}
Let $\cD\mathrm{ec}_d$ be the topologically-enriched category defined as follows. Its objects are all \emph{decorated manifolds} $(M,A,e_{1},e_{2})$ of dimension $d$, as in Definition~\ref{def:decorated-manifolds}. The space of morphisms from $M = (M,A,e_{1},e_{2})$ to $M' = (M',A',e'_{1},e'_{2})$ is the space $C^\infty_{\mathrm{dec}}(M,M')$ whose underlying set is specified in Definition~\ref{def:decorated-manifolds-morphisms} and which is topologised in Definition~\ref{def:decorated-manifolds-morphism-spaces}.

Let $\cD_{d}$ be the underlying topologically-enriched groupoid of $\cD\mathrm{ec}_d$. In other words, its objects are all decorated manifolds of dimension $d$ and its morphisms are those morphisms $(M,A,e_{1},e_{2}) \to (M',A',e'_{1},e'_{2})$ of decorated manifolds whose underlying smooth map $\varphi \colon M \to M'$ is a diffeomorphism and $\varphi(A) = A'$.

Let $\cD_{d}^{+}$ denote the topologically-enriched groupoid whose objects are decorated $d$-manifolds $(M,A,e_{1},e_{2})$ together with an orientation of $A \subset \Int{M}$ (called an \emph{orientedly} decorated $d$-manifold), and whose morphisms are maps $\varphi$ as in Definition~\ref{def:decorated-manifolds-morphisms} that are diffeomorphisms and such that the restriction $\varphi|_A \colon A \to A'$ is an orientation-preserving diffeomorphism.
\end{defi}

The morphisms of the groupoid $\cD_{d}$, i.e.~the invertible morphisms of $\cD\mathrm{ec}_d$, are called \emph{decorated diffeomorphisms}, for which we introduce the following notation.
\begin{notation}
\label{not:diffdec}
Let us write $\diffdec(M,N)$ for the space of decorated diffeomorphisms $M \to N$ and abbreviate $\diffdec(M) = \diffdec(M,M)$.
We similarly define $\diffdec^{+}(M) = \diffdec^{+}(M,M)$ for an \emph{orientedly} decorated $d$-manifold $M$, where the superscript ${}^{+}$ means that diffeomorphisms are orientation-preserving on the submanifold $A \subset M$.
\end{notation}

The space $\diffdec(M,N)$ is topologised as a subspace of $\cdec(M,N)$, which is topologised as in Definition~\ref{def:decorated-manifolds-morphism-spaces}.

\paragraph*{Semi-monoidal structure.}

The topologically-enriched categories $\cD\mathrm{ec}_d$, $\cD_{d}$ and $\cD_{d}^{+}$ are naturally endowed with semi-monoidal structures defined from the boundary connected sum, which we detail here. Let $(M,A,e_{1},e_{2})$ and $(M',A',e'_{1},e'_{2})$ be two decorated $d$-manifolds. Define
\[
M \natural M' = (M \sqcup M')/{\sim},
\]
where $\sim$ is the equivalence relation generated by $e_{2}(x,0) \sim e'_{1}(x,0)$ for all $(x,0) \in b\bB^d_{1}$. We give this a smooth structure as follows. There are obvious topological embeddings
\begin{equation}\label{eq:embedding_smooth_boundary_sum}
M \smallsetminus e_{2}(b\bB^d_{1}) \longhookrightarrow M \natural M' \qquad\text{and}\qquad M' \smallsetminus e'_{1}(b\bB^d_{1}) \longhookrightarrow M \natural M',
\end{equation}
and another topological embedding
\begin{equation}\label{eq:embedding_smooth_boundary_sum_{2}}
\bD^{d-1} \times [-1,1] \longhookrightarrow M \natural M'
\end{equation}
given by $(x,t) \mapsto e_{2}(x,-t)$ for $t \leq 0$ and $(x,t) \mapsto e'_{1}(x,t)$ for $t \geq 0$, where we have implicitly chosen representative boundary-cylinders for the boundary-cylinder germs $e_{2}$ and $e'_{1}$.

\begin{lemm}\label{lem:boundary_connected_sum}
Declaring that the topological embeddings \eqref{eq:embedding_smooth_boundary_sum} and \eqref{eq:embedding_smooth_boundary_sum_{2}} are both \textbf{smooth} embeddings defines a smooth structure on $M \natural M'$.
\end{lemm}

\begin{proof}
First, note that the smooth structures induced by the topological embeddings \eqref{eq:embedding_smooth_boundary_sum} and \eqref{eq:embedding_smooth_boundary_sum_{2}} are compatible on intersections, since boundary-cylinders are \emph{smooth} embeddings away from $\partial b\bB^d_{1}$.
Then the smooth structure of $M \natural M'$ is determined, except on $e_{2}(b\bB^d_{1}) = e'_{1}(b\bB^d_{1})$, by the smooth structures of $M$ and $M'$. The embedding of $\bD^{d-1} \times [-1,1]$ induced by (boundary-cylinders representing the boundary-cylinder germs) $e_{2}$ and $e'_{1}$ is therefore only required to extend this smooth structure to $e_{2}(b\bB^d_{1}) = e'_{1}(b\bB^d_{1})$. As a result, the smooth structure does not depend on the choice of representative boundary-cylinders, but only their germs.
Hence the smooth structure is well-defined.
\end{proof}

\begin{defi}[Boundary connected sum.]
\label{def:boundary-connected-sum}
We define the boundary connected sum of decorated $d$-manifolds by
\[
(M,A,e_{1},e_{2}) \,\natural\, (M',A',e'_{1},e'_{2}) = (M \natural M',A\sqcup A',e_{1},e'_{2}),
\]
where $M \natural M'$ is defined by Lemma~\ref{lem:boundary_connected_sum} and $A\sqcup A'$ is the disjoint union of $A$ and $A'$.
For orientedly decorated $d$-manifolds, the boundary connected sum is defined in the same way, with the orientation for $A \sqcup A'$ being induced from those of $A$ and $A'$. 
See Figure~\ref{fig:decorated-manifolds} for a schematic illustration.
\end{defi}

\begin{rmk}\label{rmk:boundary_connected_sum_comparison}
The usual definition of boundary connected sum of two smooth manifolds $M,M'$ depends on a choice of embedded disc in the boundary of each manifold, and a method of ``straightening corners'' after gluing these discs together. Up to diffeomorphism, the resulting smooth manifold $M \natural M'$ depends only on the choice of a boundary component of $M$ and of $M'$, and orientations of these if they are orientable (this is a result of Palais' \emph{Disc theorem} \cite[Th.~B and Cor.~1]{Palais1960Extendingdiffeomorphisms} and the existence of collar neighbourhoods \cite{Brown1962}). However, in order for $\natural$ to induce a well-defined monoidal structure on some category of manifolds with boundary (which we will do just below), it must be well-defined on the nose, not just up to diffeomorphism (since objects are manifolds, not diffeomorphism classes of manifolds). Here, these additional choices are \emph{built in} for decorated manifolds, and no additional choices are required in Definition~\ref{def:boundary-connected-sum} above.
\end{rmk}

\begin{rmk}
\label{rmk:cyclic-operad-structure}
The boundary connected sum of Definition \ref{def:boundary-connected-sum} is a fragment of a richer structure: a non-unital cyclic operad (see \cite{GetzlerKapranov1995}) in the category of small categories. This consists in a contravariant functor $P$ from the category of finite sets and bijections to the category of small categories, together with structure maps (functors)
\begin{equation}
\label{eq:cyclic-operad-structure-maps}
{}_x \circ_y \colon P(X) \times P(Y) \too P((X \smallsetminus \{x\}) \sqcup (Y \smallsetminus \{y\}))
\end{equation}
for each pair $(x,y) \in X \times Y$, satisfying certain associativity and equivariance axioms.
In our setting, the category $P(X)$ is defined just as $\cD\mathrm{ec}_d$, except that each manifold is equipped with $\lvert X \rvert$ boundary-cylinder germs, labelled by the elements of $X$, each bijection $f \colon X \to X'$ induces a functor $P(X) \to P(X')$ by applying $f$ to the labels and the functor \eqref{eq:cyclic-operad-structure-maps} acts by gluing together the boundary-cylinder germs labelled by $x$ and $y$.
We will only use the fragment of this structure given by $\cD\mathrm{ec}_d = P(\{1,2\})$ and \eqref{eq:cyclic-operad-structure-maps} for $X=Y=\{1,2\}$ and $(x=2,y=1)$, which is the operation described in Definition \ref{def:boundary-connected-sum}.
\end{rmk}

The boundary connected sum of Definition~\ref{def:boundary-connected-sum} induces a semi-monoidal structure on $\cD\mathrm{ec}_d$ as follows.
Let $\varphi \colon (L,A,e_{1},e_{2}) \to (L',A',e'_{1},e'_{2})$ and $\psi \colon (M,B,f_{1},f_{2}) \to (M',B',f'_{1},f'_{2})$ be morphisms in $\cD\mathrm{ec}_d$. By definition, these are smooth, proper maps $L \to L'$ and $M \to M'$ that take $A$ and $B$ into $A'$ and $B'$ respectively, and are compatible with the given boundary-cylinder germs. This compatibility implies that they glue to a well-defined, smooth, proper map $L \natural M \to L' \natural M'$ such that the image of $A\sqcup B$ is contained in $A'\sqcup B'$ and satisfying the boundary-cylinder germ conditions of Definition~\ref{def:decorated-manifolds-morphisms}: this is thus a morphism
\[
(L,A,e_{1},e_{2}) \,\natural\, (M,B,f_{1},f_{2}) \too (L',A',e'_{1},e'_{2}) \,\natural\, (M',B',f'_{1},f'_{2})
\]
of $\cD\mathrm{ec}_d$ that we denote by $\varphi\natural \psi$. This extends to morphisms of $\cD_{d}^{+}$ since the restricted diffeomorphisms $A\to A'$ and $B\to B'$ are orientation-preserving by assumption and $(\varphi\natural \psi)|_{A \sqcup B} = \varphi|_A \sqcup \psi|_B$.

That the operation $\natural$ preserves composition of morphisms is a straightforward consequence of its assignments. Also, it immediately follows from the definitions that the boundary connected sum of Definition~\ref{def:boundary-connected-sum} is associative, and that the isomorphisms defined from this associativity are natural with respect to the morphisms of $\cD\mathrm{ec}_d$. Hence, we equip $\natural$ with the associator defined by these natural isomorphisms, for which one easily checks the pentagon diagram. Hence we prove:

\begin{prop}[Semi-monoidal structure.]
\label{prop:Dd-semi-monoidal}
There is a semifunctor
\[
\natural \colon \cD\mathrm{ec}_d \times \cD\mathrm{ec}_d \too \cD\mathrm{ec}_d
\]
defined by assigning the boundary connected sum of Definition~\ref{def:boundary-connected-sum} on objects and $\varphi\natural \psi$ for two morphisms $\varphi$ and $\psi$ in $\cD\mathrm{ec}_d$. Moreover, the semifunctor $\natural$ induces a semi-monoidal structure on the topologically-enriched categories $\cD\mathrm{ec}_d$, $\cD_{d}$ and $\cD_{d}^{+}$.
\end{prop}

\begin{rmk}[Monoidal and semi-monoidal structures.]
\label{r:semi-monoidal}
The semi-monoidal structure on $\cD\mathrm{ec}_d$ defined above does not have a unit, since there is no \emph{natural} way of identifying $M \natural \bB^d_{1}$ with $M$ for all $M$ (although they are clearly non-naturally diffeomorphic). If decorated manifolds had been defined to be equipped with boundary cylinders (not just \emph{germs} of boundary-cylinders), then $\cD\mathrm{ec}_d$ would have had an obvious 
\emph{genuine} monoidal structure. However, the proof of Theorem~\ref{thm:fibre-bundle} below, which tells us that the Serre fibration hypothesis of Proposition~\ref{prop:topological-Quillen_morphisms} is satisfied for subgroupoids of $\cD_{d}$ (see Lemma~\ref{lem:Serre-fibration-condition}), depends crucially on the fact that morphisms of decorated manifolds are only required to preserve \emph{germs} of boundary-cylinders, rather than entire boundary-cylinders. Thus we are forced either to formally adjoin a unit to $\cD\mathrm{ec}_d$ (which is unnatural since $\pi_{0}(\cD\mathrm{ec}_d)$ already has a unit by Lemma~\ref{lem:monoidal-pi0}), or to work directly with semi-monoidal categories, which is the choice that we make here.
\end{rmk}

\paragraph*{Groupoids of path-components.}

Applying the functor $\pi_{0}$, we may consider the discrete groupoids $\pi_{0}(\cD_{d})$ and $\pi_{0}(\cD_{d}^{+})$. The following lemma may be checked straightforwardly from the definitions.

\begin{lemm}\label{lem:monoidal-pi0}
The discrete groupoids $\pi_{0}(\cD_{d})$ and $\pi_{0}(\cD_{d}^{+})$ inherit semi-monoidal structures from those on $\cD_{d}$ and $\cD_{d}^{+}$. Moreover, these semi-monoidal structures are monoidal, with unit object given in each case by the solid cylinder $(\bB^d_{1},\emptyset,\id,r)$, where $r \colon \bB^d_{1} \to \bB^d_{1}$ is the reflection $(x,t) \mapsto (x,1-t)$.

More generally, if $\cG \subseteq \cD_{d}$ is any subgroupoid closed under the semi-monoidal structure and containing the solid cylinder $(\bB^d_{1},\emptyset,\id,r)$, then the semi-monoidal structure inherited by $\pi_{0}(\cG)$ is monoidal.
Furthermore, if $\cM \subseteq \cD_{d}$ is a subgroupoid closed under the operation $g \natural -$ for each object $g$ of $\cG$ (thus $\cM$ is a left module over the semi-monoidal groupoid $\cG$), then $\pi_{0}(\cM)$ inherits a structure of a (genuine) left module over the monoidal groupoid $\pi_{0}(\cG)$.
Similar results hold for analagous subgroupoids $\cG \subseteq \cD_{d}^{+}$ and $\cM \subseteq \cD_{d}^{+}$.
\end{lemm}

\subsection{Split short exact sequences}
\label{ss:split-ses}

In this subsection, we establish several split homotopy fibration sequences for embedding spaces and decorated diffeomorphism groups of manifolds, whose associated split short exact sequences on $\pi_{1}$ and $\pi_{0}$ will be used in our constructions in \S\ref{s:general_construction}. In particular, the split short exact sequences \eqref{eq:split-ses-1} and \eqref{eq:split-ses-2} below will be used in the two versions of our general construction of homological representations in \S\ref{ss:homological_representation_functor_motion_groups} and \S\ref{ss:homological_representation_functor_mcg} respectively. In addition, the split short exact sequence \eqref{eq:split-ses-2} implies that any motion group, under a certain peripherality condition, is a \emph{braided mapping class group} (see Proposition~\ref{prop:braided-diff-groups}); in particular, a normal subgroup of a mapping class group.

\emph{Throughout \S\ref{ss:split-ses}, we fix a closed submanifold $Z \subset \bR^d$ and an open subgroup $\sG \leq \Diff(Z)$.}
Note that open subgroups $\sG \leq \Diff(Z)$ correspond to subgroups of $\pi_{0}(\Diff(Z))$. We recall our notational convention that, for $X$ a manifold with boundary, $\Int{X}$ denotes its interior.

\subsubsection{The first short exact sequence}\label{sss:split-ses-1}

The first split short exact sequence \eqref{eq:split-ses-1} that we construct deals with fundamental groups of certain spaces of embeddings, which we now define.

\begin{defi}
\label{def:cE}
For smooth manifolds $X$ and $Y$, let us write $\cE(X,Y) = \Emb(X,Y) / \Diff(X)$. For a subgroup $H \leq \Diff(X)$, we also write $\cE_{H}(X,Y) = \Emb(X,Y)/H$. In the special case when $H = \Diff^{+}(X)$ (when $X$ is orientable), we also write $\cE^{+}(X,Y) = \Emb(X,Y)/\Diff^{+}(X)$.
\end{defi}

We specify once and for all an identification of $\bR^d$ with the interior of the solid cylinder $\bB^d_{1}$. The choice of closed submanifold $Z \subset \bR^d$ above therefore determines a decorated manifold $(\bB^d_{1},Z) = (\bB^d_{1},Z,[\id],[r])$ whose boundary-cylinder germs are represented, respectively, by the identity and by the reflection $r$ of the solid cylinder $\bB^d_{1} = \bD^{d-1} \times [0,1]$ in its second coordinate.
To justify that this is a decorated manifold, we need to verify that there are boundary-cylinders $e_1$ and $e_2$, having the same germs as $\id$ and $r$, whose images are disjoint from $Z$. We will explain this for $e_1$, the case of $e_2$ being similar. Since $Z$ is a compact subset of the interior of $\bD^{d-1} \times [0,1]$, it is disjoint from $\bD^{d-1} \times [0,\epsilon]$ for some $\epsilon > 0$. Choose a diffeomorphism $[0,1] \cong [0,\epsilon]$ that agrees with the identity on $[0,\delta]$ for some $0 < \delta < \epsilon$. The product of this with the identity of $\bD^{d-1}$ is then a boundary-cylinder $e_1$ for $\bB^d_{1}$ whose image is disjoint from $Z$ and whose germ is the same as the identity.

For any other decorated manifold $(M,A) \in \cD_{d}$, we may therefore consider the boundary connected sum $(M,A) \natural (\bB^d_{1},Z)$.
If we fix an orientation of $Z$, we may do the same for any orientedly decorated manifold $(M,A) \in \cD_{d}^{+}$.

\begin{notation}
\label{not:union-with-Z}
To shorten notation, we write $\mbar = M \natural \bB^d_{1}$ and we denote the interior of $\mbar$ by $\Breve{M}$. Notice that $\Breve{M}$ contains (disjointly) both $\Int{M}$ and $\bR^d$ under its fixed identification with the interior of $\bB^d_{1}$.
\end{notation}

\begin{construction}
\label{construction:Theta}
It will be useful at several points (see for instance the proofs of Propositions \ref{prop:split-fibration-sequence-1} and \ref{prop:split-fibration-sequence-2}) to have a specified self-embedding $\Theta$ of $\mbar$ whose image is $M \subset \mbar$, together with a path of embeddings from $\Theta$ to the identity. Let us first choose a representative for the boundary-cylinder germ $e_2$ associated to $M$; this is an embedding $\bB^d_1 = \bD^{d-1} \times [0,1] \hookrightarrow M$, which we also denote by $e_2$. Let us moreover assume that its image $e_2(\bB^d_1)$ is disjoint from the submanifold $A \subset \mathring{M}$. The decorated manifold $\mbar$ may then be viewed as the union of $M$ and $\bD^{d-1} \times [-1,0]$ along $\bD^{d-1} \times \{0\}$. Using this viewpoint, we construct a path of self-embeddings
\[
\Theta_t \colon \mbar \longhookrightarrow \mbar
\]
such that $\Theta_0 = \id$ and $\Theta_1(\mbar) = M$. To do this, we first specify a path of embeddings of intervals $\theta_t \colon [-1,1] \hookrightarrow [-1,1]$ where $\theta_0 = \id$, each $\theta_t$ agrees with the identity near $1$ and $\theta_t([-1,1]) = [t-1,1]$. Taking the product of $\theta_t$ with the identity on $\bD^{d-1}$ and extending by the identity on the complement $\mbar \smallsetminus (\bD^{d-1} \times [-1,1])$, we obtain a self-embedding $\mbar \hookrightarrow \mbar$, which we define to be $\Theta_t$. In particular, $\Theta = \Theta_1$ is the advertised self-embedding of $\mbar$.
\end{construction}

\begin{rmk}
We emphasise that, by construction, each self-embedding $\Theta_t$ of Construction~\ref{construction:Theta} agrees with the identity on $M \smallsetminus e_2(\bB^d_1)$; in particular it agrees with the identity on $A \subset \mathring{M}$.
\end{rmk}

\begin{prop}
\label{prop:split-fibration-sequence-1}
For any decorated manifold $(M,A) \in \cD_{d}$ and an open subgroup $\sG \leq \Diff(Z)$, there is a homotopy fibration sequence
\begin{equation}
\label{eq:split-fibration-sequence-1}
\centering
\begin{split}
\begin{tikzpicture}
[x=1mm,y=1mm]
\node (l) at (0,0) {$\cE_{\sG}(Z,\Breve{M}\smallsetminus A)$};
\node (m) [anchor=west] at ($ (l.east) + (10,0) $) {$\cE_{\Diff(A) \times \sG}(A \sqcup Z,\Breve{M})$};
\node (r) [anchor=west] at ($ (m.east) + (10,0) $) {$\cE(A,\Breve{M}),$};
\draw[->] (l) to (m);
\draw[->] (m) to (r);
\draw[->,densely dashed] (r.north west) to[out=160,in=20] (m.north east);
\end{tikzpicture}
\end{split}
\end{equation}
in which the second map admits a section up to homotopy, as pictured. The same holds for objects of $\cD_{d}^{+}$, with $\Diff(A)$ replaced by $\Diff^{+}(A)$ in the middle term of \eqref{eq:split-fibration-sequence-1}.
\end{prop}
\begin{proof}
The map $\cE_{\Diff(A) \times \sG}(A \sqcup Z,\Breve{M}) \to \cE(A,\Breve{M})$ that forgets the embedding of $Z$ is equivariant with respect to the left action of the topological group $\Diff_c(\Breve{M})$ of compactly-supported diffeomorphisms of $\Breve{M}$. By \cite[Prop.~4.15]{Palmer2018HomologicalstabilitymoduliI}, the action of $\Diff_c(\Breve{M})$ on $\cE(A,\Breve{M})$ is locally retractile, i.e.~it admits local sections. Thus, by \cite[Th.~A]{Palais1960Localtrivialityof}, the map
\begin{equation}
\label{eq:restriction-to-A}
\cE_{\Diff(A) \times \sG}(A \sqcup Z,\Breve{M}) \too \cE(A,\Breve{M})
\end{equation}
is a fibre bundle, in particular a Serre fibration. Write $\mathrm{incl.}$ for the inclusion of $A$ into $\Breve{M}$ and $[\mathrm{incl.}]$ for its $\Diff(A)$-orbit; this is a natural basepoint for $\cE(A,\Breve{M})$. The point-set fibre of \eqref{eq:restriction-to-A} over $[\mathrm{incl.}] \in \cE(A,\Breve{M})$ is clearly equal to $\cE_{\sG}(Z,\Breve{M}\smallsetminus A)$, so \eqref{eq:split-fibration-sequence-1} is a homotopy fibration sequence.

A section up to homotopy for \eqref{eq:restriction-to-A} is defined by sending $[\varphi] \in \cE(A,\Breve{M})$ to $[(\Theta \circ \varphi) \sqcup \iota]$, where $\iota$ is the inclusion $Z \subset \bR^d \subset \Breve{M}$ and $\Theta$ is the self-embedding of Construction~\ref{construction:Theta}.
\end{proof}

\begin{coro}\label{coro:split-ses-1}
For any decorated manifold $(M,A) \in \cD_{d}$ there is a split short exact sequence
\begin{equation}
\label{eq:split-ses-1}
\centering
\begin{split}
\begin{tikzpicture}
[x=1mm,y=1mm]
\node (ll) at (0,0) {$1$};
\node (l) [anchor=west] at ($ (ll.east) + (5,0) $) {$\pi_{1}(\cE_{\sG}(Z,\Breve{M} \smallsetminus A))$};
\node (m) [anchor=west] at ($ (l.east) + (10,0) $) {$\pi_{1}(\cE_{\Diff(A) \times \sG}(A \sqcup Z,\Breve{M}))$};
\node (r) [anchor=west] at ($ (m.east) + (10,0) $) {$\pi_{1}(\cE(A,\Breve{M}))$};
\node (rr) [anchor=west] at ($ (r.east) + (5,0) $) {$1,$};
\draw[->] (ll) to (l);
\draw[->] (l) to (m);
\draw[->] (m) to (r);
\draw[->] (r) to (rr);
\draw[->,densely dashed] (r.north west) to[out=160,in=20] (m.north east);
\end{tikzpicture}
\end{split}
\end{equation}
where the basepoint of each space of embeddings modulo diffeomorphisms is given by the inclusion. The same holds for objects of $\cD_{d}^{+}$, with $\Diff(A)$ replaced by $\Diff^{+}(A)$ in the middle term of \eqref{eq:split-ses-1}.
\end{coro}

We will sometimes consider spaces of embeddings of $Z$ into $\mbar$, rather than its interior $\Breve{M}$. We record here the following fact, for later reference.

\begin{prop}\label{prop:inclusion-heq}
The inclusion $\cE_{\sG}(Z,\Breve{M} \smallsetminus A) \subset \cE_{\sG}(Z,\mbar \smallsetminus A)$ is a homotopy equivalence.
\end{prop}
\begin{proof}
A deformation retraction may be defined by post-composing embeddings with an isotopy of self-embeddings of $\mbar$ starting from the identity and ``shrinking'' a collar neighbourhood of its boundary.
\end{proof}

\subsubsection{The second short exact sequence}\label{sss:split-ses-2}

We now construct a second split short exact sequence \eqref{eq:split-ses-2}. We first need some notation. For a decorated manifold $(M,A) = (M,A,e_{1},e_{2}) \in \cD_{d}$, recall from Notation~\ref{not:diffdec} that we write $\diffdec(M,A)$ for its automorphism group in $\cD_{d}$, in other words, the self-diffeomorphisms of $M$ that send $A$ onto itself and that are compatible with the boundary-cylinder germs $e_{1}$ and $e_{2}$. Similarly, we write $\diffdec^{+}(M,A)$ for the automorphism group of an orientedly decorated manifold $(M,A) \in \cD_{d}^{+}$, the additional condition being that the restriction to $A$ must be orientation-preserving.

\begin{defi}\label{def:diff_dec_partitioned}
For a decorated manifold $(M,A) \in \cD_{d}$ such that $A$ decomposes as $A_{1} \sqcup\ldots \sqcup A_{k}$, we define $\diffdec(M,A_{1},\ldots,A_{k})$ to be the subgroup of $\diffdec(M,A)$ of those diffeomorphisms that preserve this decomposition, in other words, that send each $A_{i}$ onto itself for each $1\leq i\leq k$. In the setting of $\cD_{d}^{+}$, we analogously define the subgroup $\diffdec^{+}(M,A_{1},\ldots,A_{k}) \subset \diffdec^{+}(M,A)$.
\end{defi}

In the following, we continue to use the Notation~\ref{not:union-with-Z} from \S\ref{sss:split-ses-1}.

\begin{prop}\label{prop:split-fibration-sequence-2}
For any decorated manifold $(M,A) \in \cD_{d}$ and an open subgroup $\sG \leq \Diff(Z)$, there is a homotopy fibration sequence
\begin{equation}
\label{eq:split-fibration-sequence-2}
\centering
\begin{split}
\begin{tikzpicture}
[x=1mm,y=1mm]
\node (l) at (0,0) {$\diffdec(\mbar,A,Z|\sG)$};
\node (m) [anchor=west] at ($ (l.east) + (10,0) $) {$\diffdec(\mbar,A)$};
\node (r) [anchor=west] at ($ (m.east) + (10,0) $) {$\cE_{\sG}(Z,\Breve{M} \smallsetminus A),$};
\draw[->] (l) to (m);
\draw[->] (m) to (r);
\draw[->,densely dashed] (m.north west) to[out=160,in=20] (l.north east);
\end{tikzpicture}
\end{split}
\end{equation}
in which the first map admits a section up to homotopy, as pictured, and where $\diffdec(\mbar,A,Z|\sG)$ denotes the subgroup of $\varphi \in \diffdec(\mbar,A,Z)$ such that $\varphi|_Z \in \sG \leq \Diff(Z)$.

The analogous statement also holds for objects of $\cD_{d}^{+}$, with $\diffdec$ replaced by $\diffdec^{+}$.
\end{prop}

\begin{proof}
The right-hand map above is equivariant with respect to the left action of $\Diff_c(\Breve{M} \smallsetminus A)$. By \cite[Prop.~4.15]{Palmer2018HomologicalstabilitymoduliI}, since $\sG \leq \Diff(Z)$ is an open subgroup, its action on $\cE_{\sG}(Z,\Breve{M} \smallsetminus A)$ is locally retractile, so by \cite[Th.~A]{Palais1960Localtrivialityof}, the right-hand map above is a fibre bundle. Its point-set fibre over the basepoint $[\mathrm{incl.}]$ is clearly equal to $\diffdec(\mbar,A,Z|\sG)$.

We now define a section up to homotopy for the inclusion $\diffdec(\mbar,A,Z|\sG) \subset \diffdec(\mbar,A)$, using the path of self-embeddings $\Theta_t$ for $t \in [0,1]$ of Construction~\ref{construction:Theta}. Let us define a homotopy $H \colon \diffdec(\mbar,A) \times [0,1] \to \diffdec(\mbar,A)$ by $H(\varphi,t)(x) = \Theta_t(\varphi(\Theta_t^{-1}(x)))$ if $x \in \Theta_t(\mbar)$ and by $H(\varphi,t)(x) = x$ if $x \notin \Theta_t(\mbar)$. To see that this defines a section up to homotopy for the inclusion, we must check (i) that $H(\varphi,t)$ is indeed a diffeomorphism of $\mbar$, (ii) that $H$ is continuous, (iii) that $H(-,t)$ is a group homomorphism for each $t \in [0,1]$, (iv) that $H(-,0)$ is the identity and (v) that $H(\varphi,1) \in \diffdec(\mbar,A,Z|\sG)$ for each $\varphi \in \diffdec(\mbar,A)$.

By construction, the conjugate $\Theta_t \circ \varphi \circ \Theta_t^{-1}$ is a diffeomorphism of $\Theta_t(\mbar)$. Since $\varphi$ is a decorated diffeomorphism of $\mbar$, it agrees with the identity on $\bD^{d-1} \times [-1,-1+\delta]$ for some $\delta > 0$ and so $\Theta_t \circ \varphi \circ \Theta_t^{-1}$ agrees with the identity on $\bD^{d-1} \times [t-1,t-1+\delta']$ for some $\delta' > 0$. The result of extending this by the identity on $\mbar \smallsetminus \Theta_t(\mbar) = \bD^{d-1} \times [-1,t-1)$ is therefore also a diffeomorphism. This is precisely $H(\varphi,t)$, so we have verified point (i).

Composition of proper maps (such as diffeomorphisms or closed embeddings like $\Theta_t$) is continuous in the Whitney topology (see \cite[\S 2, Prop.~1]{Mather1969}), as is extending diffeomorphisms by the identity. Since the topology on decorated diffeomorphism groups is derived from the Whitney topology (see Definitions~\ref{def:decorated-manifolds-morphism-spaces} and \ref{def:Dd} for precise details), it follows that $H(-,t)$ is continuous for each $t \in [0,1]$. Continuity in $t$ follows from the fact that $t \mapsto \Theta_t$ was constructed to be a continuous path of self-embeddings; see Construction~\ref{construction:Theta}. This verifies point (ii).

Points (iii) and (iv) are obvious: conjugating by a fixed self-embedding and extending by the identity is clearly a group homomorphism, and $H(-,0)$ is the identity since $\Theta_0$ is the identity. Finally, $H(\varphi,1)$ agrees with the identity on $\bB^d_1 \subset \mbar$, by construction, so in particular it acts by the identity on $Z \subset \bB^d_1$; hence it lies in $\diffdec(\mbar,A,Z|\sG)$. This verifies point (v).

Thus $H(-,1) \colon \diffdec(\mbar,A) \to \diffdec(\mbar,A,Z|\sG)$ is a section up to homotopy for the inclusion, witnessed by the homotopy $H$.
\end{proof}

\begin{coro}\label{coro:split-ses-2}
For any decorated manifold $(M,A) \in \cD_{d}$ there is a split short exact sequence
\begin{equation}
\label{eq:split-ses-2}
\centering
\begin{split}
\begin{tikzpicture}
[x=1mm,y=1mm]
\node (ll) at (0,0) {$1$};
\node (l) [anchor=west] at ($ (ll.east) + (5,0) $) {$\pi_{1}(\cE_{\sG}(Z,\Breve{M} \smallsetminus A))$};
\node (m) [anchor=west] at ($ (l.east) + (10,0) $) {$\pi_{0}(\diffdec(\mbar,A,Z|\sG))$};
\node (r) [anchor=west] at ($ (m.east) + (10,0) $) {$\pi_{0}(\diffdec(\mbar,A))$};
\node (rr) [anchor=west] at ($ (r.east) + (5,0) $) {$1,$};
\draw[->] (ll) to (l);
\draw[->] (l) to (m);
\draw[->] (m) to (r);
\draw[->] (r) to (rr);
\draw[->,densely dashed] (r.north west) to[out=160,in=20] (m.north east);
\end{tikzpicture}
\end{split}
\end{equation}
where the basepoint of the space $\cE_{\sG}(Z,\Breve{M} \smallsetminus A)$ is the orbit of the inclusion. The analogous statement also holds for objects of $\cD_{d}^{+}$, with $\diffdec$ replaced by $\diffdec^{+}$.
\end{coro}
\begin{proof}
In contrast to Corollary~\ref{coro:split-ses-1}, where it was immediate that \eqref{eq:split-fibration-sequence-1} induces \eqref{eq:split-ses-1} by considering the long exact sequence of homotopy groups, this corollary requires a little more explanation, since it involves also $\pi_{0}$, while we aim to prove that we obtain a short exact sequence of \emph{groups} and \emph{homomorphisms}: more precisely, it is \emph{a priori} unclear that the connecting homomorphism $\pi_{1}(\cE_{\sG}(Z,\Breve{M} \smallsetminus A))\to \pi_{0}(\diffdec(\mbar,A,Z|\sG))$ of the homotopy long exact sequence is a group homomorphism. Let us write $\Gamma = \diffdec(\mbar,A)$, $H = \diffdec(\mbar,A,Z|\sG)$ and $B = \cE_{\sG}(Z,\Breve{M} \smallsetminus A)$. Since we are not interested in $\pi_{0}(B)$, we may implicitly replace $B$ with the image of the fibre bundle $\pi \colon \Gamma \to B$. This has the effect of possibly throwing away some path-components of $B$ (not containing its basepoint) and ensuring that $\pi$ is a surjective fibre bundle. It is also $\Gamma$-equivariant, where in each case $\Gamma$ acts on the left by post-composition. Lemma~\ref{lem:equivariant-fibre-bundles} below thus implies that $\pi \colon \Gamma \to B$ is isomorphic (in the category of maps out of $\Gamma$) to the quotient map $\Gamma \to \Gamma/H$, which is therefore a fibre bundle. We shall prove in a moment that it is in fact a principal $H$-bundle. This implies that the homotopy fibration sequence $H \hookrightarrow \Gamma \twoheadrightarrow \Gamma/H$ continues two steps to the right with $\Gamma/H \to \mathscr{B}H \to \mathscr{B}\Gamma$ (where $\mathscr{B}H$ and $\mathscr{B}\Gamma$ denote the classifying spaces of $H$ and $\Gamma$ respectively). Moreover, the map $\mathscr{B}H \to \mathscr{B}\Gamma$ in this sequence, which is induced by the inclusion, admits a section up to homotopy, induced by the section up to homotopy in \eqref{eq:split-fibration-sequence-2}. Considering the long exact sequence of homotopy groups associated to the split fibration sequence $\Gamma/H \to \mathscr{B}H \to \mathscr{B}\Gamma$, we obtain the split short exact sequence \eqref{eq:split-ses-2}, as claimed.

It remains to prove our assertion that $\Gamma \to \Gamma/H$ is a principal $H$-bundle. For this, we will use a refinement of \cite[Th.~A]{Palais1960Localtrivialityof}, which is remarked on page 307 of \cite{Palais1960Localtrivialityof} and written down explicitly in \cite[Prop.~4.8]{Palmer2018HomologicalstabilitymoduliI}. Namely, let $X$ be a space with a left $\Gamma_{0}$-action and a right $H$-action that commute, and such that the induced left $\Gamma_{0}$-action on $X/H$ is locally retractile. Assume that the $H$-action on $X$ is free and that the action map $x.- \colon H \to X$ is a topological embedding for each $x \in H$. Then $X \to X/H$ is a principal $H$-bundle. In our setting we have $X = \Gamma = \diffdec(\mbar,A)$ and we take $\Gamma_{0} = \Diff_c(\Breve{M} \smallsetminus A)$ acting on $\Gamma$ by extending compactly-supported diffeomorphisms of $\Breve{M} \smallsetminus A$ to $\mbar$ by the identity and post-composing. This obviously commutes with the right action of $H$ by pre-composition. The induced left action of $\Gamma_{0}$ on $\Gamma/H$ is the same as its left action on $\cE_{\sG}(Z,\Breve{M} \smallsetminus A)$ by post-composition under the homeomorphism $\Gamma/H \cong B = \cE_{\sG}(Z,\Breve{M} \smallsetminus A)$, which is locally retractile by \cite[Prop.~4.15]{Palmer2018HomologicalstabilitymoduliI}, as mentioned already in the proof of Proposition~\ref{prop:split-fibration-sequence-2}. The right $H$-action on $\Gamma$ is simply right multiplication ($H$ is a subgroup of $\Gamma$), so it is obviously free. It is also clear that the action map $g.- \colon H \to \Gamma$ is a topological embedding for each $g \in \Gamma$, since it is just the inclusion $H \subset \Gamma$ (which is a topological embedding) followed by the homeomorphism $\Gamma \cong \Gamma$ given by left multiplication by $g$. It therefore follows from \cite[Prop.~4.8]{Palmer2018HomologicalstabilitymoduliI} that $\Gamma \to \Gamma/H$ is a principal $H$-bundle, as claimed.
\end{proof}

\begin{lemm}
\label{lem:equivariant-fibre-bundles}
Let $\Gamma$ be a topological group acting on the left on a space $B$ and let $\pi \colon \Gamma \to B$ be a surjective, $\Gamma$-equivariant fibre bundle. Denote by $H$ the stabiliser of the basepoint $\pi(1_\Gamma) \in B$. Then there is a homeomorphism $B \cong \Gamma/H$ commuting with $\pi$ and the projection $\Gamma \to \Gamma/H$ that quotients by the right action of $H$ on $\Gamma$ by multiplication.
\end{lemm}
\begin{proof}
Surjective fibre bundles are quotient maps, so in order to show that the two quotient maps $\Gamma \to \Gamma/H$ and $\pi \colon \Gamma \to B$ are isomorphic, it suffices to check that their fibres give the same partition of the underlying set of $\Gamma$. The fibres of $\Gamma \to \Gamma/H$ are by definition the right $H$-cosets in $\Gamma$. Let $b \in B$ and choose $g \in \Gamma$ so that $\pi(g) = b$, which is possible since $\pi$ is surjective. The fibre $\pi^{-1}(b)$ consists of all $g' \in \Gamma$ such that $g' \pi(1_\Gamma) = \pi(g') = b = \pi(g) = g \pi(1_\Gamma)$, where we have used $\Gamma$-equivariance of $\pi$. This condition is equivalent to $g^{-1}g' \in \mathrm{Stab}(\pi(1_\Gamma)) = H$, in other words $g' \in gH$. Hence the fibres of $\pi$ are also the right $H$-cosets in $\Gamma$.
\end{proof}

\subsubsection{Connecting the short exact sequences}\label{sec:map_of_ses}

We now describe the relationship between the split short exact sequences \eqref{eq:split-ses-1} and \eqref{eq:split-ses-2}.
\begin{prop}
\label{prop:map_of_ses}
Let $(M,A) \in \cD_{d}$ be a decorated manifold, let $Z \subset \bR^d$ be a closed submanifold and let $\sG \leq \Diff(Z)$ be an open subgroup. Then there is a commutative diagram of the form:
\begin{equation}
\label{eq:3x3diagram}
\centering
\begin{split}
\begin{tikzpicture}
[x=1mm,y=1.1mm]
\node (tl) at (0,24) {$\pi_{1}(\cE_{\sG}(Z,\Breve{M}\smallsetminus A))$};
\node (tm) at (40,24) {$\pi_{1}(\cE_{\Diff(A) \times \sG}(A \sqcup Z,\Breve{M}))$};
\node (tr) at (80,24) {$\pi_{1}(\cE(A,\Breve{M}))$};
\node (ml) at (0,12) {$\pi_{1}(\cE_{\sG}(Z,\Breve{M}\smallsetminus A))$};
\node (mm) at (40,12) {$\pi_{0}(\diffdec(\mbar,A,Z|\sG))$};
\node (mr) at (80,12) {$\pi_{0}(\diffdec(\mbar,A))$};
\node (bl) at (0,0) {$1$};
\node (bm) at (40,0) {$\pi_{0}(\diffdec(\mbar,\emptyset))$};
\node (br) at (80,0) {$\pi_{0}(\diffdec(\mbar,\emptyset))$};
\node (tll) at (-25,24) {$1$};
\node (trr) at (105,24) {$1$};
\node (mll) at (-25,12) {$1$};
\node (mrr) at (105,12) {$1$};
\node (bll) at (-25,0) {$1$};
\node (brr) at (105,0) {$1$};
\draw[->] (tll) to (tl);
\draw[->] (tl) to (tm);
\draw[->] (tm) to (tr);
\draw[->] (tr) to (trr);
\draw[->] (mll) to (ml);
\draw[->] (ml) to (mm);
\draw[->] (mm) to (mr);
\draw[->] (mr) to (mrr);
\draw[->] (bll) to (bl);
\draw[->] (bl) to (bm);
\draw[->] (bm) to node[above,font=\small]{$\id$} (br);
\draw[->] (br) to (brr);
\draw[->] (tl) to node[left,font=\small]{$\id$} (ml);
\draw[->] (ml) to (bl);
\draw[->] (tm) to (mm);
\draw[->] (mm) to (bm);
\draw[->] (tr) to (mr);
\draw[->] (mr) to (br);
\draw[->,densely dashed] (tr.north west) to[out=160,in=20] (tm.north east);
\draw[->,densely dashed] (mr.north west) to[out=160,in=20] (mm.north east);
\end{tikzpicture}
\end{split}
\end{equation}
where the top and middle rows are \eqref{eq:split-ses-1} and \eqref{eq:split-ses-2} respectively and the columns are also exact. The top-right square also commutes when the solid horizontal arrows are replaced by the dotted arrows.
\end{prop}
\begin{proof}
Let us consider the following diagram:
\begin{equation}
\label{eq:3x3diagram-proof}
\centering
\begin{split}
\begin{tikzpicture}
[x=1mm,y=1mm]
\node (tl) at (0,24) {$\diffdec(\mbar,A,Z|\sG)$};
\node (tm) at (35,24) {$\diffdec(\mbar,A)$};
\node (tr) at (70,24) {$\cE_{\sG}(Z,\Breve{M}\smallsetminus A)$};
\node (ml) at (0,12) {$\diffdec(\mbar,\emptyset)$};
\node (mm) at (35,12) {$\diffdec(\mbar,\emptyset)$};
\node (bll) at (-40,0) {$\cE_{\sG}(Z,\Breve{M}\smallsetminus A)$};
\node (bl) at (0,0) {$\cE_{\Diff(A) \times \sG}(A \sqcup Z,\Breve{M})$};
\node (bm) at (35,0) {$\cE(A,\Breve{M}).$};
\draw[->>] (ml) to (bl);
\draw[->>] (mm) to (bm);
\draw[->>] (tm) to (tr);
\draw[->>] (bl) to (bm);
\incl{(tl)}{(tm)}
\incl{(tl)}{(ml)}
\incl{(tm)}{(mm)}
\incl{(bll)}{(bl)}
\draw[->] (ml) to node[above,font=\small]{$\id$} (mm);
\draw[->,densely dashed] (bm.north west) to[out=160,in=20] (bl.north east);
\draw[->,densely dashed] (tm.north west) to[out=160,in=20] (tl.north east);
\end{tikzpicture}
\end{split}
\end{equation}
The top row is \eqref{eq:split-fibration-sequence-2} and the bottom row is \eqref{eq:split-fibration-sequence-1}. The bottom vertical maps are the restriction maps to $A \sqcup Z$ and to $A$ respectively. They are fibre bundles by the same argument as in the proof of Proposition~\ref{prop:split-fibration-sequence-2}, and the top two vertical maps are the corresponding inclusions of point-set fibres. Both squares evidently commute, since restricting a diffeomorphism to $A \sqcup Z$ and then to $A$ gives the same result as restricting directly to $A$.
Taking the associated long exact sequences of homotopy groups of these homotopy fibration sequences horizontally and vertically, we obtain the diagram \eqref{eq:3x3diagram}, where for the two columns we apply the same argument as in the proof of Corollary~\ref{coro:split-ses-2} (except without the section up to homotopy). One may check directly that the top-right square of \eqref{eq:3x3diagram} commutes when the solid horizontal arrows are replaced by the dotted arrows by considering the diagram \eqref{eq:3x3diagram-proof} and using the definition of the connecting homomorphisms in the vertical long exact sequences.
\end{proof}

\begin{rmk}
We shall prove in \S\ref{sss:functorial-split-ses} below that the diagram \eqref{eq:3x3diagram-proof} is functorial with respect to morphisms of the category $\fU\cD_{d}$ defined in \S\ref{ss:Quillen_bracket_construction} below; in particular it is functorial with respect to decorated diffeomorphisms of $(M,A)$, i.e.~with respect to the morphisms of the groupoid $\cD_{d}$.
\end{rmk}

\subsubsection{Motion groups and braided mapping class groups}\label{sec:braided-mcg}

The goal of this section is to introduce \emph{braided diffeomorphism groups} (see Definition~\ref{def:braided-diffeomorphisms}) and \emph{motion groups} (see Definition~\ref{def:motion_groups}), and to show (see  Proposition~\ref{prop:braided-diff-groups}) that the groups of path-components of braided diffeomorphism groups (\emph{braided mapping class groups}) are naturally identified with the corresponding motion groups when $A \subset M$ is \emph{peripheral} in the sense of Definition~\ref{def:peripheral} below.

\begin{defi}\label{def:braided-diffeomorphisms}
The \emph{braided diffeomorphism group} $\diffdecbr(M,A)$ of a decorated manifold $(M,A)$ is the kernel of the natural homomorphism
\[
\diffdec(M,A) \too \pi_{0}(\diffdec(M,\emptyset)),
\]
in other words, it is the subgroup of diffeomorphisms of $(M,A)$ that become isotopic to the identity after forgetting $A$. We also define $\diffdecbrplus(M,A) = \diffdecbr(M,A) \cap \diffdec^{+}(M,A)$.
\end{defi}

\begin{defi}\label{def:motion_groups}
Given a closed submanifold $Y \subset \Int{M}$, the corresponding \emph{motion group} $\Mot_Y(M)$ is the fundamental group $\pi_{1}(\cE(Y,\Int{M}))$. When $Y$ is orientable, the corresponding \emph{oriented motion group} is $\Mot_Y^{+}(M) = \pi_{1}(\cE^{+}(Y,\Int{M}))$.
\end{defi}

\begin{defi}
\label{def:peripheral}
Let $(M,A)$ be a decorated manifold. The submanifold $A \subset M$ is called \emph{peripheral} if it is contained in a codimension-zero ball in $M$.
\end{defi}

\begin{rmk}
If $A \subset M$ is peripheral, the columns of the commutative diagram \eqref{eq:3x3diagram} are moreover split short exact. Indeed, in this case, one may construct sections up to homotopy for the top two vertical inclusions in \eqref{eq:3x3diagram-proof} (or alternatively consider the two columns as instances of \eqref{eq:split-fibration-sequence-2} with $Z$ replaced by $A$ or $A \sqcup Z$), which implies that the columns in \eqref{eq:3x3diagram} are split short exact.
\end{rmk}

When $A \subset M$ is peripheral (see Definition~\ref{def:peripheral}), there is a canonical identification of (oriented) motion groups with $\pi_{0}$ of (oriented) braided diffeomorphism groups. Thus -- in this case -- it follows that \emph{motion groups are braided mapping class groups} (see Remark~\ref{rmk:globality} below).

\begin{prop}
\label{prop:braided-diff-groups}
Let $(M,A)$ be a decorated manifold. There is a canonical surjection
\begin{equation}
\label{eq:from-motion-to-bmcg}
\Mot_{A}(M) = \pi_{1}(\cE(A,\Int{M})) \longtwoheadrightarrow \pi_{0}(\diffdecbr(M,A))
\end{equation}
that is an isomorphism if $A \subset M$ is peripheral. The analogous statement also holds in the oriented setting, adding a superscript ${}^{+}$ to all groups and embedding spaces involved.
\end{prop}
\begin{proof}
The restriction map $\diffdec(M,\emptyset) \to \cE(A,\mathring{M})$ is a fibre bundle by the same reasoning as in the proof of Proposition~\ref{prop:split-fibration-sequence-2}, so (using the same argument as in the proof of Corollary~\ref{coro:split-ses-2}) it induces a long exact sequence of the form
\[
\cdots \to \pi_{1}(\diffdec(M,\emptyset)) \too \Mot_{A}(M) \too \pi_{0}(\diffdec(M,A)) \too \pi_{0}(\diffdec(M,\emptyset)).
\]
By exactness, the image of the middle map is precisely $\pi_{0}(\diffdecbr(M,A))$, so we obtain the surjection \eqref{eq:from-motion-to-bmcg}.
When $A \subset M$ is peripheral, we may choose a codimension-zero ball in $M$ containing $A$ and an embedded path from a point on its boundary to a point $p \in \partial M$ that lies in one of the boundary cylinders of $M$. Since all diffeomorphisms in $\diffdec(M,\emptyset)$ are the identity in a neighbourhood of $p$, we may deform them continuously until they are the identity on the chosen ball and path; in particular on $A$. This deformation provides a section up to homotopy of the inclusion $\diffdec(M,A) \hookrightarrow \diffdec(M,\emptyset)$. Hence this inclusion induces a split-surjection on homotopy groups; by exactness, this implies that the middle map in the long exact sequence above is injective, and hence \eqref{eq:from-motion-to-bmcg} is an isomorphism.
\end{proof}

\begin{rmk}[Globality.]
\label{rmk:globality}
In particular, Proposition~\ref{prop:braided-diff-groups} implies: \emph{each motion group $\Mot_{A}(M)$, where $A \subset M$ is peripheral, is canonically isomorphic to a normal subgroup of the corresponding decorated mapping class group $\pi_{0}(\diffdec(M,A))$}. All of the examples of motion groups discussed in this paper are of this form; see \S\ref{ss:categories_for_families_of_groups}.

Although motion groups and braided mapping class groups are in general different without the peripherality condition, there is always the canonical surjection \eqref{eq:from-motion-to-bmcg}. Thus we may conclude: \emph{a representation of the discrete groupoid $\pi_{0}(\cD_{d})$ canonically induces representations of all motion groups and all mapping class groups of $d$-dimensional manifolds}.
\end{rmk}

\subsubsection{Actions of motion groups}
\label{sss:actions}

To conclude this subsection, we prove a result identifying the actions coming from the two split short exact sequences \eqref{eq:split-ses-1} and \eqref{eq:split-ses-2}. We first need the following general result about connecting homomorphisms in long exact sequences associated to quotients of topological groups.

\begin{prop}
\label{prop:equivariant-connecting-homomorphism}
Let $\Gamma$ be a topological group and $H \leq \Gamma$ a subgroup such that the quotient map $\pi \colon \Gamma \to \Gamma/H$ is a fibre bundle. There is a left $H$-action
\begin{itemizeb}
\item on $\pi_{1}(\Gamma/H)$ induced by the obvious left $H$-action on $\Gamma/H$,
\item on $\pi_{0}(H)$ by conjugation.
\end{itemizeb}
With respect to these actions, the connecting homomorphism
\begin{equation}
\label{eq:conecting-homomorphism}
\delta \colon \pi_{1}(\Gamma/H) \longrightarrow \pi_{0}(H)
\end{equation}
of the long exact sequence associated to $\pi$ is $H$-equivariant.
\end{prop}

\begin{rmk}
\label{rmk:multiplication-vs-conjugation}
Note that the obvious left action of the \emph{whole} group $\Gamma$ on $\Gamma/H$ does \emph{not} induce an action on $\pi_{1}(\Gamma/H)$, since it does not preserve the basepoint.
Also, Proposition~\ref{prop:equivariant-connecting-homomorphism} may seem counterintuitive, since the action on the domain of \eqref{eq:conecting-homomorphism} is induced by \emph{multiplication} in $\Gamma$ whereas the action on its codomain is induced by \emph{conjugation} in $H$. The reason for this is explained in Remark~\ref{rmk:why-conjugation} after the proof.
\end{rmk}

\begin{proof}[Proof of Proposition~\ref{prop:equivariant-connecting-homomorphism}]
Let $[\alpha] \in \pi_{1}(\Gamma/H)$ represented by a loop $\alpha \colon [0,1] \to \Gamma/H$ and let $h \in H$. Our goal is to show that $h.\delta([\alpha]) = \delta(h.[\alpha])$. We first compute the left-hand side of this expression. Choose a path $\widetilde{\alpha} \colon [0,1] \to \Gamma$ such that $\widetilde{\alpha}(0) = 1_\Gamma$ and $\pi \circ \widetilde{\alpha} = \alpha$ (this is possible since $\pi$ is a fibre bundle, in particular a Serre fibration). By definition of the connecting homomorphism, we have $\delta([\alpha]) = [\widetilde{\alpha}(1)]$. By definition of the left-action of $H$ on $\pi_{0}(H)$ by conjugation, we thus have
\begin{equation}
\label{eq:equivariance-left}
h.\delta([\alpha]) = [h \widetilde{\alpha}(1) h^{-1}].
\end{equation}
We now compute the right-hand side, namely $\delta(h.[\alpha])$. The element $h.[\alpha] \in \pi_{1}(\Gamma/H)$ is represented by the loop $h.\alpha \colon [0,1] \to \Gamma/H$ given by applying pointwise the left action of $H$ on $\Gamma/H$, i.e.~it sends $t \in [0,1]$ to $h.\alpha(t) \in \Gamma/H$. We lift this to $\Gamma$ by the path $h.\widetilde{\alpha}.h^{-1} \colon [0,1] \to \Gamma$, which is the composition along the top of the following diagram:
\begin{equation}
\label{eq:lifting-alpha}
\begin{tikzcd}
& \Gamma \ar[r,"h.-"] \ar[d,"\pi"] & \Gamma \ar[d,"\pi"] \ar[r,"-.h^{-1}"] & \Gamma \ar[dl,"\pi"] \\
{[0,1]} \ar[r,"\alpha"] \ar[ur,"\widetilde{\alpha}"] & \Gamma/H \ar[r,"h.-"] & \Gamma/H.
\end{tikzcd}
\end{equation}
The square commutes since $\pi$ is equivariant with respect to the left action of $H$ on $\Gamma$ by multiplication and the right-hand triangle commutes since $\pi$ quotients by the right action of $H$ on $\Gamma$ by multiplication. Hence $\beta := h.\widetilde{\alpha}.h^{-1}$ is a lift of $h.\alpha$ to $\Gamma$, which sends $0$ to $h.\widetilde{\alpha}(0).h^{-1} = h.h^{-1} = 1_\Gamma$. Hence, by definition of the connecting homomorphism, we have:
\begin{equation}
\label{eq:equivariance-right}
\delta(h.[\alpha]) = \delta([h.\alpha]) = [\beta(1)] = [h.\widetilde{\alpha}(1).h^{-1}].
\end{equation}
This is equal to \eqref{eq:equivariance-left}, so we have shown that $\delta$ is $H$-equivariant.
\end{proof}

\begin{rmk}
\label{rmk:why-conjugation}
In the above proof, rather than \eqref{eq:lifting-alpha}, the first naive guess for lifting $h.\alpha$ to $\Gamma$ might instead be the path $h.\widetilde{\alpha}$, namely the composition along the top of the following diagram:
\begin{equation}
\label{eq:lifting-alpha-naive}
\begin{tikzcd}
& \Gamma \ar[r,"h.-"] \ar[d,"\pi"] & \Gamma \ar[d,"\pi"] \\
{[0,1]} \ar[r,"\alpha"] \ar[ur,"\widetilde{\alpha}"] & \Gamma/H \ar[r,"h.-"] & \Gamma/H.
\end{tikzcd}
\end{equation}
This is indeed a lift of $h.\alpha$ due to the fact that $\pi$ is $H$-equivariant, so the square commutes. However, this path will not do for the definition of $\delta$, since it sends $0$ to $h.\widetilde{\alpha}(0) = h$ instead of $1_\Gamma$. This is why we must consider the lift \eqref{eq:lifting-alpha} rather than \eqref{eq:lifting-alpha-naive} in the proof, and why the left $H$-action on the domain of $\delta$ by \emph{multiplication} corresponds to the left $H$-action on the codomain of $\delta$ by \emph{conjugation}; see Remark~\ref{rmk:multiplication-vs-conjugation}.
\end{rmk}

\begin{coro}
\label{coro:equivariant-connecting-homomorphism}
In the setting of Proposition~\ref{prop:equivariant-connecting-homomorphism}, if we have a continuous group homomorphism $s \colon \Gamma \to H$ whose composition with the inclusion $H \subset \Gamma$ is homotopic to the identity on $\Gamma$, there is a split short exact sequence
\begin{equation}
\label{eq:split-ses-abstract}
\centering
\begin{split}
\begin{tikzpicture}
[x=1mm,y=1mm]
\node (ll) at (0,0) {$1$};
\node (l) [anchor=west] at ($ (ll.east) + (5,0) $) {$\pi_{1}(\Gamma/H)$};
\node (m) [anchor=west] at ($ (l.east) + (10,0) $) {$\pi_{0}(H)$};
\node (r) [anchor=west] at ($ (m.east) + (10,0) $) {$\pi_{0}(\Gamma)$};
\node (rr) [anchor=west] at ($ (r.east) + (5,0) $) {$1.$};
\draw[->] (ll) to (l);
\draw[->] (l) to node[above,font=\small]{$\delta$} (m);
\draw[->] (m) to (r);
\draw[->] (r) to (rr);
\draw[->,densely dashed] (r.north west) to[out=160,in=20] node[above,font=\small]{$s_*$} (m.north east);
\end{tikzpicture}
\end{split}
\end{equation}
The left action of $[g] \in \pi_{0}(\Gamma)$ on $\pi_{1}(\Gamma/H)$ arising from this split short exact sequence is precisely the action of $s(g) \in H$ on $\pi_{1}(\Gamma/H)$ induced by the obvious left action of $H$ on $\Gamma/H$.
\end{coro}
\begin{proof}
The split short exact sequence \eqref{eq:split-ses-abstract} follows immediately by taking the long exact sequence of homotopy groups associated to the fibre bundle $\Gamma \to \Gamma/H$ and noting that $s$ provides a section $s_* \colon \pi_{0}(\Gamma) \to \pi_{0}(H)$ for the inclusion-induced map $\pi_{0}(H) \to \pi_{0}(\Gamma)$.
The left action of $[g] \in \pi_{0}(\Gamma)$ on $\pi_{1}(\Gamma/H)$ arising from \eqref{eq:split-ses-abstract} is given by lifting $[g]$ to $\pi_{0}(H)$ using $s_*$ and then acting by conjugation on $\pi_{1}(\Gamma/H)$, considered as a normal subgroup via $\delta$. Hence it is given explicitly by:
\[
[g] \cdot [\alpha] = \delta^{-1} \bigl( s_*([g]) \delta([\alpha]) s_*([g])^{-1} \bigr).
\]
By the $H$-equivariance of $\delta$ deduced from Proposition~\ref{prop:equivariant-connecting-homomorphism}, this is the same as the action of $s_*([g]) = [s(g)]$ on $\pi_{1}(\Gamma/H)$ induced by the obvious left action of $H$ on $\Gamma/H$.
\end{proof}

\begin{prop}
\label{prop:two-actions}
Let $(M,A)$ be a decorated manifold. Extending decorated diffeomorphisms of $(M,A)$ by the identity on $\mbar \smallsetminus M$, there is an action of $\diffdec(M,A)$ on $\cE_{\sG}(Z,\Breve{M} \smallsetminus A)$ by post-composition. This action fixes the basepoint, since $Z$ lies in $\mbar \smallsetminus M$, and so it induces an action
\begin{equation}
\label{eq:natural-action}
\pi_{0}(\diffdec(M,A)) \;\curvearrowright\; \pi_{1}(\cE_{\sG}(Z,\Breve{M} \smallsetminus A)).
\end{equation}

\textup{(i)} The action \eqref{eq:natural-action} agrees with the action induced by the split short exact sequence \eqref{eq:split-ses-2}, pre-composed with the isomorphism
\[
\pi_{0}(\diffdec(M,A)) \cong \pi_{0}(\diffdec(\mbar,A))
\]
induced by extending decorated diffeomorphisms by the identity on $\mbar \smallsetminus M$.

\textup{(ii)} Precomposing with the canonical map (see Proposition~\ref{prop:braided-diff-groups})
\begin{equation}
\label{eq:from-motion-to-mcg}
\pi_{1}(\cE(A,\Int{M})) \longtwoheadrightarrow \pi_{0}(\diffdecbr(M,A)) \leq \pi_{0}(\diffdec(M,A)),
\end{equation}
the action \eqref{eq:natural-action} agrees with the action induced by the split short exact sequence \eqref{eq:split-ses-1}, pre-composed with the isomorphism $\pi_{1}(\cE(A,\Int{M})) \cong \pi_{1}(\cE(A,\Breve{M}))$ induced by the inclusion $\Int{M} \subset \Breve{M}$.

Both statements also hold in the oriented setting, adding a superscript ${}^{+}$ to all groups and all embedding spaces (except those already with a subscript ${}_G$) in the statement.
\end{prop}
\begin{proof}
As explained in the proof of Corollary~\ref{coro:split-ses-2}, the fibre bundle
\[
\diffdec(\mbar,A) \too \cE_{\sG}(Z,\Breve{M} \smallsetminus A)
\]
from \eqref{eq:split-fibration-sequence-2} is isomorphic (in the category of maps out of $\Gamma$) to the quotient map $\Gamma \twoheadrightarrow \Gamma/H$, where $\Gamma = \diffdec(\mbar,A)$ and $H = \diffdec(\mbar,A,Z|\sG)$ and we have possibly thrown away some non-basepoint path-components of $\cE_{\sG}(Z,\Breve{M} \smallsetminus A)$, which does not matter since we are only interested in its fundamental group. We are therefore in the setting of Corollary~\ref{coro:equivariant-connecting-homomorphism}, with $\eqref{eq:split-ses-abstract} = \eqref{eq:split-ses-2}$. It then follows from Corollary~\ref{coro:equivariant-connecting-homomorphism} that the action of $\pi_{0}(\Gamma)$ on $\pi_{1}(\cE_{\sG}(Z,\Breve{M} \smallsetminus A))$ induced by the split short exact sequence \eqref{eq:split-ses-2} is equal to the section $\pi_{0}(\Gamma) \to \pi_{0}(H)$ followed by the action of $\pi_{0}(H) = \pi_{0}(\diffdec(\mbar,A,Z|\sG))$ on $\pi_{1}(\cE_{\sG}(Z,\Breve{M} \smallsetminus A))$ by post-composition. Pre-composing the latter action by the isomorphism $\pi_{0}(\diffdec(M,A)) \cong \pi_{0}(\diffdec(\mbar,A)) = \pi_{0}(\Gamma)$, we obtain exactly the action \eqref{eq:natural-action}, because extending by the identity on $\mbar \smallsetminus M$, conjugating by the isomorphism $\Theta \colon \mbar \to M$ (see Construction~\ref{construction:Theta}) and extending by the identity again is isotopic to simply extending by the identity once. This establishes part (i) of the proposition.

Part (ii) of the proposition then follows from part (i), together with diagram \eqref{eq:3x3diagram} of Proposition~\ref{prop:map_of_ses} (in particular the fact that its top-right square commutes also when the solid horizontal arrows are replaced by the dotted arrows) and the observation that the top-right vertical map of that diagram coincides with the canonical map \eqref{eq:from-motion-to-mcg} under the two isomorphisms
\begin{equation}
\label{eq:two-isomorphisms}
\pi_{1}(\cE(A,\Int{M})) \cong \pi_{1}(\cE(A,\Breve{M})) \quad\text{ and }\quad \pi_{0}(\diffdec(M,A)) \cong \pi_{0}(\diffdec(\mbar,A)).
\end{equation}
This last observation follows because \eqref{eq:from-motion-to-mcg} is the connecting homomorphism of the restriction fibre bundle $\diffdec(M,\emptyset) \to \cE(A,\Int{M})$ (see the proof of Proposition~\ref{prop:braided-diff-groups}), the top-right vertical map of \eqref{eq:3x3diagram} is the connecting homomorphism of the restriction fibre bundle $\diffdec(\mbar,\emptyset) \to \cE(A,\Breve{M})$ (see diagram \eqref{eq:3x3diagram-proof}) and the isomorphisms \eqref{eq:two-isomorphisms} are induced by the map of fibre bundles induced by extending decorated diffeomorphisms of $M$ to $\mbar$ by the identity on $\mbar \smallsetminus M$ and post-composing embeddings into $\Int{M}$ with the inclusion $\Int{M} \subset \Breve{M}$.
\end{proof}

\subsection{A Serre fibration of decorated diffeomorphism groups}
\label{ss:fibration-condition}

This section is devoted to the proof of a technical result (see Theorem~\ref{thm:fibre-bundle}) concerning decorated diffeomorphism groups, which will play an important role in \S\ref{ss:Quillen_bracket_construction} to ensure that the topologically-enriched Quillen bracket construction behaves well with respect to the $\pi_{0}$ functor on the topologically-enriched groupoids of decorated manifolds that we shall consider; see Lemma~\ref{lem:Serre-fibration-condition} and Proposition~\ref{prop:morphism-spaces-bracket}. \textbf{From this point onwards, we assume that the dimension $d$ is not equal to $4$}; see Remark~\ref{rmk:d_neq_4_explanation} for the reason why this is necessary. We start by introducing the following notion.

\begin{defi}[Decorated embeddings.]
\label{def:embdec}
Let $M = (M,A)$ and $N = (N,B)$ be decorated manifolds. Define $\embdec(M,N)$ to be the set of smooth, proper embeddings $\varphi \colon M \hookrightarrow N$, equipped with a germ of an extension $\varphi'$ to an embedding $\bB_{1}^d \natural M \hookrightarrow N$, such that:
\begin{itemizeb}
\item $\varphi(A) \subseteq B$;
\item for some $\epsilon > 0$ we have $\varphi(e_{2}(\bB_{\epsilon}^d)) = e'_{2}(\bB_{\epsilon}^d)$ and $(e'_{2})^{-1} \circ \varphi \circ e_{2}$ is the identity map $\bB_{\epsilon}^d \to \bB_{\epsilon}^d$, where $e_{1},e_{2}$ are the boundary-cylinder germs of $(M,A)$ and $e'_{1},e'_{2}$ are those of $(N,B)$;
\item there is a decorated manifold $M'$ and diffeomorphism of decorated manifolds $\bar{\varphi} \colon M' \natural M \to N$ such that $\varphi = \bar{\varphi} \circ \iota_{M,M'}$, where $\iota_{M,M'}$ denotes the canonical embedding of $M$ into $M' \natural M$. This extension of $\varphi$ to $\bar{\varphi}$ should be compatible with the given germ of an extension $\varphi'$ of $\varphi$.
\end{itemizeb}
In other words -- modulo decorated diffeomorphisms of the codomain -- a decorated embedding is the canonical inclusion of (a germ of a neighbourhood of) the right-hand factor of a boundary connected sum. As an illustration, Figure~\ref{fig:decorated-manifolds-Mt} depicts the inclusion of a neighbourhood $M_t$ of $M$ into $L\natural M$, representing a decorated embedding of $M$ into $L\natural M$.

We now describe the topology on this set. Recall that $\mbar = \bB^d_1 \natural M = (\bD^{d-1} \times [-1,0]) \cup M$ and write $M_t = (\bD^{d-1} \times [-t,0]) \cup M$ for $t \in (0,1)$, as depicted in Figure~\ref{fig:decorated-manifolds-Mt}. For $\epsilon \in (0,1)$, write $\embdece(M_t,N)$ for the set of embeddings defined as above except that $\epsilon$ is fixed in the second bullet point and the embedding is defined on $M_t$ rather than on $M$ plus a germ of a neighbourhood. Topologise $\embdece(M_t,N)$ as a subspace of $C^\infty(M_t,N)$ equipped with the Whitney topology. We may take the colimit as $\epsilon \to 0$ (along inclusion maps) and $t \to 0$ (along restriction maps), and we define $\embdec(M,N) := o\bigl( \underset{\epsilon,t \to 0}{\mathrm{colim}}(\embdece(M_t,N)) \bigr)$, where $o(-)$ is the functor of Construction~\ref{construction:open-path-components}.
\end{defi}

\begin{notation}
\label{notation:decorated-embeddings}
For decorated manifolds $L,M,N$, we denote by $\embdec(M,N)_L$ the subspace of $\embdec(M,N)$ of those embeddings for which we may take $M' = L$ in the third point above.
Additionally, if $A$ and $B$ are oriented, we denote by $\embdec^{+}(M,N)$ the subspace of $\embdec(M,N)$ of those decorated embeddings $\varphi \colon M \hookrightarrow N$ whose restriction $\varphi|_A \colon A \hookrightarrow B$ is orientation-preserving.
\end{notation}

The composition of decorated embeddings is defined as follows.

\begin{defi}
\label{def:decorated-embeddings-composition}
Let $\varphi\in\embdec(M,N)$ and $\psi\in\embdec(N,P)$ for decorated manifolds $M$, $N$ and $P$. Choose decorated diffeomorphisms $\bar{\varphi} \colon M' \natural M \to N$ and $\bar{\psi} \colon N' \natural N \to P$ extending $\varphi$ and $\psi$ respectively; these exist by definition. We may then form the composite decorated diffeomorphism $\bar{\psi} \circ (\mathrm{id}_{N'} \natural \bar{\varphi}) \colon N' \natural M' \natural M \to P$, which restricts to the composite embedding $\psi \circ \varphi$ on $M \subset N' \natural M' \natural M$. Moreover, its restriction to an $\epsilon$-neighbourhood of $M \subset N' \natural M' \natural M$ provides a germ of an extension of $\psi \circ \varphi$ to $\bB_{1}^d \natural M$. The first two points of Definition~\ref{def:embdec} clearly hold for the composite embedding $\psi \circ \varphi$, so this embedding together with the germ of an extension is an element of $\embdec(M,P)$, which we define to be the composition of the decorated embeddings $\varphi$ and $\psi$.
\end{defi}

\begin{lemm}
\label{lem:decorated-embeddings-composition-continuous-associative}
The composition of decorated embeddings given in Definition \ref{def:decorated-embeddings-composition} defines a continuous, associative operation $\embdec(M,N) \times \embdec(N,P) \to \embdec(M,P)$.
\end{lemm}

\begin{proof}
We must first check that this operation is a well-defined map, i.e.~that the construction of Definition \ref{def:decorated-embeddings-composition} does not depend on the choice of decorated diffeomorphisms $\bar{\varphi}$ and $\bar{\psi}$ extending $\varphi$ and $\psi$. The only part of the construction that depends on these choices of extensions is the germ of an extension of $\psi \circ \varphi$ to an embedding $\bB_{1}^d \natural M \hookrightarrow P$. The fact that this germ is independent of the choice of decorated diffeomorphisms extending the decorated embeddings follows from the fact that decorated diffeomorphisms are required to be compatible with boundary-cylinder germs. Continuity of the operation follows from the fact that composition of proper, smooth maps is continuous in the Whitney topology; see \cite[\S 2, Prop.~1]{Mather1969}.

Finally, let us check associativity of this operation. For $\varphi\in\embdec(M,N)$, $\psi\in\embdec(N,P)$ and $\chi\in\embdec(P,Q)$, the two iterated compositions $(\chi \circ \psi) \circ \varphi$ and $\chi \circ (\psi \circ \varphi) \in \embdec(M,Q)$ each consist of a smooth embedding $M \hookrightarrow Q$ and a germ of an extension to $\bB_{1}^d \natural M \hookrightarrow Q$. The two smooth embeddings $M \hookrightarrow Q$ are evidently equal, by associativity of composition of maps. Moreover, the fact that the two germs of extensions to embeddings $\bB_{1}^d \natural M \hookrightarrow Q$ are equal follows from the fact that they are restrictions, to an $\epsilon$-neighbourhood of $M \subset P' \natural N' \natural M' \natural M$, of the decorated diffeomorphisms $(\bar{\chi} \circ (\mathrm{id}_{P'} \natural \bar{\psi})) \circ (\mathrm{id}_{P' \natural N'} \natural \bar{\varphi})$ and $\bar{\chi} \circ ((\mathrm{id}_{P'} \natural \bar{\psi}) \circ (\mathrm{id}_{P' \natural N'} \natural \bar{\varphi}))$, which are equal by associativity of composition of maps.
\end{proof}

\begin{lemm}
\label{lem:decomposition-disjoint-union}
The embedding space $\embdec(M,N)$ decomposes as a topological disjoint union
\begin{equation}
\label{eq:embdec-decomposition}
\embdec(M,N) \;\cong\; \bigsqcup_L \, \embdec(M,N)_L,
\end{equation}
where the disjoint union runs over representatives of isomorphism classes of decorated manifolds.
\end{lemm}
\begin{proof}
We denote by $(-)^{\mathfrak{c}}$ the function
\[
\embdec(M,N) \too \{ \text{isomorphism classes of decorated manifolds} \}
\]
that assigns to each decorated embedding the isomorphism class of the decorated manifold given by the complement of its image. We must show that this function is locally constant. Let $\varphi \in \embdec(M,N)$ be any decorated embedding; we must show that it admits a neighbourhood $\cU$ such that $\psi^{\mathfrak{c}} = \varphi^{\mathfrak{c}}$ for each $\psi \in \cU$.

Since $\varphi$ is a decorated embedding, we may fix an identification of $N$ with $L \natural M$, for some decorated manifold $L$, such that $\varphi$ is (the germ of an extension to a collar neighbourhood of) the canonical inclusion of $M$ into $L \natural M$.

We first note that there exists an open neighbourhood $\cU$ of $\varphi$ such that, for each $\psi \in \cU$, the interface $\psi(M) \cap (\overline{(L \natural M)) \smallsetminus \psi(M)})$ between the image of $\psi$ and its complement is a smoothly embedded $(d-1)$-dimensional disc $D(\psi)$ contained in the central solid cylinder $\bD^{d-1} \times [-1,1]$ along which the boundary connected sum of $L \natural M$ is formed; see Figure~\ref{fig:decorated-manifolds-Mt}. Moreover, $D(\psi)$ intersects the boundary of $\bD^{d-1} \times [-1,1]$ precisely in $\partial D(\psi)$. Such an open neighbourhood $\cU$ exists in the compact-open topology, and hence also in the topology that we are considering on $\embdec(M,N)$ (a refinement of a colimit of Whitney topologies), which is finer than the compact-open topology.

It remains to show, for a given $\psi \in \cU$, that $\psi^{\mathfrak{c}} = \varphi^{\mathfrak{c}}$. The embedded codimension-$1$ disc $D(\psi)$ cuts the solid cylinder $\bD^{d-1} \times [-1,1]$ into two pieces; let us denote by $A(\psi)$ the piece that contains $\bD^{d-1} \times \{-1\}$. We claim that $A(\psi)$ is diffeomorphic to $\bD^{d-1} \times [-1,0]$ by a diffeomorphism $\Xi$ that acts by the identity on a neighbourhood of $\bD^{d-1} \times \{-1\}$. This claim will immediately imply that $\psi^{\mathfrak{c}} = \varphi^{\mathfrak{c}}$, as desired, by extending $\Xi$ by the identity on $L \smallsetminus (\bD^{d-1} \times [-1,0])$.

Smoothing corners, we may rephrase the claim in the previous paragraph as follows. Consider the standard $d$-ball $\bD^d \subset \bR^d$, denote by $\bD^d_+$ the closed upper half-ball and fix a small $(d-1)$-dimensional subdisc $D_0$ of the boundary $\partial \bD^d$ centred at the north pole. Suppose we are given a smoothly embedded $(d-1)$-dimensional disc $D \subset \bD^d$ with $D \cap \partial \bD^d = \partial D$ disjoint from $D_0$ and denote by $A$ the closure of the component of $\bD^d \smallsetminus D$ containing $D_0$. The claim is that there exists a diffeomorphism $\Xi \colon A \cong \bD^d_+$ that acts by the identity on a neighbourhood of $D_0$. As explained above, this claim will complete the proof.

To prove this claim, we first consider the \emph{double} $2(\bD^d)$ of $\bD^d$, which is a smooth $d$-sphere. Smoothing corners again, the boundary of $A \subset \bD^d \subset 2(\bD^d)$ is a smoothly embedded $(d-1)$-sphere. The smooth Schoenflies theorem (see \cite{Schoenflies1906} for $d=2$, see \cite{Mazur1959,Morse1960,Brown1960} for $d=3$ and see \cite[Prop.~D, \S 9]{Milnor1965} for $d\geq 5$) then implies that $A$ is diffeomorphic to a $d$-ball. Smoothing corners yet again, this gives a diffeomorphism $\Xi' \colon A \cong \bD^d_+$. However, this does not necessarily act by the identity on a neighbourhood of $D_0$.

To ensure this second condition, let $e \colon \bD^d \hookrightarrow \bD^d$ be an orientation-preserving embedding so that $e(\bD^d)$ is a neighbourhood of $D_0$ and so that $e(\bD^d)$ is contained in $A \cap \bD^d_+$. Then $e$ and $\Xi' \circ e$ are two embeddings of $\bD^d$ into $\bD^d_+$. By composing $\Xi'$ with a reflection if necessary, $\Xi' \circ e$ is also orientation-preserving. Palais' disc embedding theorem \cite[Thm.~B]{Palais1960Extendingdiffeomorphisms} implies that there is a diffeomorphism $\Xi'' \colon \bD^d_+ \cong \bD^d_+$ such that $\Xi'' \circ e = \Xi' \circ e$. The composition $\Xi := (\Xi'')^{-1} \circ \Xi' \colon A \cong \bD^d_+$ is then a diffeomorphism that acts by the identity on the neighbourhood $e(\bD^d)$ of $D_0$, as desired.
\end{proof}

\begin{rmk}
\label{rmk:d_neq_4_explanation}
The proof above uses the smooth Schoenflies theorem, which is an open question (rather than a theorem) in dimension $d=4$; this is the reason for the restriction $d \neq 4$ that we make in the present subsection \S\ref{ss:fibration-condition}, as well as consequently in \S\ref{sss:quillen-bracket-categories}--\S\ref{sss:functorial-split-ses}, \S\ref{ss:homological_representation_functor_motion_groups} and \S\ref{ss:homological_representation_functor_mcg}.
\end{rmk}

Now, for the remainder of \S\ref{ss:fibration-condition}, we consider two decorated manifolds $L = (L,A,e_{1},e_{2})$ and $M = (M,B,e'_{1},e'_{2})$. There is a continuous right action of $\diffdec(L)$ on $\diffdec(L \natural M)$ given by $\varphi \cdot \psi = \varphi \circ (\psi \natural \id_M)$, and hence a quotient map
\begin{equation}\label{eq:quotient-diff}
\Psi \colon \diffdec(L \natural M) \too \diffdec(L \natural M) / \diffdec(L).
\end{equation}

\begin{thm}\label{thm:fibre-bundle}
The quotient map \eqref{eq:quotient-diff} is a Serre fibration. There is a homeomorphism between its codomain and $\embdec(M,L \natural M)_L$, induced by the restriction map.

If the submanifolds $A \subset \Int{L}$ and $B \subset \Int{M}$ are oriented and we consider the subgroups $\diffdec^{+}(-)$ of $\diffdec(-)$, then the analogue of the quotient map \eqref{eq:quotient-diff} is also a Serre fibration and there is a homeomorphism from its codomain to $\embdec^{+}(M,L \natural M)_L$, induced by the restriction map.
\end{thm}

\begin{rmk}
This is related to results of Cerf \cite[p.~294, \S II.2.2.2, Cor.~2]{Cerf1961Topologiedecertains}, Palais \cite[Th.~B]{Palais1960Localtrivialityof} and Lima \cite{Lima1963localtrivialityof}, but we were not able to find an instance of their results that covers exactly the setting that we need here. We therefore give a complete proof of Theorem~\ref{thm:fibre-bundle} below, using as an input two results of Cerf and Palais, namely Lemme II.2.1.2 (page 291) of \cite{Cerf1961Topologiedecertains} and \cite[Th.~A]{Palais1960Localtrivialityof}.
\end{rmk}

\begin{figure}[t]
\centering
\begin{tikzpicture}
[x=1mm,y=1mm,scale=1]
\begin{scope}[xshift=5mm]
\fill[blue!20] (0,0) -- (10,0) .. controls (15,0) and (15,-5) .. (20,-5) -- (21,-5) -- (29,-5) -- (30,-5) .. controls (35,-5) and (35,0) .. (40,0) -- (50,0) -- (50,10) -- (40,10) .. controls (35,10) and (35,15) .. (30,15) -- (29,15) -- (21,15) -- (20,15) .. controls (15,15) and (15,10) .. (10,10) -- (0,10) -- cycle;
\draw[thick] (0,0) -- (10,0) .. controls (15,0) and (15,-5) .. (20,-5) -- (21,-5);
\draw[thick,densely dashed] (21,-5) -- (29,-5);
\draw[thick] (29,-5) -- (30,-5) .. controls (35,-5) and (35,0) .. (40,0) -- (50,0);
\draw[thick] (50,10) -- (40,10) .. controls (35,10) and (35,15) .. (30,15) -- (29,15);
\draw[thick,densely dashed] (29,15) -- (21,15);
\draw[thick] (21,15) -- (20,15) .. controls (15,15) and (15,10) .. (10,10) -- (0,10) -- (0,0);
\draw[blue!50,thick] (18,2) circle (2);
\draw[blue!50,thick] (28,7) circle (2);
\fill[blue!50,thick] (22,9) circle (0.5);
\node at (34,5) [blue!50] {$A$};
\fill[white] (29,1) circle (1.5);
\draw[thick] (29,1) circle (1.5);
\draw (7,0) -- (7,10);
\draw (10,0) -- (10,10);
\node at (0,0) [anchor=north,font=\footnotesize] {$0$};
\node at (7,0) [anchor=north,font=\footnotesize] {$\vphantom{0}\epsilon$};
\node at (10,0) [anchor=north,font=\footnotesize] {$1$};
\node at (5,12) [anchor=south] {$e_{1}$};
\draw[decorate,decoration={brace,amplitude=5pt,raise=2pt}] (0,10) -- (10,10);
\node at (0,5) [anchor=east] {$L \natural M =$};
\end{scope}
\begin{scope}[xshift=55mm]
\fill[green!20] (0,0) -- (10,0) .. controls (15,0) and (15,-5) .. (20,-5) -- (21,-5) -- (29,-5) -- (30,-5) .. controls (35,-5) and (35,0) .. (40,0) -- (50,0) -- (50,10) -- (40,10) .. controls (35,10) and (35,15) .. (30,15) -- (29,15) -- (21,15) -- (20,15) .. controls (15,15) and (15,10) .. (10,10) -- (0,10) -- cycle;
\draw[thick] (0,0) -- (10,0) .. controls (15,0) and (15,-5) .. (20,-5) -- (21,-5);
\draw[thick,densely dashed] (21,-5) -- (29,-5);
\draw[thick] (29,-5) -- (30,-5) .. controls (35,-5) and (35,0) .. (40,0) -- (50,0) -- (50,10) -- (40,10) .. controls (35,10) and (35,15) .. (30,15) -- (29,15);
\draw[thick,densely dashed] (29,15) -- (21,15);
\draw[thick] (21,15) -- (20,15) .. controls (15,15) and (15,10) .. (10,10) -- (0,10);
\draw[green!80!black,thick] (17,7) circle (2);
\draw[green!80!black,thick] (22,1) circle (2);
\draw[green!80!black,thick] (29,-1) circle (1.5);
\fill[green!80!black,thick] (29,5) circle (0.5);
\fill[white] (24,9.5) circle (2.5);
\draw[thick] (24,9.5) circle (2.5);
\node at (34,5) [green!80!black] {$B$};
\draw (40,0) -- (40,10);
\draw (43,0) -- (43,10);
\node at (40,0) [anchor=north,font=\footnotesize] {$1$};
\node at (43,0) [anchor=north,font=\footnotesize] {$\vphantom{0}\epsilon$};
\node at (50,0) [anchor=north,font=\footnotesize] {$0$};
\node at (45,12) [anchor=south] {$e'_{2}$};
\draw[decorate,decoration={brace,amplitude=5pt,raise=2pt}] (40,10) -- (50,10);
\draw (-10,0) -- (-10,10);
\draw (-4,0) -- (-4,10);
\draw (0,0) -- (0,10);
\draw (10,0) -- (10,10);
\node at (-10,0) [anchor=north,font=\footnotesize] {$1$};
\node at (-4,0) [anchor=north,font=\footnotesize] {$\vphantom{0}t$};
\node at (0,0) [anchor=north,font=\footnotesize] {$0$};
\node at (10,0) [anchor=north,font=\footnotesize] {$1$};
\draw[decorate,decoration={brace,amplitude=5pt,raise=2pt}] (-10,10) -- (0,10);
\draw[decorate,decoration={brace,amplitude=5pt,raise=2pt}] (0,10) -- (10,10);
\node at (-5,12) [anchor=south] {$e_{2}$};
\node at (5,12) [anchor=south] {$e'_{1}$};
\draw[decorate,decoration={brace,amplitude=5pt,mirror,raise=2pt}] (-4,-7) -- (50,-7);
\node at (23,-9) [anchor=north] {$M_t$};
\end{scope}
\end{tikzpicture}
\caption{The boundary connected sum $L \natural M$ from the proof of Theorem~\ref{thm:fibre-bundle}.}
\label{fig:decorated-manifolds-Mt}
\end{figure}

\begin{proof}[Proof of Theorem~\ref{thm:fibre-bundle}]
In the proof below, we use the topology on the spaces of decorated diffeomorphisms and of decorated embeddings given by a colimit of Whitney topologies, as described in Definitions~\ref{def:decorated-manifolds-morphism-spaces} and \ref{def:embdec}, thus proving that \eqref{eq:quotient-diff} is a Serre fibration for this topology. This result automatically implies that the same result holds when using the finer topology making all path-components open; see Lemma~\ref{lem:open-path-components-Serre-fibrations}. Hence we shall make no further mention of this finer topology and work directly with the colimit of Whitney topologies.

Moreover, we detail here the proof for the first (unoriented) case of Theorem~\ref{thm:fibre-bundle}. All of the arguments below repeat mutatis mutandis in the setting where the submanifolds $A$ and $B$ are oriented and the morphisms of decorated manifolds preserve these orientations.

The decorated manifolds $L = (L,A,e_{1},e_{2})$ and $M = (M,B,e'_{1},e'_{2})$ come equipped with germs $e_{1},e_{2},e'_{1},e'_{2}$ of boundary cylinders; let us once and for all choose representative boundary cylinders for these germs, and denote them by the same symbols, by abuse of notation.

For $\epsilon \in (0,1)$, let $\diffdece(L\natural M)$ denote the group of self-diffeomorphisms of $L\natural M$ sending $A \sqcup B$ onto itself and restricting to the identity on $e_{1}(\bB_{\epsilon}^d)$ and on $e'_{2}(\bB_{\epsilon}^d)$. If we give this the Whitney topology, then we have
\begin{equation}
\label{eq:diff-colimit-identification}
\diffdec(L\natural M) \cong \underset{\epsilon\to 0}{\colim} (\diffdece(L\natural M)),
\end{equation}
by Definitions~\ref{def:decorated-manifolds-morphisms} and \ref{def:decorated-manifolds-morphism-spaces} (recalling that, after the first paragraph of the present proof, we are ignoring the refinement $o(-)$ of the topology). Similarly, for each $\epsilon,t \in (0,1)$, let $\diffdecet(L)$ denote the subgroup of $\diffdece(L\natural M)$ consisting of diffeomorphisms that restrict to the identity on the submanifold $M_t = M \cup e_{2}(\bB_t^d)$ of $L\natural M$ pictured in Figure~\ref{fig:decorated-manifolds-Mt}. We have a quotient map
\[
\Psi_{\epsilon,t} \colon \diffdece(L\natural M) \too \diffdece(L\natural M) / \diffdecet(L).
\]
For any $\epsilon,\epsilon',t,t' \in (0,1)$ with $\epsilon' \leq \epsilon$ and $t' \leq t$ there are natural maps
\[
\diffdece(L\natural M) / \diffdecet(L) \too \diffdeceprime(L\natural M) / \diffdecetprime(L),
\]
so we may take the directed colimit of the maps $\Psi_{\epsilon,t}$ to obtain
\[
\underset{\epsilon,t\to 0}{\colim}(\Psi_{\epsilon,t}) \colon \diffdec(L\natural M) \too \underset{\epsilon,t\to 0}{\colim}(\diffdece(L\natural M) / \diffdecet(L)),
\]
where we have used the identification \eqref{eq:diff-colimit-identification} in the domain. Since each $\Psi_{\epsilon,t}$ is a quotient map, using the characterisation of quotient maps as coequalisers and the fact that colimits commute, we deduce that \raisebox{0pt}[0pt][0pt]{$\underset{\epsilon,t\to 0}{\colim}(\Psi_{\epsilon,t})$} is also a quotient map. The map
\[
\Psi \colon \diffdec(L \natural M) \too \diffdec(L\natural M) / \diffdec(L),
\]
i.e.~the map \eqref{eq:quotient-diff} that we would like to show is a Serre fibration, is also a quotient map, with the same domain. Since $M_t$ is a cofinal family of neighbourhoods of $M$ in $L \natural M$, two diffeomorphisms of $\diffdec(L \natural M)$ have the same image under $\Psi$ if and only if they have the same image under
$\underset{\epsilon,t\to 0}{\colim}(\Psi_{\epsilon,t})$. As they are quotient maps with the same domain, it follows that $\Psi \cong \underset{\epsilon,t\to 0}{\colim}(\Psi_{\epsilon,t})$.

We will prove below that each $\Psi_{\epsilon,t}$ is a fibre bundle (and hence a Serre fibration), and then deduce that $\Psi$ is also a Serre fibration using the following general fact.
\begin{itemizeb}
\item[$(*)$] Any filtered colimit of based Serre fibrations between compactly-generated weak-Hausdorff spaces is again a Serre fibration.
\end{itemizeb}
For a reference for this fact, see \cite[Prop.~1.2.3.5(1)]{ToenVezzosi}, which states that a filtered colimit of fibrations is a fibration in any compactly generated model category. The classical model category of based compactly-generated weak-Hausdorff spaces, with its Quillen model structure in which the fibrations are the Serre fibrations, is compactly generated; see for example \cite[Prop.~6.3]{MMSS}.

To apply $(*)$ in our situation, first note that we are taking a directed colimit, which is in particular a filtered colimit. We then need to check that the diffeomorphism groups $\diffdece(L \natural M)$ and their quotients are compactly-generated weak-Hausdorff spaces. Recall that, when restricting to subspaces of proper maps, the Whitney topology on spaces of smooth maps coincides with the weak topology, which is second-countable and hence also first-countable; see \cite[Chap.~2.1]{Hirsch1976Differentialtopology}. Also, diffeomorphism groups of smooth manifolds are Hausdorff in the compact-open topology, and thus also in the Whitney topology, since the latter is a finer topology. Hence diffeomorphism groups of smooth manifolds, in the Whitney topology, are always first-countable and Hausdorff, and thus compactly-generated and weak-Hausdorff. Moreover, the property of being compactly-generated is preserved when taking quotients. The property of being (weak) Hausdorff is \emph{not} preserved when taking quotients; however, in the process of proving that each $\Psi_{\epsilon,t}$ is a fibre bundle below, we will also show that its target space $\diffdece(L\natural M) / \diffdecet(L)$ is Hausdorff.

It therefore remains to show that each $\Psi_{\epsilon,t}$ is a fibre bundle and that its target space is Hausdorff. Write
\[
\embdece(M_t,L\natural M)_L
\]
for the space of smooth, proper embeddings $\varphi \colon M_t \to L \natural M$ such that $\varphi(B) \subseteq A \sqcup B$, the restriction of $\varphi$ to $e'_{2}(\bB_{\epsilon}^d)$ is the identity and there exists $\bar{\varphi} \in \diffdece(L \natural M)$ such that $\varphi = \bar{\varphi} \circ \iota$, where $\iota$ is the inclusion of $M_t$ into $L\natural M$. Then we have
\begin{equation}
\label{eq:emb-colimit-identification}
\embdec(M,L\natural M)_L \cong \underset{\epsilon,t\to 0}{\colim} (\embdece(M_t,L\natural M)_L),
\end{equation}
see Definition~\ref{def:embdec}, Notation~\ref{notation:decorated-embeddings} and Lemma~\ref{lem:decomposition-disjoint-union}. There is a restriction map
\[
\Phi_{\epsilon,t} \colon \diffdece(L\natural M) \too \embdece(M_t,L\natural M)_L,
\]
which is equivariant with respect to the left action of $\diffdece(L\natural M)$ by post-composition. This factors through the quotient map $\Psi_{\epsilon,t}$, so we have an induced map
\begin{center}
\begin{tikzpicture}
[x=1mm,y=1mm]
\node (b) at (0,0) {$\diffdece(L \natural M) / \diffdecet(L)$};
\node (t) at (0,15) {$\diffdece(L \natural M)$};
\node (r) at (40,15) {$\embdece(M_t,L\natural M)_L.$};
\draw[->>] (t) to node[left,font=\small]{$\Psi_{\epsilon,t}$} (b);
\draw[->] (t) to node[above,font=\small]{$\Phi_{\epsilon,t}$} (r);
\draw[->] (b.north east) to (r);
\node at (33,5) [font=\small] {$\widehat{\Phi}_{\epsilon,t}$};
\end{tikzpicture}
\end{center}
By definition of the right-hand embedding space, the map $\Phi_{\epsilon,t}$ is surjective, and so is the induced map $\widehat{\Phi}_{\epsilon,t}$. Moreover, if two diffeomorphisms of $\diffdece(L\natural M)$ have the same image under $\Phi_{\epsilon,t}$, their difference lies in $\diffdecet(L)$, so the induced map $\widehat{\Phi}_{\epsilon,t}$ is also injective. 

We claim that it is sufficient to prove that $\Phi_{\epsilon,t}$ is a fibre bundle. Indeed, it is then a quotient map, since surjective fibre bundles are always quotient maps. Thus the induced map $\widehat{\Phi}_{\epsilon,t}$ must be a homeomorphism. This implies:
\begin{itemizeb}
\item[$\circ$] The map $\Psi_{\epsilon,t}$ is also a fibre bundle, hence a Serre fibration.
\item[$\circ$] Its target space is homeomorphic to the embedding space $\embdece(M_t,L\natural M)_L$, which we have given the Whitney topology, so it is Hausdorff.
\item[$\circ$] We also obtain the second statement of the proposition:
\begin{align*}
\diffdec(L \natural M) / \diffdec(L) &\cong
\underset{\epsilon,t\to 0}{\colim}(\diffdece(L\natural M) / \diffdecet(L)) \\
&\cong \underset{\epsilon,t\to 0}{\colim}(\embdece(M_t,L\natural M)_L) \\
&\cong \embdec(M,L\natural M)_L,
\end{align*}
by combining the identification of $\Psi \cong \underset{\epsilon,t\to 0}{\colim}(\Psi_{\epsilon,t})$ with the colimit of $\widehat{\Phi}_{\epsilon,t}$ and \eqref{eq:emb-colimit-identification}.
\end{itemizeb}

It therefore remains just to prove that $\Phi_{\epsilon,t}$ is a fibre bundle. Since it is equivariant with respect to the left action of $\diffdece(L\natural M)$, it suffices to prove that the action of $\diffdece(L\natural M)$ on $\embdece(M_t,L\natural M)_L$ is \emph{locally retractile}. This is because, by \cite[Th.~A]{Palais1960Localtrivialityof}, any $G$-equivariant map into a $G$-locally retractile space is a fibre bundle.

Thus, we have to prove the following: given an embedding $e \in \embdece(M_t,L\natural M)_L$, we may find an open neighbourhood $\cU$ of $e$ and continuous map $\gamma \colon \cU \to \diffdece(L\natural M)$ such that $\gamma(e) = \id$ and $\gamma(f) \circ e = f$ for any $f \in \cU$. Note that, since $\diffdece(L\natural M)$ acts transitively on $\embdece(M_t,L\natural M)_L$ (because $\Phi_{\epsilon,t}$ is both equivariant and surjective), it suffices to prove this for just one such $e$, which we take to be the inclusion $M_t \hookrightarrow L \natural M$.

To prove this, we apply a result of Cerf \cite[p.~291, Lem.~II.2.1.2]{Cerf1961Topologiedecertains}, which we first recall. Let $X$ be a manifold-with-corners. This means in particular that $X$ has a stratification into \emph{faces} (for example, if $X$ is a connected manifold with boundary, but no higher-codimension corners, then its set of faces is $\pi_{0}(\partial X) \sqcup \{X\}$). Each point $x \in X$ may lie in many faces, but it has a unique \emph{smallest} face (according to inclusion) in which it lies, which we denote by $\mathsf{C}_{X}(x)$. Now if $Y$ is any submanifold-with-corners of $X$, we define
\[
C^{\infty}_{\mathrm{face}}(Y,X) = \{ \text{smooth maps } \varphi \colon Y \to X \text{ such that } \mathsf{C}_{X}(\varphi(x)) = \mathsf{C}_{X}(x) \text{ for each } x \in Y \} ,
\]
equipped with the Whitney topology. The \emph{Extension Lemma} II.2.1.2 of \cite{Cerf1961Topologiedecertains} says that, if $Y$ is closed in $X$ and $V$ is any neighbourhood of $Y$ in $X$, then the restriction map
\[
C^{\infty}_{\mathrm{face}}(X,X) \too C^{\infty}_{\mathrm{face}}(Y,X)
\]
admits a section $s$ defined on an open neighbourhood $\cV$ of the inclusion in $C^{\infty}_{\mathrm{face}}(Y,X)$, such that $s(\mathrm{incl}) = \id$ and $s(f)(x) = x$ for all $f \in \cV$ and $x \in X \smallsetminus V$.

\textbf{Step 1.} Let us write $\partial_\bullet L$ for the union of all boundary components of $L$ except for the one that intersects the image of $e_{2}$; see Figure~\ref{fig:decorated-manifolds-Mt}. Note that $\partial_\bullet L$ may or may not intersect the image of $e_{1}$. Then there is a canonical identification:
\begin{equation}
\label{eq:boundary-components-identification}
\pi_{0}(\partial(L\natural M)) \cong \pi_{0}(\partial_\bullet L) \sqcup \pi_{0}(\partial M).
\end{equation}
This is necessarily asymmetric in $L$ and $M$. Each embedding $f \in \embdece(M_t,L\natural M)_L$ extends to a diffeomorphism of $L \natural M$, so it induces an injection $f_\partial \colon \pi_{0}(\partial M) \to \pi_{0}(\partial (L\natural M))$. In particular, if $f$ is the inclusion, then $f_\partial$ is also the inclusion, under the identification \eqref{eq:boundary-components-identification}. In addition, we know that $f$ sends $B$ into $A \sqcup B$, so it also induces a map $f_\sharp \colon \pi_{0}(B) \to \pi_{0}(A) \sqcup \pi_{0}(B)$, which must be an injection since $A$ and $B$ are closed manifolds and $f$ is an embedding. The function $f \mapsto (f_\partial,f_\sharp)$ is locally constant, so its fibres are open. Let $\cU'$ be the open subset of $\embdece(M_t,L\natural M)_L$ consisting of all $f$ such that $f_\partial$ is the inclusion and $f_\sharp(\pi_{0}(B)) = \pi_{0}(B)$. Note that the second condition implies that $f(B) = B$, since $f$ is an embedding and $B$ is a closed manifold.

\textbf{Step 2.} Write $M_{\epsilon,t} = e_{1}(\bB_{\epsilon}^d) \sqcup M_t$ (pictured in Figure~\ref{fig:decorated-manifolds-Mepsilont}). Let
\[
\gamma' \colon \embdece(M_t,L\natural M)_L \too C^\infty(M_{\epsilon,t},L\natural M)
\]
be the continuous map that extends a given embedding $M_t \hookrightarrow L \natural M$ to a smooth map $M_{\epsilon,t} \to L \natural M$ by defining it to be the identity on $e_{1}(\bB_{\epsilon}^d)$. Note that this may fail to be injective, so it is just a smooth map, not necessarily an embedding. Also observe that, if $f$ lies in the open subset $\cU'$ from Step 1, then $\gamma'(f)$ lies in the subspace $C_{\mathrm{face}}^\infty(M_{\epsilon,t},L\natural M)$, since it takes points of $\Int{(L\natural M)} \cap M_{\epsilon,t}$ into the interior of $L\natural M$ and, for any boundary component $P$ of $L\natural M$, it takes $P \cap M_{\epsilon,t}$ into $P$ (this uses the fact that $f_\partial = \id$). Restricting $\gamma'$ to $\cU'$, we therefore have a continuous map
\[
\gamma' \colon \cU' \too C^{\infty}_{\mathrm{face}}(M_{\epsilon,t}, L \natural M)
\]
such that $\gamma'(\mathrm{incl}) = \mathrm{incl}$ and $\gamma'(f)|_{M_t} = f$ for all $f \in \cU'$.

\textbf{Step 3.} Now set $X = L \natural M$ and $Y = M_{\epsilon,t}$ in the Extension Lemma of Cerf above, and choose $V$ to be any open neighbourhood of $M_{\epsilon,t}$ in $L \natural M$ that is disjoint from the submanifold $A \subset \Int{L}$. Composing the local section $s$ obtained from the Extension Lemma~with $\gamma'$, we have a continuous map
\[
\gamma'' = s \circ \gamma' \colon \cU'' = (\gamma')^{-1}(\cV) \too C^{\infty}_{\mathrm{face}}(L \natural M,L \natural M)
\]
such that $\gamma''(\mathrm{incl}) = \id$ and for any $f \in \cU''$ we have $\gamma''(f)|_{M_t} = f$ and $\gamma''(f)(A) = A$. Moreover, by construction, we also know that $\gamma''(f)(B) = B$ and $\gamma''(f)(x) = x$ for all $x \in e_{1}(\bB_{\epsilon}^d) \sqcup e'_{2}(\bB_{\epsilon}^d)$.

\textbf{Step 4.} Finally, note that $\Diff(L \natural M)$ is open in $C^\infty(L\natural M,L\natural M)$, so
\[
\cU = (\gamma'')^{-1}(C^{\infty}_{\mathrm{face}}(L \natural M,L \natural M) \cap \Diff(L\natural M))
\]
is an open neighbourhood of the inclusion in $\embdece(M_t,L\natural M)_L$. For each $f \in \cU$, the diffeomorphism $\gamma''(f)$ of $L \natural M$ fixes each point of $e_{1}(\bB_{\epsilon}^d) \sqcup e'_{2}(\bB_{\epsilon}^d)$ and sends $A \sqcup B$ onto itself, so it is an element of $\diffdece(L\natural M)$. So we have a continuous map
\[
\gamma = \gamma''|_{\cU} \colon \cU \too \diffdece(L\natural M)
\]
such that $\gamma(\mathrm{incl}) = \id$ and, for all $f \in \cU$, we have $\gamma(f) \circ \mathrm{incl} = \gamma(f)|_{M_t} = f$.

\textbf{Summary.} The $4$-step construction above may be summarised in the following commutative diagram:
\begin{center}
\begin{tikzpicture}
[x=1mm,y=1mm]
\node (tl) at (0,20) {$\embdece(M_t,L\natural M)$};
\node (tm) at (40,20) {$C^\infty(M_{\epsilon,t},L\natural M)$};
\node (ml) at (0,10) {$\cU'$};
\node (mm) at (40,10) {$C_{\mathrm{face}}^\infty(M_{\epsilon,t},L\natural M)$};
\node (mr) at (80,10) {$C_{\mathrm{face}}^\infty(L\natural M,L\natural M)$};
\node (bl) at (0,0) {$\cU''$};
\node (bm) at (40,0) {$\cV$};
\node (bbl) at (0,-10) {$\cU$};
\node (bbr) at (80,-10) {$C_{\mathrm{face}}^\infty(L\natural M,L\natural M) \cap \Diff(L\natural M)$};
\node (bbrr) at (120,-10) {$\Diff(L\natural M),$};
\draw[->] (tl) to node[above,font=\small]{$\gamma'$} (tm);
\draw[->] (ml) to (mm);
\draw[->] (mr) to (mm);
\draw[->] (bl) to (bm);
\draw[->] (bm) to[out=0,in=210] node[above,font=\small]{$s$} (mr);
\draw[->] (bbl) to node[above,font=\small]{$\gamma$} (bbr);
\node at ($ (bbr.east)!0.5!(bbrr.west) $) {$\subseteq$};
\incl{(bbr)}{(mr)}
\node at (0,-5) {\rotatebox{90}{$\subseteq$}};
\node at (0,5) {\rotatebox{90}{$\subseteq$}};
\node at (0,15) {\rotatebox{90}{$\subseteq$}};
\node at (40,5) {\rotatebox{90}{$\subseteq$}};
\node at (40,15) {\rotatebox{90}{$\subseteq$}};
\end{tikzpicture}
\end{center}
where the construction of $\gamma$ ensures that its image lies in $\diffdece(L\natural M) \subseteq \Diff(L\natural M)$.
\end{proof}

\begin{figure}[t]
\centering
\begin{tikzpicture}
[x=1mm,y=1mm,scale=1]
%
\begin{scope}[xshift=5mm]
\fill[green!20] (0,0) -- (7,0) -- (7,10) -- (0,10) -- cycle;
\end{scope}
\begin{scope}[xshift=55mm]
\fill[green!20] (-4,0) -- (10,0) .. controls (15,0) and (15,-5) .. (20,-5) -- (21,-5) -- (29,-5) -- (30,-5) .. controls (35,-5) and (35,0) .. (40,0) -- (50,0) -- (50,10) -- (40,10) .. controls (35,10) and (35,15) .. (30,15) -- (29,15) -- (21,15) -- (20,15) .. controls (15,15) and (15,10) .. (10,10) -- (-4,10) -- cycle;
\end{scope}
\begin{scope}[xshift=5mm]
\draw[thick] (0,0) -- (10,0) .. controls (15,0) and (15,-5) .. (20,-5) -- (21,-5);
\draw[thick,densely dashed] (21,-5) -- (29,-5);
\draw[thick] (29,-5) -- (30,-5) .. controls (35,-5) and (35,0) .. (40,0) -- (50,0);
\draw[thick] (50,10) -- (40,10) .. controls (35,10) and (35,15) .. (30,15) -- (29,15);
\draw[thick,densely dashed] (29,15) -- (21,15);
\draw[thick] (21,15) -- (20,15) .. controls (15,15) and (15,10) .. (10,10) -- (0,10) -- (0,0);
\draw[black!20,thick] (18,2) circle (2);
\draw[black!20,thick] (28,7) circle (2);
\fill[black!20,thick] (22,9) circle (0.5);
\node at (34,5) [black!20] {$A$};
\fill[white] (29,1) circle (1.5);
\draw[thick] (29,1) circle (1.5);
\draw (7,0) -- (7,10);
\draw (10,0) -- (10,10);
\node at (0,0) [anchor=north,font=\footnotesize] {$0$};
\node at (7,0) [anchor=north,font=\footnotesize] {$\vphantom{0}\epsilon$};
\node at (10,0) [anchor=north,font=\footnotesize] {$1$};
\node at (5,12) [anchor=south] {$e_{1}$};
\draw[decorate,decoration={brace,amplitude=5pt,raise=2pt}] (0,10) -- (10,10);
\node at (0,5) [anchor=east] {$L \natural M =$};
\end{scope}
\begin{scope}[xshift=55mm]
\draw[thick] (0,0) -- (10,0) .. controls (15,0) and (15,-5) .. (20,-5) -- (21,-5);
\draw[thick,densely dashed] (21,-5) -- (29,-5);
\draw[thick] (29,-5) -- (30,-5) .. controls (35,-5) and (35,0) .. (40,0) -- (50,0) -- (50,10) -- (40,10) .. controls (35,10) and (35,15) .. (30,15) -- (29,15);
\draw[thick,densely dashed] (29,15) -- (21,15);
\draw[thick] (21,15) -- (20,15) .. controls (15,15) and (15,10) .. (10,10) -- (0,10);
\draw[green!80!black,thick] (17,7) circle (2);
\draw[green!80!black,thick] (22,1) circle (2);
\draw[green!80!black,thick] (29,-1) circle (1.5);
\fill[green!80!black,thick] (29,5) circle (0.5);
\fill[white] (24,9.5) circle (2.5);
\draw[thick] (24,9.5) circle (2.5);
\node at (34,5) [green!80!black] {$B$};
\draw (40,0) -- (40,10);
\draw (43,0) -- (43,10);
\node at (40,0) [anchor=north,font=\footnotesize] {$1$};
\node at (43,0) [anchor=north,font=\footnotesize] {$\vphantom{0}\epsilon$};
\node at (50,0) [anchor=north,font=\footnotesize] {$0$};
\node at (45,12) [anchor=south] {$e'_{2}$};
\draw[decorate,decoration={brace,amplitude=5pt,raise=2pt}] (40,10) -- (50,10);
\draw (-10,0) -- (-10,10);
\draw (-4,0) -- (-4,10);
\draw (0,0) -- (0,10);
\draw (10,0) -- (10,10);
\node at (-10,0) [anchor=north,font=\footnotesize] {$1$};
\node at (-4,0) [anchor=north,font=\footnotesize] {$\vphantom{0}t$};
\node at (0,0) [anchor=north,font=\footnotesize] {$0$};
\node at (10,0) [anchor=north,font=\footnotesize] {$1$};
\draw[decorate,decoration={brace,amplitude=5pt,raise=2pt}] (-10,10) -- (0,10);
\draw[decorate,decoration={brace,amplitude=5pt,raise=2pt}] (0,10) -- (10,10);
\node at (-5,12) [anchor=south] {$e_{2}$};
\node at (5,12) [anchor=south] {$e'_{1}$};
\end{scope}
\end{tikzpicture}
\caption{The submanifold $M_{\epsilon,t}$ (shaded in green) of $L \natural M$ from the proof of Theorem~\ref{thm:fibre-bundle}.}
\label{fig:decorated-manifolds-Mepsilont}
\end{figure}

\subsection{The topologically-enriched Quillen bracket construction}\label{ss:Quillen_bracket_construction}

This section deals with the \emph{bracket construction} due to Quillen, its generalisations (see \S\ref{sss:Topological-Quillen}) and its applications for subcategories of the groupoids of decorated manifolds (see \S\ref{sss:quillen-bracket-categories}). The various instances of the Quillen bracket construction will be essential for the construction of the homological representation functors in \S\ref{s:general_construction}. First, we recall the definition of the original version of Quillen's bracket construction, which is a particular case of a more general construction described in \cite[p.219]{graysonQuillen}; see also \cite[\S 1.1]{RWW}.

\begin{defi}[Original Quillen bracket construction.]
\label{def:non-enriched-Quillen}
We fix a (discrete) monoidal groupoid $(\cG_{\circ},\natural,\zero)$ and a (discrete) left module $(\cM_{\circ},\natural)$ over $\cG_{\circ}$.
The \emph{Quillen bracket construction} $\langle \cG_{\circ},\cM_{\circ}\rangle $ is the category with the same objects as $\cM_{\circ}$ and whose morphisms are given by the colimit
\begin{equation}
\label{eq:Quillen-bracket-colimit}
\langle \cG_{\circ},\cM_{\circ}\rangle(X,Y)=\underset{\cG_{\circ}}{\colim}[\cM_{\circ}(-\natural X,Y)].
\end{equation}
\end{defi}

\begin{eg}\label{eg:UG}
As an example, a monoidal groupoid $\cG_{\circ}$ has a canonical left action on itself given by its monoidal structure, so we may always take $\cM_{\circ} = \cG_{\circ}$ and consider the category $\langle \cG_{\circ},\cG_{\circ} \rangle$. As an abbreviation, we denote the category $\langle \cG_{\circ},\cG_{\circ} \rangle$ by $\fU\cG_{\circ}$.
\end{eg}

\begin{notation}
\label{not:Quillen-bracket-construction}
By \eqref{eq:Quillen-bracket-colimit}, a morphism in $\langle \cG_{\circ},\cM_{\circ}\rangle$ from $X$ to $Y$ corresponds to the equivalence class of a morphism $\varphi \in \cM_{\circ}(A \natural X,Y)$ for an object $A$ of $\cG_{\circ}$, which we denote by $[A,\varphi]$.
\end{notation}

There is a canonical faithful functor $c_{\langle \cG_{\circ},\cM_{\circ}\rangle }\colon\cM_{\circ}\hookrightarrow\langle \cG_{\circ},\cM_{\circ}\rangle $ defined as the identity on objects and sending each morphism $\varphi$ of $\cM_{\circ}(X,Y)$ to $[\zero,\varphi]$.
Assuming in addition that $\cM_{\circ}$ is a groupoid (as for all the examples discussed in this paper, see \S\ref{ss:categories_for_families_of_groups}), following mutatis mutandis \cite[Prop.~1.7]{RWW}, if $(\cG_{\circ},\natural,\zero)$ has \emph{no zero divisors} -- meaning that $A\natural B\cong \zero$ if and only if $A\cong B\cong \zero$ for all objects $A$ and $B$ of $\cG_{\circ}$ -- and if $\Aut_{\cG_{\circ}}(\zero)=\{\id_{\zero}\}$, then the above canonical functor $c_{\langle \cG_{\circ},\cM_{\circ}\rangle }$ is an isomorphism from $\cM_{\circ}$ onto the maximal subgroupoid of $\langle \cG_{\circ},\cM_{\circ} \rangle$ (i.e.~the subcategory with the same objects as $\langle \cG_{\circ},\cM_{\circ}\rangle$ and whose morphisms are the isomorphisms of $\langle \cG_{\circ},\cM_{\circ}\rangle$).

\subsubsection{A topological enrichment of the Quillen bracket construction}\label{sss:Topological-Quillen}

We now explain how to generalise the Quillen bracket construction to \emph{semicategories} and \emph{topologically-enriched semicategories}.

\begin{prop}
\label{prop:topological-Quillen}
Let $(\cG,\natural)$ be a semi-monoidal groupoid and let $\cM$ be a category with a left action of $\cG$ also denoted by $\natural$. There is a semicategory $\langle \cG,\cM \rangle$ with the same objects as $\cM$, given by assigning $\langle \cG,\cM \rangle(X,Y)$, for objects $X,Y$ of $\cM$, to be the quotient set
\[
\Biggl( \bigsqcup_{\obj(\cG)} \cM(- \natural X,Y) \Biggr) / {\sim},
\]
where $\sim$ is the equivalence relation given by $(A,\varphi) \sim (A',\varphi')$ if and only if $\varphi = \varphi' \circ (\sigma \natural \id_{X})$ for some $\sigma \in \cG(A,A')$. For two morphisms $[A,\varphi] \colon X\rightarrow Y$ and $[B,\psi] \colon Y\rightarrow Z$ in $\langle \cG,\cM\rangle $, their composition is defined to be $[B,\psi]\circ[A,\varphi]=[B\natural A,\psi\circ(\id_{B}\natural\varphi)]$.

Moreover, if the categories $\cG$ and $\cM$ are topologically-enriched and the left-action of $\cG$ on $\cM$ is continuous, then performing the same constructions on spaces rather than sets produces a topologically-enriched semicategory $\langle \cG,\cM \rangle$.
\end{prop}
\begin{proof}
It straightforwardly follows from the assignments that $\sim$ is an equivalence relation, and so $\langle \cG,\cM \rangle(X,Y)$ is well-defined.
Hence the associativity of the composition is the only point to check in order to prove that we define a semicategory $\langle \cG,\cM \rangle$.
Let $[A,\varphi]:X\rightarrow Y$, $[B,\psi] \colon Y\rightarrow Z$ and $[C,\chi] \colon Z\rightarrow W$ be morphisms in $\langle \cG,\cM\rangle $. Recall that the associator of $\cG$ provides an isomorphism $\alpha_{C,B,A}\colon C\natural(B\natural A)\cong (C\natural B)\natural A$. Then the functoriality of the left action $\cG \times \cM \to \cM$ implies that $\id_{C}\natural(\psi\circ(\id_{B}\natural\varphi))=
(\id_{C}\natural\psi)\circ(\id_{C}\natural(\id_{B}\natural\varphi))$. We also note from the naturality of the associator that $\id_{C}\natural(\id_{B}\natural\varphi) = (\alpha_{C,B,A}^{-1}\natural\id_{X}) \circ(\id_{C\natural B}\natural\varphi)\circ (\alpha_{C,B,A}\natural\id_{X})$.
Therefore, we deduce from the associativity of the composition in $\cM$ that
\[
\chi\circ(\id_{C}\natural(\psi\circ(\id_{B}\natural\varphi))) = ((\chi\circ(\id_{C}\natural\psi))\circ(\alpha_{C,B,A}^{-1}\natural\id_{X} )\circ(\id_{C\natural B}\natural\varphi))\circ (\alpha_{C,B,A}\natural\id_{X})
\]
and so $[C\natural(B\natural A),\chi\circ(\id_{C}\natural(\psi\circ(\id_{B}\natural\varphi)))] = [(C\natural B)\natural A,(\chi\circ(\id_{C}\natural\psi))\circ(\alpha_{C,B,A}^{-1}\natural\id_{X} )\circ(\id_{C\natural B}\natural\varphi)]$.
Hence, it follows from the assignment for the composition in $\langle \cG,\cM \rangle$ that $[C,\chi]\circ([B,\psi]\circ[A,\varphi])=([C,\chi]\circ[B,\psi])\circ[A,\varphi]$, which proves associativity of the composition.

In the topologically-enriched setting, the definition of the hom-set $\langle \cG,\cM \rangle(X,Y)$ in the statement of the lemma equips it with an evident topology, namely the quotient topology induced from the topological disjoint union of the hom-spaces $\cM(A \natural X,Y)$ of the topologically-enriched category $\cM$. Continuity of the composition with respect to this topology straightforwardly follows from that of $\cM$, since morphisms of $\langle \cG,\cM \rangle$ are equivalence classes of morphisms of $\cM$.
\end{proof}

\begin{rmk}
\label{rmk:open-path-components-Quillen-bracket-construction}
If the morphism spaces of $\cM$ have all of their path-components open, then the same holds for $\langle \cG , \cM \rangle$. This follows from Lemma~\ref{lem:open-path-components}, since the morphism spaces of $\langle \cG , \cM \rangle$ are constructed as quotients of disjoint unions of morphism spaces of $\cM$.
\end{rmk}

When $\cG$ and $\cM$ are both \emph{discrete} categories with a \emph{genuine} monoidal structure and left $\natural$-module structure respectively, then the version of the Quillen bracket construction defined by Proposition~\ref{prop:topological-Quillen} recovers that of Definition~\ref{def:non-enriched-Quillen}.
Just as for the non-enriched case, there is a canonical continuous faithful functor $c_{\langle \cG,\cM\rangle }\colon\cM\hookrightarrow\langle \cG,\cM\rangle$ defined as the identity on objects and sending $\phi\in \cM(X,Y)$ to $[\zero,\phi]$. Furthermore, Example~\ref{eg:UG} (i.e.~taking $\cM = \cG$ and defining the category $\fU\cG:=\langle \cG,\cG \rangle$) repeats verbatim for $\cG$ a \emph{topologically-enriched} monoidal groupoid, since the canonical left action of $\cG$ on itself given by its monoidal structure is continuous.

\begin{rmk}\label{rmk:small_categories_Quillen}
From now on, we assume that all the categories we consider are \emph{small}.
The point is that the morphisms in the (topologically-enriched) Quillen bracket construction quantify over all objects of the groupoid that we start with, so if that groupoid is not small, then the Quillen bracket construction will fail even to be locally small. Hence it is important that we do this \emph{before} applying the Quillen bracket construction, so that the categories arising from this construction are small.
In particular, following Convention~\ref{convention:small_categories}, we stress here that we may assume that the decorated manifold category $\cD\mathrm{ec}_d$ (see Definition~\ref{def:Dd}), and hence all of its subcategories that we consider, are small. These are the categories to which we will apply the Quillen bracket construction; see \S\ref{sss:quillen-bracket-categories}.
\end{rmk}

\begin{rmk}
A topological version of Quillen's bracket construction is mentioned briefly in \cite[Rem.~2.10]{Krannich2017Homologicalstabilitytopological}, although there the categories are \emph{topological} in the sense of being categories internal to the category of topological spaces, rather than topologically-enriched categories. Proposition~\ref{prop:topological-Quillen} is stated and proved for topologically-enriched categories, but it is likely that it has an analogue for categories internal to the category of topological spaces, in which case \cite[Lem.~2.10]{Krannich2017Homologicalstabilitytopological} would be a particular case of this analogue.
\end{rmk}

The construction defined from Proposition~\ref{prop:topological-Quillen} is clearly \emph{functorial} in $\cG$ and $\cM$ in an appropriate sense. We detail here some properties of this functoriality that we will need.

\begin{lemm}
\label{lem:restriction}
Let $(\cD,\natural)$ be a topologically-enriched semi-monoidal groupoid and let $\cG_{1} \subseteq \cG_{2} \subseteq \cD$ and $\cM_{1} \subseteq \cM_{2} \subseteq \cD$ be subgroupoids such that, for $i \in \{1,2\}$, $\cG_{i}$ is closed under $- \natural -$ and $\cM_{i}$ is closed under $g \natural -$ for each object $g$ of $\cG_{i}$. Then there is a canonical continuous semifunctor $\langle \cG_{1} , \cM_{1} \rangle \to \langle \cG_{2} , \cM_{2} \rangle .$
Moreover,
\begin{itemizeb}
\item if $\cG_{1} = \cG_{2} = \cG$, the semifunctor $\langle \cG , \cM_{1} \rangle \to \langle \cG , \cM_{2} \rangle$ is an inclusion of a subsemicategory, which is full if the inclusion $\cM_{1} \subseteq \cM_{2}$ is full;
\item if $\cM_{1} = \cM_{2} = \cM$, the semifunctor $\langle \cG_{1} , \cM \rangle \to \langle \cG_{2} , \cM \rangle$ is the identity on objects, and is faithful -- thus an inclusion of a subsemicategory -- if the inclusion $\cG_{1} \subseteq \cG_{2}$ is full.
\end{itemizeb}
In particular, there is a canonical semifunctor $\langle \cG_{1} , \cM_{1} \rangle \to \fU\cD$, which is an inclusion of a subsemicategory if the inclusion $\cG_{1} \subseteq \cD$ is full.

All these results repeat mutatis mutandis in the unital setting, considering genuine (i.e.~unital) monoidal $\natural$ and left $\natural$-module structures, so that the Quillen bracket construction defines categories, and then functors between these categories.
\end{lemm}
\begin{proof}
The objects of $\langle \cG_{i},\cM_{i} \rangle$ are the objects of $\cM_{i}$, so we define the semifunctor on objects as the class inclusion $\obj(\cM_{1}) \hookrightarrow \obj(\cM_{2})$. A morphism in $\langle \cG_{i},\cM_{i} \rangle$ from $X$ to $Y$ is represented by a choice of object $A$ of $\cG_{i}$ and a morphism $A \natural X \to Y$ of $\cM_{i}$. We may therefore send such a morphism, for $i=1$, to the morphism, for $i=2$, represented by the same data, since $\cG_{1} \subseteq \cG_{2}$ and $\cM_{1} \subseteq \cM_{2}$. It is straightforward to check that this assignment respects the defining equivalence relation, so induces a continuous map of morphism spaces, and that it also respects composition (and identities, if we are in the unital setting). The statements in the two bullet points may straightforwardly be verified by unwinding the definition of morphisms in $\langle \cG_{i},\cM_{i} \rangle$. In particular, that the semifunctor $\langle \cG_{1} , \cM \rangle \to \langle \cG_{2} , \cM \rangle$ is faithful when the inclusion $\cG_{1} \subseteq \cG_{2}$ is full is proven as follows. By definition (see Proposition~\ref{prop:topological-Quillen}), for two morphisms $[A,\varphi]$ and $[A',\varphi']$ of $\langle \cG_{1} , \cM \rangle(X,Y)$, their images under $\langle \cG_{1} , \cM \rangle \to \langle \cG_{2} , \cM \rangle$ are equal if and only if there exists $\sigma \in \cG_{2}(A,A')$ such that $\varphi = \varphi' \circ (\sigma \natural \id_{X})$. Since $\cG_{1} \subseteq \cG_{2}$ is full, there exists $\tilde{\sigma} \in \cG_{1}(A,A')$ such that $\varphi = \varphi' \circ (\tilde{\sigma} \natural \id_{X})$, and hence, by definition, $[A,\varphi]=[A',\varphi']$.
\end{proof}

Finally, we study the interaction between the topologically-enriched Quillen bracket construction and the path-component functor $\pi_{0}$ in Proposition~\ref{prop:topological-Quillen_morphisms}. As an ingredient for its proof, we need the following Lemma~\ref{lem:extract_proof_topological_Quillen}, which will also be used later in the proof of Lemma~\ref{lem:Serre-fibration-condition}.

Let $(\cG,\natural)$ be a topologically-enriched semi-monoidal groupoid and $\cM$ a topologically-enriched category with a continuous left action of $\cG$ also denoted by $\natural$. Since $\cG$ is a groupoid, setting that two objects $A,A'$ of $\cG$ are equivalent if and only if there is a morphism $A \to A'$ in $\cG$ defines an equivalence relation; we generically denote each one of its associated equivalence classes by $\cO_{\alpha}$ for some index $\alpha$. In particular, the equivalence classes $\{\cO_{\alpha}\}_{\alpha}$ partition $\obj(\cG)$. Set
\[
\Phi \;=\!\! \bigsqcup_{A \in \obj(\cG)} \cM(A \natural X,Y) \,\,\text{ and }\,\, \Phi_{\alpha} \;=\! \bigsqcup_{A \in \cO_{\alpha}} \cM(A\natural X,Y)
\]
for each $\alpha$ indexing an equivalence class of $\obj(\cG)$. It is straightforward to check that we define an equivalence relation $\sim_{t}$ on $\Phi$ by assigning $(A,\varphi) \sim_{t} (A',\varphi')$ if and only if there is a morphism $\sigma \in \cG(A,A')$ such that $\varphi = \varphi' \circ (\sigma\natural\id_{A})$. We denote by $\varrho \colon \Phi \twoheadrightarrow \Phi / {\sim_{t}}$ the associated quotient map. 

\begin{lemm}\label{lem:extract_proof_topological_Quillen}
For each pair of objects $X,Y$ of $\cM$ and picking an object $A_{\alpha} \in \cO_{\alpha}$ for each $\alpha$, there is a commutative diagram
\begin{equation}\label{eq:diagram_decomposition_disjoint_union}
\begin{tikzcd}
\coprod_{\alpha}\Phi_{\alpha} \ar[rr,"\cong"] \ar[d,two heads] && \Phi \ar[d,two heads,"\varrho"] \\
\coprod_{\alpha}\cM(A_{\alpha} \natural X,Y)/\Aut_{\cG}(A_{\alpha}) \ar[rr,"\cong"] && \Phi/{\sim_{t}}
\end{tikzcd}
\end{equation}
where the horizontal maps are homeomorphisms.
\end{lemm}
\begin{proof}
First, the equivalence classes $\cO_{\alpha}$ provide a partition of $\obj(\cG)$, which defines the top horizontal homeomorphism $\bigsqcup_{\alpha} \Phi_{\alpha}\cong \Phi$ of \eqref{eq:diagram_decomposition_disjoint_union}.
The equivalence relation $\sim_{t}$ on $\Phi$ clearly preserves the topological disjoint union $\bigsqcup_{\alpha} \Phi_{\alpha}$, so this homeomorphism descends to $\bigsqcup_{\alpha} (\Phi_{\alpha} / {\sim_{t}}) \cong \Phi/{\sim_{t}}$.
Moreover, note that, for a fixed $\alpha$ and for any two objects $A,A' \in \cO_{\alpha}$, we have $\varrho(\cM(A\natural X,Y))=\varrho(\cM(A'\natural X,Y))$.
Hence we have a decomposition of $\Phi/{\sim_{t}}$ as a topological disjoint union:
\begin{equation}\label{eq:decomposition_disjoint_union_Phi}
\Phi/{\sim_{t}} \;\cong\; \bigsqcup_{\alpha} \; \varrho(\cM(A_{\alpha} \natural X,Y)).
\end{equation}
Thus we have established the analogue of the commutative diagram \eqref{eq:diagram_decomposition_disjoint_union} where the bottom-left corner is $\bigsqcup_{\alpha} \varrho(\cM(A_{\alpha} \natural X,Y))$.

For each $\alpha$ and $A \in \cO_{\alpha}$, we note that two elements $\varphi,\varphi' \in \cM(A\natural X,Y)$ have the same image under $\varrho$ if and only if they are $\sim_{t}$-equivalent, which is equivalent to saying that they lie in the same orbit of the $\Aut_{\cG}(A)$-action on $\cM(A\natural X,Y)$. Hence the restriction
\begin{equation}\label{eq:restriction_of_q}
\varrho_{A} \colon \cM(A \natural X,Y) \too \varrho(\cM(A\natural X,Y))
\end{equation}
of $\varrho$ is isomorphic to
\begin{equation}\label{eq:Hom_Aut_surjection}
\cM(A \natural X,Y) \too \cM(A \natural X,Y) / \Aut_{\cG}(A),
\end{equation}
at least on underlying sets. It now just remains to prove that \eqref{eq:restriction_of_q} and \eqref{eq:Hom_Aut_surjection} are isomorphic also as continuous maps of spaces. Since \eqref{eq:restriction_of_q} and \eqref{eq:Hom_Aut_surjection} are surjective continuous maps with the same domain and the same point-fibres, and we know moreover that \eqref{eq:Hom_Aut_surjection} is a quotient map, it suffices to prove that \eqref{eq:restriction_of_q} is also a quotient map.

Let $U \subseteq \varrho(\cM(A\natural X,Y))$ be a subset such that $\varrho_{A}^{-1}(U)$ is open in $\cM(A\natural X,Y)$. We need to show that $U$ is open in $\varrho(\cM(A\natural X,Y))$. The fact that the equivalence relation $\sim_t$ preserves the decomposition of $\Phi$ into $\Phi_{\alpha}$ implies that the restriction
\[
\varrho_{\alpha} = \varrho|_{\Phi_{\alpha}} \colon \Phi_{\alpha} \too \varrho(\Phi_{\alpha}) = \varrho(\cM(A\natural X,Y))
\]
is a quotient map. So it suffices to show that $\varrho_{\alpha}^{-1}(U)$ is open in $\Phi_{\alpha}$. Now, from the definitions, we observe the following description of the subset
\[
\varrho_{\alpha}^{-1}(U) \;\subseteq\; \bigsqcup_{A' \in \cO_{\alpha}} \cM(A'\natural X,Y).
\]
For each object $A' \in \cO_{\alpha}$, choose an isomorphism $\sigma_{A'} \colon A' \to A$ in $\cG$. This induces a homeomorphism
\[
\Upsilon_{A'} \;=\; -\circ (\sigma_{A'} \natural \id) \colon \cM(A\natural X,Y) \too \cM(A'\natural X,Y).
\]
Then we have
\[
\varrho_{\alpha}^{-1}(U) \;=\! \bigsqcup_{A' \in \cO_{\alpha}} \Upsilon_{A'}(\varrho_A^{-1}(U)).
\]
Since $\varrho_A^{-1}(U)$ is open in $\cM(A\natural X,Y)$, it follows that $\Upsilon_{A'}(\varrho_A^{-1}(U))$ is open in $\cM(A'\natural X,Y)$ for each $A' \in \cO_{\alpha}$. Thus $\varrho_{\alpha}^{-1}(U)$ is open in $\Phi_{\alpha}$, as required.
\end{proof}

\begin{prop}\label{prop:topological-Quillen_morphisms}
Let $(\cG,\natural)$ be a topologically-enriched semi-monoidal groupoid. Let $\cM$ be a topologically-enriched category with a continuous left action of $\cG$ also denoted by $\natural$. Assume that, for each object $A$ of $\cG$ and each pair of objects $X,Y$ of $\cM$, the quotient map
\begin{equation}\label{eq:Quillen-bracket-lemma-condition}
\cM(A \natural X,Y) \too \cM(A \natural X,Y) / \Aut_{\cG}(A)
\end{equation}
is a Serre fibration. Then there is a canonical isomorphism of semicategories
\begin{equation}\label{eq:Quillen-bracket-lemma-conclusion}
\pi_{0}(\langle \cG,\cM \rangle) \cong \langle \pi_{0}(\cG),\pi_{0}(\cM) \rangle .
\end{equation}
Moreover, if the semi-monoidal structure induced from $(\cG,\natural)$ on the groupoid $\pi_{0}(\cG)$ admits a unit (and thus upgrades to a genuine monoidal structure) that acts by the identity on $\pi_{0}(\cM)$ (so that $\pi_{0}(\cM)$ is a genuine left $\natural$-module over $\pi_{0}(\cG)$), then the isomorphism \eqref{eq:Quillen-bracket-lemma-conclusion} upgrades to an isomorphism of categories.
\end{prop}
\begin{proof}
First note that $\pi_{0}(\langle \cG,\cM \rangle)$ and $\langle \pi_{0}(\cG),\pi_{0}(\cM) \rangle$ have the same class of objects, by the definition of the discrete and topologically-enriched Quillen bracket constructions, and the functor $\pi_{0}$. Specifically, their common class of objects is $\obj(\cM)$. It therefore remains to show that, for objects $X$ and $Y$ of $\cM$, there is a natural bijection between $\pi_{0}(\langle \cG,\cM \rangle(X,Y))$ and $\langle \pi_{0}(\cG),\pi_{0}(\cM) \rangle(X,Y)$. We shall use the above notation $\Phi = \bigsqcup_{A \in \obj(\cG)} \cM(A\natural X,Y)$ and the equivalence relation $\sim_{t}$ from Lemma~\ref{lem:extract_proof_topological_Quillen}.
Unravelling the definitions, what we need to prove is that there is a natural bijection
\[
\pi_{0}(\Phi/{\sim_{t}}) \cong \pi_{0}(\Phi)/{\sim_h},
\]
where $\sim_h$ is the equivalence relation given by $(A,[\varphi]) \sim_h (A',[\varphi'])$ if and only if there is a morphism $\sigma \in \cG(A,A')$ such that $\varphi \simeq \varphi' \circ (\sigma\natural\id_{A})$. We note that the only difference between the definitions of $\sim_t$ and $\sim_h$ is that the equality is replaced by a homotopy in the definition of $\sim_h$.
As sets, these are both quotients of (the underlying set of) $\Phi$, so we just need to show that, given two elements $(A,\varphi)$ and $(A',\varphi')$ of $\Phi$, they have the same image in $\pi_{0}(\Phi/{\sim_{t}})$ if and only if they have the same image in $\pi_{0}(\Phi)/{\sim_h}$.

(a) Suppose first that $(A,\varphi)$ and $(A',\varphi')$ have the same image in $\pi_{0}(\Phi)/{\sim_h}$. This means that there is a morphism $\sigma \in \cG(A,A')$ and a path $\gamma \colon [0,1] \too \cM(A\natural X,Y) \subseteq \Phi$ with $\gamma(0) = (A,\varphi)$ and $\gamma(1) = (A,\varphi' \circ (\sigma\natural\id_{A}))$. Composing with the projection $\varrho \colon \Phi \twoheadrightarrow \Phi/{\sim_{t}}$ and writing $[-]_t$ for the equivalence classes with respect to $\sim_{t}$, we obtain a path in $\Phi/{\sim_{t}}$ from $[(A,\varphi)]_t$ to $[(A,\varphi' \circ (\sigma\natural\id_{A}))]_t = [(A',\varphi')]_t$. Hence $(A,\varphi)$ and $(A',\varphi')$ have the same image in $\pi_{0}(\Phi/{\sim_{t}})$.

(b) To prove the converse, we first show that the quotient map $\varrho \colon \Phi \twoheadrightarrow \Phi / {\sim_{t}}$ is a Serre fibration. Directly from the definition, one may straightforwardly verify the following two facts:
\begin{itemizeb}
\item $\bigsqcup_{i} f_{i} \colon \bigsqcup_{i} E_{i} \to B$ is a Serre fibration if and only if each $f_{i} \colon E_{i} \to B$ is a Serre fibration.
\item $f \colon E \to B$ is a Serre fibration if and only if $f(E)$ is a union of path-components of $B$ and $f \colon E \to f(E)$ is a Serre fibration.
\end{itemizeb}
Hence, since we know by hypothesis that \eqref{eq:Quillen-bracket-lemma-condition} is a Serre fibration, so is the left-hand vertical map of \eqref{eq:diagram_decomposition_disjoint_union}, and thus so is $\varrho \colon \Phi \twoheadrightarrow \Phi / {\sim_{t}}$ by Lemma~\ref{lem:extract_proof_topological_Quillen}.

Now assume that $(A,\varphi)$ and $(A',\varphi')$ have the same image in $\pi_{0}(\Phi/{\sim_{t}})$, so there is a path $\delta \colon [0,1] \to \Phi/{\sim_{t}}$ with $\delta(0) = [(A,\varphi)]_t$ and $\delta(1) = [(A',\varphi')]_t$. Since $\varrho$ is a Serre fibration, we may lift this to a path $\varepsilon \colon [0,1] \to \Phi$ with $\varepsilon(0) = (A,\varphi)$ and $\varepsilon(1) \sim_{t} (A',\varphi')$. Its image $\varepsilon([0,1])$ is path-connected, so it must lie in $\cM(A\natural X,Y) \subseteq \Phi$. Hence we have a path $\varepsilon \colon [0,1] \too \cM(A\natural X,Y)$ with $\varepsilon(0) = (A,\varphi)$ and $\varepsilon(1) = (A,\varphi'') \sim_{t} (A',\varphi')$, for some $\varphi'' \in \cM(A\natural X,Y)$. The relation $(A,\varphi'') \sim_{t} (A',\varphi')$ means that there is a morphism $\sigma \in \cG(A,A')$ such that $\varphi'' = \varphi' \circ (\sigma\natural\id_{A})$. Hence $\varepsilon$ is a homotopy witnessing that $\varphi \simeq \varphi' \circ (\sigma\natural\id_{A})$, so we have shown that $(A,[\varphi]) \sim_h (A',[\varphi'])$, in other words, $(A,\varphi)$ and $(A',\varphi')$ have the same image in $\pi_{0}(\Phi)/{\sim_h}$.

Finally, if the semi-monoidal structure induced from $(\cG,\natural)$ on $\pi_{0}(\cG)$ admits a unit, then $\langle \pi_{0}(\cG),\pi_{0}(\cM) \rangle$ is a genuine category via the original Quillen bracket construction of Definition~\ref{def:non-enriched-Quillen}. The isomorphism \eqref{eq:Quillen-bracket-lemma-conclusion} automatically preserves identities, since the image of each identity morphism of $\langle \pi_{0}(\cG),\pi_{0}(\cM) \rangle$ satisfies the relations of an identity and thus defines the identity morphism of its associated object, making $\pi_{0}(\langle \cG,\cM \rangle)$ a genuine category.
\end{proof}

\subsubsection{Quillen bracket categories of manifolds}
\label{sss:quillen-bracket-categories}

We now focus on applications of the Quillen bracket construction to (subgroupoids of) the groupoids of decorated manifolds $\cD_{d}$ and $\cD_{d}^{+}$, studying in particular their path components (see Lemma~\ref{lem:Serre-fibration-condition}) and describing their morphism spaces in terms of embedding spaces (see Proposition~\ref{prop:morphism-spaces-bracket}). The results that we shall need for our construction of homological representation functors are summarised in Corollary~\ref{cor:description_UD_{d}}. Since we make use of the results of \S\ref{ss:fibration-condition}, \textbf{we resume the assumption that the dimension $d$ is not equal to $4$}; cf.~Remark~\ref{rmk:d_neq_4_explanation}.
To begin with, we introduce the notion of \emph{$0$-full subcategory}, which is weaker than that of full subcategory:

\begin{defi}\label{def:zero-full}
An inclusion of topologically-enriched categories $\cC \subseteq \cD$ is called \emph{$0$-full} if, for each pair of objects $c,c'$ of $\cC$, the subspace $\cC(c,c') \subseteq \cD(c,c')$ is a union of path-components.
\end{defi}

\begin{rmk}
We recall that a \emph{full} inclusion of topologically-enriched categories $\cC \subseteq \cD$ is one where, for each pair of objects $c,c'$ of $\cC$, the inclusion $\cC(c,c') \subseteq \cD(c,c')$ is an equality.
We also note that a $0$-full inclusion $\cC \subseteq \cD$ induces an inclusion $\pi_{0}(\cC) \subseteq \pi_{0}(\cD)$. Indeed, subcategories of $\pi_{0}(\cD)$ correspond bijectively to $0$-full subcategories of $\cD$.
\end{rmk}

For the remainder of this section, we place ourselves in the following standard framework that will be further used at many points in the paper:

\begin{hypothesis}\label{hypo:standard_framework}
We fix an integer $d \geq 2$. Let $\cG$ and $\cM$ be subgroupoids of either $\cD_{d}$ or $\cD_{d}^{+}$ (see \S\ref{ss:topological_groupoids_of_decorated_manifolds}) such that $\cG$ is full and closed under $\natural$, while $\cM$ is $0$-full and closed under the (continuous) left action of $\cG$ via $\natural$.
\end{hypothesis}

Hypothesis~\ref{hypo:standard_framework} implies that $\cG$ is a semi-monoidal groupoid with a left action on $\cM$; we may thus form the Quillen bracket construction $\langle \cG , \cM \rangle$, which is a topologically-enriched semicategory; see \S\ref{sss:Topological-Quillen}. It is called a \emph{Quillen bracket category of manifolds}.

\begin{rmk}
\label{rmk:framework-open-path-components}
If $\cG$ and $\cM$ are groupoids as in Hypothesis~\ref{hypo:standard_framework}, then the morphism spaces of the Quillen bracket category of manifolds $\langle \cG , \cM \rangle$ have all of their path-components open, i.e.~the projection $\langle \cG , \cM \rangle \to \pi_{0}(\langle \cG , \cM \rangle)$ is continuous. Indeed, we ensured that the morphism spaces of $\cD_{d}$ and $\cD_{d}^{+}$ have all of their path-components open in Definition~\ref{def:decorated-manifolds-morphism-spaces}, so it follows that the same is true for $\cM$ since it is $0$-full, and hence also for $\langle \cG , \cM \rangle$ by Remark~\ref{rmk:open-path-components-Quillen-bracket-construction}.
\end{rmk}

\paragraph*{Path-components.}

The following result shows that the topologically-enriched Quillen bracket construction commutes with the path-component functor $\pi_{0}$ under the above framework.

\begin{lemm}\label{lem:Serre-fibration-condition}
The Serre fibration condition of Proposition~\ref{prop:topological-Quillen_morphisms} is satisfied for $\cG$ and $\cM$ as in Hypothesis~\ref{hypo:standard_framework}, and hence we have an isomorphism of semicategories
\[
\pi_{0}(\langle \cG , \cM \rangle) \cong \langle \pi_{0}(\cG) , \pi_{0}(\cM) \rangle ,
\]
which upgrades to an isomorphism of categories if the semi-monoidal structure induced from $(\cG,\natural)$ on $\pi_{0}(\cG)$ admits a unit that acts by the identity on $\pi_{0}(\cM)$.
\end{lemm}
\begin{proof}
Let $X,Y$ be objects of $\cM$ and let $A$ be an object of $\cG$; we must show that \eqref{eq:Quillen-bracket-lemma-condition} is a Serre fibration. If $A\natural X$ is not isomorphic to $Y$ in $\cM$ then it is the map $\emptyset \to \emptyset$, which is vacuously a Serre fibration. If $A\natural X$ is isomorphic to $Y$ in $\cM$ then we may use this isomorphism to replace $\cM(A\natural X,Y)$ with $\cM(A\natural X,A\natural X)$; thus we may assume that $Y = A\natural X$. The fact that \eqref{eq:Quillen-bracket-lemma-condition} is a Serre fibration in this case follows directly from Theorem~\ref{thm:fibre-bundle} if $\cM$ is a \emph{full} subgroupoid of $\cD_{d}^{(+)}$. If $\cM$ only satisfies the weaker property of being a \emph{$0$-full} subgroupoid of $\cD_{d}^{(+)}$, then it follows from Theorem~\ref{thm:fibre-bundle} together with Lemma~\ref{lem:Serre-fibration-path-components} below.
\end{proof}

\begin{lemm}\label{lem:Serre-fibration-path-components}
Let $X$ be a space with a continuous right action of a topological group $G$ such that the projection $X \to X/G$ is a Serre fibration. Let $X_{0} \subseteq X$ be a union of path-components such that the $G$-action sends $X_{0}$ into itself. Then the projection $X_{0} \to X_{0}/G$ is also a Serre fibration.
\end{lemm}
\begin{proof}
More generally, by considering lifting diagrams, one may prove that, in the following commutative square:
\begin{center}
\begin{tikzpicture}
[x=1mm,y=1mm]
\node (tl) at (0,12) {$A$};
\node (tr) at (30,12) {$B$};
\node (bl) at (0,0) {$C$};
\node (br) at (30,0) {$D$};
\draw[->] (tl) to node[above,font=\small]{$a$} (tr);
\draw[->] (bl) to node[below,font=\small]{$b$} (br);
\draw[->] (tl) to node[left,font=\small]{$g$} (bl);
\draw[->] (tr) to node[right,font=\small]{$f$} (br);
\end{tikzpicture}
\end{center}
if $f$ is a Serre fibration, $a$ is an inclusion of a union of path-components, and $b$ is injective, then $g$ is also a Serre fibration. In our setting, $a$ is the inclusion $X_{0} \hookrightarrow X$, which is assumed to be a union of path-components, and $b$ is the induced map $X_{0}/G \to X/G$, which is injective.
\end{proof}

\paragraph*{Morphism spaces.}

For decorated manifolds $M = (M,A)$ and $N = (N,B)$, recall from Notation~\ref{not:diffdec} and Definition~\ref{def:embdec} the topological group of decorated diffeomorphisms $\diffdec(M)$ and the space of decorated embeddings $\embdec(M,N)$ respectively, as well as their oriented versions $\diffdec^{+}(M)$ and $\embdec^{+}(M,N)$.

\begin{defi}[Decorated diffeomorphisms and embeddings, revisited.]
\label{def:diffdec-2}

If $M$ lies in a subgroupoid $\cH \subseteq \cD_{d}$, define $\Diff_{\cH}(M) \subseteq \diffdec(M)$ to be the subgroup of automorphisms of $M$ in $\cH$.

If $M$ and $N$ lie in $\cM \subseteq \cD_{d}$, define $\embgm(M,N) \subseteq \embdec(M,N)$ to be the subspace where, in the third condition of Definition~\ref{def:embdec}, the decorated manifold $M'$ lies in $\cG$ and the decorated diffeomorphism $\bar{\varphi}$ lies in $\cM$.
For an object $L$ of $\cG$, we also write $\embgm(M,N)_L$ for the subspace of $\embgm(M,N)$ where we may take $M' = L$ in Definition~\ref{def:embdec} (and $\bar{\varphi}$ lies in $\cM$).

For subgroupoids of $\cD_{d}^{+}$, we make the analogous definitions, with a superscript ${}^{+}$ in all of the notation.
\end{defi}

\begin{rmk}
\label{rmk:decorated-embeddings-gm-composition}
It is easy to check that composition of decorated embeddings $\embdec(-,-)$, as described in Definition~\ref{def:decorated-embeddings-composition}, restricts to a well-defined composition on $\embgm(-,-)$, using the facts that $\obj(\cG)$ is closed under $\natural$ and $\cM$ is closed under composition and the action of $\cG$ via $\natural$.
\end{rmk}

\begin{rmk}
We note that the inclusion $\Diff_{\cH}(M) \subseteq \diffdec(M)$ when $M$ lies in a subgroupoid $\cH \subseteq \cD_{d}$ is an equality if $\cH \subseteq \cD_{d}$ is full.
Also, the condition in the definition of $\embgm(M,N)$ that the decorated diffeomorphism $\bar{\varphi}$ lies in $\cM$ is automatic if $\cM \subseteq \cD_{d}$ is full.

Furthermore, similarly to the decomposition \eqref{eq:embdec-decomposition}, the space $\embgm(M,N)$ decomposes as a topological disjoint union:
\begin{equation}
\label{eq:embgm-decomposition}
\embgm(M,N) \;\cong\; \bigsqcup_L \, \embgm(M,N)_L,
\end{equation}
where the disjoint union runs over representatives of isomorphism classes of objects $L$ of $\cG$. This follows directly from Lemma~\ref{lem:decomposition-disjoint-union}, since $\embgm(M,N) \subseteq \embdec(M,N)$ is equipped with the subspace topology and $\embgm(M,N)_L = \embdec(M,N)_L \cap \embgm(M,N)$.
\end{rmk}

\begin{prop}
\label{prop:morphism-spaces-bracket}
Let $\cG$ and $\cM$ be subgroupoids of $\cD_{d}$ satisfying the conditions of Hypothesis~\ref{hypo:standard_framework}. Then the morphism spaces of the Quillen bracket category of manifolds $\langle \cG , \cM \rangle$ may be identified as follows:
\begin{equation}
\label{eq:morphism-spaces-bracket}
\langle \cG , \cM \rangle (M,N) \;\cong\; \embgm(M,N).
\end{equation}
Moreover, composition in the semicategory $\langle \cG , \cM \rangle$ on the left-hand side corresponds to composition of decorated embeddings on the right-hand side, as described in Definition~\ref{def:decorated-embeddings-composition} (see Remark~\ref{rmk:decorated-embeddings-gm-composition}).
The same also holds when $\cG$ and $\cM$ are subgroupoids of $\cD_{d}^{+}$.
\end{prop}

\begin{rmk}
When $\cM \subseteq \cD_{d}$ is full, this says that morphisms in $\langle \cG , \cM \rangle$ are embeddings of decorated manifolds such that the complement of the image of the embedding is an object of $\cG$. In particular, the morphism space $\langle \cG , \cM \rangle (M,N)$ is non-empty if and only if there exists an object $M'$ of $\cG$ such that $M' \natural M$ is diffeomorphic to $N$ as decorated manifolds.
\end{rmk}

\begin{proof}[Proof of Proposition~\ref{prop:morphism-spaces-bracket}]
The space $\langle \cG , \cM \rangle (M,N)$ is described in Proposition~\ref{prop:topological-Quillen}: using the notation of Lemma~\ref{lem:extract_proof_topological_Quillen}, it is the quotient space $\Phi/{\sim_{t}}$. Lemma~\ref{lem:extract_proof_topological_Quillen} then implies that we have a homeomorphism
\[
\langle \cG , \cM \rangle (M,N) \;=\; \Phi/{\sim_{t}} \;\cong\; \bigsqcup_L \, \cM(L\natural M,N) / \Aut_{\cG}(L),
\]
where the disjoint union runs over representatives $L$ of isomorphism classes of objects of $\cG$.
Since $\cM$ is a groupoid, the space $\cM(L\natural M,N) / \Aut_{\cG}(L)$ is empty unless $L \natural M$ is isomorphic to $N$ in $\cM$, in which case we may rewrite it as $\Aut_\cM(L\natural M) / \Aut_{\cG}(L) = \Diff_{\cM}(L\natural M) / \diffdec(L)$, using the notation of Definition~\ref{def:diffdec-2} (recall that $\cG \subseteq \cD_{d}$ is full). By Theorem~\ref{thm:fibre-bundle}, the restriction map induces a homeomorphism
\begin{equation}
\label{eq:restriction-map-identification}
\diffdec(L\natural M) / \diffdec(L) \;\cong\; \embdec(M,L\natural M)_L,
\end{equation}
and one may easily see that this sends the subspace $\Diff_{\cM}(L\natural M) / \diffdec(L)$ homeomorphically onto the subspace $\embgm(M,L\natural M)_L$. Putting this all together, we have a homeomorphism
\begin{equation}
\label{eq:gm-decomposition}
\langle \cG , \cM \rangle (M,N) \;\cong\; \bigsqcup_L \, \embgm(M,L\natural M)_L,
\end{equation}
where the disjoint union runs over representatives $L$ of isomorphism classes of objects of $\cG$ such that $L \natural M$ is isomorphic to $N$ in $\cM$.

Next, consider the topological decomposition \eqref{eq:embgm-decomposition}, where the disjoint union is indexed by representatives $L$ of isomorphism classes of all objects of $\cG$. If $L \natural M \not\cong N$ in $\cM$, the corresponding term is empty, whereas if $L \natural M \cong N$ in $\cM$, we may rewrite the corresponding term by replacing $N$ with $L \natural M$, to obtain:
\begin{equation}
\label{eq:embgm-decomposition-2}
\embgm(M,N) \;\cong\; \bigsqcup_L \, \embgm(M,L \natural M)_L,
\end{equation}
where the disjoint union now runs over representatives $L$ of isomorphism classes of objects of $\cG$ such that $L \natural M$ is isomorphic to $N$ in $\cM$. Combining \eqref{eq:gm-decomposition} and \eqref{eq:embgm-decomposition-2}, we obtain the desired identification \eqref{eq:morphism-spaces-bracket}.

To prove the second statement of the proposition, we first note that the isomorphism \eqref{eq:morphism-spaces-bracket} is given concretely as follows. A morphism $f$ on the left-hand side is represented by an object $L$ of $\cG$ and a decorated diffeomorphism $L \natural M \to N$ in $\cM \subseteq \cD_{d}$. Restricting to an $\epsilon$-neighbourhood of $M \subset L \natural M$ we obtain a decorated embedding of $M$ into $N$. This is the image of $f$ under \eqref{eq:morphism-spaces-bracket}. (This is ultimately because the identification \eqref{eq:restriction-map-identification} from Theorem~\ref{thm:fibre-bundle} is given by the restriction map.) From this description, it is evident that the square formed by (two instances of) \eqref{eq:morphism-spaces-bracket}, composition in $\langle \cG , \cM \rangle$ and composition of decorated embeddings is commutative.

Finally, all of the above proof adapts identically to the setting where $\cD_{d}$ is replaced by $\cD_{d}^{+}$.
\end{proof}

\begin{coro}\label{cor:cancellation1}
Suppose that \emph{cancellation} holds for our chosen subgroupoids $\cG$ and $\cM$ --- meaning that $L \natural M \cong L' \natural M$ in $\cM$ implies $L \cong L'$ in $\cG$ for objects $L,L'$ of $\cG$ and $M$ of $\cM$. Then we have:
\[
\langle \cG , \cM \rangle (M,N) \;\cong\; \begin{cases}
\embgm(M,L\natural M)_L & \text{if } \exists L \in \cG \text{ such that } L \natural M \cong N \text{ in } \cM ; \\
\emptyset & \text{otherwise.}
\end{cases}
\]
\end{coro}
\begin{proof}
This follows from the fact that the disjoint union \eqref{eq:gm-decomposition} is taken either over the empty set or a set of size one in this situation.
\end{proof}

\begin{rmk}
\label{rmk:cancellation1}
All of our examples of $\cG$ and $\cM$ in \S\ref{ss:categories_for_families_of_groups} satisfy cancellation. This follows from the classification of compact surfaces in all of our examples in dimension $2$ (see \S\ref{sss:category_surface_braid_groups} and \S\ref{sss:category_mcg}). In the examples of \S\ref{sss:category_loop_braid_groups}, it follows from the fact that the groupoid of finite sets and bijections under disjoint union satisfies cancellation, since the objects of the groupoids in the examples of \S\ref{sss:category_loop_braid_groups} are determined up to isomorphism by the number of copies of $\bS^{1}$ in the submanifold $A \subset M$.
\end{rmk}

As a conclusion of this subsection, we have the following corollary, which gives an explicit description of the semicategory $\fU\cD_{d}$ defined by applying the topologically-enriched Quillen bracket construction (see Proposition~\ref{prop:topological-Quillen}) to the semi-monoidal groupoid $\cD_{d}$ (see Definition~\ref{def:Dd}).

\begin{coro}\label{cor:description_UD_{d}}
For each $d\geq 2$, the semicategory $\fU\cD_{d}$ is isomorphic to the semicategory whose objects are the decorated manifolds $M=(M,A,e_{1},e_{2})$ of dimension $d$ and whose spaces of morphisms are given by the embedding spaces $\embdec(M,N)$ introduced in Definition~\ref{def:embdec} with composition as described in Definition~\ref{def:decorated-embeddings-composition}.
For subgroupoids $\cG$ and $\cM$ of $\cD_{d}$ satisfying the conditions of Hypothesis~\ref{hypo:standard_framework}, the subsemicategory $\langle \cG , \cM \rangle$ corresponds under this isomorphism to the semicategory whose objects are the objects of $\cM$ and whose morphisms are given by the embedding spaces $\embgm(M,N)$ of Definition \ref{def:diffdec-2}.

In addition, we have isomorphisms of discrete semicategories $\pi_{0}(\fU\cD_{d}) \cong \fU(\pi_{0}(\cD_{d}))$ and more generally $\pi_{0}(\langle \cG , \cM \rangle) \cong \langle \pi_{0}(\cG) , \pi_{0}(\cM) \rangle$ for $\cG$ and $\cM$ as above. 
In particular, $\pi_{0}(\cD_{d})$ and $\pi_{0}(\cM)$ are the underlying groupoids of $\pi_{0}(\fU\cD_{d})$ and $\pi_{0}(\langle \cG , \cM \rangle)$ respectively. 

All of the above holds equally for $\fU\cD_{d}^{+}$ and subgroupoids $\cG , \cM \subseteq \cD_{d}^{+}$ as in Hypothesis~\ref{hypo:standard_framework}.
\end{coro}

\begin{proof}
Let us denote by $\fU\cD_{d}'$ the semicategory described in the first statement of the corollary. To begin with, we note that this is a well-defined (topologically-enriched) semicategory: its composition is defined in Definition~\ref{def:decorated-embeddings-composition} and is verified to be continuous and associative in Lemma~\ref{lem:decorated-embeddings-composition-continuous-associative}.
An isomorphism of semicategories between $\fU\cD_{d}$ and $\fU\cD_{d}'$ is then given by the identity on objects and by Proposition~\ref{prop:morphism-spaces-bracket} (with $\cG = \cM = \cD_{d}$) on morphisms. Moreover, the second statement of the corollary also follows directly from Proposition~\ref{prop:morphism-spaces-bracket} (in its general setting). The third statement is Lemma~\ref{lem:Serre-fibration-condition}.
\end{proof}

\subsubsection{Functorial split short exact sequences}
\label{sss:functorial-split-ses}

As an immediate application of the description in Corollary~\ref{cor:description_UD_{d}} of $\fU\cD_{d}$ in terms of decorated embeddings of manifolds, we prove that the split short exact sequences \eqref{eq:split-ses-1} and \eqref{eq:split-ses-2}  -- and indeed the whole diagram \eqref{eq:3x3diagram} -- constructed in \S\ref{ss:split-ses} are functorial on $\fU\cD_{d}$. Recall that we are continuing to assume that $d \neq 4$ (cf.~Remark~\ref{rmk:d_neq_4_explanation}).

\begin{prop}
\label{prop:split-ses-functoriality}
Each morphism $\varphi \colon (M,A) \to (N,B)$ in $\fU\cD_{d}$ induces a map of diagrams of the form \eqref{eq:3x3diagram-proof}, preserving basepoints of the embedding spaces and respecting the group structure of the diffeomorphism groups. As a result, we obtain a map of diagrams of the form \eqref{eq:3x3diagram}. In particular, we obtain maps of split short exact sequences of the form \eqref{eq:split-ses-1} and \eqref{eq:split-ses-2}.
\end{prop}
\begin{proof}
By Corollary~\ref{cor:description_UD_{d}}, $\varphi$ is a decorated embedding $(M,A) \hookrightarrow (N,B)$, which means by definition (see Definition~\ref{def:embdec}) that there is a diffeomorphism of decorated manifolds $\bar{\varphi} \colon (M',A') \natural (M,A) \to (N,B)$, for some decorated manifold $(M',A')$, so that $\varphi$ is equal to its pre-composition with the inclusion of (an arbitrarily small neighbourhood of) $(M,A)$ into $(M',A') \natural (M,A)$. Let us fix such a diffeomorphism $\bar{\varphi}$ and extend it by the identity on $\bB^d_{1}$ to $\widehat{\varphi} \colon (M',A') \natural (\mbar,A) \to (\nbar,B)$, where $\mbar = M \natural \bB^d_{1}$ and $\nbar = N \natural \bB^d_{1}$ as in Notation~\ref{not:union-with-Z}.

On each of the diffeomorphism groups of \eqref{eq:3x3diagram-proof}, the induced map is defined by sending a decorated diffeomorphism $\psi$ of $\mbar$ to the decorated diffeomorphism $E_\varphi(\psi) = \widehat{\varphi} \circ (\id_{M'} \natural \psi) \circ \widehat{\varphi}^{-1}$ of $\nbar$. It is straightforward to check that if $\psi$ sends $A$ onto itself then $E_\varphi(\psi)$ sends $B$ onto itself and that if $\psi$ sends $Z$ onto itself then $E_\varphi(\psi)$ sends $Z$ onto itself by the same diffeomorphism. It is also clear that $E_\varphi(-)$ is a group homomorphism.

On the embedding spaces of \eqref{eq:3x3diagram-proof} with domain $Z$, the induced map is defined simply by post-composing with $\widehat{\varphi}$. On the embedding spaces of \eqref{eq:3x3diagram-proof} with domain $A$, the induced map is defined by sending an orbit of embeddings $[e]$ to $[\widehat{\varphi} \circ (\mathrm{incl}_{A'} \sqcup e) \circ (\widehat{\varphi}^{-1})|_B]$, where $\mathrm{incl}_{A'}$ denotes the inclusion $A' \subset \Int{M}'$. On the embedding spaces of \eqref{eq:3x3diagram-proof} with domain $A \sqcup Z$, the induced map is defined separately on $Z$ and on $A$ as in the previous two sentences. These maps preserve basepoints because $\widehat{\varphi}$ sends $A' \sqcup A$ onto $B$ and because it fixes $Z$ pointwise, since $Z \subset \bB^d_{1} \subset \mbar$.

It remains to check that these induced maps commute with the maps of diagram \eqref{eq:3x3diagram-proof}. For the (solid) maps between two diffeomorphism groups or between two embedding spaces, this is obvious from the construction. For the three maps that go from a diffeomorphism group to an embedding space, it is also straightforward to check, the only subtlety being the observation that $E_\varphi(\psi)|_Z = \widehat{\varphi} \circ \psi|_Z$, since $\widehat{\varphi}$ restricts to the identity on $Z \subset \bB^d_{1}$. Finally, for the dotted maps of \eqref{eq:3x3diagram-proof}, their commutation with the maps induced by $\varphi$ may also be checked directly, and follows essentially because the dotted maps are defined by ``squashing'' the second (right-hand) boundary-cylinder of $M$ or $N$ using the self-embedding $\Theta$ (see Construction~\ref{construction:Theta}), whereas $\varphi$ is defined by taking the boundary connected sum with $M'$ or $N'$ along the first (left-hand) boundary-cylinder of $M$ or $N$ and then applying the diffeomorphism $\bar{\varphi}$; the two operations thus have ``disjoint support''.
\end{proof}

\subsubsection{Skeleta of Quillen bracket categories}
\label{sss:skeleta}

It is often convenient to be able to pass to a skeleton of a category. For a plain (small) category this is easy: one simply has to pick one object from each isomorphism class (using the axiom of choice, if necessary) and consider the full subcategory on these objects. (We recall that all categories considered in this paper are small; see Convention~\ref{convention:small_categories} and Remark~\ref{rmk:small_categories_Quillen}.) For categories with additional structure, one has to be a little more careful to preserve the structure. The authors learned of the following folklore fact from \cite{MO_31599}.

\begin{prop}
\label{prop:monoidal-skeleton}
\textup{\textbf{(i)}} Let $\cC$ be a monoidal category and let $\cC_{0}$ be any skeleton of $\cC$.
Then there is a monoidal structure on $\cC_{0}$ such that the inclusion $\cC_{0} \subset \cC$ may be enhanced to a monoidal equivalence.

\textup{\textbf{(ii)}} Let $\cM$ be a category equipped with a left action of a skeletal monoidal category $\cC_{0}$ and let $\cM_{0}$ be any skeleton of $\cM$. Then there is a left action of $\cC_{0}$ on $\cM_{0}$ such that the inclusion $\cM_{0} \subset \cM$ may be enhanced to an equivalence of categories equipped with left $\cC_{0}$-actions.
\end{prop}
\begin{proof}[Sketch proof]
Let us first consider part (i). By \cite[Thm.~1, \S IV.4]{MacLane1}, the inclusion $\cC_{0} \subset \cC$ is part of an adjoint equivalence, namely an adjunction whose unit and counit are natural isomorphisms. Choosing such an adjoint equivalence, one may transfer the monoidal structure of $\cC$ along it to obtain a monoidal structure on $\cC_{0}$. The coherence isomorphisms needed to enhance the inclusion $\cC_{0} \subset \cC$ to a (strong) monoidal equivalence may then be constructed from the units and counits of the adjoint equivalence. The same strategy also deals with part (ii).
\end{proof}

\begin{rmk}
In general, the monoidal structure induced on $\cC_{0}$ and the enhancement of the inclusion $\cC_{0} \subset \cC$ to a monoidal equivalence depend on the choice of adjoint inverse to the inclusion in the proof above. We also note that Proposition~\ref{prop:monoidal-skeleton} is not true if one insists on the inclusion $\cC_{0} \subset \cC$ being a \emph{strict} monoidal equivalence.
\end{rmk}

As a consequence, we may always find a skeleton of a Quillen bracket category that is itself a Quillen bracket category:

\begin{coro}
\label{coro:QBC-skeleton}
Let $\cC$ be a monoidal category and let $\cM$ be a category with a left action of $\cC$. Then there are skeleta $\cC_{0} \subset \cC$ and $\cM_{0} \subset \cM$, together with a monoidal structure on $\cC_{0}$ and a left $\cC_{0}$-action on $\cM_{0}$, such that the inclusion $\langle \cC_{0} , \cM_{0} \rangle \subset \langle \cC , \cM \rangle$ is an equivalence. Moreover, the subcategory $\langle \cC_{0} , \cM_{0} \rangle$ is skeletal.
\end{coro}
\begin{proof}
We may choose $\cC_{0} \subset \cC$ and $\cM_{0} \subset \cM$ to be any skeleta.
Proposition~\ref{prop:monoidal-skeleton} then provides us with a monoidal structure on $\cC_{0}$ and a left $\cC_{0}$-action on $\cM_{0}$ together with enhancements of $\cC_{0} \subset \cC$ to a monoidal equivalence and of $\cM_{0} \subset \cM$ to an equivalence of categories equipped with left $\cC_{0}$-actions. These induce an inverse equivalence for the inclusion $\langle \cC_{0} , \cM_{0} \rangle \subset \langle \cC , \cM \rangle$. Finally, $\langle \cC_{0} , \cM_{0} \rangle$ is skeletal since its underlying groupoid is $\cM_{0}$, which is skeletal.
\end{proof}

\section{Construction of homological representation functors}\label{s:general_construction}

In this section, we construct homological representation functors for mapping class groups and motion groups. We consider the topologically-enriched category $\fU\cD_{d}$ introduced in \S\ref{s:categorical_framework}, which contains all motion groups and all mapping class groups of $d$-dimensional manifolds, as explained in Remark~\ref{rmk:globality}.
The construction consists of two main parts:
\begin{enumerate}[noitemsep,label=(\arabic*)]
\item\label{construction:step1} constructing continuous semifunctors $\fU\cD_{d} \to \covr$;
\item\label{construction:step2} using the output of step \ref{construction:step1} to construct functors $\pi_{0}(\langle \cG , \cM \rangle) \to \modr$, by restricting to some subcategory $\langle \cG , \cM \rangle \subset \fU\cD_{d}$ and applying the procedure summarised in diagram \eqref{eq:construction}.
\end{enumerate}
Here, $\covr$ and $\modr$ denote categories of spaces equipped, respectively, with coverings and with right modules over rings; see Definitions~\ref{def:bicoverings} and \ref{def:bimodule_category}.

In \S\ref{ss:ingredients_construction_representations}, we deal with step \ref{construction:step2} in a general setting: for any topologically-enriched semicategory $\cC$ and continuous semifunctor $F \colon \cC \to \covr$, as well as a continuous semifunctor $V \colon \cC \to \modlr$ that is ``compatible'' in a precise way with $F$, we produce a (semi)functor $L_{i}(F;V) \colon \pi_{0}(\cC) \to \modr$ for each $i\geq 0$.
We first give precise definitions of the categories involved in \S\ref{sss:categories}. The three steps of the construction are then described in \S\ref{sss:lift}--\S\ref{sss:twisted-homology}, and are put together in \S\ref{sss:construction}.

In \S\ref{ss:homological_representation_functor_motion_groups} and \S\ref{ss:homological_representation_functor_mcg} we then focus on step \ref{construction:step1}. This is the more geometric, or topological, part of our construction of representations, whereas step \ref{construction:step2} is the more formal and algebraic part. We note that the semifunctors $\cF \colon \fU\cD_{d} \to \covr$ obtained as outputs of step \ref{construction:step1} restrict, after passing to $\pi_0$, to all motion groups and all mapping class groups in dimension $d$. For this reason, we refer to them as \emph{global} semifunctors.

In contrast, we then restrict each of these semifunctors $\cF$ to a subcategory $\langle \cG , \cM \rangle \subset \fU\cD_{d}$, where $\cG$ and $\cM$ are a pair of groupoids as in Hypothesis~\ref{hypo:standard_framework}, satisfying further mild properties; see \S\ref{sss:col_coeff_unwtisted_rep}. This is in order to construct in \S\ref{sss:col_coeff_unwtisted_rep} a continuous semifunctor $V = V_{\col}(\cF) \colon \langle \cG , \cM \rangle \to \modlr$ that is compatible with $\cF$ in the above sense. It is constructed so that, under certain conditions, the functor $L_{i}(\cF;V_{\col}(\cF)) \colon \pi_{0}(\langle \cG , \cM \rangle) \to \modr$ produced in step \ref{construction:step2} has image contained in the subcategory $\modr[R] \subset \modr$ for a fixed ring $R$, in other words it is an \emph{untwisted} functorial representation.

\subsection{Elements of the construction}\label{ss:ingredients_construction_representations}

This section presents the elements of the construction of functorial homological representations, in the sense of step \ref{construction:step2} above, and how to assemble them.

A minor subtlety is that the framework in \S\ref{ss:ingredients_construction_representations} is set up in order to construct homological representation functors starting from a topologically-enriched \emph{semi}category $\cC$, in order to be able to consider $\cC := \fU\cD_{d}$ defined in \S\ref{s:categorical_framework}. This is the reason why the construction produces a \emph{semi}functor. However, in our applications of the construction in \S\ref{ss:homological_representation_functor_motion_groups} and \S\ref{ss:homological_representation_functor_mcg}, the discrete semicategory $\pi_{0}(\cC)$ will always be a \emph{genuine} category (i.e.~with identities) and the induced semifunctor from $\pi_{0}(\cC)$ will always be a \emph{genuine} functor (i.e.~preserving identities); see Remark~\ref{rmk:homological_representation_semifunctor_genuine}.

\subsubsection{Categories}\label{sss:categories}

In this subsection, we define the different categories involved in the construction.

\begin{rmk}
The topological spaces involved in the definitions of \S\ref{ss:ingredients_construction_representations} will always be \emph{locally path-connected} and \emph{semi-locally simply-connected}.
\end{rmk}

\begin{defi}[The category of bicoverings.]
\label{def:bicoverings}
The topologically-enriched category $\covlr$ of \emph{spaces equipped with bicoverings} is defined as follows. An object of $\covlr$ consists of a pair of groups $Q_{1},Q_{2}$ and a path-connected, based space $X$ admitting a universal covering (i.e.~locally path-connected and semi-locally simply-connected), equipped with a surjective homomorphism $\phi\colon \pi_{1}(X) \twoheadrightarrow Q_{1} \times Q_{2}$. For $i\in\{1,2\}$, we denote by $\phi_{i}$ the further projection of $\phi$ onto $Q_{i}$. Via the correspondence between connected, regular coverings of $X$ and normal subgroups of $\pi_{1}(X)$, this is the same as specifying a pullback square of connected, regular coverings:
\begin{center}
\begin{tikzpicture}
[x=1mm,y=1mm]
\node (tl) at (0,12) {$X^{\phi}$};
\node (tr) at (30,12) {$X^{\phi_{2}}$};
\node (bl) at (0,0) {$X^{\phi_{1}}$};
\node (br) at (30,0) {$X$};
\draw[->] (tl) to (tr);
\draw[->] (bl) to (br);
\draw[->] (tl) to (bl);
\draw[->] (tr) to (br);
\node at (4,8) {$\lrcorner$};
\end{tikzpicture}
\end{center}
A morphism in $\covlr$ from $(Q_{1},Q_{2},X,\phi)$ to $(Q'_{1},Q'_{2},X',\phi')$ is a based, continuous map $f \colon X \to X'$ such that the induced homomorphism $\pi_{1}(f)$ sends $\ker(\phi_{1})$ into $\ker(\phi'_{1})$ and $\ker(\phi_{2})$ into $\ker(\phi'_{2})$. Equivalently, this condition may be phrased as saying that there exist (uniquely defined) homomorphisms $\alpha_{1} \colon Q_{1} \to Q'_{1}$ and $\alpha_{2} \colon Q_{2} \to Q'_{2}$ such that
\[
\phi'_{1} \circ \pi_{1}(f) = \alpha_{1} \circ \phi_{1} \qquad\text{and}\qquad \phi'_{2} \circ \pi_{1}(f) = \alpha_{2} \circ \phi_{2} .
\]
The morphism space of $\covlr$ from $(Q_{1},Q_{2},X,\phi)$ to $(Q'_{1},Q'_{2},X',\phi')$ is topologised as a subspace of the space of continuous maps $X \to X'$ in the compact-open topology.

If $G$ is a group, the category $\covlr[G]$ is the subcategory of $\covlr$ on those objects $(X,\phi)$ such that $Q_{1} = G$ and those morphisms $f$ such that the induced homomorphism $\alpha_{1}\colon Q_{1} \to Q'_{1}$ is equal to $\id_G$. Similarly, we have a subcategory $\covlr[\bullet][G]$. If $G$ is the trivial group, we drop it from the notation, and write $\covr$ and $\covl$ respectively for these subcategories of $\covlr$. We note that $\covl$ and $\covr$ are abstractly isomorphic, but not equal as subcategories of $\covlr$. As another variant, we write $\covr[G]^{\tw} \subset \covr$ for the \emph{full} subcategory on those objects $(X,\phi)$ such that $Q_{1} = G$. (The superscript ${}^{\tw}$ indicates that the morphisms are permitted to act on the transformation group $G$, i.e.~we are considering \emph{twisted} covering spaces.)
We define the subcategories $\covlr^{\pr}$, $\covr^{\pr}$, $\covl^{\pr}$ and $\covr[G]^{\pr,\tw}$ of $\covlr$, $\covr$, $\covl$ and $\covr[G]^{\tw}$, respectively, to have the same objects but only those morphisms that are \emph{proper} as maps between spaces.
\end{defi}

We recall that all of the rings that we consider are assumed to be \emph{associative} and \emph{unital} and that their morphisms preserve units.

\begin{defi}[The category of bimodules.]\label{def:bimodule_category}
The category $\modlr$ is the category of \emph{bimodules over rings}. An object of $\modlr$ is a pair of rings $(R,S)$ together with an $(R,S)$-bimodule $V$. A morphism from $(R,S,V)$ to $(R',S',V')$ is a pair of ring homomorphisms $\alpha \colon R \to R'$ and $\beta \colon S \to S'$, together with a morphism of $(R,S)$-bimodules $\theta \colon V \to (\alpha,\beta)^{*} (V')$.

The category $\modlr$ is topologically-enriched by equipping each morphism space with the discrete topology.

For any ring $R$, we have a subcategory $\modlr[R] \subset \modlr$ on those objects $(R',S',V')$ where $R'=R$ and those morphisms $(\alpha,\beta,\theta)$ where $\alpha = \id_R$ (note that this is generally not a full subcategory). If moreover $R=\bZ$, we drop it from the notation, and write $\modr$ for $\modlr[R]$. Similarly, we have $\modlr[\bullet][S]$ and $\modlr[R][S]$ for rings $R$ and $S$. We note that $\modlr[R][S]$ is just the category of $(R,S)$-bimodules, as usually defined. Furthermore, for a fixed ring $S$, we write $\modr[S]^{\tw} \subset \modr$ for the \emph{full} subcategory on those objects $(S',V')$ where $S'=S$. (The superscript ${}^{\tw}$ indicates that the morphisms are permitted to act on the underlying ring $S$, i.e.~we are considering \emph{twisted} modules.)
\end{defi}

\begin{defi}[Bundles of $G$-sets and of $(R,S)$-bimodules.]
\label{d:bundlesII}
For any space $X$ and group $G$, a bundle of left $G$-sets over $X$ is a functor $\pi_{\leq 1}(X) \to \setlr[G][]$, where $\pi_{\leq 1}(X)$ is the fundamental groupoid of $X$ and $\setlr[G][]$ is the category of left $G$-sets. Similarly, for a pair of rings $(R,S)$, a bundle of $(R,S)$-bimodules over $X$ is a functor $\xi \colon \pi_{\leq 1}(X)\to \modlr[R][S]$.

A morphism $[\tau] = [f,\alpha,\beta,\tau] \colon \xi \to \xi'$ of bundles of bimodules from $\xi\colon \pi_{\leq 1}(X)\to \modlr[R][S]$ to $\xi'\colon\pi_{\leq 1}(X')\to \modlr[R'][S']$ is a continuous map $f \colon X \to X'$, two ring homomorphisms $\alpha \colon R \to R'$ and $\beta \colon S \to S'$, and a natural transformation $\tau \colon \xi \Rightarrow (\alpha,\beta)^{*} \circ \xi' \circ \pi_{\leq 1}(f)$, where $(\alpha,\beta)^{*} \colon {}_{R'}\modr[S'] \to \modlr[R][S]$ denotes the restriction functor induced by $\alpha$ and $\beta$.
\end{defi}

\begin{lemm}\label{lem:composition_bundle_bimodules}
In notation of Definition~\ref{d:bundlesII}, there is a well-defined, associative composition of bundles of bimodules, as follows. Let $[\tau_{1}] = [f_{1},\alpha_{1},\beta_{1},\tau_{1}] \colon \xi \to \xi'$ and $[\tau_{2}] = [f_{2},\alpha_{2},\beta_{2},\tau_{2}] \colon \xi' \to \xi''$ be bundles of bimodules. Their composition $[\tau_{2}]\circ[\tau_{1}]\colon \xi\to \xi''$ is the natural transformation
\begin{equation}\label{eq:composition_bundle_bimodules}
[(\alpha_{1},\beta_{1})^{*} \circ \tau_{2} \circ \pi_{\leq 1}(f_{1})] \circ \tau_{1} \colon \xi \Rightarrow (\alpha_{2} \circ \alpha_{1} , \beta_{2} \circ \beta_{1})^{*} \circ \xi'' \circ \pi_{\leq 1}(f_{2} \circ f_{1})
\end{equation}
given by pasting together the natural transformations $\tau_{1}$ and $\tau_{2}$ in the diagram:
\begin{equation}
\label{eq:Top-composition-diagram}
\begin{split}
\begin{tikzpicture}
[x=1mm,y=1mm]
\node (tl) at (0,15) {$\pi_{\leq 1}(X)$};
\node (tm) at (30,15) {$\modlr[R][S]$};
\node (ml) at (0,0) {$\pi_{\leq 1}(X')$};
\node (mm) at (30,0) {$\modlr[R'][S']$};
\node (bl) at (0,-15) {$\pi_{\leq 1}(X'')$};
\node (bm) at (30,-15) {$\modlr[R''][S''].$};
\draw[->] (tl) to node[above,font=\small]{$\xi$} (tm);
\draw[->] (ml) to node[above,font=\small]{$\xi'$} (mm);
\draw[->] (bl) to node[above,font=\small]{$\xi''$} (bm);
\draw[->] (tl) to node[left,font=\small]{$\pi_{\leq 1}(f_{1})$} (ml);
\draw[->] (ml) to node[left,font=\small]{$\pi_{\leq 1}(f_{2})$} (bl);
\draw[->] (mm) to node[right,font=\small]{$(\alpha_{1},\beta_{1})^{*}$} (tm);
\draw[->] (bm) to node[right,font=\small]{$(\alpha_{2},\beta_{2})^{*}$} (mm);
\node (tau1) at (15,9) {\rotatebox{270}{$\Rightarrow$}};
\node (tau2) at (15,-6) {\rotatebox{270}{$\Rightarrow$}};
\node at (tau1.east) [anchor=west,font=\small] {$\tau_{1}$};
\node at (tau2.east) [anchor=west,font=\small] {$\tau_{2}$};
\end{tikzpicture}
\end{split}
\end{equation}
\end{lemm}
\begin{proof}
Since $(\alpha_{1},\beta_{1})^{*} \circ \tau_{2} \circ \pi_{\leq 1}(f_{1})$ is a natural transformation (pre-composing or post-composing a natural transformation by a functor produces a natural transformation), the assignment \eqref{eq:composition_bundle_bimodules} is well-defined as it is a composition of natural transformations.
Associativity follows from associativity of composition of natural transformations, as well as functoriality of $\pi_{\leq 1}(-) \colon \topo \to \cat$ and of $\modlr[-][-] \colon \Rings \times \Rings \to \cat$.
\end{proof}

\begin{rmk}
\label{r:bundles}
For a path-connected space $X$, a bundle of left $G$-sets over $X$ is more often defined as a fibre bundle over $X$ with fibre a left $G$-set $T$ and structure group the automorphism group $\Aut_G(T)$ of $T$ as a left $G$-set; see \cite[\S 2]{Steenrod1951} for example. With this definition, by unique path-lifting, a bundle of left $G$-sets over a path-connected space $X$ determines a functor $\pi_{\leq 1}(X) \to \setlr[G][]$. Moreover, whenever $X$ is locally path-connected and semi-locally simply-connected, this gives an identification between isomorphism classes of bundles of left $G$-sets over $X$ and isomorphism classes of functors $\pi_{\leq 1}(X) \to \setlr[G][]$. In all of our examples, the base space $X$ has these properties, and it is more convenient to use Definition~\ref{d:bundlesII} rather than this definition in terms of fibre bundles.
\end{rmk}

\begin{defi}[The category of bundles of bimodules.]
\label{def:bundles_of_bimodules}
The topologically-enriched category of \emph{bundles of bimodules} $\toplr$ is defined as follows. An object of $\toplr$ is a space $X$ together with a pair of rings $(R,S)$ and a bundle $\xi$ of $(R,S)$-bimodules over $X$ in the sense of Definition~\ref{d:bundlesII}. A morphism from $(X,R,S,\xi)$ to $(X',R',S',\xi')$ is a morphism of bundles of bimodules as described in Definition~\ref{d:bundlesII}.
The composition of morphisms and its associativity follow from Lemma~\ref{lem:composition_bundle_bimodules}.
The category $\toplr$ is then topologically-enriched by topologising each morphism space from $(X,R,S,\xi)$ to $(X',R',S',\xi')$ as follows: a morphism is given by a four-tuple $[f,\alpha,\beta,\tau]$ (see Definition~\ref{d:bundlesII}); we use the compact-open topology for the continuous map $f$ and the discrete topology for the other three components of this tuple.

For any ring $R$, we have a subcategory $\toplr[R] \subset \toplr$ on those objects $(X',R',S',\xi')$ where $R'=R$ and those morphisms $(f,\alpha,\beta,\tau)$ where $\alpha = \id_R$. If moreover $R=\bZ$, we drop it from the notation, and write $\topr$ for $\toplr[R]$. Similarly, we have $\toplr[\bullet][S]$ and $\toplr[R][S]$ for rings $R$ and $S$. In addition, for a fixed ring $S$, we write $\topr[S]^{\tw} \subset \topr$ for the \emph{full} subcategory on those objects $(X',S',\xi')$ where $S'=S$. (The superscript ${}^{\tw}$ indicates that the morphisms are permitted to act on the ground ring $S$, i.e.~we are considering \emph{twisted} bundles of bimodules.)

We define the subcategories $\toplr^{\pr}$, $\topr^{\pr}$, $\topl^{\pr}$ and $\topr[S]^{\pr,\tw}$ of $\toplr$, $\topr$, $\topl$ and $\topr[S]^{\tw}$, respectively, to have the same objects but only those morphisms whose underlying map of spaces is a \emph{proper} map.
\end{defi}

Note that $\modlr$ is equivalent to the full subcategory of $\toplr$ on those objects whose underlying space is a point.

\begin{notation}
\label{not:lr_forgetful}
Writing $\Rings$ for the category of (associative, unital) rings and $\groups$ for the category of groups, there are obvious forgetful functors
\begin{center}
\begin{tikzpicture}
[x=1mm,y=1mm]
\node (l1) at (0,0) [anchor=east] {$\covlr$};
\node (r1) at (10,0) [anchor=west] {$\groups$};
\node (l2) at (70,0) [anchor=east] {$\modlr \subset \toplr$};
\node (r2) at (80,0) [anchor=west] {$\Rings$};
\draw[->,transform canvas={yshift=0.4ex}] (l1) to node[above,font=\small]{$\sL$} (r1);
\draw[->,transform canvas={yshift=-0.4ex}] (l1) to node[below,font=\small]{$\sR$} (r1);
\draw[->,transform canvas={yshift=0.4ex}] (l2) to node[above,font=\small]{$\sL$} (r2);
\draw[->,transform canvas={yshift=-0.4ex}] (l2) to node[below,font=\small]{$\sR$} (r2);
\end{tikzpicture}
\end{center}
that remember just the left (respectively right) underlying group or ring. For example, an object $(X,R,S,\xi)$ of $\toplr$ is sent under $\sL$ to $R$ and under $\sR$ to $S$.
\end{notation}

\begin{defi}\label{def:pulbacktopmod}
Write $\topmodlr$ for the pullback of the forgetful functors $\sL \colon \modlr \to \Rings$ and $\sR \colon \toplr \to \Rings$:
\begin{equation}
\label{eq:topmod}
\centering
\begin{split}
\begin{tikzpicture}
[x=1mm,y=1mm]
\node (tl) at (0,12) {$\topmodlr$};
\node (tr) at (30,12) {$\modlr$};
\node (bl) at (0,0) {$\toplr$};
\node (br) at (30,0) {$\Rings .$};
\draw[->] (tl) to (tr);
\draw[->] (bl) to node[below,font=\small]{$\sR$} (br);
\draw[->] (tl) to (bl);
\draw[->] (tr) to node[right,font=\small]{$\sL$} (br);
\node at (4,8) {$\lrcorner$};
\end{tikzpicture}
\end{split}
\end{equation}
Explicitly, an object of $\topmodlr$ consists of a space $X$, a triple of rings $(R,S,T)$, a bundle of $(R,S)$-bimodules $\xi$ over $X$ and an $(S,T)$-bimodule $V$. A morphism from $(X,R,S,T,\xi,V)$ to $(X',R',S',T',\xi',V')$ consists of a morphism of bundles of bimodules as described in Definition~\ref{d:bundlesII} (a continuous map $f \colon X \to X'$, a pair of ring homomorphisms $(\alpha \colon R \to R , \beta \colon S \to S')$ and a natural transformation $\tau \colon \xi \Rightarrow (\alpha,\beta)^{*} \circ \xi' \circ \pi_{\leq 1}(f)$), together with another ring homomorphism $\gamma \colon T \to T'$ and a morphism of $(S,T)$-bimodules $\theta \colon V \to (\beta,\gamma)^{*} (V')$. The category $\topmodlr$ is then topologically-enriched by topologising each morphism space from $(X,R,S,\xi)$ to $(X',R',S',\xi')$ as follows: a morphism is given by a $6$-tuple $[f,\alpha,\beta,\gamma,\tau,\theta]$; we use the compact-open topology for the continuous map $f$ and the discrete topology for the other five components of the tuple.
We define the subcategory $\topmodlrpr$ via the analogue of \eqref{eq:topmod}, using the category $\toplr^{\pr}$ instead of $\toplr$.
\end{defi}

\begin{rmk}
The categories of bimodules $\modlr$ (see Definition~\ref{def:bimodule_category}), of bundles of bimodules $\toplr$ (see Definition~\ref{def:bundles_of_bimodules}) and $\topmodlr$ of Definition~\ref{def:pulbacktopmod} may easily be generalised using a fixed ring $\bA$ rather than $\bZ$ as ground ring, $\bA$-algebras $(R,S)$ rather than rings and the category $\Aalg$ rather than $\Rings$. The work of \S\ref{sss:lift}--\S\ref{sss:construction} repeats verbatim for this more general framework.
\end{rmk}

\subsubsection{From bicoverings to bundles of bimodules}\label{sss:lift}

In this subsection, we introduce the \emph{linearisation} continuous functor $\covlr \to \toplr$.

\begin{prop}
\label{prop:lift}
There is a natural continuous functor $\covlr \to \toplr$, taking bicoverings to bundles of bimodules, such that the squares
\begin{equation}
\label{eq:lift-condition}
\centering
\begin{split}
\begin{tikzpicture}
[x=1mm,y=1mm]
\node (tl) at (0,15) {$\covlr$};
\node (tr) at (30,15) {$\toplr$};
\node (bl) at (0,0) {$\groups$};
\node (br) at (30,0) {$\Rings$};
\draw[dashed,->] (tl) to node[above,font=\small]{} (tr);
\draw[->] (bl) to node[above,font=\small]{$\bZ[(-)^{\opp}]$} (br);
\draw[->] (tl) to node[left,font=\small]{$\sL$} (bl);
\draw[->] (tr) to node[right,font=\small]{$\sL$} (br);
\begin{scope}[xshift=60mm]
\node (tla) at (0,15) {$\covlr$};
\node (tra) at (30,15) {$\toplr$};
\node (bla) at (0,0) {$\groups$};
\node (bra) at (30,0) {$\Rings$};
\draw[dashed,->] (tla) to node[above,font=\small]{} (tra);
\draw[->] (bla) to node[above,font=\small]{$\bZ[-]$} (bra);
\draw[->] (tla) to node[left,font=\small]{$\sR$} (bla);
\draw[->] (tra) to node[right,font=\small]{$\sR$} (bra);
\end{scope}
\end{tikzpicture}
\end{split}
\end{equation}
commute, where $\bZ[-] \colon \groups \to \Rings$ is the group ring functor and $(-)^{\opp} \colon \groups \to \groups$ takes a group to its opposite. By abuse of notation, we denote this functor $\covlr \to \toplr$ also by $\bZ[-]$.
Furthermore, the restriction of the linearisation functor $\bZ[-]$ to the subcategory $\covr^{\pr} \subset \covr$ takes values in the subcategory $\topr^{\pr} \subset \topr$.
\end{prop}

\begin{proof}
On objects, this is defined as follows. Let $(Q_{1},Q_{2},X,x_{0},\phi \colon \pi_{1}(X,x_{0}) \twoheadrightarrow Q_{1} \times Q_{2})$ be an object of $\covlr$. The normal subgroup
\begin{equation}
\label{eq:kernel-intersection}
K \coloneqq \ker(\phi_{1}) \cap \ker(\phi_{2}) = \ker(\phi) \;\triangleleft\; \pi_{1}(X,x_{0})
\end{equation}
corresponds to a regular covering of $X$ with deck transformation group $Q \coloneqq Q_{1} \times Q_{2}$. In order to specify a particular regular covering of $X$, rather than just an \emph{isomorphism class} of such, we will be slightly more careful. Start with the universal covering $\widetilde{X}$ of $X$; more specifically, the standard model for $\widetilde{X}$ whose underlying set consists of endpoint-preserving homotopy classes of paths in $X$ starting at $x_{0}$. This is equipped with an action of $\pi_{1}(X,x_{0})$; take the quotient of $\widetilde{X}$ by the action of the subgroup $K$. We denote this regular covering by $\xi_{\phi} \colon X^\phi = \widetilde{X}/K \rightarrow X$.
Whether deck transformations act on the left or the right is an arbitrary convention. We will consider them to act on the \emph{right}, since this agrees with the typical convention that the structure group of a principal bundle (of which regular coverings are examples) acts on the total space on the \emph{right}.
This is a bundle of right $Q$-sets over $X$ (whose fibres all happen to be isomorphic to $Q$ itself). Since $Q$ is the product $Q_{1} \times Q_{2}$, we may equally view $\xi_{\phi}$ as a bundle of $(Q_{1}^{\opp},Q_{2})$-bisets over $X$, where a $(G,H)$-biset is a set equipped with a left $G$-action and a compatible right $H$-action.

Now replace each fibre of $\xi_{\phi}$ with the free $\bZ$-module generated by that fibre; this forms a bundle of $(\bZ[Q_{1}^{\opp}],\bZ[Q_{2}])$-bimodules over $X$. This operation of \emph{taking free $\bZ$-modules fibrewise} is simple to describe when viewing bundles of $(G,H)$-bisets over $X$ as functors $\pi_{\leq 1}(X) \to \setlr[G][H]$ to the category of $(G,H)$-bisets and bundles of $(R,S)$-bimodules over $X$ as functors $\pi_{\leq 1}(X) \to \modlr[R][S]$ to the category of $(R,S)$-bimodules. Namely, the operation is given by post-composition with the functor $\bZ[-] \colon \setlr[Q_{1}^{\opp}][Q_{2}] \rightarrow \modlr[{\bZ[Q_{1}^{\opp}]}][{\bZ[Q_{2}]}]$ defined by assigning the free $(\bZ[Q_{1}^{\opp}],\bZ[Q_{2}])$-bimodule $\bZ[S]$ to each $(Q_{1}^{\opp},Q_{2})$-biset $S$.
We denote the resulting bundle of bimodules by $\bZ_{\fib}[\xi_{\phi}] \colon \bZ_{\fib}[X^\phi] \rightarrow X$.
This defines the linearisation functor $\bZ[-]\colon \covlr \to \toplr$ on objects:
\[
\bZ[-](Q_{1},Q_{2},X,x_{0},\phi) \;=\; (X,\bZ[Q_{1}^{\opp}],\bZ[Q_{2}],\bZ_{\fib}[\xi_{\phi}]).
\]

In order to define the linearisation functor $\bZ[-]\colon \covlr \to \toplr$ on morphisms, we first note that, although we did not need it to define the functor on objects, the regular covering $\xi_{\phi} \colon X^\phi \to X$ associated to $(Q_{1},Q_{2},X,x_{0},\phi)$ comes equipped with a particular choice of basepoint of $X^\phi$, covering the basepoint $x_{0}$ of $X$. This is because the standard construction of the universal cover $\widetilde{X}$, which we used above, has a canonical basepoint (namely the constant path at $x_{0}$), and therefore so does its quotient $X^\phi$. Let us denote this basepoint by $\widetilde{x}_{0} \in X^\phi$.

Now suppose we have a morphism $(Q_{1},Q_{2},X,x_{0},\phi) \to (Q'_{1},Q'_{2},Y,y_{0},\phi')$ in $\covlr$, that is, a continuous, based map $f \colon X \to Y$ such that $f_*(\ker(\phi_{1})) \subseteq \ker(\phi'_{1})$ and $f_*(\ker(\phi_{2})) \subseteq \ker(\phi'_{2})$. We recall from Definition~\ref{def:bicoverings} that this determines certain homomorphisms $\alpha_{1} \colon Q_{1} \to Q'_{1}$ and $\alpha_{2} \colon Q_{2} \to Q'_{2}$, which determine ring homomorphisms $\bZ[\alpha_{1}^{\opp}] \colon \bZ[Q_{1}^{\opp}] \to \bZ[(Q'_{1})^{\opp}]$ and $\bZ[\alpha_{2}] \colon \bZ[Q_{2}] \to \bZ[Q'_{2}]$.
Let us write $\xi_{\phi} \colon X^\phi \rightarrow X$ and $\xi_{\phi'} \colon Y^{\phi'} \rightarrow Y$ for the regular covering spaces corresponding to $\phi$ and $\phi'$ respectively. By covering space theory, for each point $\widetilde{y} \in \xi_{\phi'}^{-1}(y_{0})$, there is a unique continuous map $X^\phi \to Y^{\phi'}$ that lifts the composition $f \circ \xi_{\phi} \colon X^\phi \to Y$ and that takes $\widetilde{x}_{0}$ to $\widetilde{y}$. We therefore obtain a uniquely-determined lift $\widetilde{f} \colon X^\phi \to Y^{\phi'}$ by requiring that $\widetilde{f}(\widetilde{x}_{0}) = \widetilde{y}_{0}$. Extending this map $\bZ$-linearly in each fibre results in a map
\[
\bZ_{\fib} [\widetilde{f}] \colon \bZ_{\fib}[X^\phi] \longrightarrow \bZ_{\fib}[Y^{\phi'}]
\]
of bundles of $\bZ$-modules. Finally, one may check that $\bZ_{\fib} [\widetilde{f}]$ is a map of bundles of $(\bZ[Q_{1}^{\opp}],\bZ[Q_{2}])$-bimodules, covering $f$, from $\bZ_{\fib}[\xi_{\phi}]$ to $\bZ_{\fib}[\xi_{\phi'}]$, where the latter is given the structure of a bundle of $(\bZ[Q_{1}^{\opp}],\bZ[Q_{2}])$-bimodules via $\bZ[\alpha_{1}^{\opp}]$ and $\bZ[\alpha_{2}]$. This uses the interpretation of morphisms of $\toplr$ from Remark~\ref{r:bundles}. Hence we may define
\[
\bZ[-](f) = (f,\bZ[\alpha_{1}^{\opp}],\bZ[\alpha_{2}],\bZ_{\fib}[\widetilde{f}]).
\]
The assignment for $\bZ[-]$ on morphisms obviously satisfies the identity axiom, and one may straightforwardly check using the unique path lifting property that the composition axiom holds for the fourth component $\bZ_{\fib}[\widetilde{f}]$ of $\bZ[-](f)$, this axiom being evident for the other components.
Finally, the continuity of the linearisation functor $\bZ[-] \colon \covlr \too \toplr$ is immediate from the fact that it acts by the identity on the part of the input data given by the continuous map $f \colon X \to X'$.
\end{proof}

\subsubsection{Fibrewise tensor product}\label{sss:fibrewise_tensor}

In this subsection, we define the \emph{fibrewise tensor product} continuous functor
\begin{equation}
\label{eq:fibrewise-tensor-product}
\otimes \colon \topmodlr \too \toplr
\end{equation}
via the following assignments:
\begin{itemizeb}
\item (\emph{On objects.})
We recall that an object of $\topmodlr$ consists of a space $X$, three rings $R,S,T$, a bundle $\xi \colon \pi_{\leq 1}(X) \to \modlr[R][S]$ of $(R,S)$-bimodules over $X$ and an $(S,T)$-bimodule $V$. Define its image under \eqref{eq:fibrewise-tensor-product} to be the following object of $\toplr$:
\[
(X,R,T, \pi_{\leq 1}(X) \xrightarrow{\;\xi\;} \modlr[R][S] \xrightarrow{- \otimes_S V} \modlr[R][T]).
\]
\item (\emph{On morphisms.})
A morphism of $\topmodlr$ from the object $(X,R,S,T,\xi,V)$ to the object $(X',R',S',T',\xi',V')$ consists of a continuous map $f \colon X \to X'$, ring homomorphisms $\alpha \colon R \to R'$, $\beta \colon S \to S'$ and $\gamma \colon T \to T'$, a natural transformation $\tau \colon \xi \Rightarrow (\alpha,\beta)^{*} \circ \xi'\circ \pi_{\leq 1}(f)$ and a homomorphism $\theta \colon V \to (\beta,\gamma)^{*}V'$ of $(S,T)$-bimodules, where $(\alpha,\beta)^{*}$ and $(\beta,\gamma)^{*}$ denote the restriction functors $\modlr[R'][S'] \to \modlr[R][S]$ and $\modlr[S'][T'] \to \modlr[S][T]$ respectively. Define its image under \eqref{eq:fibrewise-tensor-product} to be the morphism
\[
(f,\alpha,\gamma,\hat{\tau})
\]
of $\toplr$, where $\hat{\tau}$ is the natural transformation given by pasting the following diagram:
\begin{equation}
\label{eq:TopMod-pasting-diagram}
\begin{split}
\begin{tikzpicture}
[x=1mm,y=1mm]
\node (tl) at (0,15) {$\pi_{\leq 1}(X)$};
\node (tm) at (30,15) {$\modlr[R][S]$};
\node (tr) at (60,15) {$\modlr[R][T]$};
\node (bl) at (0,0) {$\pi_{\leq 1}(X')$};
\node (bm) at (30,0) {$\modlr[R'][S']$};
\node (br) at (60,0) {$\modlr[R'][T']$.};
\draw[->] (tl) to node[above,font=\small]{$\xi$} (tm);
\draw[->] (tm) to node[above,font=\small]{$- \otimes_S V$} (tr);
\draw[->] (bl) to node[below,font=\small]{$\xi'$} (bm);
\draw[->] (bm) to node[below,font=\small]{$- \otimes_{S'} V'$} (br);
\draw[->] (tl) to node[left,font=\small]{$\pi_{\leq 1}(f)$} (bl);
\draw[->] (bm) to node[right,font=\small]{$(\alpha,\beta)^{*}$} (tm);
\draw[->] (br) to node[right,font=\small]{$(\alpha,\gamma)^{*}$} (tr);
\node at (15,7.5) {\rotatebox{270}{$\Rightarrow$}};
\node at (45,7.5) {\rotatebox{-45}{$\Rightarrow$}};
\end{tikzpicture}
\end{split}
\end{equation}
Here, the left-hand natural transformation is $\tau$ and the right-hand natural transformation, for each $(R',S')$-bimodule $W$, is given by the $(R,T)$-bimodule homomorphism
\begin{equation}
\label{eq:RT-bimodule-hom}
\begin{split}
\begin{tikzpicture}
[x=1mm,y=1mm]
\node (tl) at (tr.west) [anchor=east] {$(\alpha,\beta)^*(W) \otimes_S V \;=\;$};
\node (tr) at (0,0) [anchor=west] {${}_R R' \otimes_{R'} W \otimes_{S'} S' \otimes_S V_T$};
\node (mr) at (0,-15) [anchor=west] {${}_R R' \otimes_{R'} W \otimes_{S'} S' \otimes_S S' \otimes_{S'} V' \otimes_{T'} T'_T$};
\node (br) at (0,-30) [anchor=west] {${}_R R' \otimes_{R'} W \otimes_{S'} S' \otimes_{S'} V' \otimes_{T'} T'_T$};
\node (bbr) at (0,-40) [anchor=west] {${}_R R' \otimes_{R'} W \otimes_{S'} V' \otimes_{T'} T'_T,$};
\node (bbl) at (bbr.west) [anchor=east] {$(\alpha,\gamma)^*(W \otimes_{S'} V') \;=\;$};
\draw[->] (22,-3) to node[right,font=\small]{$({}_R R' \otimes_{R'} W \otimes_{S'} S') \otimes_S \theta$} (22,-12);
\draw[->] (22,-18) to node[right,font=\small]{$({}_R R' \otimes_{R'} W) \otimes_{S'} \mu \otimes_{S'} (V' \otimes_{T'} T'_T)$} (22,-27);
\node at (22,-35) {\rotatebox{270}{$\cong$}};
\end{tikzpicture}
\end{split}
\end{equation}
where the $(S',S')$-bimodule homomorphism $\mu \colon S' \otimes_S S' \to S'$ is given by multiplication in the ring $S'$.
\end{itemizeb}

\begin{lemm}
The above assignments for \eqref{eq:fibrewise-tensor-product} define a continuous functor $\topmodlr \to \toplr$. Its restriction to the subcategory $\topmodlrpr \subset \topmodlr$ takes values in the subcategory $\topr^{\pr} \subset \topr$.
\end{lemm}
\begin{proof}
The diagram \eqref{eq:TopMod-pasting-diagram} for the identity morphism of $(X,R,S,T,\xi,V)$ is the identity natural transformation, so \eqref{eq:fibrewise-tensor-product} satisfies the identity axiom.
The fact that \eqref{eq:fibrewise-tensor-product} satisfies the composition axiom is clear, since composition in $\toplr$ is given by vertical pasting of diagrams of the form \eqref{eq:Top-composition-diagram}, and composition in $\topmodlr$ is given by a similar vertical pasting of diagrams.
The continuity of the functor $\otimes$ is obvious since it acts by the identity on the part of the input data given by the continuous map $f \colon X \to X'$.
The statement about the restriction to $\topmodlrpr \subset \topmodlr$ holds since the underlying continuous map is not changed by \eqref{eq:fibrewise-tensor-product}.
\end{proof}

\begin{notation}
\label{r:tensor-notation}
There are forgetful functors
\[
\topmodlr \too \toplr \qquad\text{and}\qquad \topmodlr \too \modlr
\]
coming from the pullback square \eqref{eq:topmod}. For a continuous functor $F \colon \cC \to \topmodlr$, we denote its compositions with these two forgetful functors by
\[
F_{1} \colon \cC \too \toplr \qquad\text{and}\qquad F_{2} \colon \cC \too \modlr
\]
respectively. In this notation, we denote the composite functor $\otimes \circ F = \eqref{eq:fibrewise-tensor-product} \circ F$ by $F_{1} \otimes F_{2}$.
\end{notation}

\subsubsection{Twisted homology functor}\label{sss:twisted-homology}

Over a fixed ring $R$, we view local coefficient systems on a space $X$ as bundles of right $R$-modules over $X$. From this viewpoint, homology with local coefficients (defined over $R$ and in degree $i\geq 0$) is a continuous functor
\begin{equation}
\label{eq:twisted-homology-R}
H_{i} \colon \topr[R] \too \modr[R].
\end{equation}
It is defined on objects by sending $(X,\xi)$ to the $R$-module $H_{i}(X;\cL_{\xi})$, where $\cL_{\xi}$ is the local system defined by the bundle of $R$-modules $\xi\colon \pi_{\leq 1}(X)\to \modlr[][R]$. See, for example, \cite[\S 5.4]{daviskirk} or \cite[\S 5.1]{Palmer2018Twistedhomologicalstability}. The following fact may be proven by directly generalising the usual construction of singular twisted homology, keeping careful track of the variable bimodule structure.

\begin{prop}
\label{prop:twisted-homology}
In any degree $i\geq 0$, homology with local coefficients extends to a continuous functor
\begin{equation}
\label{eq:twisted-homology}
H_{i} \colon \toplr \too \modlr
\end{equation}
such that the square
\begin{equation}
\label{eq:twisted-homology-compatibility}
\centering
\begin{split}
\begin{tikzpicture}
[x=1mm,y=1mm]
\node (tl) at (0,15) {$\toplr$};
\node (tr) at (30,15) {$\modlr$};
\node (bl) at (0,0) {$\topr[R]$};
\node (br) at (30,0) {$\modr[R]$};
\incl{(bl)}{(tl)}
\incl{(br)}{(tr)}
\draw[->] (tl) to node[above,font=\small]{$H_{i}$} (tr);
\draw[->] (bl) to node[above,font=\small]{$H_{i}$} (br);
\end{tikzpicture}
\end{split}
\end{equation}
commutes for any ring $R$. These statements repeat mutatis mutandis for Borel-Moore homology with local coefficients, providing a continuous functor $H^{\BM}_{i} \colon \toplr^{\pr} \to \modlr^{\pr}$ in any degree $i\geq 0$.
\end{prop}

\begin{proof}
The continuous functor $H_{i}$ is defined on objects by sending $(X,R,S,\xi)$ to $(R,S,H_{i}(X;\cL_{\xi}))$, where $\cL_{\xi}$ is the local system corresponding to the bundle of $(R,S)$-bimodules $\xi$ over $X$. For a morphism $(f,\alpha,\beta,\tau)$ with $\alpha=\id_R$ and $\beta=\id_S$, there is a well-defined map $H_{i}(X;\cL_{\xi}) \to H_{i}(X';\cL_{\xi'})$ of $(R,S)$-bimodules, induced by $f$ for the spaces and by $\tau$ for the local coefficients, constructed in \cite[Th.~5.11~and~Ex.~83]{daviskirk}. For general $\alpha$ and $\beta$, the same construction produces a well-defined map $H_{i}(X;\cL_{\xi}) \to (\alpha,\beta)^* H_{i}(X';\cL_{\xi'})$ of $(R,S)$-bimodules that is compatible with the canonical map $(\alpha,\beta)^* H_{i}(X';\cL_{\xi'}) \to H_{i}(X';\cL_{(\alpha,\beta)^* \xi'})$ of the universal coefficient theorem.
It is straightforward to check from this construction that it respects composition of morphisms.
For Borel-Moore homology with local coefficients, the construction is the straightforward analogue of this, generalising \cite[Chap.~V, \S 4]{bredonsheaf}.
\end{proof}

\begin{rmk}[Twisted homology as homology of covering spaces.]
\label{rmk:covering-spaces}
A regular covering $\rho \colon \hat{X} \to X$ with deck transformation group $G$ induces a bundle of $\bZ[G]$-modules $\xi(\rho) \colon \pi_{\leq 1}(X) \to \modr[{\bZ[G]}]$ by unique path lifting, which in turn determines a local system $\cL_{\xi(\rho)}$ on $X$. By Shapiro's lemma for covering spaces, there is a canonical isomorphism $H_*(X;\cL_{\xi(\rho)}) \cong H_*(\hat{X};\bZ)$. If $\rho$ is a \emph{finite} covering (not in general if it is infinite), we also have $H_*^{\BM}(X;\cL_{\xi(\rho)}) \cong H_*^{\BM}(\hat{X};\bZ)$. In the applications (see \S\ref{s:applications}) of our general construction of homological representations, the relevant local systems will typically arise from infinite regular coverings, so when we use ordinary homology (but not when we use Borel-Moore homology) we may interpret them as actions on the untwisted homology of covering spaces.
\end{rmk}

\subsubsection{The general construction}\label{sss:construction}

We now put together the linearisation functor (\S\ref{sss:lift}), the fibrewise tensor product functor (\S\ref{sss:fibrewise_tensor}) and the twisted homology functor (\S\ref{sss:twisted-homology}) to give a general construction (taking two continuous semifunctors as input) of homological representation functors, described in Theorem~\ref{thm:construction} and Definition~\ref{def:construction} below.

\emph{Throughout \S\ref{sss:construction}, we assume that the projection semifunctor $\cC \to \pi_{0}(\cC)$ is \textbf{continuous}, in other words that the path-components of the morphism spaces of $\cC$ are open.} In particular, this assumption is always satisfied when considering $\cC=\langle \cG, \cM\rangle$ as in Hypothesis~\ref{hypo:standard_framework}; see Remark~\ref{rmk:framework-open-path-components}.
We first record the following observation, saying that continuous semifunctors with discrete codomain factor through $\pi_{0}$ of their domain.

\begin{lemm}\label{lem:factorisation_pi_{0}}
Let $\cC$ be a topologically-enriched semicategory and $\cD$ a discrete semicategory.
Any continuous semifunctor $M\colon \cC \to \cD$ factors uniquely into the projection semifunctor $\cC \to \pi_{0}(\cC)$ followed by a semifunctor $\pi_{0}(M)\colon\pi_{0}(\cC) \to \cD$.
\end{lemm}
\begin{proof}
The categories $\cC$ and $\pi_{0}(\cC)$ have the same objects, so we may define $\pi_{0}(M)(c):=M(c)$ for each $c\in\obj(\cC)$. For each pair of objects $(x,y)$ of $\cC$ and morphisms $\phi$ and $\phi'$ in $\cC(x,y)$ such that $\pi_{0}(\phi)=\pi_{0}(\phi')$, we have $M(\phi)=M(\phi')$ since $M$ is continuous and $\cD$ is discrete, so in particular the path-components of its morphism spaces are points. We may therefore define $\pi_{0}(M)$ on morphisms by assigning $\pi_{0}(M)([\phi])=M(\phi)$ for each morphism $[\phi]\in\pi_{0}(\cC)(x,y)$ with $\phi\in\cC(x,y)$ such that $\pi_{0}(\phi)=[\phi]$.
Composition of morphisms and the associativity axiom for $\pi_{0}(M)$ follow from those for $M$.
This defines a semifunctor $\pi_{0}(M)$ with the required property, which is, moreover, uniquely determined by $M$.
\end{proof}

We recall (see Notation~\ref{not:lr_forgetful}) that $\sL$ and $\sR$ denote the semifunctors
\begin{center}
\begin{tikzpicture}
[x=1mm,y=1mm]
\node (l1) at (0,0) [anchor=east] {$\covlr$};
\node (r1) at (10,0) [anchor=west] {$\groups$};
\node (l2) at (70,0) [anchor=east] {$\modlr \subset \toplr$};
\node (r2) at (80,0) [anchor=west] {$\Rings$};
\draw[->,transform canvas={yshift=0.4ex}] (l1) to node[above,font=\small]{$\sL$} (r1);
\draw[->,transform canvas={yshift=-0.4ex}] (l1) to node[below,font=\small]{$\sR$} (r1);
\draw[->,transform canvas={yshift=0.4ex}] (l2) to node[above,font=\small]{$\sL$} (r2);
\draw[->,transform canvas={yshift=-0.4ex}] (l2) to node[below,font=\small]{$\sR$} (r2);
\end{tikzpicture}
\end{center}
that remember just the first (respectively second) underlying group or ring.
We use the same notation for their restrictions to the subcategories $\covlr^{\pr}$ and $\toplr^{\pr}$.

We are now ready to describe the general construction of homological representation functors. The input for the construction consists of:
\begin{itemizeb}
\item[$\circ$] a topologically-enriched (small) semicategory $\cC$,
\item[$\circ$] a continuous semifunctor $F \colon \cC \to \covr$,
\item[$\circ$] a continuous semifunctor $V \colon \cC \to \modlr$,
\item[$\circ$] a positive integer $i\geq 0$,
\end{itemizeb}
where $F$ and $V$ satisfy the following condition.

\begin{condition}
\label{condition:input}
The continuous semifunctors $F$ and $V$ are required to be compatible in the sense that we have
\[
\bZ[\sR\circ F] = \sL \circ V,
\]
where $\bZ[-]$ denotes the group ring functor $\groups \to \Rings$.
\end{condition}

\begin{eg}\label{eg:trivial_fibrewise_tensor_product}
For fixed $F \colon \cC \to \covr$, a choice of $V \colon \cC \to \modlr$ for which Condition~\ref{condition:input} is automatically satisfied, is the \emph{trivial} coefficient system $V = V_{\triv}(F)$, defined as follows.
We denote by $\triv \colon \Rings \to \modlr$ the functor that sends $R$ to the bimodule $(R,R,R)$. The \emph{trivial} coefficient system $V_{\triv}(F)$ is the continuous semifunctor defined as the composition $\triv \circ \bZ[-] \circ \sR \circ F$.
In particular, using the notational convention of Notation~\ref{r:tensor-notation}, we have $(\bZ[-] \circ F) \otimes V_{\triv}(F) \cong \bZ[-] \circ F$; this is the sense in which the choice $V = V_{\triv}(F)$ is \emph{trivial} for $F$.
It follows from this definition that Condition~\ref{condition:input} is satisfied for $F$ and $V = V_{\triv}(F)$.

Another important (and non-trivial) choice of $V \colon \cC \to \modlr$ constructed from $F \colon \cC \to \covr$ (together with a subsemicategory of $\cC$) is the \emph{colimit coefficient system} $V_{\col}(F)$ introduced in \S\ref{sss:col_coeff_unwtisted_rep}. In contrast with $V_{\triv}(F)$, it is not defined for any category $\cC$: it requires one to consider $\cC = \langle \cG , \cM \rangle$, where $\cG$ and $\cM$ are a pair of groupoids as in Hypothesis~\ref{hypo:standard_framework}, satisfying further mild properties; see \S\ref{sss:col_coeff_unwtisted_rep}.
\end{eg}

\begin{thm}\label{thm:construction}
Continuous semifunctors $F$ and $V$ satisfying Condition~\ref{condition:input}  induce a continuous semifunctor $\cC \to \topmodr$. The composition of this with the fibrewise tensor product functor $\otimes$ \textup{(\S\ref{sss:fibrewise_tensor})} and the twisted homology functor $H_{i}$ \textup{(\S\ref{sss:twisted-homology})} then induces a semifunctor $\pi_{0}(\cC) \to \modr$.

If $F$ takes values in the subcategory $\covr^{\pr}$, then the above procedure repeats verbatim using the twisted Borel-Moore homology functor $H^{\BM}_{i}$ in place of $H_i$.
\end{thm}
\begin{proof}
Condition~\ref{condition:input} implies the identity $\sL \circ V = \sR \circ \bZ[-] \circ F$, where $\bZ[-]$ denotes the linearisation functor of Proposition~\ref{prop:lift}. Hence the continuous semifunctors $\bZ[-] \circ F$ and $V$ jointly determine a continuous semifunctor $\cC \to \topmodr$ by the universal property of the pullback, proving the first statement. This may then be composed with the fibrewise tensor product $\otimes$ and twisted homology $H_{i}$ to obtain a continuous semifunctor of the form $\cC \to \modr$. Since $\modr$ is a discrete category, Lemma~\ref{lem:factorisation_pi_{0}} implies that this factors uniquely through a semifunctor $\pi_{0}(\cC) \to \modr$, which is the output of the construction. This may be summarised diagrammatically as follows, where the input is red and the output is blue. (We note that the right-hand square of diagram \eqref{eq:lift-condition} appears as a subdiagram here, where the top $\bZ[-]$ is the linearisation functor introduced in Proposition~\ref{prop:lift} while the bottom $\bZ[-]$ denotes the group ring functor.)
\begin{equation}
\label{eq:construction}
\centering
\begin{split}
\begin{tikzpicture}
[x=1mm,y=1mm]
\node [red] (tl) at (0,24) {$\cC$};
\node [red] (tm) at (25,24) {$\covr$};
\node (tr) at (50,24) {$\topr$};
\node [red] (ml) at (0,12) {$\modlr$};
\node (mm) at (25,12) {$\groups$};
\node (mr) at (50,12) {$\Rings$};
\node (bl) at (75,0) {$\topmodr$};
\node (bm) at (100,0) {$\topr$};
\node [blue] (br) at (120,0) {$\modr$};
\node [blue] (bb) at (0,-12) {$\pi_{0}(\cC)$};
\draw[->,red] (tl) to node[above,font=\small]{$F$} (tm);
\draw[->] (tm) to node[above,font=\small]{$\bZ[-]$} (tr);
\draw[->,red] (tl) to node[left,font=\small]{$V$} (ml);
\draw[->] (tm) to node[left,font=\small]{$\sR$} (mm);
\draw[->] (tr) to node[right,font=\small]{$\sR$} (mr);
\draw[->] (mm) to node[above,font=\small]{$\bZ[-]$} (mr);
\draw[->] (ml) to[out=-20,in=200] node[above,font=\small]{$\sL$} (mr);
\draw[->] (bl) to[out=180,in=-25] (ml);
\draw[->] (bl) to[out=90,in=-10] (tr);
\draw[->] (bl) to node[above,font=\small]{$\otimes$} (bm);
\draw[->] (bm) to node[above,font=\small]{$H_{i}$} (br);
\draw[->] (tl.west) -- (-8,24) -- (-8,-12) -- (bb.west);
\draw[->,blue] (bb) to[out=0,in=190] node[below,font=\small]{$L_{i}(F;V)$} (br);
\node at (67,4) {$\ulcorner$};
\end{tikzpicture}
\end{split}
\end{equation}
If $F$ takes values in the subcategory $\covr^{\pr}$, the alternative construction is similar, restricting both the linearisation functor $\bZ[-]$ and the fibrewise tensor product functor $\otimes$ to the subcategories of their (co)domains with superscript ${}^{\pr}$ and using $H^{\BM}_{i}$ in place of $H_i$.
\end{proof}

\begin{defi}[The general construction.]
\label{def:construction}
Given continuous semifunctors $F$ and $V$ satisfying Condition~\ref{condition:input}, define $L_{i}(F;V)$ to be the semifunctor given by Theorem~\ref{thm:construction}, called the \emph{homological representation semifunctor} associated to $F$ and $V$. Using the notational convention described in Notation~\ref{r:tensor-notation}, this may be written symbolically as
\[
L_{i}(F;V) = H_{i} \circ ((\bZ[-] \circ F) \otimes V).
\]
Similarly, the homological representation semifunctor $H^{\BM}_{i} \circ ((\bZ[-] \circ F) \otimes V)$ given by Theorem~\ref{thm:construction} using Borel-Moore homology is denoted by $L^{\BM}_{i}(F;V)$.
\end{defi}

Given a continuous semifunctor $F \colon \cC \to \covr$, we may apply Theorem~\ref{thm:construction} to $F$, and a fortiori Definition~\ref{def:construction}, along with the coefficient systems of Example~\ref{eg:trivial_fibrewise_tensor_product}: hence the key ingredient that we must construct is the semifunctor $F \colon \cC \to \covr$, which is the aim of \S\ref{ss:homological_representation_functor_motion_groups}--\S\ref{ss:homological_representation_functor_mcg}.

\begin{rmk}\label{rmk:homological_representation_semifunctor_genuine}
Although the output $L_{i}(F;V)$ of Definition~\ref{def:construction} is a \emph{semifunctor}, we check in \S\ref{sss:F-properties} that it is a \emph{genuine} functor under some standard conditions; see Corollary~\ref{cor:homological_rep_genuine_functor}. In particular, this will be the case in all of the situations that we address in our applications of \S\ref{s:applications}.
\end{rmk}

\subsubsection{Functorial quotients of groups}\label{sss:functorial_quotient_group}

The geometric input for the general construction of \S\ref{sss:construction} is the semifunctor $F \colon \cC \to \covr$; this will be the subject of \S\ref{ss:homological_representation_functor_motion_groups} and \S\ref{ss:homological_representation_functor_mcg}. The definition of $F$ will depend on a number of parameters, one of which is a \emph{functorial quotient of groups}. Here, we introduce this notion and study its behaviour with respect to split short exact sequences.
Recall that $\groups$ denotes the category of groups.

\begin{defi}\label{def:functorial_quotient_group}
A \emph{functorial quotient of groups} is a functor $Q \colon \groups \to \groups$ which is equipped with a natural transformation $q \colon \id_{\groups} \Rightarrow Q$ such that $q(G) \colon G \to Q(G)$ is a quotient homomorphism for each group $G$.
\end{defi}

\begin{eg}[Canonical examples.]
\label{eg:examples_FQG}
The trivial functor $0$ (constant at the trivial group) and the identity functor $\id_{\groups}$ provide obvious examples of functorial quotients of groups.
Some more sophisticated examples of functorial quotients of groups naturally arise in families. A \emph{functorial descending normal series} is a sequence indexed by integers $\ell\geq 1$ of endofunctors $s_{\ell} \colon \groups \to \groups$ with $s_{1} = \id_\groups$, such that $s_{\ell + 1}(G) \subseteq s_{\ell}(G)$ and these inclusions assemble to a natural transformation $s_{\ell + 1} \Rightarrow s_{\ell}$ for each $\ell\geq 1$, and each inclusion $s_{\ell}(G) \subseteq G$ is normal.
For each $\ell$, we then define the associated functorial quotient of groups $Q_{s_{\ell}} \colon \groups \to \groups$ by the natural transformation $G \mapsto G/s_{\ell}(G)$.
For instance, the lower central series $\{ \LCS_{\ell}(-)\} _{\ell\geq 1}$ and the derived series are functorial descending normal series. In particular, we will use functorial quotients of groups defined from the lower central series in the applications of \S\ref{s:applications}.
\end{eg}

We now prove two crucial properties of functorial quotients of groups, concerning their interactions with split short exact sequences. They will have a decisive role in the underlying structure of the homological representation functors of \S\ref{ss:homological_representation_functor_motion_groups} and \S\ref{ss:homological_representation_functor_mcg}.

\begin{lemm}
\label{lem:6-term-sequence}
For any functorial quotient of groups $Q \colon \groups \to \groups$ and split short exact sequence
\begin{equation}\label{eq:SES_functorial_quotient_groups}
\begin{split}
\begin{tikzpicture}
[x=1mm,y=1mm]
\node (t1) at (5,12) {$1$};
\node (t2) at (20,12) {$A$};
\node (t3) at (40,12) {$B$};
\node (t4) at (60,12) {$C$};
\node (t5) at (75,12) {$1,$};
\draw[->] (t1) to (t2);
\draw[->] (t2) to node[above,font=\small]{$f$} (t3);
\draw[->] (t3) to node[below,font=\small]{$g$} (t4);
\draw[->] (t4) to (t5);
\draw[->,densely dotted] (t4) to[out=150,in=30] (t3);
\end{tikzpicture}
\end{split}
\end{equation}
we have $\im(q(B) \circ f \colon A \to B \twoheadrightarrow Q(B)) = \ker(Q(g) \colon Q(B) \to Q(C))$. Denoting this group by $Q^{*}(A)$, we have an induced commutative diagram
\begin{equation}\label{eq:diagram_functorial_quotient_groups}
\begin{split}
\begin{tikzpicture}
[x=1mm,y=1mm]
\node (t1) at (5,12) {$1$};
\node (t2) at (20,12) {$A$};
\node (t3) at (40,12) {$B$};
\node (t4) at (60,12) {$C$};
\node (t5) at (75,12) {$1$};
\node (b1) at (5,0) {$1$};
\node (b2) at (20,0) {$Q^{*}(A)$};
\node (b3) at (40,0) {$Q(B)$};
\node (b4) at (60,0) {$Q(C)$};
\node (b5) at (75,0) {$1$};
\draw[->] (t1) to (t2);
\draw[->] (t2) to node[above,font=\small]{$f$} (t3);
\draw[->] (t3) to node[below,font=\small]{$g$} (t4);
\draw[->] (t4) to (t5);
\draw[->] (b1) to (b2);
\draw[->] (b2) to (b3);
\draw[->] (b3) to (b4);
\draw[->] (b4) to (b5);
\draw[->>] (t2) to (b2);
\draw[->>] (t3) to (b3);
\draw[->>] (t4) to (b4);
\draw[->,densely dotted] (t4) to[out=150,in=30] (t3);
\draw[->,densely dotted] (b4) to[out=150,in=30] (b3);
\end{tikzpicture}
\end{split}
\end{equation}
in which both rows are split short exact sequences and the middle and right-hand vertical maps are given by the natural transformation $q$. Furthermore, there is an action of $C$ on the quotient $Q^{*}(A)$ induced by the conjugation action of the semi-direct product defined by \eqref{eq:SES_functorial_quotient_groups}.
\end{lemm}
\begin{proof}
The second statement is clear, given the first one: we define the right-hand square by applying the functor $Q$ to the map $g$ and its given section and by applying the natural transformation $q$ to the groups $B$ and $C$. We then fill in the left-hand square by factoring $q(B) \circ f$ uniquely as a surjection followed by an injection. The final thing to check is exactness in the middle of the bottom row, which is precisely the first statement of the lemma.

To prove the first statement, first note that the inclusion $\im(q(B) \circ f) \subseteq \ker(Q(g))$ follows immediately from exactness of the top row. To prove the opposite inclusion, let $x \in Q(B)$ with $Q(g)(x) = 1$; we need to find a lift $y \in B$ of $x$ such that $y \in f(A)$. To do this, first pick any lift $y' \in B$ of $x$ and set $z = g(y') \in C$. Denoting the given section of $g$ by $s$, note that $s(z)$ projects to $1 \in Q(B)$, since $z$ projects to $1 \in Q(C)$. Thus $y = y' \cdot s(z)^{-1} \in B$ is another lift of $x$, and moreover $g(y) = g(y') \cdot z^{-1} = z \cdot z^{-1} = 1$, so $y \in f(A)$ by exactness of the top row.

Finally, the action of $C$ on $Q^{*}(A)$ is formally defined either by lifting elements of $Q^{*}(A)$ along the vertical map $A\twoheadrightarrow Q^{*}(A)$ and using the conjugation action of the semi-direct product on the top row of  \eqref{eq:diagram_functorial_quotient_groups}, or equivalently by projecting along the vertical map $C\twoheadrightarrow Q(C)$ to the bottom-right group of \eqref{eq:diagram_functorial_quotient_groups} and then using the conjugation action of the semi-direct product on the bottom row.
\end{proof}

\begin{notation}
We denote by $Q^{*}(A)_{C}$ the \emph{coinvariants} of the group $Q^{*}(A)$ under the action of $C$ shown in Lemma~\ref{lem:6-term-sequence}, i.e.~the largest quotient of $Q^{*}(A)$ that collapses the orbits of the $C$-action.
\end{notation}

Let us use the notation $[A,B,C]$ for a split short exact sequence of the form \eqref{eq:SES_functorial_quotient_groups}. Recall that a morphism of split short exact sequences $[A,B,C]\to [A',B',C']$ is a triple $\Psi=[\Psi_{A},\Psi_{B},\Psi_{C}]$ of group homomorphisms $\Psi_{A}\colon A\to A'$, $\Psi_{B}\colon B\to B'$ and $\Psi_{C}\colon C\to C'$ making the obvious diagram commute. In particular, we denote by $[q^{*}(A),q(B),q(C)]$ the morphism of split short exact sequences defined by diagram \eqref{eq:diagram_functorial_quotient_groups}. In this notation, we have the following naturality property of the construction of Lemma~\ref{lem:6-term-sequence} illustrated in diagram \eqref{eq:diagram_functorial_quotient_groups}:

\begin{lemm}\label{lem:naturality_functorial_quotient_group}
For a morphism of split short exact sequences $\Psi\colon [A,B,C]\to [A',B',C']$, there exists a unique group homomorphism $\phi\colon Q^{*}(A)\to Q^{*}(A')$ such that there is a morphism of split short exact sequences $Q(\Psi):=[\phi,Q(\Psi_{B}),Q(\Psi_{C})]$ making the following diagram commute:
\begin{equation}
\label{eq:naturality_functorial_quotient_group}
\centering
\begin{split}
\begin{tikzpicture}
[x=1mm,y=1mm]
\node (tl) at (0,15) {$[A,B,C]$};
\node (tr) at (60,15) {$[A',B',C']$};
\node (bl) at (0,0) {$[Q^{*}(A),Q(B),Q(C)]$};
\node (br) at (60,0) {$[Q^{*}(A'),Q(B'),Q(C')]$.};
\draw[->] (tl) to node[above,font=\small]{$\Psi$} (tr);
\draw[->] (bl) to node[above,font=\small]{$Q(\Psi)$} (br);
\draw[->>] (tl) to node[left,font=\small]{$[q^{*}(A),q(B),q(C)]$} (bl);
\draw[->>] (tr) to node[right,font=\small]{$[q^{*}(A'),q(B'),q(C')]$} (br);
\end{tikzpicture}
\end{split}
\end{equation}
Furthermore, there exists a unique group homomorphism $\phi_{\mathrm{coinv}} \colon Q^{*}(A)_{C}\to Q^{*}(A')_{C'}$ such that the following diagram (where the vertical maps are the canonical projections onto the coinvariants) is commutative:
\begin{equation}
\label{eq:naturality_functorial_quotient_group_untwisted}
\centering
\begin{split}
\begin{tikzpicture}
[x=1mm,y=1mm]
\node (tl) at (0,13) {$Q^{*}(A)$};
\node (tr) at (30,13) {$Q^{*}(A')$};
\node (bl) at (0,0) {$Q^{*}(A)_{C}$};
\node (br) at (30,0) {$Q^{*}(A')_{C'}$.};
\draw[->] (tl) to node[above,font=\small]{$\phi$} (tr);
\draw[->] (bl) to node[above,font=\small]{$\phi_{\mathrm{coinv}}$} (br);
\draw[->>] (tl) to node[left,font=\small]{} (bl);
\draw[->>] (tr) to node[right,font=\small]{} (br);
\end{tikzpicture}
\end{split}
\end{equation}
\end{lemm}
\begin{proof}
Since $q \colon \id \Rightarrow Q$ is a natural transformation, each of the morphisms $g\colon B\to C$, $g'\colon B'\to C'$, $\Psi_{B}\colon B\to B'$ and $\Psi_{C}\colon C\to C'$ induces a commutative square of the form of the right-hand side of \eqref{eq:diagram_functorial_quotient_groups}. We then obtain a cubical commutative diagram with these four commutative squares along with the one defined by $(g,g',\Psi_{B},\Psi_{C})$ (the top of the cube) and its image under the functor $Q$ (the bottom of the cube), which is of the form
\begin{equation}
\label{eq:naturality_functorial_quotient_group_proof}
\centering
\begin{split}
\begin{tikzpicture}
[x=1mm,y=1mm]
\node (tl) at (0,12) {$Q(B)$};
\node (tr) at (24,12) {$Q(C)$};
\node (bl) at (0,0) {$Q(B')$};
\node (br) at (24,0) {$Q(C'),$};
\draw[->>] (tl) to node[above,font=\small]{$Q(g)$} (tr);
\draw[->>] (bl) to node[above,font=\small]{$Q(g')$} (br);
\draw[->] (tl) to node[left,font=\small]{$Q(\Psi_{B})$} (bl);
\draw[->] (tr) to node[right,font=\small]{$Q(\Psi_{C})$} (br);
\draw[->,densely dotted] (tr) to[out=140,in=40] (tl);
\draw[->,densely dotted] (br) to[out=140,in=40] (bl);
\end{tikzpicture}
\end{split}
\end{equation}
where the dotted arrows are the images under $Q$ of the given sections of $g$ and $g'$. Commutativity of \eqref{eq:naturality_functorial_quotient_group_proof} along with the universal property of $Q^{*}(A')$ as the kernel of $Q(g')$ implies the existence of a unique group homomorphism $\phi \colon Q^{*}(A)\to Q^{*}(A')$ such that $Q(\Psi):=[\phi,Q(\Psi_{B}),Q(\Psi_{C})]$ is a well-defined morphism of split short exact sequences $[Q^{*}(A),Q(B),Q(C)] \to [Q^{*}(A'),Q(B'),Q(C')]$.
Moreover, since the cube diagram of which \eqref{eq:naturality_functorial_quotient_group_proof} is the bottom face is commutative, the universal property of $Q^{*}(A')$ as the kernel of $Q(g')$ ensures the existence of a unique group homomorphism $\phi'\colon A\to Q^{*}(A')$ such that
\begin{equation}
\label{eq:morphism_of_split_ses}
[\phi',q(B')\circ \Psi_{B},q(C')\circ \Psi_{C}]=[\phi',Q(\Psi_{B})\circ q(B),Q(\Psi_{C})\circ q(C)]
\end{equation}
is a well-defined morphism of short exact sequences. The assignment $\phi' = q^*(A') \circ \Psi_{A}$ makes the left-hand side of \eqref{eq:morphism_of_split_ses} a well-defined morphism of short exact sequences, whereas the assignment $\phi' = \phi \circ q^*(A)$ makes the right-hand side a well-defined morphism of short exact sequences. By uniqueness of $\phi'$, we conclude that $q^*(A')\circ \Psi_{A}=\phi'=\phi\circ q^*(A)$, so diagram \eqref{eq:naturality_functorial_quotient_group} is commutative.

For a semi-direct product $H \rtimes G$, we denote the conjugation action $G\times H\to H$ by $\mathrm{conj}_{G}(H)$. By the commutativity of diagram \eqref{eq:naturality_functorial_quotient_group}, we have $\phi\circ\mathrm{conj}_{C}(Q^{*}(A)) = \mathrm{conj}_{C'}(Q^{*}(A'))\circ (\Psi_{C},\phi)$. The universal property of the quotient $Q^{*}(A)_{C}$ then provides the unique group homomorphism $\phi_{\mathrm{coinv}} \colon Q^{*}(A)_{C}\to Q^{*}(A')_{C'}$ making the diagram \eqref{eq:naturality_functorial_quotient_group_untwisted} commute.
\end{proof}

\subsection{First constructions of homological representation functors}\label{ss:homological_representation_functor_motion_groups}

In this section, we describe the first version of our recipe for constructing homological representation functors on the category $\fU\pi_{0}(\cD_{d}) \cong \pi_{0}(\fU\cD_{d})$ (see Lemma~\ref{lem:Serre-fibration-condition}) via the general construction of \S\ref{sss:construction} above. This consists in specifying the input $\cF \colon \fU\cD_{d} \to \covr$ for the general construction, which is done in Theorem~\ref{thm:global_functor_motion_groups} below, and depends on a closed submanifold $Z \subset \bR^d$ and open subgroup $\sG \leq \Diff(Z)$, as well as a functorial quotient of groups $Q$ (see \S\ref{sss:functorial_quotient_group}).

Throughout \S\ref{ss:homological_representation_functor_motion_groups}, we fix an integer $d \geq 2$, $d \neq 4$ (cf.~Remark~\ref{rmk:d_neq_4_explanation}) and consider the Quillen bracket semicategories of manifolds $\fU\cD_{d}$ and $\fU\cD_{d}^{+}$. In addition, we consider groupoids $\cG$ and $\cM$ satisfying the conditions of Hypothesis~\ref{hypo:standard_framework}. By Lemma~\ref{lem:restriction}, there is a natural inclusion of semicategories
\begin{equation}
\label{eq:restriction_of_global_functor}
\langle \cG , \cM \rangle \longhookrightarrow \fU\cD_{d},
\end{equation}
or $\langle \cG , \cM \rangle \hookrightarrow \fU\cD_{d}^{+}$ when considering \emph{orientedly} decorated $d$-manifolds.
We recall that this is given by the inclusion of the objects of $\cM$ into the objects of $\cD_{d}$ (or $\cD_{d}^{+}$), while on morphisms it is an inclusion of embedding spaces, under the identification of Corollary~\ref{cor:description_UD_{d}}. In this section, we construct \emph{global} semifunctors on the categories $\fU\cD_{d}$ or $\fU\cD_{d}^{+}$. We may then restrict them to subcategories of the form $\langle \cG , \cM \rangle$ along \eqref{eq:restriction_of_global_functor} in order to consider the homological representations for a specific family of groups; see \S\ref{ss:applications_motion_groups}.

We fix a closed submanifold $Z \subset \bR^d$ (that is orientable, if we are working with $\fU\cD_{d}^{+}$). We also consider an open subgroup $\sG$ of $\Diff(Z)$, which must lie in the subgroup $\Diff^{+}(Z) \subset \Diff(Z)$ if we are working with $\fU\cD_{d}^{+}$.
We then construct, in Theorem~\ref{thm:global_functor_motion_groups}, continuous semifunctors associated to each functorial quotient of groups $Q$:
\begin{equation}
\label{eq:global_functor_motion_groups}
\Int{\cF}_{(Z,\sG,Q)} \,\text{ and }\, \cF_{(Z,\sG,Q)} \colon \fU\cD_{d} \too \covr ,
\end{equation}
or of the form $\fU\cD_{d}^{+} \to \covr$ when considering \emph{orientedly} decorated $d$-manifolds.
We also construct variants $\Int{\cF}^{\unt}_{(Z,\sG,Q)}$ and $\cF^{\unt}_{(Z,\sG,Q)}$ by taking \emph{coinvariants} at a certain point during the construction (see \S\ref{sss:F-construction}).
These continuous semifunctors are the essential geometric input in our general construction of \S\ref{sss:construction} (see Definition~\ref{def:construction}).
Their construction occupies \S\ref{sss:F-construction} and their elementary properties are established in \S\ref{sss:F-properties}.

Restricting along \eqref{eq:restriction_of_global_functor} to any subcategory $\langle \cG , \cM \rangle$, we construct in \S\ref{sss:col_coeff_unwtisted_rep} a \emph{colimit} coefficient system $V \colon \langle \cG , \cM \rangle \to \modlr$ adapted to each semifunctor \eqref{eq:global_functor_motion_groups}. These are devised so that, under certain conditions, the construction of \S\ref{sss:construction} leads to \emph{untwisted} functorial homological representations; see Proposition~\ref{prop:factoringthroughtopQ}.

\paragraph*{The homological representation functors.}

Having constructed the continuous semifunctors \eqref{eq:global_functor_motion_groups}, the construction of homological representation functors is then immediate. Namely, we apply the machinery of \S\ref{ss:ingredients_construction_representations} (summarised in Definition~\ref{def:construction}) to the input $F = \eqref{eq:global_functor_motion_groups}$ or one of its variants, possibly restricted along \eqref{eq:restriction_of_global_functor} to a subcategory $\langle \cG , \cM \rangle \subset \cU\cD_{d}$, together with a \emph{coefficient system} $V$ adapted to $F$, for instance the \emph{trivial} coefficient system $V_{\triv}(F)$ of Example~\ref{eg:trivial_fibrewise_tensor_product} or the \emph{colimit} coefficient system $V_{\col}(F)$ introduced in \S\ref{sss:col_coeff_unwtisted_rep} below.
Before stating the core construction of Theorem~\ref{thm:global_functor_motion_groups}, we record this consequence.

\begin{coro}\label{coro:def_homological_rep_functor_motion_groups}
For groupoids $\cG$ and $\cM$ as in Hypothesis~\ref{hypo:standard_framework}, the semifunctors of Theorem~\ref{thm:global_functor_motion_groups} provide semifunctors of the form
\begin{equation}
\label{eq:output_of_general_construction_motion_groups}
\langle \pi_{0}(\cG) , \pi_{0}(\cM) \rangle \too \modr,
\end{equation}
where $V$ is any continuous semifunctor satisfying Condition~\ref{condition:input}, such as $V_{\triv}(F)$ or $V_{\col}(F)$. We denote them by $L_{i}(\Int{\cF}_{(Z,\sG,Q)};V)$, $L_{i}(\Int{\cF}^{\unt}_{(Z,\sG,Q)};V)$, $L_{i}(\cF_{(Z,\sG,Q)};V)$ and $L_{i}(\cF^{\unt}_{(Z,\sG,Q)};V)$. If the semifunctor $\pi_{0}(V)$ induced from $V$ by Lemma~\ref{lem:factorisation_pi_{0}} is a genuine functor, then the semifunctors \eqref{eq:output_of_general_construction_motion_groups} upgrade to genuine functors.
\end{coro}
\begin{proof}
The construction of the semifunctors \eqref{eq:output_of_general_construction_motion_groups} is a direct application of Definition~\ref{def:construction}, where we note that $\pi_{0}(\langle \cG , \cM \rangle)$ is identified canonically with $\langle \pi_{0}(\cG) , \pi_{0}(\cM) \rangle$ by Lemma~\ref{lem:Serre-fibration-condition}.
The upgrade of the resulting homological representation semifunctors into genuine functors is a consequence of Corollary~\ref{cor:homological_rep_genuine_functor} below.
\end{proof}

\subsubsection{The construction of the functors \texorpdfstring{$\Int{\cF}$}{F} and \texorpdfstring{$\cF$}{F}}
\label{sss:F-construction}

The goal of this section is to prove the following result, defining the input for the construction of homological representation functors:

\begin{thm}
\label{thm:global_functor_motion_groups}
For any integer $d \geq 2$, closed submanifold $Z \subset \bR^d$ and open subgroup $\sG$ of $\Diff(Z)$, each functorial quotient of groups $Q$ determines continuous semifunctors
\begin{equation}
\label{eq:global-functor2}
\Int{\cF}_{(Z,\sG,Q)},\,\,\Int{\cF}^{\unt}_{(Z,\sG,Q)},\,\,\cF_{(Z,\sG,Q)}\, \text{ and }\, \cF^{\unt}_{(Z,\sG,Q)} \colon \fU\cD_{d} \too \covr.
\end{equation}
Similar continuous semifunctors are defined on the semicategory $\fU\cD_{d}^{+}$ if $Z$ is orientable and $\sG$ is contained in $\Diff^{+}(Z)$.
\end{thm}

The proof of Theorem~\ref{thm:global_functor_motion_groups} occupies the remainder of \S\ref{sss:F-construction}. We first describe the construction of the semifunctors of Theorem~\ref{thm:global_functor_motion_groups} in its $\Int{\cF}_{(Z,\sG,Q)}$ and $\Int{\cF}^{\unt}_{(Z,\sG,Q)}$ (``open'') variants. After this, we summarise the modifications involved in defining the ``closed'' variants $\cF_{(Z,\sG,Q)}$ and $\cF^{\unt}_{(Z,\sG,Q)}$ of \eqref{eq:global-functor2}. Throughout \S\ref{sss:F-construction}, we focus on the functors defined on the category $\fU\cD_{d}$, with the oriented setting (functors defined on $\fU\cD_{d}^{+}$) being exactly analogous, as we record here:

\begin{prop}
All the work of \S\ref{sss:F-construction} below repeats mutatis mutandis considering $\cD_{d}^{+}$ instead of $\cD_{d}$ and assuming that $Z$ and $\sG$ are orientable and orientation-preserving respectively.
\end{prop}

\paragraph*{The construction of the semifunctor \texorpdfstring{$\Int{\cF}$}{F} on objects.}

We define the continuous semifunctors $\Int{\cF}_{(Z,\sG,Q)}$ and $\Int{\cF}^{\unt}_{(Z,\sG,Q)}$ on the objects of $\fU\cD_{d}$, which are the same as the objects of $\cD_{d}$, namely decorated $d$-dimensional manifolds. Let $(M,A) \in \obj(\cD_{d})$. We shall associate to this:
\begin{enumerate}[noitemsep,label=(\roman*)]
\item\label{objects-step-1} a based, path-connected space $X_{(Z,\sG)}(M,A)$ that admits a universal cover,
\item\label{objects-step-2} surjective homomorphisms $\phi_{(Z,\sG,Q)}(M,A) \colon \pi_{1}(X_{(Z,\sG)}(M,A)) \twoheadrightarrow \cQ_{(Z,\sG,Q)}(M,A)$ and $\phi^{\unt}_{(Z,\sG,Q)}(M,A) \colon \pi_{1}(X_{(Z,\sG)}(M,A)) \twoheadrightarrow \cQ^{\unt}_{(Z,\sG,Q)}(M,A)$.
\end{enumerate}
Together, these data determine two objects of $\covr$ (corresponding to $\Int{\cF}_{(Z,\sG,Q)}$ and $\Int{\cF}^{\unt}_{(Z,\sG,Q)}$ respectively). To simplify the notation, since the choice of $(Z,\sG,Q)$ is fixed throughout this construction, we will drop the subscripts, denoting the space by $X(M,A)$ and the surjective homomorphisms by $\phi(M,A) \colon \pi_{1}(X(M,A)) \to \cQ(M,A)$ and $\phi^{\unt}(M,A) \colon \pi_{1}(X(M,A)) \to \cQ^{\unt}(M,A)$.

\paragraph*{The space.}

Recall from Notation~\ref{not:union-with-Z} that we write $\Breve{M}$ for the interior of $\mbar = M \natural \bB_{1}^d$, where $\natural$ denotes the boundary connected sum along boundary-cylinder-germs (i.e.~the semi-monoidal structure of $\cD_{d}$), and that we have fixed an identification of the interior of $\bB_{1}^d$ with $\bR^d$, so that there is a preferred embedding $\bR^d \hookrightarrow \Breve{M}$. Its image is disjoint from $\Int{M}$, hence in particular disjoint from $A \subset \Int{M}$. Thus, restricting this embedding to $Z \subset \bR^d$, we obtain a preferred embedding $Z \hookrightarrow \Breve{M} \smallsetminus A$, which determines a basepoint of the relative embedding space
\[
\cE_{\sG}(Z , \Breve{M} \smallsetminus A) = \Emb(Z , \Breve{M} \smallsetminus A) / \sG.
\]

\begin{defi}\label{def:space_{X}_M_A}
We define $X(M,A)$ to be the path-component of the space $\cE_{\sG}(Z , \Breve{M} \smallsetminus A)$ that contains the basepoint.
\end{defi}

To complete step \ref{objects-step-1} of the construction, we have to show that $X(M,A)$ admits a universal cover: this is the result of Corollary~\ref{coro:universal-cover} below. To prove this, we need the following point-set topological result:

\begin{lemm}
\label{lem:semi-locally-1-connected}
Let $f \colon X \to Y$ be a surjective fibre bundle and suppose that $X$ is semi-locally simply-connected. Then $Y$ is also semi-locally simply-connected.
\end{lemm}
\begin{proof}
Let $y \in Y$ and let $U$ be an open neighbourhood of $y$ in $Y$. We need to find a smaller open neighbourhood $V \subseteq U$ of $y$ such that any loop in $V$ based at $y$ is nullhomotopic in $U$. First, choose a smaller open neighbourhood $U' \subseteq U$ such that $f$ is trivialisable over $U'$, and choose a trivialisation $\varphi \colon f^{-1}(U') \cong U' \times F$ (where $F = f^{-1}(y)$). Also, since $f$ is surjective, we can choose a point $z \in F$. Since $X$ is semi-locally simply-connected, we may find an open neighbourhood $W \subseteq f^{-1}(U')$ of $\tilde{y} = \varphi^{-1}(y,z)$ such that any loop in $W$ based at $\tilde{y}$ is nullhomotopic in $f^{-1}(U')$. By the definition of the product topology, we may then find open subsets $V \subseteq U'$ and $F' \subseteq F$ such that $y \in V$, $z \in F'$ and $\varphi^{-1}(V \times F') \subseteq W$. Now let $\gamma$ be any loop in $V$ based at $y$. Then $\tilde{\gamma} = \varphi^{-1} \circ (\gamma \times \{z\})$ is a loop in $W$ based at $\tilde{y}$. By above, we may find a nullhomotopy of $\tilde{\gamma}$ in $f^{-1}(U')$. Composing this nullhomotopy with $f$, it becomes a nullhomotopy of $\gamma$ in $U' \subseteq U$.
\end{proof}

\begin{prop}\label{prop:X_locally_path_connected_semi-locally_simply_connected}
The space $X(M,A)$ is locally path-connected and semi-locally simply-connected.
\end{prop}
\begin{proof}
First, note that it suffices to show that the space $\cE_{\sG}(Z,\Breve{M}\smallsetminus A)$ is locally path-connected and semi-locally simply-connected, since $X(M,A)$ is one path-component of this space. There is a quotient map
\begin{equation}
\label{eq:quotient-G}
\Emb(Z,\Breve{M}\smallsetminus A) \longtwoheadrightarrow \cE_{\sG}(Z,\Breve{M}\smallsetminus A).
\end{equation}
Since $Z$ is compact (being a closed submanifold of $\bR^d$), the embedding space $\Emb(Z,\Breve{M}\smallsetminus A)$, equipped with the Whitney topology, is locally contractible; see \cite[p.~281, Coro.~of Prop.~$4'$]{Cerf1961Topologiedecertains}. Thus in particular it is locally path-connected and semi-locally simply-connected.
Local path-connectedness is preserved under taking quotients, so it follows that $\cE_{\sG}(Z,\Breve{M}\smallsetminus A)$ is also locally path-connected.

Semi-local simply-connectedness is not preserved under taking quotients in general, so we shall instead invoke Lemma~\ref{lem:semi-locally-1-connected} above.
The space $\cE_{\sG}(Z,\Breve{M}\smallsetminus A)$ is locally retractile with respect to the left action of the group $\Diff_c(\Breve{M}\smallsetminus A)$ of compactly-supported diffeomorphisms, by \cite[Prop.~4.15]{Palmer2018HomologicalstabilitymoduliI}. (Here we use the assumption that $\sG$ is an \emph{open} subgroup of $\Diff(Z)$.) Moreover, the map \eqref{eq:quotient-G} is equivariant with respect to the left action of $\Diff_c(\Breve{M}\smallsetminus A)$, so by \cite[Th.~A]{Palais1960Localtrivialityof}, the map \eqref{eq:quotient-G} is a fibre bundle. Therefore, Lemma~\ref{lem:semi-locally-1-connected} implies that its target $\cE_{\sG}(Z,\Breve{M}\smallsetminus A)$ is also semi-locally simply-connected.
\end{proof}

\begin{coro}
\label{coro:universal-cover}
The space $X(M,A)$ admits a universal cover.
\end{coro}
\begin{proof}
By classical covering space theory (see \cite[Chap.~1, \S 3]{hatcheralgebraic} for instance), this is equivalent to $X(M,A)$ being path-connected (which is true by definition) and locally path-connected and semi-locally simply-connected (which is true by Proposition~\ref{prop:X_locally_path_connected_semi-locally_simply_connected}).
\end{proof}

\paragraph*{The surjective homomorphism.}

To complete the definition of the functors $\Int{\cF}_{(Z,\sG,Q)}$ and $\Int{\cF}^{\unt}_{(Z,\sG,Q)}$ on objects, we need to specify two quotients $\phi(M,A)$ and $\phi^{\unt}(M,A)$ of the fundamental group $\pi_{1}(X(M,A)) = \pi_{1}(\cE_{\sG}(Z,\Breve{M} \smallsetminus A))$. To do this, we use the split homotopy fibration sequence \eqref{eq:split-fibration-sequence-1} of Proposition~\ref{prop:split-fibration-sequence-1}, which induces the split short exact sequence \eqref{eq:split-ses-1}; this is the top row of diagram \eqref{eq:split-short-exact-sequence-quotient} below. We then apply the natural transformation $q \colon \id \Rightarrow Q$ (which is part of the data of the functorial quotient of groups $Q$) to the middle and right-hand groups, applying Lemma~\ref{lem:6-term-sequence} to obtain the following $6$-term commutative diagram, in which the two rows are split short exact sequences. We define $\phi(M,A)$ to be the left-hand vertical map of this diagram, denoting its target $Q^{*}(\pi_{1}(\cE_{\sG}(Z,\Breve{M}\smallsetminus A)))$ by $\cQ(M,A)$ for brevity.
\begin{equation}
\label{eq:split-short-exact-sequence-quotient}
\centering
\begin{split}
\begin{tikzpicture}
[x=1mm,y=1mm,font=\small]
\node at (10,8) {$\pi_{1}(X(M,A))$};
\node at (10,4) {\rotatebox{270}{$\cong$}};
\begin{scope}
\node (ll) at (-15,0) {$1$};
\node (l) at (10,0) {$\pi_{1}(\cE_{\sG}(Z,\Breve{M}\smallsetminus A))$};
\node (m) at (57,0) {$\pi_{1}(\cE_{\Diff(A) \times \sG}(A \sqcup Z , \Breve{M}))$};
\node (r) at (107,0) {$\pi_{1}(\cE(A,\Breve{M}))$};
\node (rr) at (131,0) {$1$};
\draw[->] (ll) to (l);
\draw[->] (l) to (m);
\draw[->] (m) to (r);
\draw[->] (r) to (rr);
\draw[->,densely dashed] ($ (r.west) + (0,2)  $) to[out=160,in=20] ($ (m.east) + (0,2) $);
\end{scope}
\begin{scope}[yshift=-20mm]
\node (lla) at (-15,0) {$1$};
\node (la) at (10,0) {$\cQ(M,A)$};
\node (ma) at (57,0) {$Q\left(\pi_{1}(\cE_{\Diff(A)\times \sG}(A\sqcup Z,\Breve{M}))\right)$};
\node (ra) at (107,0) {$Q\left(\pi_{1}(\cE(A,\Breve{M}))\right)$};
\node (rra) at (131,0) {$1$};
\draw[->] (lla) to (la);
\draw[->] (la) to (ma);
\draw[->] (ma) to (ra);
\draw[->] (ra) to (rra);
\draw[->,densely dashed] ($ (ra.west) + (0,2)  $) to[out=160,in=20] ($ (ma.east) + (0,2) $);
\end{scope}
\draw[->>] (l) to node[left,font=\footnotesize]{$\phi(M,A)$} (la);
\draw[->>] (m) to node[left,font=\footnotesize]{$q(M,A)$} (ma);
\draw[->>] (r) to node[left,font=\footnotesize]{$\bar{q}(M,A)$} (ra);
\end{tikzpicture}
\end{split}
\end{equation}
The vertical morphisms in \eqref{eq:split-short-exact-sequence-quotient} depend upon $(Z,\sG,Q)$; this has been elided from the notation to avoid cluttering the diagram.

The functor $\Int{\cF}_{(Z,\sG,Q)}$ is then defined on objects by:
\begin{equation}\label{eq:def_Int(F)_objects}
(M,A) \longmapsto (X(M,A),\phi(M,A)).
\end{equation}

\begin{eg}\label{eg:trivial_universal_covers}
When $Q = \id_{\groups}$ we have $\phi(M,A) = \id$, corresponding to the universal covering, while when $Q = 0$ we have $\cQ(M,A) = \{\id\}$, corresponding to the trivial covering.
\end{eg}

\paragraph*{Actions and untwisted quotients.}

By Proposition~\ref{prop:split-ses-functoriality} and Lemma~\ref{lem:naturality_functorial_quotient_group}, the entire commutative diagram \eqref{eq:split-short-exact-sequence-quotient} is functorial in the input $(M,A)$ as an object of $\fU\cD_{d}$. Restricting to the automorphism group $\diffdec(M,A)$ of a single object $(M,A)$ and looking just at the bottom-left group $\cQ(M,A)$ of \eqref{eq:split-short-exact-sequence-quotient}, we obtain an action of $\diffdec(M,A)$ on $\cQ(M,A)$. Concretely, this action is induced via the quotient $\phi(M,A)$ from the action on $\pi_{1}(\cE_{\sG}(Z,\Breve{M}\smallsetminus A))$ induced by post-composition with diffeomorphisms of $(M,A)$ extended by the identity on $\mbar \smallsetminus M$. Since the group $\cQ(M,A)$ is discrete, this action factors through an action of $\pi_{0}(\diffdec(M,A))$ on $\cQ(M,A)$. Pre-composing with the canonical map \eqref{eq:from-motion-to-mcg}, this gives an action of $\pi_{1}(\cE(A,\Int{M}))$ on $\cQ(M,A)$.

We also have an action of $\pi_{1}(\cE(A,\Breve{M}))$ (the top-right group of \eqref{eq:split-short-exact-sequence-quotient}) on $\cQ(M,A)$ (the bottom-left group of \eqref{eq:split-short-exact-sequence-quotient}) described in Lemma~\ref{lem:6-term-sequence}. Concretely, this action is induced by the splittings in diagram \eqref{eq:split-short-exact-sequence-quotient}. Pre-composing with the isomorphism $\pi_{1}(\cE(A,\Int{M})) \cong \pi_{1}(\cE(A,\Breve{M}))$ induced by the inclusion $\Int{M} \subset \Breve{M}$, we thus have another action of $\pi_{1}(\cE(A,\Int{M}))$ on $\cQ(M,A)$.

\begin{prop}
\label{prop:two-actions-on-Q}
The two actions of $\pi_{1}(\cE(A,\Int{M}))$ on $\cQ(M,A)$ described above are equal.
\end{prop}
\begin{proof}
Since both of these actions are induced via the quotient $\phi(M,A)$ from actions defined on $\pi_{1}(\cE_{\sG}(Z,\Breve{M}\smallsetminus A))$, it suffices to prove that the two actions agree on this group, before passing to the quotient. This is precisely part (ii) of Proposition~\ref{prop:two-actions}.
\end{proof}

\begin{notation}\label{nota:untwisted}
We denote by $\cQ^{\unt}(M,A)$ the coinvariants $\cQ(M,A)_{\pi_{1}(\cE(A,\Int{M}))}$ associated to this action, and by $\phi^{\unt}(M,A)$ the composite $\pi_{1}(X(M,A)) \twoheadrightarrow \cQ(M,A) \twoheadrightarrow \cQ^{\unt}(M,A)$, where the first arrow is $\phi(M,A)$ while the second one is the canonical projection onto the coinvariants.
\end{notation}

We then define $\Int{\cF}^{\unt}_{(Z,\sG,Q)}$ on objects by:
\begin{equation}\label{eq:def_Int(F)_unt_objects}
(M,A) \longmapsto (X(M,A),\phi^{\unt}(M,A)).
\end{equation}

\begin{rmk}\label{rmk:justification_untwisted}
By construction, any quotient of $\cQ(M,A)$ on which the group $\pi_{1}(\cE(A,\Int{M}))$ acts trivially factors through the quotient onto $\cQ^{\unt}(M,A)$. In particular, the covering spaces encoded by $\Int{\cF}^{\unt}_{(Z,\sG,Q)}$ have a trivial action of these motion groups.
In this sense they are \emph{untwisted} (in a universal way), whence the superscript ${}^{\unt}$ in the notation for this variant.
\end{rmk}

\paragraph*{The construction of the semifunctor $\Int{\cF}$ on morphisms.}

We now define the semifunctors $\Int{\cF}_{(Z,\sG,Q)}$ and $\Int{\cF}^{\unt}_{(Z,\sG,Q)}$ on morphisms of $\fU\cD_{d}$ and check the composition axiom. Recall from  Corollary~\ref{cor:description_UD_{d}} that the morphisms of $\fU\cD_{d}$ may be described as embeddings of manifolds satisfying the three properties of Definition~\ref{def:embdec}, with composition induced by the evident composition of embeddings (see Definition~\ref{def:decorated-embeddings-composition}); this description will be used throughout this paragraph.
By Proposition~\ref{prop:split-ses-functoriality}, each morphism $\varphi \colon (M,A) \to (N,B)$ in $\fU\cD_{d}$ induces a map of split homotopy fibration sequences $\cE_{\varphi}\colon\eqref{eq:split-fibration-sequence-1}_{(M,A)} \to \eqref{eq:split-fibration-sequence-1}_{(N,B)}$, in particular a map
\[
f_{\varphi} \colon \cE_{\sG}(Z,\Breve{M}\smallsetminus A) \too \cE_{\sG}(Z,\Breve{N}\smallsetminus B).
\]
Notation of the form $\eqref{eq:split-fibration-sequence-1}_{(N,B)}$ means the sequence \eqref{eq:split-fibration-sequence-1} with each instance of $(M,A)$ replaced by $(N,B)$.
This map preserves basepoints and therefore restricts to a based map
\begin{equation}\label{eq:def_f_varphi_final}
f^{X}_{\varphi} \colon X(M,A) \too X(N,B).
\end{equation}

This defines the actions of $\Int{\cF}_{(Z,\sG,Q)}$ and $\Int{\cF}^{\unt}_{(Z,\sG,Q)}$ on the morphism $\varphi$ at the level of spaces. To verify that \eqref{eq:def_f_varphi_final} determines well-defined morphisms in $\covr$ (see Definition~\ref{def:bicoverings}), we must verify that the homomorphism $\pi_{1}(f^X_\varphi)$ descends along the quotients $\phi(M,A)$ and $\phi(N,B)$ and also along the quotients $\phi^{\unt}(M,A)$ and $\phi^{\unt}(N,B)$.

The map of split homotopy fibration sequences $\cE_{\varphi}\colon\eqref{eq:split-fibration-sequence-1}_{(M,A)} \to \eqref{eq:split-fibration-sequence-1}_{(N,B)}$ induces a map of split short exact sequences $\Psi_{\varphi}\colon\eqref{eq:split-ses-1}_{(M,A)} \to \eqref{eq:split-ses-1}_{(N,B)}$. It then follows from Lemma~\ref{lem:naturality_functorial_quotient_group} that $\Psi_{\varphi}$ induces a map of diagrams $\eqref{eq:split-short-exact-sequence-quotient}_{(M,A)} \to \eqref{eq:split-short-exact-sequence-quotient}_{(N,B)}$, in particular homomorphisms
\begin{equation}\label{eq:assignment_morphisms_theta}
\theta_{\varphi} \colon \cQ(M,A) \too \cQ(N,B) \quad \text{ and } \quad \theta^{\unt}_{\varphi} \colon \cQ^{\unt}(M,A) \too \cQ^{\unt}(N,B).    
\end{equation}
In particular, this means that
\[
\theta_{\varphi} \circ \phi(M,A) = \phi(N,B) \circ \pi_{1}(f^{X}_{\varphi})\quad \text{ and } \quad \theta^{\unt}_{\varphi} \circ \phi^{\unt}(M,A) = \phi^{\unt}(N,B) \circ \pi_{1}(f^{X}_{\varphi}),
\]
i.e.~that $\pi_{1}(f^X_\varphi)$ descends along the quotients $\phi(M,A)$ and $\phi(N,B)$ and along the quotients $\phi^{\unt}(M,A)$ and $\phi^{\unt}(N,B)$. Equivalently, it sends $\ker(\phi(M,A))$ into $\ker(\phi(N,B))$ and $\ker(\phi^{\unt}(M,A))$ into $\ker(\phi^{\unt}(N,B))$. So $f^{X}_{\varphi}$ is a morphism in $\covr$ from $(X(M,A),\phi(M,A))$ to $(X(N,B),\phi(N,B))$ and also a morphism in $\covr$ from $(X(M,A),\phi^{\unt}(M,A))$ to $(X(N,B),\phi^{\unt}(N,B))$.

In summary, we define both $\Int{\cF}_{(Z,\sG,Q)}$ and $\Int{\cF}^{\unt}_{(Z,\sG,Q)}$ on morphisms by $\varphi \mapsto f^{X}_{\varphi}$, interpreted as a morphism in $\covr$ in these two different ways.

\paragraph*{Continuity and composition axioms.}

We first prove the continuity of the construction at the level of hom-spaces.

\begin{lemm}\label{lem:continuity_semifunctors}
The function
\[
\fU\cD_{d}((M,A),(N,B)) \too \covr(\Int{\cF}_{(Z,\sG,Q)}(M,A), \Int{\cF}_{(Z,\sG,Q)}(N,B))
\]
and its analogue for $\Int{\cF}^{\unt}$ defined above are continuous functions of hom-spaces.
\end{lemm}
\begin{proof}
The function $\fU\cD_{d}((M,A),(N,B)) \to \Map_*(\cE_{\sG}(Z,\Breve{M}\smallsetminus A), \cE_{\sG}(Z,\Breve{N}\smallsetminus B))$ given by $\varphi \mapsto f_{\varphi}$ is continuous, since it is the adjoint of the function $\embdec(M,N) \times \cE_{\sG}(Z,\Breve{M}\smallsetminus A) \to \cE_{\sG}(Z,\Breve{N}\smallsetminus B)$, which is a quotient of a restriction of the composition map $\Emb((\mbar,A),(\nbar,B)) \times \Emb(Z,\Breve{M}\smallsetminus A) \to \Emb(Z,\Breve{N}\smallsetminus B)$, which is continuous in the Whitney topology since $Z$ is compact (see \cite[\S 2, Prop.~1]{Mather1969}).
Restriction to a subspace is always continuous in the compact-open topology, so it follows that the function $\fU\cD_{d}((M,A),(N,B)) \to \Map_*(X(M,A),X(N,B))$ given by $\varphi \mapsto f^{X}_{\varphi}$ is also continuous.

Now, the hom-space $\covr(\Int{\cF}_{(Z,\sG,Q)}(M,A), \Int{\cF}_{(Z,\sG,Q)}(N,B))$ is by definition the subspace of $\Map_*(X(M,A),X(N,B))$ of those maps whose induced homomorphism on $\pi_{1}$ sends $\ker(\phi(M,A))$ into $\ker(\phi(N,B))$ and we know by construction that the function $\varphi \mapsto f^{X}_{\varphi}$ takes values in this subspace. Hence, by restricting the codomain, we conclude that the function
\[
\fU\cD_{d}((M,A),(N,B)) \too \covr(\Int{\cF}_{(Z,\sG,Q)}(M,A), \Int{\cF}_{(Z,\sG,Q)}(N,B))
\]
is continuous, as claimed. The proof for $\Int{\cF}^{\unt}$ is identical; we just restrict in the last step of the proof to the subspace of $\Map_*(X(M,A),X(N,B))$ of those maps whose induced homomorphism on $\pi_{1}$ sends $\ker(\phi^{\unt}(M,A))$ into $\ker(\phi^{\unt}(N,B))$.
\end{proof}

The composition axiom thus remains the only point to check in order to prove that $\Int{\cF}_{(Z,\sG,Q)}$ and $\Int{\cF}^{\unt}_{(Z,\sG,Q)}$ are well-defined continuous semifunctors. Let us consider morphisms $\varphi \colon (M,A) \to (N,B)$ and $\varphi' \colon (N,B) \to (P,C)$ in $\fU\cD_{d}$, viewed as embeddings via Corollary~\ref{cor:description_UD_{d}}.

\begin{lemm}\label{lem:composition_axiom}
There is an equality $f^{X}_{\varphi'\circ\varphi}=f^{X}_{\varphi'}\circ f^{X}_{\varphi}$.
\end{lemm}
\begin{proof}
The composite embedding $\varphi'\circ\varphi$ clearly induces the same map of split homotopy fibration sequences $\eqref{eq:split-fibration-sequence-1}_{(M,A)} \to \eqref{eq:split-fibration-sequence-1}_{(P,C)}$ as that obtained by concatenating the maps induced by the two embeddings $\varphi$ and $\varphi'$ separately. In particular, there is an equality $f_{\varphi'\circ\varphi}=f_{\varphi'} \circ f_{\varphi}$ of maps $\cE_{\sG}(Z,\Breve{M}\smallsetminus A) \to \cE_{\sG}(Z,\Breve{P}\smallsetminus C)$. This identity is preserved when restricting to path-components of embedding spaces, so we obtain the claimed equality.
\end{proof}

By the definition above, Lemma~\ref{lem:composition_axiom} tells us that $\Int{\cF}_{(Z,\sG,Q)}(\varphi'\circ\varphi)=\Int{\cF}_{(Z,\sG,Q)}(\varphi')\circ \Int{\cF}_{(Z,\sG,Q)}(\varphi)$ and $\Int{\cF}^{\unt}_{(Z,\sG,Q)}(\varphi'\circ\varphi)=\Int{\cF}^{\unt}_{(Z,\sG,Q)}(\varphi')\circ\Int{\cF}^{\unt}_{(Z,\sG,Q)}(\varphi)$; we thus have well-defined semifunctors.

This concludes the proof of Theorem~\ref{thm:global_functor_motion_groups} for the ``open'' variants $\Int{\cF}_{(Z,\sG,Q)}$ and $\Int{\cF}^{\unt}_{(Z,\sG,Q)}$.

\paragraph*{The closed variant $\cF$.}

The construction of the ``closed'' variants $\cF_{(Z,\sG,Q)}$ and $\cF^{\unt}_{(Z,\sG,Q)}$ of Theorem~\ref{thm:global_functor_motion_groups} is a slight modification of the construction above of the ``open'' variants $\Int{\cF}_{(Z,\sG,Q)}$ and $\Int{\cF}^{\unt}_{(Z,\sG,Q)}$.
To do this, we simply replace $\cE_{\sG}(Z,\Breve{M} \smallsetminus A)$ with the slightly larger embedding space $\cE_{\sG}(Z,\mbar \smallsetminus A)$.
This has the effect of replacing $X(M,A)$ with another space $X'(M,A)$ that is homotopy equivalent to $X(M,A)$ by Proposition~\ref{prop:inclusion-heq}. In particular, we have a canonical isomorphism $\pi_{1}(X(M,A)) \cong \pi_{1}(X'(M,A))$.

\begin{rmk}
\label{rmk:not-proper-homotopy-equivalence}
We note that the spaces $X'(M,A)$ and $X(M,A)$, although homotopy equivalent, are not \emph{proper} homotopy equivalent, which is why the open and closed variants of \eqref{eq:global_functor_motion_groups} give rise to different homological representations when using Borel-Moore homology (see Remark~\ref{rmk:invariance_F_open_closed}).
\end{rmk}

The only small technical difficulty arises when checking that $X'(M,A)$ has good local properties, so it admits a universal cover (see Corollary~\ref{coro:universal-cover}): the argument is exactly analogous, except that one has to use a variant of \cite[Prop.~4.15]{Palmer2018HomologicalstabilitymoduliI} where the target manifold is allowed to have non-empty boundary. This may be proved by a small adaptation of the proof of Proposition~4.15 of \cite{Palmer2018HomologicalstabilitymoduliI}, using ideas of \cite{Cerf1961Topologiedecertains} (where all manifolds are allowed to have corners of any codimension).

The quotient $\phi'(M,A) \colon \pi_{1}(X'(M,A)) \twoheadrightarrow \cQ(M,A)$ is defined to be equal to the quotient $\phi(M,A) \colon \pi_{1}(X(M,A)) \twoheadrightarrow \cQ(M,A)$ under the isomorphism $\pi_{1}(X(M,A)) \cong \pi_{1}(X'(M,A))$ induced by the homotopy equivalence of Proposition~\ref{prop:inclusion-heq}.
The semifunctor $\cF_{(Z,\sG,Q)}$ is defined on objects by $(M,A) \mapsto (X'(M,A),\phi'(M,A))$.

To define the \emph{untwisted} version, we quotient further onto the coinvariants of $\cQ(M,A)$ by the action of $\pi_{1}(\cE(A,M))$. This is identical to the action of $\pi_{1}(\cE(A,\Int{M}))$ under the isomorphism $\pi_{1}(\cE(A,\Int{M})) \cong \pi_{1}(\cE(A,M))$ induced by the canonical inclusion $\cE(A,\Int{M}) \subset \cE(A,M)$, so this further quotient is precisely $\cQ^{\unt}(M,A)$. Let us write $\phi^{\unt,\prime}(M,A)$ for the quotient map $\phi^{\unt}(M,A)$ pre-composed with the inverse of the isomorphism $\pi_{1}(X(M,A)) \cong \pi_{1}(X'(M,A))$.
The semifunctor $\cF^{\unt}_{(Z,\sG,Q)}$ is then defined on objects by $(M,A) \mapsto (X'(M,A),\phi^{\unt,\prime}(M,A))$.

On morphisms, the construction is exactly the same as the construction in the ``open'' setting above, using the fact that the actions of embeddings of decorated manifolds by post-composition on $\cE_{\sG}(Z,\Breve{M} \smallsetminus A)$ and on $\cE_{\sG}(Z,\mbar \smallsetminus A)$ commute with the inclusion $\cE_{\sG}(Z,\Breve{M} \smallsetminus A) \subset \cE_{\sG}(Z,\mbar \smallsetminus A)$.
This concludes the proof of Theorem~\ref{thm:global_functor_motion_groups}.

To finish \S\ref{sss:F-construction}, we note the following relationship between the ``open'' and ``closed'' semifunctors that we have constructed. We note that a \emph{homotopy} in $\covr$ between two morphisms $f_{1},f_{2} \colon (X,\phi \colon \pi_{1}(X) \twoheadrightarrow Q) \to (X',\phi' \colon \pi_{1}(X') \twoheadrightarrow Q')$ is a homotopy of based maps $X \to X'$ from $f_{1}$ to $f_{2}$; in particular it can only exist if the homomorphisms $Q \to Q'$ induced by $f_{1}$ and $f_{2}$ are equal. Then, a \emph{homotopy equivalence} in $\covr$ from $(X,\phi \colon \pi_{1}(X) \twoheadrightarrow Q)$ to $(X',\phi' \colon \pi_{1}(X') \twoheadrightarrow Q')$ is a based homotopy equivalence $X \to X'$ such that the induced homomorphism $Q \to Q'$ is an isomorphism. Finally, a \emph{natural homotopy equivalence} is a natural transformation between functors with $\covr$ as target category that is componentwise a homotopy equivalence.

\begin{lemm}
\label{lem:natural-homotopy-equivalence}
There are natural homotopy equivalences between the semifunctors of Theorem~\ref{thm:global_functor_motion_groups} $\Int{\cF}_{(Z,\sG,Q)} \Rightarrow \cF_{(Z,\sG,Q)}$ and $\Int{\cF}^{\unt}_{(Z,\sG,Q)}\Rightarrow \cF^{\unt}_{(Z,\sG,Q)}$.
\end{lemm}

\begin{rmk}
\label{rmk:invariance_F_open_closed}
These natural homotopy equivalences imply that, when using \emph{ordinary} twisted homology (or any other homotopy invariant flavour of twisted homology) in our general construction, it does not matter whether we use the \emph{open} or the \emph{closed} variants of the semifunctors. However, this remark does not apply if we use \emph{Borel-Moore} homology; see Remark~\ref{rmk:not-proper-homotopy-equivalence}.
\end{rmk}

\begin{proof}[Proof of Lemma~\ref{lem:natural-homotopy-equivalence}]
Our goal is to construct based homotopy equivalences $X(M,A) \to X'(M,A)$, naturally in $(M,A) \in \fU\cD_{d}$, so that the two induced endomorphisms $\cQ(M,A) \to \cQ(M,A)$ and $\cQ^{\unt}(M,A) \to \cQ^{\unt}(M,A)$ are isomorphisms (in fact they will both be the identity).

By Proposition~\ref{prop:inclusion-heq}, the inclusion $\cE_{\sG}(Z,\Breve{M} \smallsetminus A) \subset \cE_{\sG}(Z,\mbar \smallsetminus A)$ is a based homotopy equivalence. It is clearly natural with respect to embeddings of decorated manifolds, i.e.~natural in $(M,A) \in \fU\cD_{d}$.
Restricting to path-components of the basepoints, we obtain the desired natural based homotopy equivalence $X(M,A) \to X'(M,A)$. The fact that it induces the identity endomorphism of $\cQ(M,A)$ is immediate, since the quotient $\phi'(M,A)$ was defined to coincide with the quotient $\phi(M,A)$ under the identification $\pi_{1}(X(M,A)) \cong \pi_{1}(X'(M,A))$ induced by this homotopy equivalence. For the same reason ($\phi^{\unt,\prime}(M,A)$ was defined to coincide with $\phi^{\unt}(M,A)$ under this identification), the induced endomorphism of $\cQ^{\unt}(M,A)$ is also the identity.
\end{proof}

\subsubsection{Elementary properties}
\label{sss:F-properties}

We now record some fundamental properties of the continuous semifunctors of Theorem~\ref{thm:global_functor_motion_groups}.

\paragraph*{Uniqueness up to isotopy.}

The continuous semifunctors constructed in Theorem~\ref{thm:global_functor_motion_groups} depend only on the isotopy class of the submanifold $Z \subset \bR^d$, in the following precise sense.

\begin{defi}
\label{def:isotopic}
Let $Z,Z' \subset \bR^d$ be two closed submanifolds and let $\sG \leq \Diff(Z)$ and $\sG' \leq \Diff(Z')$ be open subgroups. We say that $(Z,\sG)$ and $(Z',\sG')$ are \emph{isotopic} if there is a diffeomorphism $Z \cong Z'$ such that the induced isomorphism $\Diff(Z) \cong \Diff(Z')$ sends $\sG$ onto $\sG'$ and the induced bijection of path-components $\pi_{0}(\Emb(Z,\bR^d)/G) \cong \pi_{0}(\Emb(Z',\bR^d)/G')$ sends $[\mathrm{incl}_Z]$ to $[\mathrm{incl}_{Z'}]$.
\end{defi}

\begin{prop}
\label{prop:isotopic}
If $(Z,\sG)$ is isotopic to $(Z',\sG')$, there is a natural isomorphism of continuous semifunctors $\Int{\cF}_{(Z,\sG,Q)} \cong \Int{\cF}_{(Z',\sG',Q)}$, and similarly for each of the other variants of the construction.
\end{prop}
\begin{proof}
The assumption that $(Z,\sG)$ is isotopic to $(Z',\sG')$ implies that we have a homeomorphism $\eta \colon \cE_{\sG}(Z,\bR^d) \cong \cE_{\sG'}(Z',\bR^d)$ that sends the basepoint $[\mathrm{incl}_Z]$ of the domain into the path-component containing the basepoint $[\mathrm{incl}_{Z'}]$ of the codomain. Let us choose a path in $\cE_{\sG'}(Z',\bR^d)$ from $[\mathrm{incl}_{Z'}]$ to $\eta([\mathrm{incl}_Z])$. Since the restriction map $\Diff_c(\bR^d) \to \cE_{\sG'}(Z',\bR^d)$ is a fibre bundle (by the same argument as in the proof of Proposition~\ref{prop:split-fibration-sequence-2}), we may lift this to a path in $\Diff_c(\bR^d)$ from the identity to a diffeomorphism $\varphi \in \Diff_c(\bR^d)$ whose pre-composition with $\eta([\mathrm{incl}_Z])$ is $[\mathrm{incl}_{Z'}]$.
For each $(M,A) \in \obj(\fU\cD_{d})$, pre-composition with (the inverse of) the diffeomorphism $Z \cong Z'$ also defines a homeomorphism $\eta_{(M,A)} \colon \cE_{\sG}(Z,\Breve{M} \smallsetminus A) \cong \cE_{\sG'}(Z',\Breve{M} \smallsetminus A)$, which agrees with $\eta$ on $[\mathrm{incl}_Z]$ under the preferred embedding $\bR^d \subset \Breve{M} \smallsetminus A$ (see Notation~\ref{not:union-with-Z}).
Denoting by $\varphi_{(M,A)}$ the extension of $\varphi \in \Diff_c(\bR^d)$ to a diffeomorphism of $\Breve{M} \smallsetminus A$ by setting it to be the identity outside of $\bR^d \subset \Breve{M} \smallsetminus A$, we therefore have a composite homeomorphism
\[
\begin{tikzcd}
\cE_{\sG}(Z,\Breve{M} \smallsetminus A) \ar[rr,"\eta_{(M,A)}"] && \cE_{\sG'}(Z',\Breve{M} \smallsetminus A) \ar[rr,"\varphi_{(M,A)} \circ -"] && \cE_{\sG'}(Z',\Breve{M} \smallsetminus A)
\end{tikzcd}
\]
that sends $[\mathrm{incl}_Z]$ to $[\mathrm{incl}_{Z'}]$, i.e.~is basepoint-preserving. We thus obtain a based homeomorphism $\zeta_{(M,A)} \colon X_{(Z,\sG)}(M,A) \cong X_{(Z',\sG')}(M,A)$ by restricting to the path-components of the basepoints. This induces an isomorphism on $\pi_{1}$ that extends to an isomorphism of split short exact sequences of the form \eqref{eq:split-ses-1}. By Lemma~\ref{lem:naturality_functorial_quotient_group}, this induces an isomorphism of commutative diagrams of the form \eqref{eq:split-short-exact-sequence-quotient}, which implies that $\pi_{1}(\zeta_{(M,A)})$ descends to an isomorphism $\cQ_{(Z,\sG,Q)}(M,A) \cong \cQ_{(Z',\sG',Q)}(M,A)$. Thus $\zeta_{(M,A)}$ is an isomorphism $\Int{\cF}_{(Z,\sG,Q)}(M,A) \cong \Int{\cF}_{(Z',\sG',Q)}(M,A)$ in $\covr$. This isomorphism is moreover natural with respect to morphisms of $\fU\cD_{d}$, since they act by post-composition with embeddings of decorated manifolds, whereas $\zeta_{(M,A)}$ acts by pre-composition with (the inverse of) a diffeomorphism $Z \cong Z'$ and post-composition with a diffeomorphism $\varphi$ supported in $\bR^d \subset \bB^d_{1}$ (where all embeddings of decorated manifolds act by the identity).
\end{proof}

\paragraph*{Functors of path-components.}

Note that Theorem~\ref{thm:global_functor_motion_groups} defines only \emph{semifunctors}, since $\fU\cD_{d}$ and $\fU\cD_{d}^{+}$ are only \emph{semicategories}; see \S\ref{sss:quillen-bracket-categories}. However, this technical issue does not matter once we pass to $\pi_{0}$:

\begin{prop}
\label{prop:functor-on-pi0}
The continuous semifunctors of Theorem~\ref{thm:global_functor_motion_groups} induce genuine functors on $\pi_{0}$.
\end{prop}
\begin{proof}
We consider any one of the semifunctors \eqref{eq:global-functor2} of Theorem~\ref{thm:global_functor_motion_groups}, which we denote by $\cF$; the proof repeats verbatim for the semifunctors defined on the category $\fU\cD_{d}^{+}$.
Let $(M,A)$ be an object of $\pi_{0}(\fU\cD_{d})$ (i.e.~an object of $\cD_{d}$, a decorated manifold). We have to show that $\pi_{0}(\cF)$ sends $\id_{(M,A)}$ to an identity morphism of $\pi_{0}(\covr)$. To do this, we first have to \emph{identify} the identity morphism of $(M,A)$ in $\pi_{0}(\fU\cD_{d})$.

Let $e_{1} \colon \bD^{d-1} \times [0,1] = \bB_{1}^d \hookrightarrow M$ denote one of the boundary cylinders that $(M,A)$ is equipped with (more precisely, a representative of one of the \emph{germs of} boundary cylinders that $(M,A)$ is equipped with). The boundary connected sum $(\bB_{1}^d,\emptyset) \natural (M,A)$ may be viewed as the union of $\bD^{d-1} \times [-1,1]$ with $M$ along $\bD^{d-1} \times [0,1]$ via the embedding $e_{1}$. Choose a diffeomorphism $[-1,1] \to [0,1]$ that is given by $t \mapsto t+1$ on $[-1,-1+\epsilon]$ and by $t \mapsto t$ on $[1-\epsilon,1]$ for some $\epsilon > 0$. Taking the product with $\bD^{d-1}$ (acting by the identity on this factor of the product) and extending by the identity over $M \smallsetminus \im(e_{1})$, this determines an isomorphism of decorated manifolds
\begin{equation}
\label{eq:self-embedding}
\upsilon_{(M,A)} \colon (\bB_{1}^d , \emptyset) \natural (M,A) \too (M,A).
\end{equation}
Consider the endomorphism of $(M,A)$ in $\fU\cD_{d}$ given by $\varUpsilon_{(M,A)} = [(\bB_{1}^d , \emptyset) , \upsilon_{(M,A)}]$.
One may check that, for any endomorphism $\varphi$ of $(M,A)$ in $\fU\cD_{d}$, the compositions $\varUpsilon_{(M,A)} \circ \varphi$ and $\varphi \circ \varUpsilon_{(M,A)}$ are both isotopic to $\varphi$. Hence $\varUpsilon_{(M,A)}$ is the identity of $(M,A)$ in the category $\pi_{0}(\fU\cD_{d})$. Under the identification of Corollary~\ref{cor:description_UD_{d}}, this corresponds to the self-embedding of $(M,A)$ given by restricting the diffeomorphism \eqref{eq:self-embedding} to the submanifold $(M,A) \subset (\bB_{1}^d,\emptyset) \natural (M,A)$. Since this self-embedding is isotopic to the identity, the induced self-map of embedding spaces $X(M,A) \to X(M,A)$ is homotopic to the identity, and hence is an identity morphism in $\pi_{0}(\covr)$.
\end{proof}

Recall that, apart from the continuous semifunctors provided by Theorem~\ref{thm:global_functor_motion_groups}, the other input for the general construction of \S\ref{sss:construction} is a continuous semifunctor $V \colon \fU\cD_{d} \to \modlr$ satisfying Condition~\ref{condition:input}. By Lemma~\ref{lem:factorisation_pi_{0}}, it induces an (abstract) semifunctor $\pi_{0}(V) \colon \pi_{0}(\fU\cD_{d}) \to \modlr$.

\begin{coro}\label{cor:homological_rep_genuine_functor}
If the semifunctor $\pi_{0}(V)$ is a genuine functor, then the associated homological representation semifunctor \eqref{eq:output_of_general_construction_motion_groups} of Corollary~\ref{coro:def_homological_rep_functor_motion_groups} upgrades to a genuine functor.
\end{coro}
\begin{proof}
We consider any one of the semifunctors of Theorem~\ref{thm:global_functor_motion_groups}, which we denote by $\cF$. First, it follows from Lemma~\ref{lem:factorisation_pi_{0}} that the linearisation (genuine) functor $\bZ[-]$ of \S\ref{sss:lift}, the fibrewise tensor product (genuine) functor $\otimes$ of \S\ref{sss:fibrewise_tensor} and the twisted homology (genuine) functor $H_{i}$ of \S\ref{sss:twisted-homology} induce (genuine) functors $\pi_{0}(\bZ[-])$, $\pi_{0}(\otimes)$ and $\pi_{0}(H_{i})$ on path-components, which commute with the canonical projections $\cC \to \pi_{0}(\cC)$ for the categories $\cC$ involved. Also, by Proposition~\ref{prop:functor-on-pi0}, the composite $\pi_{0}(\bZ[-])\circ \pi_{0}(\cF)$ is a genuine functor.
Now, consider the commutative diagram \eqref{eq:construction} defining the homological representation semifunctor $L_{i}(\cF;V)$. The functor $\pi_{0}(\bZ[-])\circ \pi_{0}(\cF)$ and the semifunctor $\pi_{0}(V)$ determine a semifunctor $\pi_{0}(\cF,V) \colon \pi_{0}(\fU\cD_{d}) \to \pi_{0}(\topmodr)$ by the universal property of the pullback. But then, since we assume that $\pi_{0}(V)$ is a genuine functor, so is $\pi_{0}(\cF,V)$. Finally, $L_{i}(\cF;V)$ is defined by the composite of $\pi_{0}(\cF,V)$ with the genuine functors $\pi_{0}(\otimes)$ and $\pi_{0}(H_{i})$.
\end{proof}

\paragraph*{Borel-Moore homology variant.}

We finally show that the \emph{closed} variants of the semifunctors of Theorem~\ref{thm:global_functor_motion_groups} take values in the subcategory $\covr^{\pr} \subset \covr$; see Lemma~\ref{lem:F-proper}. By Theorem~\ref{thm:construction}, these semifunctors therefore work equally well as an input for our general construction when we use twisted \emph{Borel-Moore} homology, which is functorial only with respect to proper maps of spaces.

\begin{lemm}\label{lem:F-proper}
The continuous semifunctors $\cF_{(Z,\sG,Q)}$ and $\cF^{\unt}_{(Z,\sG,Q)}$ of Theorem~\ref{thm:global_functor_motion_groups} and their analogues defined on $\fU\cD_{d}^{+}$ take values in the subcategory $\covr^{\pr} \subset \covr$.
\end{lemm}
\begin{proof}
Whichever \emph{closed} semifunctor of Theorem~\ref{thm:global_functor_motion_groups} we consider, the proof goes as follows.
We have to show that, for any morphism $\varphi \colon (M,A) \to (N,B)$ of decorated manifolds, the induced map of spaces $X'(M,A) \to X'(N,B)$ (in the \emph{closed} variant of the construction described at the end of \S\ref{sss:F-construction}) is a proper map (preimages of compact subspaces are compact). This map is (a restriction to particular path-components of) the inclusion of embedding spaces $\cE_{\sG}(Z , \mbar \smallsetminus A) \hookrightarrow \cE_{\sG}(Z , \nbar \smallsetminus B)$ induced by an embedding of decorated manifolds $(M,A) \hookrightarrow (N,B)$ satisfying the three properties of Definition~\ref{def:embdec} (by Proposition~\ref{prop:morphism-spaces-bracket}). In particular, the third property implies that the inclusion $\mbar \smallsetminus A \hookrightarrow \nbar \smallsetminus B$ has closed image, which implies that the inclusion of embedding spaces above also has closed image (this holds for the compact-open topology, hence also for the Whitney topology, which is finer); any closed inclusion is a proper map.
\end{proof}

\begin{rmk}
The continuous semifunctor $\Int{\cF}_{(Z,\sG,Q)}$ of Theorem~\ref{thm:global_functor_motion_groups} does \emph{not} take values in the subcategory $\covr^{\pr} \subset \covr$, since proof of Lemma~\ref{lem:F-proper} breaks down in this setting as the inclusion of $\Breve{M} \smallsetminus A$ into $\Breve{N} \smallsetminus B$ does not have closed image.
However, the restriction of the semifunctor $\Int{\cF}_{(Z,\sG,Q)}$ to the underlying groupoid $\cD_{d}$ of $\fU\cD_{d}$ \emph{does} take values in $\covr^{\pr}$, since the underlying groupoid of $\covr$ is contained in $\covr^{\pr}$ (because homeomorphisms are proper maps). Therefore, using Borel-Moore homology together with the semifunctor $\Int{\cF}_{(Z,\sG,Q)}$ in our general construction is not fully functorial: namely, it is functorial only for the \emph{isomorphisms} $\cD_{d}$ of $\fU\cD_{d}$, so we just obtain representations of the individual groups in this case.

Furthermore, the natural homotopy equivalences of Lemma~\ref{lem:natural-homotopy-equivalence} are not \emph{proper}. Thus, when using \emph{Borel-Moore} twisted homology -- which is invariant only under proper homotopy equivalences -- in the general construction of \S\ref{sss:construction}, the two choices $\cF_{(Z,\sG,Q)}$ and $\Int{\cF}_{(Z,\sG,Q)}$ will lead to \emph{a priori} different homological representations. As noted above, the homological representations obtained using Borel-Moore homology are defined on all of $\fU\cD_{d}$ when using $\cF_{(Z,\sG,Q)}$, but only on $\cD_{d}$ when using $\Int{\cF}_{(Z,\sG,Q)}$, so it is only on this subgroupoid that we may compare them.
\end{rmk}

\subsubsection{Colimit coefficient systems and untwisted representations}
\label{sss:col_coeff_unwtisted_rep}

Throughout this section, we fix a pair of groupoids $\cG$ and $\cM$ as in Hypothesis~\ref{hypo:standard_framework}.
Let us consider any semifunctor $F \colon \langle \cG , \cM \rangle \to \covr$ that induces a genuine functor on $\pi_{0}$, for example (see Proposition~\ref{prop:functor-on-pi0}) one of those constructed in Theorem~\ref{thm:global_functor_motion_groups} or its restriction along \eqref{eq:restriction_of_global_functor} to a subcategory $\langle \cG , \cM \rangle \subseteq \fU\cD_{d}$. In this section we construct an associated \emph{colimit} coefficient system $V_{\col}(F) \colon \langle \cG , \cM \rangle \to \modlr$. Under certain conditions on the subcategory $\langle \cG , \cM \rangle$, if $F$ is the restriction to $\langle \cG , \cM \rangle \subseteq \fU\cD_{d}$ of one of the two semifunctors $\Int{\cF}^{\unt}_{(Z,\sG,Q)}$ or $\cF^{\unt}_{(Z,\sG,Q)}$ of Theorem~\ref{thm:global_functor_motion_groups}, the homological representation functor $L_i(F;V_{\col}(F)) \colon \langle \pi_{0}(\cG) , \pi_{0}(\cM) \rangle \to \modr$ resulting from the general construction of \S\ref{sss:construction} takes values in a subcategory of $\modr$ of the form $\modr[{\bZ[\cQ]}]$ for a fixed group $\cQ$; see Proposition~\ref{prop:factoringthroughtopQ} below.

The idea to construct $V_{\col}(F)$ involves ``stabilising'' along the poset of canonical morphisms in the category $\langle \pi_{0}(\cG) , \pi_{0}(\cM) \rangle$, so we first fix notation for this poset.

\begin{defi}
\label{def:underlying-poset}
A morphism $\varphi \colon M \to N$ of $\langle \pi_{0}(\cG) , \pi_{0}(\cM) \rangle$ is called \emph{canonical} if $N = L \natural M$ and $\varphi = [L,\id_{L \natural M}]$ for some object $L$ of $\pi_{0}(\cG)$. Write $\cP(\cG,\cM)$ for the subcategory of $\langle \pi_{0}(\cG) , \pi_{0}(\cM) \rangle$ with the same objects and only the canonical morphisms. There is at most one canonical morphism between any pair of objects, so this subcategory is a poset.
For an object $M$, we also write $\cP(\cG,\cM)_{\geq M}$ for the subposet on all objects $N \geq M$, in other words all objects of the form $L \natural M$. Finally, we write $(- \natural M) \colon \cP(\cG,\cM) \cong \cP(\cG,\cM)_{\geq M}$ for the canonical isomorphism of posets.
\end{defi}

\begin{rmk}
The fact that the subcategory $\cP(\cG,\cM)$ of Definition~\ref{def:underlying-poset} is a poset is equivalent to the statement that if $L \natural M = L' \natural M$ then $L = L'$, which is tautologically true by construction. To clarify, we note that the corresponding statement with $\cong$ in place of $=$ is false in general, but that is not relevant here.
We also recall from Remark~\ref{rmk:small_categories_Quillen} that $\langle \pi_{0}(\cG) , \pi_{0}(\cM) \rangle$ is a \emph{small} category, so $\cP(\cG,\cM)$ is really a set and not a proper class.
\end{rmk}

The basic construction is the following operation on functors $\langle \pi_{0}(\cG) , \pi_{0}(\cM) \rangle \to \groups$.
For a functor $T \colon \langle \pi_{0}(\cG) , \pi_{0}(\cM) \rangle \to \groups$, since $\langle \pi_{0}(\cG) , \pi_{0}(\cM) \rangle$ is a \emph{small} category and the category of groups $\groups$ is cocomplete, we may define
\begin{equation}
\label{eq:Tcol}
T_{\col}(M) := \underset{\cP(\cG,\cM)}{\colim}(T \circ (- \natural M)).
\end{equation}
Each morphism $\varphi \colon M \to N$ induces a natural transformation $T \circ (- \natural M) \Rightarrow T \circ (- \natural N)$ given by $T(\id_{L} \natural \varphi)$ on each object $L$ of $\cP(\cG,\cM)$. Natural transformations of diagrams induce morphisms of colimits, so we may define $T_{\col}(\varphi)$ to be the morphism $T_{\col}(M) \to T_{\col}(N)$ induced by this natural transformation. The identity and composition axioms are straightforward to check from the definition, and so we prove:

\begin{lemm}
\label{lem:stabilised-functor-to-groups}
For a functor $T \colon \langle \pi_{0}(\cG) , \pi_{0}(\cM) \rangle \to \groups$, the above assignments define a functor $T_{\col} \colon \langle \pi_{0}(\cG) , \pi_{0}(\cM) \rangle \to \groups$.
\end{lemm}

Moreover, there is a natural transformation $T \Rightarrow T_{\col}$ given on each object $M$ by the morphism $T(M) \to T_{\col}(M)$ induced by the morphisms $T([L,\id_{L \natural M}]) \colon T(M) \to T(L \natural M)$ for all $L$.

An important observation is the following.

\begin{lemm}
\label{lem:isomorphic-to-constant}
If the poset $\cP(\cG,\cM)$ is a directed set, then for any functor $T \colon \langle \pi_{0}(\cG) , \pi_{0}(\cM) \rangle \to \groups$, the associated functor $T_{\col}$ is isomorphic to a functor that is constant on $\cP(\cG,\cM) \subset \langle \pi_{0}(\cG) , \pi_{0}(\cM) \rangle$, in other words it sends every object to a fixed group $\cT$ and it sends every \emph{canonical} morphism to $\id_\cT$.
\end{lemm}
\begin{proof}
This follows immediately from the fact that, when $\cP(\cG,\cM)$ is a directed set, all of the subposets $\cP(\cG,\cM)_{\geq M}$ are cofinal in each other, so all of the colimits \eqref{eq:Tcol} are canonically isomorphic to each other; see \cite[\href{https://stacks.math.columbia.edu/tag/04E7}{Lemma 04E7}]{stacks-project}. Thus they are all canonically isomorphic to a fixed choice of group $\cT$ representing this isomorphism class. Moreover, since the colimit in \eqref{eq:Tcol} is taken over all canonical morphisms, the canonical morphisms are sent to $\id_\cT$ under these identifications. Hence these canonical isomorphisms assemble into an isomorphism of functors from $T_{\col}$ to one having the properties described.
\end{proof}

\begin{notation}
\label{notation-Qcol}
When the poset $\cP(\cG,\cM)$ is a directed set, the group \eqref{eq:Tcol} (up to isomorphism) is independent of $M$, by the proof of Lemma~\ref{lem:isomorphic-to-constant}, and we denote it by $\cQ_{\col}(T)$.
\end{notation}

The technical result that we shall need is the following strengthening of Lemma~\ref{lem:isomorphic-to-constant} under an additional condition.

\begin{defi}
\label{def:motion-groupoid}
A subgroupoid $\cH \subseteq \cD_{d}$ (respectively $\cH' \subseteq \cD_{d}^{+}$) is called a \emph{motion groupoid} if, for each object $(M,A)$ of $\cH$, the subgroup $\Aut_\cH(M,A) \subseteq \diffdec(M,A)$ (resp.~$\Aut_{\cH'}(M,A) \subseteq \diffdec^{+}(M,A)$) lies in $\diffdecbr(M,A)$ (resp.~in $\diffdecbrplus(M,A)$). (The motivation for this terminology is Proposition~\ref{prop:braided-diff-groups}.)
\end{defi}

\begin{prop}
\label{prop:equivalent-to-constant}
Suppose that we have groupoids $\cG$ and $\cM$ as in Hypothesis~\ref{hypo:standard_framework} and let $F = \cF^{\unt}$ be one of the untwisted semifunctors of Theorem~\ref{thm:global_functor_motion_groups}, restricted along the inclusion \eqref{eq:restriction_of_global_functor}. Assume that the poset $\cP(\cG,\cM)$ is a directed set and that $\cM$ is a motion groupoid. Then the functor $(\sR \circ \cF^{\unt})_{\col} \colon \langle \pi_{0}(\cG) , \pi_{0}(\cM) \rangle \to \groups$ is equivalent to a constant functor.
\end{prop}
\begin{proof}
By Corollary~\ref{coro:QBC-skeleton}, we may find a skeleton of $\langle \pi_{0}(\cG) , \pi_{0}(\cM) \rangle$ that is itself a Quillen bracket category. That is, we may find skeleta $\cG_{0} \subset \pi_{0}(\cG)$ and $\cM_{0} \subset \pi_{0}(\cM)$, together with a monoidal structure on $\cG_{0}$ and a left action of $\cG_{0}$ on $\cM_{0}$, such that $\langle \cG_{0} , \cM_{0} \rangle$ is a skeleton of $\langle \pi_{0}(\cG) , \pi_{0}(\cM) \rangle$, i.e.~the inclusion $\mathsf{i} \colon \langle \cG_{0} , \cM_{0} \rangle \hookrightarrow \langle \pi_{0}(\cG) , \pi_{0}(\cM) \rangle$ is an equivalence. By construction, every morphism in $\langle \cG_{0} , \cM_{0} \rangle$ is the composition of a canonical morphism and an isomorphism of $\cM_{0}$ (via the canonical embedding of $\cM_{0}$ as the underlying groupoid of $\langle \cG_{0} , \cM_{0} \rangle$); since $\cM_{0}$ is skeletal, this isomorphism must be an \emph{automorphism}.

We will prove below that each automorphism of $\pi_{0}(\cM)$ (thus in particular of $\cM_{0}$) is sent by $T = \sR \circ \cF^{\unt}$ to an identity morphism. This implies that $T_{\col}$ also sends each automorphism of $\pi_{0}(\cM)$ to an identity morphism, since $T_{\col}(\varphi)$ is induced by a natural transformation whose component morphisms are of the form $T(\id_L \natural \varphi)$, and $\id_L \natural \varphi$ is an automorphism if $\varphi$ is an automorphism.

By Lemma~\ref{lem:isomorphic-to-constant}, $T_{\col}$ is isomorphic to a functor $T'_{\col}$ that sends every object to a fixed group $\cT$ and every canonical morphism to $\id_\cT$. By the previous paragraph, $T_{\col}$ sends every automorphism of $\pi_{0}(\cM)$ to an identity morphism; it follows that $T'_{\col}$ sends every automorphism of $\pi_{0}(\cM)$ to $\id_\cT$. Thus $T'_{\col}$ is constant on objects, canonical morphisms and automorphisms of $\pi_{0}(\cM)$ (in particular automorphisms of $\cM_{0}$). As noted above, every morphism in $\langle \cG_{0} , \cM_{0} \rangle$ is the composition of a canonical morphism and an automorphism of $\cM_{0}$, so $T'_{\col}$ is constant on $\langle \cG_{0} , \cM_{0} \rangle$.

Let $\mathsf{r} \colon \langle \pi_{0}(\cG) , \pi_{0}(\cM) \rangle \to \langle \cG_{0} , \cM_{0} \rangle$ be an inverse equivalence for the inclusion $\mathsf{i}$. Then $T_{\col}$ is equivalent to $T_{\col} \circ \mathsf{i} \circ \mathsf{r}$. By above, $T_{\col}$ is isomorphic to $T'_{\col}$, which has the property that $T'_{\col} \circ \mathsf{i}$ is constant. It thus follows that $T_{\col}$ is equivalent to $T'_{\col} \circ \mathsf{i} \circ \mathsf{r}$, which is constant.

It remains to prove the claim above that $T$ sends every automorphism $\varphi \colon (M,A) \to (M,A)$ of $\pi_{0}(\cM)$ to an identity morphism. We first recall that, by the definition of $\cF^{\unt}$ in \S\ref{sss:F-construction}, the functor $T = \sR \circ \cF^{\unt}$ sends the object $(M,A)$ to $\cQ^{\unt}(M,A)$, which is the coinvariants of $\cQ(M,A)$ under the action of $\pi_{1}(\cE(A,\Breve{M}))$ induced by the splittings in diagram \eqref{eq:split-short-exact-sequence-quotient}. It also sends the automorphism $\varphi \colon (M,A) \to (M,A)$ to the automorphism $\theta^{\unt}_\varphi \colon \cQ^{\unt}(M,A) \to \cQ^{\unt}(M,A)$ of \eqref{eq:assignment_morphisms_theta}, which is the map on coinvariants induced by the map $\theta_\varphi \colon \cQ(M,A) \to \cQ(M,A)$, which is part of the map of diagrams of the form \eqref{eq:split-short-exact-sequence-quotient} induced by $\varphi$, as explained in \S\ref{sss:F-construction}. Our goal is thus to prove that $\theta^{\unt}_\varphi$ is the identity.

By the assumption that $\cM$ is a motion groupoid, the element $\varphi$ lies in the subgroup $\pi_{0}(\diffdecbr(M,A)) \subseteq \pi_{0}(\diffdec(M,A))$. We may therefore choose a lift $\widetilde{\varphi} \in \pi_{1}(\cE(A,\Int{M}))$ of $\varphi$ under the surjection \eqref{eq:from-motion-to-bmcg} of Proposition~\ref{prop:braided-diff-groups}. Via the isomorphism $\pi_{1}(\cE(A,\Int{M})) \cong \pi_{1}(\cE(A,\Breve{M}))$ induced by the inclusion $\Int{M} \subset \Breve{M}$ we may view $\widetilde{\varphi}$ as an element of $\pi_{1}(\cE(A,\Breve{M}))$. Via the splittings in diagram \eqref{eq:split-short-exact-sequence-quotient}, it induces an automorphism $\widetilde{\varphi}_\sharp$ of $\cQ(M,A)$. Proposition~\ref{prop:two-actions-on-Q} tells us that $\theta_\varphi = \widetilde{\varphi}_\sharp$. By definition of $\cQ^{\unt}(M,A)$ as coinvariants, the automorphism of $\cQ^{\unt}(M,A)$ induced by $\widetilde{\varphi}_\sharp$ is the identity. But $\theta^{\unt}_\varphi$ is precisely the automorphism of $\cQ^{\unt}(M,A)$ induced by $\theta_\varphi = \widetilde{\varphi}_\sharp$, so it is the identity.
\end{proof}

\paragraph*{The colimit coefficient system $V_{\col}(F)$.}

Given any semifunctor $F \colon \langle \cG , \cM \rangle \to \covr$ that induces a genuine functor on $\pi_{0}$, we now define the \emph{colimit} coefficient system semifunctor; see Definition~\ref{def:colimit-coefficient-system}.
Composing with the forgetful functor $\sR \colon \covr \to \groups$ of Notation~\ref{not:lr_forgetful}, we have a semifunctor $\sR \circ F \colon \langle \cG , \cM \rangle \to \groups$, which induces by Lemma~\ref{lem:factorisation_pi_{0}} a semifunctor $\pi_{0}(\langle \cG , \cM \rangle) \to \groups$. Since $F$ induces a genuine functor on $\pi_{0}$ (by assumption), this is a genuine functor. Composing with the isomorphism of Lemma~\ref{lem:Serre-fibration-condition}, we have a functor $\langle \pi_{0}(\cG) , \pi_{0}(\cM) \rangle \to \groups$, which we also denote by $\sR \circ F$. Applying Lemma~\ref{lem:stabilised-functor-to-groups}, we may therefore consider the functor $(\sR \circ F)_{\col}$ and the natural transformation $\sR \circ F \Rightarrow (\sR \circ F)_{\col}$. Composing with the group ring functor $\bZ[-] \colon \groups \to \Rings$, this is a natural transformation of functors of the form $\langle \pi_{0}(\cG) , \pi_{0}(\cM) \rangle \to \Rings$.
Hence, we obtain:

\begin{defi}
\label{def:colimit-coefficient-system}
For $F \colon \langle \cG , \cM \rangle \to \covr$ a semifunctor that induces a genuine functor on $\pi_{0}$, we define the \emph{colimit} coefficient system semifunctor:
\begin{equation}
\label{eq:colim_coefficient_system}
V_{\col}(F) \colon \langle \cG,\cM \rangle \too \modlr
\end{equation}
to be the composite $(\bZ[\sR \circ F] \times \bZ[(\sR \circ F)_{\col}] \times \bZ[(\sR \circ F)_{\col}]) \circ \Delta_{3} \circ \pi_{0}$, where $\Delta_{3}$ is the $3$-ary diagonal functor and $\pi_{0} \colon \langle \cG , \cM \rangle \to \langle \pi_{0}(\cG) , \pi_{0}(\cM) \rangle$ is the projection onto path-components composed with the isomorphism of Lemma~\ref{lem:Serre-fibration-condition}. Concretely, on objects, $V_{\col}(F)$ sends $M$ to the group ring $\bZ[(\sR \circ F)_{\col}(M)]$ considered as a right module over itself and as a left module over $\bZ[\sR \circ F(M)]$ via the canonical group homomorphism $\sR \circ F(M) \to (\sR \circ F)_{\col}(M)$.
\end{defi}

\begin{lemm}\label{lem:colim_coeff_system}
The assignments for \eqref{eq:colim_coefficient_system} determine a continuous semifunctor, such that Condition~\ref{condition:input} is satisfied and the semifunctor $\pi_{0}(V_{\col}(F))$ is a genuine functor.
\end{lemm}
\begin{proof}
For each object $M$ of $\langle \cG,\cM \rangle$, the group ring $\bZ[(\sR \circ F)_{\col}(M)]$ has a canonical bimodule structure over the pair of rings $(\bZ[\sR \circ F(M)],\bZ[(\sR \circ F)_{\col}(M)])$: the action of $\bZ[\sR \circ F(M)]$ is induced by the canonical group homomorphism $\sR \circ F(M) \to (\sR \circ F)_{\col}(M)$ and left multiplication, while the action of $\bZ[(\sR \circ F)_{\col}(M)]$ is right multiplication. Thus $V_{\col}(F)(M)$ is an object of $\modlr$. For each morphism $\varphi \colon M \to N$ of $\langle \cG,\cM \rangle$, the ring homomorphisms $\bZ[\sR\circ F(\varphi)]$ and $\bZ[(\sR \circ F)_{\col}(\varphi)]$ are well-defined and the latter is moreover also a bimodule homomorphism $\bZ[(\sR \circ F)_{\col}(M)] \to (\bZ[\sR\circ F(\varphi)],\bZ[(\sR \circ F)_{\col}(\varphi)])^{*} (\bZ[(\sR \circ F)_{\col}(N)])$. Thus $V_{\col}(F)(\varphi)$ is a morphism of $\modlr$.

It is routine to verify that $V_{\col}(F)$ respects composition, so it is a well-defined semifunctor. It is continuous since the projection of $\langle \cG , \cM \rangle$ onto its path-components is continuous (see Remark~\ref{rmk:framework-open-path-components}) and this is then post-composed with functors between discrete categories to obtain $V_{\col}(F)$.

Applying the forgetful functor $\sL$ of Notation~\ref{not:lr_forgetful}, it follows from the definition of $V_{\col}(F)$ that $\sL \circ V_{\col}(F)$ is equal to $\bZ[\sR\circ F]$, thus Condition~\ref{condition:input} automatically holds.

Let $\Upsilon_{M}$ be an endomorphism of $M$ in $\langle \cG,\cM \rangle$ whose image in $\pi_{0}(\langle \cG,\cM \rangle)$ is the identity of $M$. (For example, we could take the explicit endomorphism defined in the proof of Proposition~\ref{prop:functor-on-pi0}.) By Proposition~\ref{prop:functor-on-pi0}, we have $F(\Upsilon_{M}) = \id_{F(M)}$ and it follows from our definition of $V_{\col}(F)$ on morphisms that $V_{\col}(F)(\Upsilon_{M})$ is an identity morphism. Thus $\pi_{0}(V_{\col}(F))$ preserves identities, i.e.~it is a genuine functor.
\end{proof}

It follows that each semifunctor $F \colon \langle \cG , \cM \rangle \to \covr$ that induces a genuine functor on $\pi_{0}$ automatically induces a homological representation functor $L_i(F;V_{\col}(F)) \colon\allowbreak \langle \pi_{0}(\cG) , \pi_{0}(\cM) \rangle \to \modr$ for each $i\geq 0$ by the construction of \S\ref{sss:construction} applied to $F$ and its associated colimit coefficient system \eqref{eq:colim_coefficient_system}.
The fundamental ``untwisting'' results of this section are the following two propositions:

\begin{prop}
\label{prop:factoringthroughtopQtw}
Suppose that we have groupoids $\cG$ and $\cM$ as in Hypothesis~\ref{hypo:standard_framework} and assume that the poset $\cP(\cG,\cM)$ is a directed set. Then, for each $i\geq 0$, the homological representation functor $L_i(F;V_{\col}(F))$ is isomorphic to a functor taking values in the subcategory $\modr[\bZ[\cQ]]^{\tw} \subset \modr$ for a fixed group $\cQ$.
\end{prop}
\begin{proof}
By Lemma~\ref{lem:isomorphic-to-constant}, the functor $(\sR \circ F)_{\col}$ is isomorphic to a functor that sends all objects to a fixed group $\cQ$. By construction it follows that, up to isomorphism, $V_{\col}(F)$ takes values in the subcategory $\modlr[\bullet][{\bZ[\cQ]}]^{\tw} \subset \modlr$. It then follows directly from the construction of \S\ref{sss:construction} that the resulting homological representation functor $L_i(F;V_{\col}(F))$, up to isomorphism, takes values in the subcategory $\modr[{\bZ[\cQ]}]^{\tw} \subset \modr$.
\end{proof}

\begin{prop}
\label{prop:factoringthroughtopQ}
Suppose that we have groupoids $\cG$ and $\cM$ as in Hypothesis~\ref{hypo:standard_framework} and let $F = \cF^{\unt}$ be one of the untwisted semifunctors of Theorem~\ref{thm:global_functor_motion_groups}, restricted along the inclusion \eqref{eq:restriction_of_global_functor}. Assume that the poset $\cP(\cG,\cM)$ is a directed set and that $\cM$ is a motion groupoid. Then, for each $i\geq 0$, the homological representation functor $L_i(\cF^{\unt};V_{\col}(\cF^{\unt}))$ is isomorphic to a functor taking values in the subcategory $\modr[\bZ[\cQ]] \subset \modr$ for a fixed group $\cQ$.
\end{prop}
\begin{proof}
By Proposition~\ref{prop:equivalent-to-constant}, the functor $(\sR \circ \cF^{\unt})_{\col}$ is isomorphic to the constant functor at a fixed group $\cQ$. By construction it follows that, up to isomorphism, $V_{\col}(\cF^{\unt})$ takes values in the subcategory $\modlr[\bullet][{\bZ[\cQ]}] \subset \modlr$. It then follows directly from the construction of \S\ref{sss:construction} that the resulting homological representation functor $L_i(\cF^{\unt};V_{\col}(\cF^{\unt}))$, up to isomorphism, takes values in the subcategory $\modr[{\bZ[\cQ]}] \subset \modr$.
\end{proof}

\begin{eg}[Low-degree lower central series.]
\label{eg:LCS}
We note that, in some examples of interest, we have $\cF = \cF^{\unt}$ for the semifunctors of Theorem~\ref{thm:global_functor_motion_groups}, which means that Proposition~\ref{prop:factoringthroughtopQ} applies to the homological representation functors $L_i(\cF;V_{\col}(\cF))$ for \emph{any} of the semifunctors $\cF$ of Theorem~\ref{thm:global_functor_motion_groups}.

This is true exactly when, for each object $(M,A)$, the $\pi_{1}(\cE(A,\Int{M})$-action on $\cQ(M,A)$ is trivial, so that $\cQ^{\unt}(M,A) = \cQ(M,A)$; see Notation~\ref{nota:untwisted}. An obvious example of this is when $Q$ is the quotient onto the trivial group, since $\cQ(M,A)$ is then always trivial. A less obvious but important example is when $Q$ is the abelianisation functor: in this case, all groups on the bottom row of diagram \eqref{eq:split-short-exact-sequence-quotient} are abelian, so the semi-direct product is a direct product, and so the $\pi_{1}(\cE(A,\Int{M})$-action on $\cQ(M,A)$ is trivial.
These are precisely the $\ell\in\{1,2\}$ cases of Example~\ref{eg:examples_FQG} applied to the lower central series, since $Q_{\LCS_{1}}$ is the quotient onto the trivial group and $Q_{\LCS_{2}}$ is the abelianisation functor.
\end{eg}

\subsection{Further constructions of homological representation functors}\label{ss:homological_representation_functor_mcg}

The construction of the functors $\cF \colon \fU\cD_{d} \to \modr$ in \S\ref{ss:homological_representation_functor_motion_groups} (see Theorem~\ref{thm:global_functor_motion_groups}) is well-adapted for motion groups in the sense that, if one restricts to a subcategory of $\fU\cD_{d}$ of the form $\langle \cG,\cM \rangle$, where $\cM$ is a \emph{motion groupoid} (see Definition~\ref{def:motion-groupoid} and Proposition~\ref{prop:braided-diff-groups}), then Proposition~\ref{prop:factoringthroughtopQ} implies that the induced homological representation functor $L_i(\cF;V_{\col}(\cF))$ takes values in the subcategory $\modr[\bZ[\cQ]] \subset \modr$ for a fixed group $\cQ$, as long as the underlying poset $\cP(\cG,\cM)$ is directed and if $\cF$ is one of the two \emph{untwisted} semifunctors constructed in Theorem~\ref{thm:global_functor_motion_groups}.
Thus the construction of \S\ref{ss:homological_representation_functor_motion_groups} is effective for producing families of untwisted representations of motion groups (or braided mapping class groups; see Remark~\ref{rmk:globality}). However, for the full mapping class groups, the construction of \S\ref{ss:homological_representation_functor_motion_groups} typically produces twisted representations.

In this subsection, we describe a variant of the construction of \S\ref{ss:homological_representation_functor_motion_groups}, based on the split short exact sequence \eqref{eq:split-ses-2} instead of \eqref{eq:split-ses-1}, which has good properties for the \emph{full} mapping class groups. More precisely, this variant satisfies an analogue of Proposition~\ref{prop:factoringthroughtopQ} where one does \emph{not} have to assume that $\cM$ is a motion groupoid; see Proposition~\ref{prop:factoringthroughcovQ-mcg} below.

Since the constructions and proofs are very similar to those of \S\ref{ss:homological_representation_functor_motion_groups}, we review the adaptation of the results of \S\ref{ss:homological_representation_functor_motion_groups}, pointing out and focusing on the differences.

\subsubsection{Construction of the homological representation functors}\label{ss:construction_semifunctor_mcg}

As in \S\ref{ss:homological_representation_functor_motion_groups}, we fix an integer $d \geq 2$, $d \neq 4$ (cf.~Remark~\ref{rmk:d_neq_4_explanation}), a closed submanifold $Z \subset \bR^d$, an open subgroup $\sG$ of $\Diff(Z)$ and a functorial quotient of groups $Q$.
We explain in the following paragraphs the construction of continuous semifunctors analogous to \eqref{eq:global-functor2}, which we denote by $\Int{\fF}_{(Z,\sG,Q)}$, $\Int{\fF}^{\unt}_{(Z,\sG,Q)}$, $\fF_{(Z,\sG,Q)}\, $ and $\fF^{\unt}_{(Z,\sG,Q)}$ to distinguish them from those of Theorem~\ref{thm:global_functor_motion_groups}. The result of this construction is summarised in Theorem~\ref{thm:global_functor_mcg} below.

\paragraph*{The semifunctor $\Int{\fF}$.}

The construction of the semifunctor $\Int{\fF}_{(Z,\sG,Q)}$ is similar to that of $\Int{\cF}_{(Z,\sG,Q)}$ in \S\ref{sss:F-construction}.
On objects, we follow the construction of \S\ref{sss:F-construction}; in particular we use the same space $X(M,A)$ (see Definition~\ref{def:space_{X}_M_A}).
However, the quotient $\phi(M,A) \colon \pi_{1}(X(M,A)) \twoheadrightarrow \cQ(M,A)$ is defined using the following 6-term commutative diagram instead of \eqref{eq:split-short-exact-sequence-quotient}.
\begin{equation}
\label{eq:split-short-exact-sequence-quotient-mcg}
\centering
\begin{split}
\begin{tikzpicture}
[x=1mm,y=1mm,font=\small]
\node at (10,8) {$\pi_{1}(X(M,A))$};
\node at (10,4) {\rotatebox{270}{$\cong$}};
\begin{scope}
\node (ll) at (-15,0) {$1$};
\node (l) at (10,0) {$\pi_{1}(\cE_{\sG}(Z,\Breve{M}\smallsetminus A))$};
\node (m) at (55,0) {$\pi_{0}(\diffdec(\mbar,A,Z|\sG))$};
\node (r) at (105,0) {$\pi_{0}(\diffdec(\mbar,A))$};
\node (rr) at (131,0) {$1$};
\draw[->] (ll) to (l);
\draw[->] (l) to (m);
\draw[->] (m) to (r);
\draw[->] (r) to (rr);
\draw[->,densely dashed] ($ (r.west) + (0,2)  $) to[out=160,in=20] ($ (m.east) + (0,2) $);
\end{scope}
\begin{scope}[yshift=-20mm]
\node (lla) at (-15,0) {$1$};
\node (la) at (10,0) {$\cQ(M,A)$};
\node (ma) at (55,0) {$Q(\pi_{0}(\diffdec(\mbar,A,Z|\sG)))$};
\node (ra) at (105,0) {$Q(\pi_{0}(\diffdec(\mbar,A)))$};
\node (rra) at (131,0) {$1$};
\draw[->] (lla) to (la);
\draw[->] (la) to (ma);
\draw[->] (ma) to (ra);
\draw[->] (ra) to (rra);
\draw[->,densely dashed] ($ (ra.west) + (0,2)  $) to[out=160,in=20] ($ (ma.east) + (0,2) $);
\end{scope}
\draw[->>] (l) to node[left,font=\footnotesize]{$\phi(M,A)$} (la);
\draw[->>] (m) to node[left,font=\footnotesize]{$q(M,A)$} (ma);
\draw[->>] (r) to node[left,font=\footnotesize]{$\bar{q}(M,A)$} (ra);
\end{tikzpicture}
\end{split}
\end{equation}
The top row is the split short exact sequence \eqref{eq:split-ses-2} of Corollary~\ref{coro:split-ses-2}. The rest of the diagram is induced from this just as in \S\ref{sss:F-construction}, using Lemma~\ref{lem:6-term-sequence}.
The construction of $\Int{\fF}_{(Z,\sG,Q)}$ on morphisms is then exactly as in \S\ref{sss:F-construction}, using the fact that the commutative diagram \eqref{eq:split-short-exact-sequence-quotient-mcg} is functorial in the morphisms of $\fU\cD_{d}$ thanks to Proposition~\ref{prop:split-ses-functoriality} and Lemma~\ref{lem:naturality_functorial_quotient_group}. In particular, Lemmas~\ref{lem:continuity_semifunctors} and \ref{lem:composition_axiom} and their proofs repeat verbatim.

\paragraph*{Untwisted and closed variants.}

By Proposition~\ref{prop:split-ses-functoriality} and Lemma~\ref{lem:naturality_functorial_quotient_group}, the commutative diagram \eqref{eq:split-short-exact-sequence-quotient-mcg} is functorial in the morphisms of $\fU\cD_{d}$, so there is an action of $\diffdec(M,A) = \Aut_{\fU\cD_{d}}((M,A))$ on the bottom-left group $\cQ(M,A)$. Since $\cQ(M,A)$ is discrete, this factors through an action of $\pi_{0}(\diffdec(M,A))$ on $\cQ(M,A)$.
On the other hand, the top-left group $\pi_{0}(\diffdec(\mbar,A))$ of \eqref{eq:split-short-exact-sequence-quotient-mcg} acts on $\cQ(M,A)$ either by lifting elements of $\cQ(M,A)$ along $\phi(M,A)$ and using the conjugation action of the semi-direct product on the top row, or equivalently by projecting along $\bar{q}(M,A)$ to the bottom-right group of \eqref{eq:split-short-exact-sequence-quotient-mcg} and then using the conjugation action of the semi-direct product on the bottom row. Via the canonical isomorphism $\pi_{0}(\diffdec(M,A)) \cong \pi_{0}(\diffdec(\mbar,A))$, induced by extending diffeomorphisms by the identity on $\mbar \smallsetminus M$, this gives another action of $\pi_{0}(\diffdec(M,A))$ on $\cQ(M,A)$.
The following is the analogue of Proposition~\ref{prop:two-actions-on-Q}.

\begin{prop}
\label{prop:two-actions-on-Q-mcg}
The two actions of $\pi_{0}(\diffdec(M,A))$ on $\cQ(M,A)$ described above are equal.
\end{prop}
\begin{proof}
Since both of these actions are induced via the quotient $\phi(M,A)$ from actions defined on $\pi_{1}(\cE_{\sG}(Z,\Breve{M}\smallsetminus A))$, it suffices to prove that the two actions agree on this group, before passing to the quotient. This is precisely part (i) of Proposition~\ref{prop:two-actions}.
\end{proof}

Analogously to Notation~\ref{nota:untwisted} in \S\ref{ss:homological_representation_functor_motion_groups}, we denote by $\cQ^{\unt}(M,A)$ the coinvariants $\cQ(M,A)_{\pi_{0}(\diffdec(M,A))}$ associated to this action, and by $\phi^{\unt}(M,A)$ the composite $X(M,A) \twoheadrightarrow \cQ(M,A) \twoheadrightarrow \cQ^{\unt}(M,A)$, where the first arrow is $\phi(M,A)$ and the second is the canonical projection onto the coinvariants. As in Remark~\ref{rmk:justification_untwisted}, this quotient group $\cQ^{\unt}(M,A)$ is the \emph{universal} one for the associated homological representations to be untwisted; see Proposition~\ref{prop:factoringthroughcovQ-mcg} below.
We thus define the untwisted variant $\Int{\fF}^{\unt}_{(Z,\sG,Q)}$ by $(M,A)\mapsto (X(M,A),\phi^{\unt}(M,A))$ on objects; on morphisms it is defined just as above, using functoriality of the diagram \eqref{eq:split-short-exact-sequence-quotient-mcg} by Proposition~\ref{prop:split-ses-functoriality} and Lemma~\ref{lem:naturality_functorial_quotient_group}.

Finally, the closed variants $\fF_{(Z,\sG,Q)}$ and $\fF^{\unt}_{(Z,\sG,Q)}$ are constructed by making the same modification as at the end of \S\ref{sss:F-construction}, using the (homotopy equivalent but not proper homotopy equivalent) space $X'(M,A)$ in place of $X(M,A)$.

The discussion above proves the following analogue of Theorem~\ref{thm:global_functor_motion_groups}.

\begin{thm}
\label{thm:global_functor_mcg}
For any integer $d \geq 2$, closed submanifold $Z \subset \bR^d$ and open subgroup $\sG$ of $\Diff(Z)$, each functorial quotient of groups $Q$ determines, via \eqref{eq:split-short-exact-sequence-quotient-mcg}, continuous semifunctors
\begin{equation}
\label{eq:global_functor_mcg}
\Int{\fF}_{(Z,\sG,Q)},\,\,\Int{\fF}^{\unt}_{(Z,\sG,Q)},\,\fF_{(Z,\sG,Q)}\, \text{ and }\, \fF^{\unt}_{(Z,\sG,Q)} \colon \fU\cD_{d} \too \covr.
\end{equation}
Similar continuous semifunctors are defined on the semicategory $\fU\cD_{d}^{+}$ if $Z$ is orientable and $\sG$ is contained in $\Diff^{+}(Z)$.
\end{thm}

\begin{rmk}
There is a natural morphism of diagrams $\eqref{eq:split-short-exact-sequence-quotient} \to \eqref{eq:split-short-exact-sequence-quotient-mcg}$ induced by Lemma~\ref{lem:naturality_functorial_quotient_group} from the map of split short exact sequences from the top row to the middle row in diagram \eqref{eq:3x3diagram} of Proposition~\ref{prop:map_of_ses}. This induces a natural transformation of semifunctors $\eqref{eq:global-functor2} \Rightarrow \eqref{eq:global_functor_mcg}$ for each of the four versions and for each fixed $(Z,\sG,Q)$.
\end{rmk}

\begin{coro}\label{coro:def_homological_rep_functor_mcg}
The semifunctors of Theorem~\ref{thm:global_functor_mcg} provide semifunctors of the form
\begin{equation}
\label{eq:output_of_general_construction_mcg}
\langle \pi_{0}(\cG) , \pi_{0}(\cM) \rangle \too \modr,
\end{equation}
which we denote by $L_{i}(\Int{\fF}_{(Z,\sG,Q)};V)$, $L_{i}(\Int{\fF}^{\unt}_{(Z,\sG,Q)};V)$, $L_{i}(\fF_{(Z,\sG,Q)};V)$ and $L_{i}(\fF^{\unt}_{(Z,\sG,Q)};V)$, where $V$ is any continuous semifunctor satisfying Condition~\ref{condition:input}. If the semifunctor $\pi_{0}(V)$ defined from $V$ by Lemma~\ref{lem:factorisation_pi_{0}} is a genuine functor, then the semifunctors \eqref{eq:output_of_general_construction_mcg} upgrade to genuine functors.
\end{coro}
\begin{proof}
This is a direct application of the general construction of \S\ref{sss:construction} (see Definition~\ref{def:construction}), together with the identification of $\pi_{0}(\langle \cG , \cM \rangle)$ with $\langle \pi_{0}(\cG) , \pi_{0}(\cM) \rangle$ by Lemma~\ref{lem:Serre-fibration-condition}.
The upgrade of \eqref{eq:output_of_general_construction_mcg} to genuine functors follows from the analogues of Proposition~\ref{prop:functor-on-pi0} and Corollary~\ref{cor:homological_rep_genuine_functor}, whose proofs are identical.
\end{proof}

\subsubsection{Elementary properties}\label{sss:elementary_properties_mcg}

The properties of \S\ref{sss:F-properties} also hold for the semifunctors of Theorem~\ref{thm:global_functor_mcg}. In particular, the proof of Proposition~\ref{prop:isotopic} repeats verbatim to prove that, for fixed $Q$, each of the four versions of \eqref{eq:global_functor_mcg} depends (up to natural isomorphism) on $(Z,\sG)$ only up to isotopy in the sense of Definition~\ref{def:isotopic}. In addition, it is routine to adapt the proofs of \S\ref{sss:F-properties} to show that the continuous semifunctors \eqref{eq:global_functor_mcg} induce genuine functors on $\pi_{0}$ (viz.~Proposition~\ref{prop:functor-on-pi0}), that the closed variants $\fF_{(Z,\sG,Q)}$ and $\fF^{\unt}_{(Z,\sG,Q)}$ take values in $\covr^{\pr}$ (viz.~Lemma~\ref{lem:F-proper}) and that there are natural homotopy equivalences $\Int{\fF}_{(Z,\sG,Q)} \Rightarrow \fF_{(Z,\sG,Q)}$ and $\Int{\fF}^{\unt}_{(Z,\sG,Q)} \Rightarrow \fF^{\unt}_{(Z,\sG,Q)}$ (viz.~Lemma~\ref{lem:natural-homotopy-equivalence}).
Indeed, the only difference between the constructions of the semifunctors of Theorem~\ref{thm:global_functor_motion_groups} and those of Theorem~\ref{thm:global_functor_mcg} is that, in the former, the groups $\cQ(M,A)$ (and the homomorphisms $\cQ(\varphi)$) are those obtained from diagram \eqref{eq:split-short-exact-sequence-quotient} (and its functoriality), whereas in the latter they are obtained from diagram \eqref{eq:split-short-exact-sequence-quotient-mcg} (and its functoriality).

\subsubsection{Colimit coefficient systems and untwisted representations}\label{sss:col_coeff_unwtisted_rep_mcg}

Recall from \S\ref{sss:col_coeff_unwtisted_rep} that each semifunctor $F \colon \langle \cG , \cM \rangle \to \covr$ that induces a genuine functor on $\pi_{0}$ has an associated \emph{colimit} coefficient system (semifunctor) $V_{\col}(F) \colon \langle \cG , \cM \rangle \to \modlr$; see Definition~\ref{def:colimit-coefficient-system} and Lemma~\ref{lem:colim_coeff_system}. Since the semifunctors \eqref{eq:global_functor_mcg} induce genuine functors on $\pi_{0}$ (see \S\ref{sss:elementary_properties_mcg}), we may apply this construction to $F = \eqref{eq:global_functor_mcg}$. Proposition~\ref{prop:factoringthroughtopQtw} applies directly to this setting, since it makes no further assumptions on $F$. The analogue of Proposition~\ref{prop:factoringthroughtopQ} is the following, in which we notably drop the assumption that $\cM$ is a motion groupoid.

\begin{prop}
\label{prop:factoringthroughcovQ-mcg}
Suppose that we have groupoids $\cG$ and $\cM$ as in Hypothesis~\ref{hypo:standard_framework} and let $F = \fF^{\unt}$ be one of the untwisted semifunctors of Theorem~\ref{thm:global_functor_mcg}, restricted along the inclusion \eqref{eq:restriction_of_global_functor}. Assume that the poset $\cP(\cG,\cM)$ is a directed set. Then, for each $i\geq 0$, the homological representation functor $L_i(\fF^{\unt};V_{\col}(\fF^{\unt}))$ is isomorphic to a functor taking values in the subcategory $\modr[\bZ[\cQ]] \subset \modr$ for a fixed group $\cQ$.
\end{prop}
\begin{proof}
Just as in the proof of Proposition~\ref{prop:factoringthroughtopQ}, this will follow immediately from the definition of the colimit coefficient system $V_{\col}(\fF^{\unt})$ and the general construction of \S\ref{sss:construction}, once we have proven the analogue of Proposition~\ref{prop:equivalent-to-constant}, in which we consider $\fF^{\unt}$ in place of $\cF^{\unt}$ (and we do not assume that $\cM$ is a motion groupoid). The proof of this analogue is identical to that of Proposition~\ref{prop:equivalent-to-constant} itself, the only difference being that we use Proposition~\ref{prop:two-actions-on-Q-mcg} instead of Proposition~\ref{prop:two-actions-on-Q} to understand the action of automorphisms of $\pi_{0}(\cM)$ on $\cQ(M,A)$. Since Proposition~\ref{prop:two-actions-on-Q-mcg} applies to the whole mapping class group $\pi_{0}(\diffdec(M,A))$ (whereas Proposition~\ref{prop:two-actions-on-Q} applies only to the subgroup $\pi_{0}(\diffdecbr(M,A))$, the braided mapping class group), we do not need to assume that $\cM$ is a motion groupoid, as we had to for Proposition~\ref{prop:factoringthroughtopQ}.
\end{proof}

\section{Applications for mapping class groups and motion groups}\label{s:applications}

In this section, we study in detail the representations obtained as outputs of the general construction of \S\ref{sss:construction} (see Definition~\ref{def:construction}), applied following the methods of \S\ref{ss:homological_representation_functor_motion_groups}--\S\ref{ss:homological_representation_functor_mcg} to some important families of mapping class groups and motion groups, namely the families of \emph{classical braid groups}, \emph{braid groups on surfaces} and \emph{loop braid groups} in \S\ref{ss:applications_motion_groups}, and the families of \emph{mapping class groups of surfaces} in \S\ref{ss:mcg_construction}. Beforehand, we detail the categorical framework coming from \S\ref{s:categorical_framework} for these families of groups in \S\ref{ss:categories_for_families_of_groups}.

\subsection{Categories for mapping class groups and motion groups}\label{ss:categories_for_families_of_groups}

The mapping class groups and motion groups that we shall study in \S\ref{ss:applications_motion_groups} and \S\ref{ss:mcg_construction} generally come in \emph{families}, in the sense that they naturally form a set $\{G_{n}\}_{n\in\mathbb{N}}$ along with canonical maps $g_{n}\colon G_{n}\to G_{n+1}$. A good way to treat these objects systematically, extracting their essential structure, is to package them using the Quillen bracket construction framework of \S\ref{ss:Quillen_bracket_construction}, as follows.

\begin{construction}\label{const:category_families_of_groups}
Let $\{g_{n}\colon G_{n}\to G_{n+1}\}_{n\in\mathbb{N}}$ be a family of groups of the following form.
\begin{enumerate}[noitemsep,label=(\roman*)]
\item\label{construction_input1} There is a fixed dimension $d\geq 2$ such that, for each $n\in\bN$, the group $G_{n}$ is isomorphic either to the mapping class group $\pi_{0}(\diffdec(M_{n},A_{n}))$ (see Notation~\ref{not:diffdec}) or to the braided mapping class group $\pi_{0}(\diffdecbr(M_{n},A_{n}))$ (see Definition~\ref{def:braided-diffeomorphisms}), where $(M_{n},A_{n})$ is a decorated manifold (see Definition~\ref{def:decorated-manifolds}); in particular, $M_{n}$ is a smooth $d$-manifold with non-empty boundary and $A_{n}$ is a closed submanifold of its interior $\Int{M}_{n}$. Moreover, the homomorphism $g_{n}$ is induced by an embedding in $\embdec((M_{n},A_{n}),(M_{n+1},A_{n+1}))$ (see Definition~\ref{def:embdec}), sending an isotopy class of diffeomorphisms of $(M_{n},A_{n})$ to the isotopy class of diffeomorphisms of $(M_{n+1},A_{n+1})$ given by applying this embedding and extending it by the identity on the complement of its image.
\item\label{construction_input2} There is a decorated manifold $(M,A)\in\cD_{d}$ such that, for each $n\in\bN$, the decorated manifold $(M_{n+1},A_{n+1})$ is isomorphic to the boundary connected sum $(M,A) \natural (M_{n},A_{n})$ and, under this isomorphism, the preferred embedding of $(M_{n},A_{n})$ into $(M_{n+1},A_{n+1})$ corresponds to the inclusion into $(M,A) \natural (M_{n},A_{n})$.
\item\label{construction_input3} In the case when $G_n$ is a \emph{braided} mapping class group in point \ref{construction_input1}, we also assume that the mapping class group $\pi_0(\diffdec((M,\emptyset)^{\natural n}))$ is trivial for all $n \in \bN$. This implies that we have an equality $\pi_{0}(\diffdecbr((M,A)^{\natural n})) = \pi_{0}(\diffdec((M,A)^{\natural n}))$.
\end{enumerate}
In order to encode the family of groups $\{g_{n}\colon G_{n}\to G_{n+1}\}_{n\in\mathbb{N}}$, we carry out the following steps.
\begin{enumerate}[noitemsep,label=(\arabic*)]
\item\label{construction_step1} We first define $\cM$ to be a certain subgroupoid of $\cD_{d}$ whose objects are all decorated manifolds that are isomorphic to $(M_{n},A_{n})$ for some $n \in \bN$. In the case when $G_n$ is a \emph{full} mapping class group, we define $\cM$ to be the full subgroupoid on these objects. In the case when $G_n$ is a \emph{braided} mapping class group, we must be slightly more careful. For each decorated manifold $(N,\emptyset)$ in the isomorphism class of $(M_{n},\emptyset)$, let us fix an isomorphism $i_{(N,\emptyset)} \colon (N,\emptyset) \cong (M_{n},\emptyset)$; in the case $(N,\emptyset) = (M_{n},\emptyset)$ we take it to be the identity. We then define the morphisms of $\cM$ to be those isomorphisms $(N,B) \cong (N',B')$ that become isotopic to $i_{(N',\emptyset)}^{-1} \circ i_{(N,\emptyset)}$ after forgetting the submanifolds $B$ and $B'$. (We call these \emph{braided isomorphisms}, generalising Definition~\ref{def:braided-diffeomorphisms} to isomorphisms between distinct objects.) This ensures that $\cM$ is a well-defined subgroupoid of $\cD_d$ and that its automorphism groups are braided diffeomorphism groups.
\item\label{construction_step2} We then define $\cG$ to be the full subgroupoid of $\cD_{d}$ on all objects that are isomorphic to $(M,A)^{\natural n}$ for some $n\in\mathbb{N}$. This is evidently closed under the semi-monoidal structure $\natural$ of $\cD_{d}$ (see Definition~\ref{prop:Dd-semi-monoidal}).
Moreover, the subgroupoid $\cM \subseteq \cD_d$ is closed under the left action of $\cG$ induced by $\natural$; in the case when the groups $G_n$ are \emph{full} mapping class groups (where we defined $\cM$ to be a full subgroupoid in step \ref{construction_step1}) this is clear, whereas in the case when they are \emph{braided} mapping class groups it follows from the assumption in point \ref{construction_input3}.
\item\label{construction_step3} Then, by Corollaries~\ref{cor:description_UD_{d}} and \ref{coro:QBC-skeleton}, there exists a skeleton $\langle \pi_{0}(\cG), \pi_{0}(\cM)\rangle_{0}$ of the Quillen bracket construction $\langle \pi_{0}(\cG), \pi_{0}(\cM)\rangle$; this encodes the family of groups $\{g_{n}\colon G_{n}\to G_{n+1}\}_{n\in\mathbb{N}}$.
\end{enumerate}
\end{construction}

\begin{rmk}\label{rmk:small_categories_pi_{0}}
Recall from Convention~\ref{convention:small_categories} that the topologically-enriched groupoids $\cD_{d}$ and $\cD_{d}^{+}$ are essentially small, and that we have implicitly replaced them with equivalent small groupoids. All of the subgroupoids of $\cD_{d}$ and $\cD_d^{+}$ that we shall consider are therefore also small, as are the corresponding discrete groupoids after applying the functor $\pi_{0}$ and the corresponding (discrete or topologically-enriched) Quillen bracket categories.
\end{rmk}

In the remainder of this subsection, we describe and discuss some properties of the subgroupoids $\cG$ and $\cM$ defined as in Construction~\ref{const:category_families_of_groups}, corresponding to mapping class groups of surfaces in \S\ref{sss:category_mcg}, to surface braid groups in \S\ref{sss:category_surface_braid_groups} and to loop braid groups and extended loop braid groups in \S\ref{sss:category_loop_braid_groups}. In \S\ref{ss:examples} we then apply the general construction of \S\ref{s:general_construction} to these examples to obtain coherent representations of these families of groups.

\subsubsection{Mapping class groups of surfaces}\label{sss:category_mcg}

We denote by $\bD^{2}$ the closed unit $2$-disc. Let $\Sigma_{0,1}^{1}$ denote the cylinder $\bS^{1}\times[0,1]$ (this notation is more commonly used to denote a once-punctured disc, which is consistent with the fact that the non-distinguished boundary component of $\Sigma_{0,1}^{1}$ may be moved freely by diffeomorphisms), let $\Sigma_{1,1}$ denote the torus with one boundary component $(\bS^{1}\times\bS^{1}) \smallsetminus\Int{\bD}^{2}$ and let $\N{}_{1,1}$ denote a Möbius band.
For $g\geq 0$, we denote by $\Sigma_{g,1}^{s}$ the boundary connected sum $(\natural_{s}\Sigma_{0,1}^{1})\natural(\natural_{g}\Sigma_{1,1})$ and by $\MCGo_{g,1}^{s}$ the mapping class group $\pi_{0}(\Diff_{\partial}(\Sigma_{g,1}^{s}))$.
For $h\geq 0$, we denote by $\N_{h,1}^{s}$ the boundary connected sum $(\natural_{s}\Sigma_{0,1}^{1})\natural(\natural_{h}\N{}_{1,1})$ and by $\MCGno_{h,1}^{s}$ the mapping class group $\pi_{0}(\Diff_{\partial}(\N_{h,1}^{s}))$.
When $s=0$, we omit it from the notation.
Let $\cM_{2}^{+,\mt}$ and $\cM_{2}^{-,\mt}$ be the full subgroupoids of $\cD_{2}$ on decorated surfaces of the form $(S,\emptyset,e_{1},e_{2})$, where $S$ is homeomorphic to $\Sigma_{g,1}$ for some $g\geq 0$ (for $\cM_{2}^{+,\mt}$) or to $\N_{h,1}$ for some $h\geq 0$ (for $\cM_{2}^{-,\mt}$) and the boundary-cylinders $e_{1}$ and $e_{2}$ both lie on the same boundary component of $S$.
The groupoids $\cM_{2}^{+,\mt}$ and $\cM_{2}^{-,\mt}$ clearly both inherit the semi-monoidal structure $\natural$ from $\cD_{2}$ (see Definition~\ref{prop:Dd-semi-monoidal}).
The choices $\cG = \cM = \cM_{2}^{+,\mt}$ and $\cG = \cM = \cM_{2}^{-,\mt}$ thus each fit into Hypothesis~\ref{hypo:standard_framework}.

We denote by $\cM_{2}^{+}$ and $\cM_{2}^{-}$ the path-components $\pi_{0}(\cM_{2}^{+,\mt})$ and $\pi_{0}(\cM_{2}^{-,\mt})$ respectively.
Since  $\cM_{2}^{+,\mt}$ and $\cM_{2}^{-,\mt}$ contain the solid cylinder $(\bB^2_{1},\emptyset,\id,r)$, the groupoids $\cM_{2}^{+}$ and $\cM_{2}^{-}$ have, by Lemma~\ref{lem:monoidal-pi0}, a monoidal structure induced by the boundary connected sum $\natural$ of Definition~\ref{def:boundary-connected-sum}.
Hence we may apply Quillen's bracket construction and by Lemma~\ref{lem:Serre-fibration-condition} we have isomorphisms of categories $\pi_{0}(\fU\cM_{2}^{+,\mt}) \cong \fU\cM_{2}^{+}$ and $\pi_{0}(\fU\cM_{2}^{-,\mt}) \cong \fU\cM_{2}^{-}$.
By Proposition~\ref{prop:monoidal-skeleton}, we may pass to monoidal skeleta of $(\cM_{2}^{+},\natural)$ and $(\cM_{2}^{-},\natural)$ and the categories $\fU\cM_{2}^{+}$ and $\fU\cM_{2}^{-}$ are thus also skeletal by Corollary~\ref{coro:QBC-skeleton}.
For diffeomorphisms of surfaces, the condition of fixing pointwise two (neighbourhoods of) intervals in a boundary component is equivalent to fixing pointwise the whole boundary component. So the morphisms in $\cM_{2}^{+}$ and $\cM_{2}^{-}$ may be identified with the isotopy classes of diffeomorphisms of the surfaces $S$ that that fix pointwise the preferred boundary component of $S$, i.e.~their automorphism groups are the \emph{mapping class groups} $\pi_{0}(\Diff_{\partial}(S))$ of $S$.

\subsubsection{Surface braid groups}\label{sss:category_surface_braid_groups}

Let $S$ be a compact, connected, smooth surface with a chosen boundary component $\partial_{0} S$. For each non-negative integer $n$, we denote by $\underline{n}$ a closed submanifold of $S$ consisting of $n$ distinct points in the interior of $S$.
Let $\cB r^S$ be the subgroupoid of $\cD_{2}$ with objects all decorated surfaces of the form $(S',\underline{n}',e_{1},e_{2})$, where the embedded $1$-discs $e_{1}(b\bB^2_{1})$ and $e_{2}(b\bB^2_{1})$ lie on the same boundary component $\partial_{0} S'$, $n$ is any non-negative integer and there is a diffeomorphism $S \cong S'$ taking $\partial_{0} S$ onto $\partial_{0} S'$ and $\underline{n}$ onto $\underline{n}'$. The morphisms of $\cB r^S$ are the \emph{braided isomorphisms} of $\cD_{2}$, as described in step \ref{construction_step1} of Construction~\ref{const:category_families_of_groups}; in particular its automorphism groups are given by the braided diffeomorphism groups of Definition~\ref{def:braided-diffeomorphisms}. Since the condition of being a \emph{braided} isomorphism is invariant under deformations along paths in the space of all isomorphisms in $\cD_{2}$, the groupoid $\cB r^S$ is $0$-full in $\cD_{2}$. If $S=\bD^{2}$, this collection of objects is closed under the semi-monoidal structure $\natural$ of Definition~\ref{prop:Dd-semi-monoidal}, so in this case $\cB r^{\bD^{2}}$ is a semi-monoidal subgroupoid of $\cD_{2}$. Moreover, in this case, $\cB r^{\bD^{2}}$ is \emph{full} (not just $0$-full) in $\cD_{2}$, since the mapping class group $\pi_{0}(\Diff_{I \sqcup I}(\bD^{2}))$ is trivial by Lemma \ref{lem:trivial_mcg_of_disc_with_two_intervals} below. In general, the groupoid $\cB r^S$ is closed under the left action of $\cB r^{\bD^{2}}$ via $\natural$.
Thus the choice of $\cG = \cB r^{\bD^{2}}$ and $\cM = \cB r^S$ fits into Hypothesis~\ref{hypo:standard_framework}.

\begin{lemm}
\label{lem:trivial_mcg_of_disc_with_two_intervals}
Let $I \sqcup I$ be a pair of disjoint closed intervals in the boundary of the disc $\bD^{2}$ and write $\Diff_{I \sqcup I}(\bD^{2})$ for the group of diffeomorphisms of $\bD^{2}$ that fix $I \sqcup I$ pointwise. Then $\Diff_{I \sqcup I}(\bD^{2})$ is weakly contractible, hence path-connected, so the mapping class group $\pi_{0}(\Diff_{I \sqcup I}(\bD^{2}))$ is trivial.
\end{lemm}
\begin{proof}
By \cite[\S II.2.2.2, Cor.~2]{Cerf1961Topologiedecertains} (see also \cite{Palais1960Localtrivialityof,Lima1963localtrivialityof}), the map
\[
\Diff_{I \sqcup I}(\bD^{2}) \too \Diff_\partial(I \sqcup I) \cong (\Diff_\partial(I))^2
\]
that remembers just the action of a diffeomorphism restricted to the \emph{complementary} pair of intervals $\partial \bD^{2} \smallsetminus (I \sqcup I) \cong I \sqcup I$ is a fibre bundle. Its fibre over the identity is $\Diff_\partial(\bD^{2})$, the group of diffeomorphisms of $\bD^{2}$ fixing all of $\partial \bD^{2}$ pointwise. The diffeomorphism group $\Diff_\partial(I)$ is easily seen to be contractible, and the diffeomorphism group $\Diff_\partial(\bD^{2})$ was shown to be contractible by Smale \cite{Smale1959}. The long exact sequence of the above fibre bundle then implies that $\Diff_{I \sqcup I}(\bD^{2})$ is weakly contractible.
\end{proof}

Let $\Beta^{S}$ be the groupoid $\pi_{0}(\cB r^{S})$. For each non-negative integer $n$, the automorphism group of its object $(S,\underline{n})$ is the surface braid group $\pi_{0}(\diffdecbr(S,\underline{n}))$ of $S$, which we denote by $\B_{n}(S)$. When $S=\bD^{2}$, we abbreviate $\Beta^{\bD^{2}}$ to $\Beta$ for brevity.
Since $\cB r^{\bD^{2}}$ contains the solid cylinder $(\bB^2_{1},\emptyset,\id,r)$, Lemma~\ref{lem:monoidal-pi0} implies that the groupoid $\Beta$ inherits a monoidal structure from the boundary connected sum $\natural$ of $\cD_{2}$, which also induces on $\Beta^{S}$ a left module structure over the monoidal groupoid $\Beta$.
Then, the Quillen bracket construction defines categories $\fU\cB r^{\bD^{2}}$, $\langle\cB r^{\bD^{2}},\cB r^{S}\rangle$, $\fU\Beta$ and $\langle \Beta,\Beta^{S}\rangle $. By Lemma~\ref{lem:Serre-fibration-condition}, we have isomorphisms of categories $\pi_{0}(\fU\cB r^{\bD^{2}}) \cong \fU\Beta$ and $\pi_{0}(\langle\cB r^{\bD^{2}},\cB r^{S}\rangle) \cong \langle \Beta,\Beta^{S}\rangle$.
By Proposition~\ref{prop:monoidal-skeleton}, we may pass to monoidal skeleta of $(\Beta,\natural)$ and $(\Beta^{S},\natural)$; after doing this, the groupoid $\Beta$ becomes equivalent to the well-known \emph{braid groupoid}; see \cite[Chap.~XI, \S 4]{MacLane1} for instance. The categories $\fU\Beta$ and $\langle \Beta,\Beta^{S}\rangle$ are then also skeletal by Corollary~\ref{coro:QBC-skeleton}.

An alternative viewpoint on surface braid groups is as fundamental groups of configuration spaces. Fix a non-negative integer $n$ and a compact, connected surface $S$. Then the embedding space $\Emb(\underline{n},\Int{S})$ is the \emph{ordered configuration space} of $n$ points in $\Int{S}$, denoted by $F_n(\Int{S})$. For a partition $\lambda = \boldsymbol{\{} \lambda_{1};\ldots;\lambda_{r} \boldsymbol{\}} \vdash n$, the quotient space $F_n(\Int{S})/{\Sym_{\lambda}}$, induced by the natural action of $\Sym_{\lambda}:= \Sym_{\lambda_{1}} \times \cdots \times \Sym_{\lambda_{r}} \subseteq \Sym_n$ on the coordinates, is the \emph{$\lambda$-partitioned configuration space} of $n$ points in $\Int{S}$, denoted by $C_{\lambda}(\Int{S})$. When $\lambda$ is the trivial partition $\boldsymbol{\{} n \boldsymbol{\}}$ (which we simply denote by $n$), this is the \emph{unordered configuration space} and is denoted by $C_n(\Int{S})$.
The fundamental group of $C_n(\Int{S})$ is isomorphic to $\B_n(S)$ (this is a classical fact, or Proposition~\ref{prop:braided-diff-groups} with $M=S$ and $Z=\underline{n}$). More generally, following Definition~\ref{def:diff_dec_partitioned}, the fundamental group $\pi_{1}(C_\lambda(\Int{S}))$, which is isomorphic to $\pi_{0}(\diffdecbr(S,\lambda_{1},\dots,\lambda_{r}))$ by Proposition~\ref{prop:braided-diff-groups}, is called the \emph{$\lambda$-partitioned braid group} $\B_{\lambda}(S)$.

\subsubsection{Loop braid groups}\label{sss:category_loop_braid_groups}

We now focus on the families of extended and non-extended loop braid groups. Their definitions are recalled here and we refer to \cite{Damianijourney} for a complete introduction to these groups. Let us denote by $\bD^{3}$ the closed $3$-disc.
For each non-negative integer $n$, we denote by $\underline{n}\bS^{1}$ a closed submanifold of $\bD^{3}$ consisting of a collection of $n$ disjoint circles forming a trivial link of $n$ components in the interior of $\bD^{3}$.
The notation $\underline{n}\bS_{+}^{1}$ denotes the same unlink as $\underline{n}\bS^{1}$, but it indicates that we shall require any diffeomorphism on this unlink to be orientation-preserving; otherwise orientation-reversing diffeomorphisms are permitted.
Choosing also two germs of boundary-cylinders (which we elide from the notation), the pair $(\ensuremath{\bD^{3}},\underline{n}\bS^{1})$ forms a decorated manifold, which we sometimes denote by $\bD^{3}_{n}$ for simplicity.
Let $\Diff_{\partial}(\bD^{3},\underline{n}\bS^{1})$ be the group of self-diffeomorphisms of $\bD^{3}$ that fix $\partial\bD^{3}$ pointwise and $\underline{n}\bS^{1}$ setwise. The \emph{extended loop braid group} $\lB_{n}'$ is the group of isotopy classes $\pi_{0}(\Diff_{\partial}(\bD^{3},\underline{n}\bS^{1}))$.
Let $\Diff_{\partial}\bigl(\bD^{3},\underline{n}\bS_{+}^{1}\bigr)$ denote the subgroup of diffeomorphisms that also preserve the orientation of $\underline{n}\bS^{1}$. The (non-extended) \emph{loop braid group} $\lB_{n}$ is the group of isotopy classes $\pi_{0}\bigl(\Diff_{\partial}\bigl(\bD^{3},\underline{n}\bS_{+}^{1}\bigr)\bigr)$.

We now set up a categorical framework for handling these families of groups.
Let $\cL\cB'$ (respectively $\cL\cB$) be the full subgroupoid of $\cD_3$ (respectively $\cD_3^{+}$) on all decorated manifolds $(M,A,e_{1},e_{2})$ such that the pair $(M,A)$ is diffeomorphic to the $3$-disc relative to an embedded $n$-component unlink. Both $\cL\cB'$ and $\cL\cB$ are closed under the semi-monoidal structure $\natural$ of Definition~\ref{prop:Dd-semi-monoidal}, so they inherit a semi-monoidal structure from $\cD_3$ and $\cD_3^{+}$.
In particular, either choice $\cG = \cM = \cL\cB$ or $\cG = \cM = \cL\cB'$ fits into Hypothesis~\ref{hypo:standard_framework}.
We denote $\pi_{0}(\cL\cB')$ and $\pi_{0}(\cL\cB)$ by $\cL\Beta'$ and $\cL\Beta$ respectively. Since $\cL\cB'$ and $\cL\cB$ both contain the solid cylinder $(\bB^3_{1},\emptyset,\id,r)$, the groupoids $\cL\Beta'$ and $\cL\Beta$ both inherit a monoidal structure from the semi-monoidal structures of $\cD_3$ and $\cD_3^{+}$ by Lemma~\ref{lem:monoidal-pi0}.
We may pass to monoidal skeleta of $(\cL\Beta',\natural)$ and $(\cL\Beta,\natural)$ by Proposition~\ref{prop:monoidal-skeleton}.
Hence Quillen's bracket construction defines categories $\fU\cL\Beta'$ and $\fU\cL\Beta$, which are skeletal by Corollary~\ref{coro:QBC-skeleton}. By Lemma~\ref{lem:Serre-fibration-condition}, we have isomorphisms of categories $\pi_{0}(\fU\cL\cB') \cong \fU\cL\Beta'$ and $\pi_{0}(\fU\cL\cB) \cong \fU\cL\Beta$.

Let us show that there are isomorphisms $\Aut_{\cL\Beta'}(\bD^{3}_{n}) \cong \lB_{n}'$ and $\Aut_{\cL\Beta}(\bD^{3}_{n}) \cong \lB_{n}$. We show this for the first case, the other one following by an identical argument. The automorphism group of $\bD^{3}_{n}$ in $\cL\Beta'$ is $\pi_{0}(\Diff_{\mathrm{dec}}(\bD^{3}_{n}))$, where $\Diff_{\mathrm{dec}}(\bD^{3}_{n})$ is the topological group of diffeomorphisms of $\bD^{3}$ that send the embedded $n$-component unlink onto itself and that restrict to the identity on a neighbourhood of two disjoint $2$-discs in $\partial \bD^{3}$. The condition of fixing \emph{a neighbourhood of} two discs in the boundary is equivalent (up to homotopy equivalence, so in particular on $\pi_{0}$) to fixing just the two discs themselves. On the other hand, $\lB_{n}'$ has the same description except that diffeomorphisms must fix the \emph{whole} boundary $\partial \bD^{3}$. To show that these two groups are isomorphic, it therefore suffices to show the following.

\begin{lemm}\label{lem:loop-braid}
Let $M$ be a $3$-manifold with a boundary component $\partial_{0} M\cong \bS^{2}$. For isotopy classes of diffeomorphisms of $M$, fixing two disjoint $2$-discs in $\partial_{0} M$ is equivalent to fixing all of $\partial_{0} M$.
\end{lemm}
\begin{proof}
Let $\Diff(M,\partial_{0} M)$ be the group of diffeomorphisms of $M$ that send $\partial_{0} M$ onto itself. The restriction map
$\Diff(M,\partial_{0} M) \rightarrow \Diff(\partial_{0} M) \cong \Diff(\bS^2)$
is a fibre bundle, by \cite[p.~294, \S II.2.2.2, Cor.~2]{Cerf1961Topologiedecertains}. Hence its restriction $\Diff_{\bD^2 \sqcup \bD^2}(M) \too \Diff_{\bD^2 \sqcup \bD^2}(\bS^2) \cong \Diff_{\partial C}(C)$ is also a fibre bundle, where the subscript ${}_{\bD^2 \sqcup \bD^2}$ means that diffeomorphisms must restrict to the identity on a given pair of disjoint discs in $\partial_{0} M \cong \bS^2$ and $C$ denotes the $2$-dimensional cylinder $\bS^1 \times [0,1]$. The fibre is $\Diff_{\partial_{0} M}(M)$ and we obtain an exact sequence
\[
\cdots \to \pi_{1}(\Diff_{\partial C}(C)) \too \pi_{0}(\Diff_{\partial_{0} M}(M)) \xrightarrow{\;(*)\;} \pi_{0}(\Diff_{\bD^2 \sqcup \bD^2}(M)) \too \pi_{0}(\Diff_{\partial C}(C)).
\]
By \cite[Th{\'e}or{\`e}me 1]{gramain1973type}, $\Diff_{\partial C}(C)$ is contractible, and hence $(*)$ is a bijection.
\end{proof}

Since their automorphism groups are the (extended) loop braid groups, we call $\cL\Beta'$ and $\cL\Beta$ respectively the \emph{extended loop braid groupoid} and the (non-extended) \emph{loop braid groupoid}. Also, recall (see Definition~\ref{def:braided-diffeomorphisms}) the braided diffeomorphism group $\diffdecbr(M,A)$, consisting of all of those (decorated) diffeomorphisms that become isotopic to the identity after forgetting $A$. Since $\pi_{0}(\diffdec(\bD^{3},\emptyset))$ is trivial by \cite{Cerf1968} (and Lemma~\ref{lem:loop-braid}), this condition is vacuous when $M = \bD^3$, so we have $\lB_{n}' \cong \Aut_{\cL\Beta'}(\bD^{3}_{n}) = \pi_{0}(\diffdec(\bD^3_n)) = \pi_{0}(\diffdecbr(\bD^3_n))$ and $\lB_{n} \cong \Aut_{\cL\Beta}(\bD^{3}_{n}) = \pi_{0}(\diffdec^{+}(\bD^3_n)) = \pi_{0}(\diffdecbrplus(\bD^3_n))$. Hence $\lB_{n}'$ and $\lB_{n}$ may be thought of as \emph{braided mapping class groups} (although the adjective \emph{braided} is redundant in this case).

Finally, we introduce here further generalisations of loop braid groups. We fix non-negative integers $n$ and $k$ and a partition $\lambda = \boldsymbol{\{} \lambda_{1};\ldots;\lambda_{r} \boldsymbol{\}} \vdash k$.
Now, consider the embedding spaces $F_{k}(\bD^{3}_{n}) = \Emb(\underline{k},\Int{\bD}^{3}\smallsetminus\underline{n}\bS^{1})$ and $\overline{U}_{k}(\bD^{3}_{n}) = \Emb^{\mathrm{unl}}(\underline{k}\bS^{1},\Int{\bD}^{3}\smallsetminus\underline{n}\bS^{1})$, where the superscript ${}^{\mathrm{unl}}$ means the path-component of the embedding space corresponding to $(n+k)$-component unlinks. We set
\begin{equation}\label{eq:configurations_loop_braid_groups}
C_{\lambda}(\bD^{3}_{n}) = F_{k}(\bD^{3}_{n}) / \Sym_{\lambda},\,\,\,
U_{\lambda}^{+}(\bD^{3}_{n}) = \overline{U}_{k}(\bD^{3}_{n}) / \Diff^{+}(\underline{\lambda}\bS^1),\,\,\, U_{\lambda}(\bD^{3}_{n}) = \overline{U}_{k}(\bD^{3}_{n}) / \Diff(\underline{\lambda}\bS^1),
\end{equation}
where $\underline{\lambda}\bS^1=(\underline{\lambda}_{1}\bS^1,\ldots,\underline{\lambda}_{r}\bS^1)$.
The first is the $\lambda$-partitioned configuration space of $k$ points in the unlink-complement $\Int{\bD}^{3}\smallsetminus\underline{n}\bS^{1}$. The middle space is the space of $\lambda$-partitioned oriented $k$-component unlinks in $\Int{\bD}^{3}\smallsetminus\underline{n}\bS^{1}$ such that the resulting $(n+k)$-component link is again trivial. The right-hand space is similar, except that the $k$-component unlinks are unoriented. By Proposition~\ref{prop:braided-diff-groups}, we have isomorphisms $\pi_{1}(U_{n}^{+}(\bD^{3})) \cong \lB_{n}$ and $\pi_{1}(U_{n}(\bD^{3})) \cong \lB'_{n}$ and, more generally, the groups $\pi_{1}(U_{\lambda}^{+}(\bD^{3}))$ and $\pi_{1}(U_{\lambda}(\bD^{3}))$ are isomorphic to the partitioned versions $\lB_{\lambda}$ and $\lB_{\lambda}'$ respectively.

\subsection{Examples of homological representation functors}
\label{ss:examples}

Recall from \S\ref{ss:homological_representation_functor_motion_groups}--\S\ref{ss:homological_representation_functor_mcg} that the input to construct a homological representation functor is a continuous semifunctor of the form $\Int{\tF}_{(Z,\sG,Q)}$, for $Z$ a closed submanifold of $\bR^d$, $\sG$ an open subgroup of $\Diff(Z)$ and $Q$ a functorial quotient of groups. (Here, $\tF$ refers either to the semifunctors denoted by $\cF$ in Theorem~\ref{thm:global_functor_motion_groups} or the semifunctors denoted by $\fF$ in Theorem~\ref{thm:global_functor_mcg}.)
\emph{Throughout \S\ref{ss:examples}, we consider an integer $k\geq 1$, which will determine the closed submanifold $Z \subset \bR^d$, as well as a partition $\lambda = \boldsymbol{\{} \lambda_{1};\ldots;\lambda_{r} \boldsymbol{\}} \vdash k$, which will determine the group $\sG$. We denote by $r'$ the number of indices $1\leq i \leq r$ such that $\lambda_i \geq 2$}.

Furthermore, we confine our study throughout \S\ref{ss:examples} to the homological representation functors with the following restrictions on the inputs of the construction of Definition~\ref{def:construction}.
\begin{itemizeb}
    \item For concision, we mainly present the ``open'' versions $\Int{\tF}_{(Z,\sG,Q)}$ rather than the ``closed'' variants $\tF_{(Z,\sG,Q)}$, although everything in \S\ref{ss:examples} repeats verbatim using the closed variants $\tF_{(Z,\sG,Q)}$.
    \item We always consider the lower central series functorial quotients of groups $Q_{\LCS_{\ell}}\colon \groups \to \groups$ (see Example~\ref{eg:examples_FQG}) for the parameter $Q$, and we simplify the indexing notation of $\Int{\tF}$ using $\ell$ instead of $Q_{\LCS_{\ell}}$. We choose to focus on constructions with these functorial quotients of groups since the lower central series of (partitioned) surface braid groups and loop braid groups are well-understood (see for example \cite{DPS}), which helps us to understand the specific examples of homological representations fitting into this framework. We recall from Example~\ref{eg:LCS} that the \emph{untwisted} variant semifunctor $\Int{\tF}^{\unt}_{(Z,\sG,\ell)}$ is equal to $\Int{\tF}_{(Z,\sG,\ell)}$ for the parameter $\ell\in\{1,2\}$.
    \item For simplicity, we always take the continuous semifunctor $V \colon \cC \to \modlr$ in the inputs for Definition~\ref{def:construction} to be equal to the \emph{colimit coefficient system} $V = V_{\col}(\Int{\tF}_{(Z,\sG,\ell)})$ associated to $\Int{\tF}_{(Z,\sG,\ell)}$ (see \S\ref{sss:col_coeff_unwtisted_rep} and \S\ref{sss:col_coeff_unwtisted_rep_mcg}). This allows us to apply Propositions~\ref{prop:factoringthroughtopQtw}, \ref{prop:factoringthroughtopQ} and \ref{prop:factoringthroughcovQ-mcg} to deduce that the homological representation functors that we construct take values in $\modr[{\bZ[\cQ]}]^{\tw}$ or $\modr[{\bZ[\cQ]}]$ for a fixed transformation group $\cQ$. We recall that the colimit coefficient system $V_{\col}(\Int{\tF}_{(Z,\sG,\ell)})$ and the transformation group $\cQ$ are determined by the semifunctor $\Int{\tF}_{(Z,\sG,\ell)}$ and the choice of subcategory $\langle \cG , \cM \rangle \subseteq \fU\cD_d$ to which we have restricted it.
    \item The sequence of construction in each setting is thus as follows: we select the parameters $k \geq 1$, $\lambda \vdash k$ and $\ell \geq 1$, determining the data $Z$, $\sG$ and $Q$ respectively; we restrict the continuous semifunctor $\Int{\tF} = \Int{\tF}_{(Z,\sG,\ell)}$ of Theorem~\ref{thm:global_functor_motion_groups} or Theorem~\ref{thm:global_functor_mcg} to a subcategory $\langle \cG , \cM \rangle \subseteq \fU\cD_d$; this determines the colimit coefficient system $V_{\col}(\Int{\tF})$ by \S\ref{sss:col_coeff_unwtisted_rep} or \S\ref{sss:col_coeff_unwtisted_rep_mcg}; the data of $\Int{\tF}$ and $V_{\col}(\Int{\tF})$ then determine a homological representation functor $L_i(\Int{\tF};V_{\col}(\Int{\tF}))$ for each $i\geq 0$.
\end{itemizeb}

\subsubsection{Motion groups}\label{ss:applications_motion_groups}

We apply the general construction of \S\ref{ss:homological_representation_functor_motion_groups} to some families of motion groups (see Definition~\ref{def:motion_groups}).
Namely, we study examples of the homological representation functors introduced in \S\ref{ss:homological_representation_functor_motion_groups}, defined on the subcategories of $\fU\cD_{2}$ and $\fU\cD_3$ that are relevant for \emph{classical braid groups}, \emph{surface braid groups} and \emph{loop braid groups}. These subcategories are introduced in \S\ref{sss:category_surface_braid_groups} and \S\ref{sss:category_loop_braid_groups}.

\paragraph*{Classical braid groups.}

We denote by $\Int{\cF}_{(\underline{k},\Sym_{\lambda},\ell)}(\bD)$ and $\Int{\cF}^{\unt}_{(\underline{k},\Sym_{\lambda},\ell)}(\bD)$ the restrictions to the subcategory $\fU\cB r^{\bD^{2}}\subset \fU\cD_{2}$ of the continuous semifunctors of Theorem~\ref{thm:global_functor_motion_groups}, where we take $Z:=\underline{k}$ a finite set of $k\geq 1$ points, $\sG:=\Sym_{\lambda}$ and $Q:=Q_{\LCS_{\ell}}$ for some $\ell\geq 1$.
Recall from \S\ref{sss:category_surface_braid_groups} that the automorphism groups of $\cB r^{\bD^{2}}$ are the braided diffeomorphism groups $\diffdecbr(\bD^{2},\underline{n})$ (see Definition~\ref{def:braided-diffeomorphisms}), so it is a motion groupoid (see Definition~\ref{def:motion-groupoid}). Also, by definition of the braid groupoid $\Beta$ (see \S\ref{sss:category_surface_braid_groups}), the poset $\cP(\Beta,\Beta)$ (see Definition~\ref{def:underlying-poset}) is isomorphic to the poset of natural numbers $\bN$ with its usual order, so in particular it is a directed set.
Therefore, using the colimit coefficient systems \eqref{eq:colim_coefficient_system} associated to the semifunctors under consideration, we obtain from Corollary~\ref{coro:def_homological_rep_functor_motion_groups}, Proposition~\ref{prop:factoringthroughtopQtw} and Proposition~\ref{prop:factoringthroughtopQ} the following homological representation functors for all $i\geq 0$:
\begin{equation}
\label{eq:output_of_general_construction_classical_braids}
\begin{split}
L_{i}(\Int{\cF}_{(\underline{k},\Sym_{\lambda},\ell)}(\bD)) \colon \fU\Beta &\too \modr[\bZ[\cQ_{(\lambda,\ell)}(\bD)]]^{\tw}; \\
L_{i}(\Int{\cF}^{\unt}_{(\underline{k},\Sym_{\lambda},\ell)}(\bD)) \colon \fU\Beta &\too \modr[\bZ[\cQ^{\unt}_{(\lambda,\ell)}(\bD)]].
\end{split}
\end{equation}
Here, $\cQ_{(\lambda,\ell)}(\bD)$ and $\cQ^{\unt}_{(\lambda,\ell)}(\bD)$ denote the colimit group $\cQ_{\col}(\sR \circ \Int{\cF})$ (Notation~\ref{notation-Qcol}) associated to $\Int{\cF} = \Int{\cF}_{(\underline{k},\Sym_{\lambda},\ell)}(\bD)$ and $\Int{\cF} = \Int{\cF}^{\unt}_{(\underline{k},\Sym_{\lambda},\ell)}(\bD)$ respectively.
For $\ell = 2$ and the trivial partition $\lambda = \boldsymbol{\{} k \boldsymbol{\}}$, the representations arising from the functors \eqref{eq:output_of_general_construction_classical_braids} are related to the well-known families of \emph{Lawrence-Bigelow representations}, originally introduced by Lawrence~\cite{Lawrence1} as representations of Hecke algebras and recovered by Bigelow \cite[\S 2]{bigelow2001braid} following a more geometric method. Namely, they define a $\B_{n}$-representation $\LB_{k}(n)$ for each $k\geq 1$ and $n\geq 0$, called the $k$-th Lawrence-Bigelow representation.
The most famous among these are the Burau representations, first introduced in \cite{burau}, and the \emph{Lawrence-Krammer-Bigelow} representations, which Bigelow \cite[Th.~1.1]{bigelow2001braid} and Krammer \cite{KrammerLK} independently proved to be faithful. The following result may be straightforwardly verified by unwinding the construction of \S\ref{ss:homological_representation_functor_motion_groups} and comparing it to that of \cite[\S 2]{bigelow2001braid}.

\begin{thm}[The Lawrence-Bigelow representations.]
\label{thm:LawrenceBigelowFunctors}
Let $k\geq 1$ and $n\geq 0$ be integers. The $\B_{n}$-representation encoded by the functor $L_{k}(\Int{\cF}_{(\underline{k},\Sym_{k},2)}(\bD))$, i.e.~\eqref{eq:output_of_general_construction_classical_braids} with $i=k=\lambda$ and $\ell = 2$, is isomorphic to the Lawrence-Bigelow representation $\LB_{k}(n)$.
\end{thm}

In addition, we compute that $\cQ_{(\lambda,2)}(\bD) \cong \bZ^{r'} \times \bZ^{r(r-1)/2} \times \bZ^{r}$; see Lemma~\ref{lem:transformation_groups_ab_quotient_surface_braid_groups}.
As far as the authors know, there are no representations of the braid groups in the literature whose ground rings are of the form $\bZ[\bZ^{r'} \times \bZ^{r(r-1)/2} \times \bZ^{r}]$ for $r\geq 2$. Hence the functors \eqref{eq:output_of_general_construction_classical_braids} for $\ell = 2$ appear to define new representations of the braid groups.

Furthermore, if $\lambda$ is a partition of the form $\boldsymbol{\{}2;\lambda'\boldsymbol{\}}$ or $\boldsymbol{\{}1;1;1;\lambda''\boldsymbol{\}}$ where $\lambda'$ and $\lambda''$ are any partitions of $k-2$ and $k-3$ respectively, we prove in the sequel \cite{PSIN} that $L_{i}(\Int{\cF}_{(\underline{k},\Sym_{\lambda},\ell)}(\bD))\neq L_{i}(\Int{\cF}_{(\underline{k},\Sym_{\lambda},\ell + 1)}(\bD))$ for each $\ell \geq 2$; see \cite[Table~2]{PSIN}. The transformation group $\cQ_{(\boldsymbol{\{}2;\lambda'\boldsymbol{\}},\ell)}(\bD)$ is computed in Proposition~\ref{prop:transformation_groups_{2}_nilpotent_quotient_surface_braid_groups} for $\ell\geq 2$ when the partition $\lambda' \vdash (k-2)$ is such that $\lambda'_{l}\geq 3$ for all $1\leq l\leq r-1$.
Also, one may easily deduce from \cite[Th.~1.1]{bigelow2001braid} that each functor $L_{k}(\Int{\cF}_{(\underline{k},\Sym_{\boldsymbol{\{}2;\lambda'\boldsymbol{\}}},\ell)}(\bD))$ for $\ell\geq 2$ encodes a \emph{faithful} representation of $\B_{n}$ for each $n$; see \cite[Rem.~4.8]{PSIN}.

\paragraph*{Surface braid groups.}

We consider the continuous semifunctors $\Int{\cF}_{(\underline{k},\Sym_{\lambda},\ell)}$ and $\Int{\cF}^{\unt}_{(\underline{k},\Sym_{\lambda},\ell)}$ of Theorem~\ref{thm:global_functor_motion_groups} where we take $Z:=\underline{k}$ a finite set of $k\geq 1$ points, $\sG:=\Sym_{\lambda}$ and $Q:=Q_{\LCS_{\ell}}$ for some $\ell\geq 1$.
Let $S$ be a compact, connected surface with boundary, different from the $2$-disc: it is therefore homeomorphic to either $\Sigma_{g,1}$ or $\N_{h,1}$ for $g \text{ or } h\geq 1$.
We denote by $\Int{\cF}_{(\underline{k},\Sym_{\lambda},\ell)}(S)$ and $\Int{\cF}^{\unt}_{(\underline{k},\Sym_{\lambda},\ell)}(S)$ the restrictions of these semifunctors to the subcategory $\langle\cB r^{\bD^{2}},\cB r^{S}\rangle \subset \fU\cD_{2}$.
Recall from \S\ref{sss:category_surface_braid_groups} that the automorphism groups of $\Beta^{S}$ are the braided diffeomorphism groups $\diffdecbr(S,\underline{n})$ (see Definition~\ref{def:braided-diffeomorphisms}), so it is a motion groupoid (see Definition~\ref{def:motion-groupoid}).
Also, by definition of the groupoid $\Beta^{S}$ (see \S\ref{sss:category_surface_braid_groups}), the poset $\cP(\Beta,\Beta^{S})$ (see Definition~\ref{def:underlying-poset}) is isomorphic to the poset of natural numbers $\bN$ with its usual order, so in particular it is a directed set.
By Corollary~\ref{coro:def_homological_rep_functor_motion_groups}, Proposition~\ref{prop:factoringthroughtopQtw} and Proposition~\ref{prop:factoringthroughtopQ}, using the colimit coefficient systems \eqref{eq:colim_coefficient_system} associated to these semifunctors, we have the following homological representation functors for all $i\geq 0$:
\begin{equation}
\label{eq:output_of_general_construction_surface_braids}
\begin{split}
L_{i}(\Int{\cF}_{(\underline{k},\Sym_{\lambda},\ell)}(S)) \colon \langle \Beta,\Beta^{S}\rangle &\too \modr[\bZ[\cQ_{(\lambda,\ell)}(S)]]^{\tw}; \\
L_{i}(\Int{\cF}^{\unt}_{(\underline{k},\Sym_{\lambda},\ell)}(S)) \colon \langle \Beta,\Beta^{S}\rangle &\too \modr[\bZ[\cQ^{\unt}_{(\lambda,\ell)}(S)]].
\end{split}
\end{equation}
Here, $\cQ_{(\lambda,\ell)}(S)$ and $\cQ^{\unt}_{(\lambda,\ell)}(S)$ denote the colimit group $\cQ_{\col}(\sR \circ \Int{\cF})$ (Notation~\ref{notation-Qcol}) associated to $\Int{\cF} = \Int{\cF}_{(\underline{k},\Sym_{\lambda},\ell)}(S)$ and $\Int{\cF} = \Int{\cF}^{\unt}_{(\underline{k},\Sym_{\lambda},\ell)}(S)$ respectively.
For $\ell\in \{2,3\}$, we compute that $\cQ_{(\lambda,2)}(S) \cong (\bZ/2)^{r'} \times H_{1}(S;\bZ)^{r}$ (see Lemma~\ref{lem:transformation_groups_ab_quotient_surface_braid_groups}), while, when the partition $\lambda \vdash k$ is such that $\lambda_{l}\geq 3$ for all $1\leq l\leq r$, the transformation groups $\cQ_{(\lambda,3)}(\Sigma_{g,1})$, $\cQ_{(\lambda,3)}^{\unt}(\Sigma_{g,1})$, $\cQ_{(\lambda,3)}(\N_{h,1})$ and $\cQ_{(\lambda,3)}^{\unt}(\N_{h,1})$ are calculated in Proposition~\ref{prop:transformation_groups_{2}_nilpotent_quotient_surface_braid_groups} (see \eqref{eq:transformation_group_orientable_partitioned} and \eqref{eq:transformation_group_non-orientable_partitioned}).
Apart from one specific case when $\ell = 3$ detailed in Example~\ref{eg:An-Ko} below, there are no representations of surface braid groups in the literature whose ground rings coincide with those of the functors \eqref{eq:output_of_general_construction_surface_braids} for $\ell\in \{2,3\}$, which thus appear to be new.

\begin{eg}[The An-Ko representations.]
\label{eg:An-Ko}
Let us focus on the closed variant $\cF_{(\underline{k},\Sym_{k},3)}(\Sigma_{g,1})$ for orientable surfaces. Since it takes values in $\covr^{\pr}$ by Lemma~\ref{lem:F-proper}, we may use Borel-Moore homology in the construction of Definition~\ref{def:construction} and we denote by $L^{\BM}_{k}(\cF_{(\underline{k},\Sym_{k},3)}(\Sigma_{g,1}))$ the corresponding homological representation functor. This functor defines representations of $\B_{n}(\Sigma_{g,1})$ for any $k\geq 1$ and $n\geq 0$.
The procedure to construct these representations may be seen as a reinterpretation following Bellingeri, Godelle and Guaschi \cite{BellingeriGodelleGuaschi} of the work of An and Ko \cite{AnKo}, who extend some homological representations from the classical braid groups to the surface braid groups.
Namely, in \cite[\S 3.A]{AnKo}, the transformation group $\cQ_{(k,3)}(\Sigma_{g,1})$ and the quotient $\B_{k,n}(\Sigma_{g,1})/\LCS_{3}$ are introduced in a completely different way, as abstract groups satisfying certain technical homological constraints. Then, for $k\geq 3$, \cite[\S 4]{BellingeriGodelleGuaschi} redefine these groups via the third lower central quotient of $\B_{k,n}(\Sigma_{g,1})$ and prove that they are isomorphic to those of \cite[\S 3.A]{AnKo}.

For $k\geq 3$, one may straightforwardly check from the definitions that the $k$-th An-Ko representation \cite[Th.~3.2]{AnKo} of $\B_{n}(\Sigma_{g,1})$ is given by the tensor product $L^{\BM}_{k}(\cF_{(\underline{k},\Sym_{k},3)}(\Sigma_{g,1}))(n) \otimes_{\bZ[\cQ_{(k,3)}(\Sigma_{g,1})]} \bZ[\B_{k,n}(\Sigma_{g,1})/\LCS_{3}]$. The case of $k\leq2$ is trickier: a careful analysis of the definitions of \cite[\S 3.A]{AnKo} shows, using the canonical map $\cQ_{(k,3)}(\Sigma_{g,1})\to \cQ_{(3,3)}(\Sigma_{g,1})\to \B_{3,n}(\Sigma_{g,1})/\LCS_{3}$, that for $k\leq 2$ the $k$-th An-Ko representation \cite[Th.~3.2]{AnKo} corresponds to $L^{\BM}_{k}(\cF_{(\underline{k},\Sym_{k},3)}(\Sigma_{g,1}))(n) \otimes_{\bZ[\cQ_{(k,3)}(\Sigma_{g,1})]} \bZ[\B_{3,n}(\Sigma_{g,1})/\LCS_{3}]$.
The general method applied in this section thus elucidates the geometric origins of these groups, moreover proving that the use of the third lower central quotient is a key tool in defining the homological representations. Our framework also gives an alternative to the technical result \cite[Lem.~3.1]{AnKo} to justify that the representations are well-defined.
In contrast, the representations encoded by the functors $L_{i}^{\BM}(\cF_{(\underline{k},\Sym_{\lambda},3)}(S))$ for any non-trivial partition $\lambda\vdash k$, as well as all of the untwisted versions $L_{i}^{\BM}(\cF^{\unt}_{(\underline{k},\Sym_{\lambda},3)}(S))$, appear to be new.
\end{eg}

In addition, if $\lambda$ is a partition of the form $\boldsymbol{\{}1;\lambda'\boldsymbol{\}}$ or $\boldsymbol{\{}2;\lambda''\boldsymbol{\}}$ where $\lambda'$ and $\lambda''$ are any partitions of $k-1$ and $k-2$ respectively, and for $S\in\{\Sigma_{g,1}, \N_{h,1} \mid g\geq 1,h\geq 2\}$, the sequel \cite{PSIN} proves that $L_{i}(\Int{\cF}_{(\underline{k},\Sym_{\lambda},\ell)}(S)) \neq L_{i}(\Int{\cF}_{(\underline{k},\Sym_{\lambda},\ell + 1)}(S))$ for each $\ell \geq 2$; see \cite[Prop.~5.2]{PSIN}. Therefore, for the aforementioned reasons, the functors \eqref{eq:output_of_general_construction_classical_braids} for $\ell\geq 2$ appear, in general, to define new representations of the braid groups on the surfaces $\Sigma_{g,1}$ and $\N_{h,1}$.

\paragraph*{Loop braid groups.}

To apply the construction of \S\ref{ss:homological_representation_functor_motion_groups} to extended loop braid groups, we consider restrictions of the continuous semifunctor $\Int{\cF}_{(Z,\sG,\ell)}\colon \fU\cD_3 \rightarrow \covr$ of Theorem~\ref{thm:global_functor_motion_groups}, with $Q:=Q_{\LCS_{\ell}}$ for some $\ell\geq 1$, to the subcategory $\fU\cL\cB' \subset \fU\cD_3$ (see \S\ref{sss:category_loop_braid_groups}). For non-extended loop braid groups, we instead consider restrictions to $\fU\cL\cB \subset \fU\cD_3^{+}$ (see \S\ref{sss:category_loop_braid_groups}) of the continuous semifunctor $\Int{\cF}_{(Z,\sG,\ell)}\colon \fU\cD^{+}_3 \to \covr$; see Theorem~\ref{thm:global_functor_motion_groups}.
The automorphism groups of $\cL\cB'$ and $\cL\cB$ are the (braided) diffeomorphism groups $\diffdecbr(\bD^{3}_{n})$ and $\diffdecbrplus(\bD^{3}_{n})$ respectively; see \S\ref{sss:category_loop_braid_groups}. Hence $\cL\cB'$ and $\cL\cB$ are motion groupoids (see Definition~\ref{def:motion-groupoid}).
It also follows from the definitions that the posets $\cP(\cL\cB',\cL\cB')$ and $\cP(\cL\cB,\cL\cB)$ (see Definition~\ref{def:underlying-poset}) are both isomorphic to the poset of natural numbers $\bN$ with its usual order, so in particular they are directed sets.
Now, two choices for the submanifold $Z \subset \bR^3$ naturally arise as relevant inputs to construct homological representations: a set of points or an unlink.

\paragraph*{Using configurations of points.}

We take $Z=\underline{k}$ a set of $k\geq 1$ points in $\bR^{3}$ and $\sG=\Sym_{\lambda}$.
We consider the restrictions of the continuous semifunctors of Theorem~\ref{thm:global_functor_motion_groups} to the subcategories $\fU\cL\cB\subset \fU\cD^{+}_{3}$ and $\fU\cL\cB'\subset \fU\cD_{3}$.
Using the colimit coefficient systems \eqref{eq:colim_coefficient_system} associated to these semifunctors, by Corollary~\ref{coro:def_homological_rep_functor_motion_groups}, Proposition~\ref{prop:factoringthroughtopQtw} and Proposition~\ref{prop:factoringthroughtopQ}, we have the following homological representation functors for all $i\geq 0$:
\begin{equation}
\label{eq:output_of_general_construction_loop_braids_points}
\begin{split}
L_{i}(\Int{\cF}_{(\underline{k},\Sym_{\lambda},\ell)}(\bD^{3}))\colon\fU\cL\Beta &\too \modr[\bZ[\cQ_{(P,\lambda,\ell)}(\bD^{3})]]^{\tw}; \\
L_{i}(\Int{\cF}^{\unt}_{(\underline{k},\Sym_{\lambda},\ell)}(\bD^{3}))\colon\fU\cL\Beta &\too \modr[\bZ[\cQ^{\unt}_{(P,\lambda,\ell)}(\bD^{3})]]; \\
L_{i}(\Int{\cF}'_{(\underline{k},\Sym_{\lambda},\ell)}(\bD^{3}))\colon\fU\cL\Beta' &\too \modr[\bZ[\cQ'_{(P,\lambda,\ell)}(\bD^{3})]]^{\tw}; \\
L_{i}(\Int{\cF}'^{\unt}_{(\underline{k},\Sym_{\lambda},\ell)}(\bD^{3}))\colon\fU\cL\Beta' &\too \modr[\bZ[\cQ'^{\unt}_{(P,\lambda,\ell)}(\bD^{3})]].
\end{split}
\end{equation}
Here, $\cQ_{(P,\lambda,\ell)}(\bD^{3})$ and $\cQ^{\unt}_{(P,\lambda,\ell)}(\bD^{3})$ denote the colimit group $\cQ_{\col}(\sR \circ \Int{\cF})$ (Notation~\ref{notation-Qcol}) associated to $\Int{\cF} = \Int{\cF}_{(\underline{k},\Sym_{\lambda},\ell)}(\bD^{3})$ and $\Int{\cF} = \Int{\cF}^{\unt}_{(\underline{k},\Sym_{\lambda},\ell)}(\bD^{3})$ respectively, and similarly for the versions with $'$. The symbol `$P$' here simply indicates that we are in the setting where $Z$ is a set of \emph{points}.
For $\ell = 2$, we compute that $\cQ_{(P,\lambda,2)}(\bD^{3})\cong \bZ^{r}\times (\bZ/2)^{r'}$ and $\cQ'_{(P,\lambda,2)}(\bD^{3})\cong (\bZ/2)^{r+r'}$; see Lemma~\ref{lem:transformation_groups_ab_quotient_loop_braid_groups}. Among the representations introduced by \eqref{eq:output_of_general_construction_loop_braids_points}, the only ones that we currently understand in detail are those defined with the parameters $i=k=1$ and $\ell = 2$:

\begin{eg}[The loop Burau representations.]
\label{eg:loop_Burau}
In \cite{PS0}, we explicitly compute the matrices of the representations encoded by the functors $L_{1}(\Int{\cF}_{(\underline{1},0,2)}(\bD^{3}))$ and $L_{1}(\Int{\cF}'_{(\underline{1},0,2)}(\bD^{3}))$: these extend the \emph{Burau representations} of the classical braid groups to $\lB_{n}$ and $\lB'_{n}$ respectively. To describe these calculations, we recall that the loop braid group $\lB_{n}$ admits a presentation given by generators $\{ \sigma_{i},\tau_{i}\mid 1\leq i\leq n-1\}$, where $\{\sigma_{1},\ldots,\sigma_{n-1}\}$ satisfy the relations of the classical braid group $\B_{n}$ and $\{\tau_{1},\ldots,\tau_{n-1}\}$ satisfy those of the symmetric group $\Sym_{n}$, together with three additional mixed relations; see \cite[Prop.~3.14 and 3.16]{Damianijourney} for details. We show in \cite{PS0} that the matrices of the representations $L_{1}(\Int{\cF}_{(\underline{1},0,2)}(\bD^{3}))(n)\colon\lB_{n}\rightarrow\Aut_{\bZ[\bZ]}(\bZ[\bZ]^{\oplus n-1})$ are those of:
\begin{itemizeb}
\item the Burau representation $\LB_{1}(n)$ (see \S\ref{ss:applications_motion_groups}) for the generators $\{\sigma_{1},\ldots,\sigma_{n-1}\}$;
\item the standard representation of the symmetric group $\Sym_{n}$ for the generators $\{\tau_{1},\ldots,\tau_{n-1}\}$.
\end{itemizeb}
The matrices for the representations $L_{1}(\Int{\cF}'_{(\underline{1},0,2)}(\bD^{3}))(n)$ of the extended loop braid groups $\lB'_{n}$ over $R = \bZ[\bZ/2] \cong \bZ[t^{\pm 1}] / (t^2 - 1)$ are more subtle, since the underlying $R$-module is not free in this case: it is $R^{\oplus n-1} \oplus R/(t-1)$. They can however be computed; see \cite[Table~1]{PS0}.
\end{eg}

On the other hand, apart from the setting of Example~\ref{eg:loop_Burau}, the representations of the (extended and non-extended) loop braid groups encoded by the functors \eqref{eq:output_of_general_construction_loop_braids_points} appear to be new because of their ground rings.
Furthermore, for any partition $\lambda = \boldsymbol{\{}2;\lambda'\boldsymbol{\}}$, we prove in \cite[Prop.~6.2]{PSIN} that $L_{i}(\Int{\cF}_{(\underline{k},\Sym_{\lambda},\ell)}(\bD^{3}))\neq L_{i}(\Int{\cF}_{(\underline{k},\Sym_{\lambda},\ell + 1)}(\bD^{3}))$ for each $\ell \geq 2$.

\paragraph*{Using configurations of unlinks.}

We now set $Z=\underline{k}\bS^{1}$ a $k$-component unlink in $\bR^{3}$. For each partition $\lambda \vdash k$, we focus on two choices for the group $\sG$: the group $\Diff(\underline{\lambda}\bS^{1})$ of diffeomorphisms of $\underline{k}\bS^{1}$ preserving the partition $\lambda$ or its subgroup $\Diff^{+}(\underline{\lambda}\bS^{1})$ of orientation-preserving diffeomorphisms.
Then we consider the restrictions of the continuous semifunctors of Theorem~\ref{thm:global_functor_motion_groups} to the full subcategories $\fU\cL\cB\subset \fU\cD^{+}_{3}$ and $\fU\cL\cB'\subset \fU\cD_{3}$ along with their associated colimit coefficient systems \eqref{eq:colim_coefficient_system}.

\textbf{\emph{The oriented version.}}
We take $\sG=\Diff^{+}(\underline{\lambda}\bS^{1})$. By Corollary~\ref{coro:def_homological_rep_functor_motion_groups}, Proposition~\ref{prop:factoringthroughtopQtw} and Proposition~\ref{prop:factoringthroughtopQ}, we have the following homological representation functors for all $i\geq 0$:
\begin{equation}
\label{eq:output_of_general_construction_loop_braids_unlinks}
\begin{split}
L_{i}(\Int{\cF}_{(\underline{k}\bS^{1},\Diff^{+}(\underline{\lambda}\bS^{1}),\ell)}(\bD^{3}))\colon\fU\cL\Beta \too \modr[\bZ[\cQ_{(\bS_{+},\lambda,\ell)}(\bD^{3})]]^{\tw};\\
L_{i}(\Int{\cF}^{\unt}_{(\underline{k}\bS^{1},\Diff^{+}(\underline{\lambda}\bS^{1}),\ell)}(\bD^{3}))\colon\fU\cL\Beta \too \modr[\bZ[\cQ^{\unt}_{(\bS_{+},\lambda,\ell)}(\bD^{3})]];\\
L_{i}(\Int{\cF}'_{(\underline{k}\bS^{1},\Diff^{+}(\underline{\lambda}\bS^{1}),\ell)}(\bD^{3}))\colon\fU\cL\Beta' \too \modr[\bZ[\cQ'_{(\bS_{+},\lambda,\ell)}(\bD^{3})]]^{\tw};\\
L_{i}(\Int{\cF}'^{\unt}_{(\underline{k}\bS^{1},\Diff^{+}(\underline{\lambda}\bS^{1}),\ell)}(\bD^{3}))\colon\fU\cL\Beta' \too \modr[\bZ[\cQ'^{\unt}_{(\bS_{+},\lambda,\ell)}(\bD^{3})]].
\end{split}
\end{equation}
The group $\cQ_{(\bS_{+},\lambda,\ell)}(\bD^{3})$ and its variants with ${}^{\unt}$ or $'$ are defined similarly to above, except that the symbol `$\bS_{+}$' indicates that we are now in the setting where $Z$ is an \emph{oriented unlink}.
For $\ell = 2$, we compute that $\cQ_{(\bS_{+},\lambda,2)}(\bD^{3})\cong \bZ^{r^2+r+r'}\times (\bZ/2)^{r'}$ and $\cQ'_{(\bS_{+},\lambda,2)}(\bD^{3})\cong \bZ^{r^2+r'}\times (\bZ/2)^{r+r'}$; see Lemma~\ref{lem:transformation_groups_ab_quotient_loop_braid_groups}. Moreover, for partitions of the form $\lambda =\boldsymbol{\{}b;\lambda'\boldsymbol{\}}$ with $b\in\{2,3,\boldsymbol{\{}1,1\boldsymbol{\}}\}$ and generically denoting both $\Int{\cF}_{(\underline{k}\bS^{1},\Diff^{+}(\underline{\lambda}\bS^{1}),\ell)}$ and $\Int{\cF}'_{(\underline{k}\bS^{1},\Diff^{+}(\underline{\lambda}\bS^{1}),\ell)}$ by $\cF_{(\lambda,\ell)}$, we prove in \cite[Prop.~6.2]{PSIN} that $L_{i}(\cF_{(\lambda,\ell)}(\bD^{3}))\neq L_{i}(\cF_{(\lambda,\ell + 1)}(\bD^{3}))$ for each $\ell \geq 2$.

\textbf{\emph{The unoriented version.}}
We take $\sG=\Diff(\underline{\lambda}\bS^{1})$. By Corollary~\ref{coro:def_homological_rep_functor_motion_groups}, Proposition~\ref{prop:factoringthroughtopQtw} and Proposition~\ref{prop:factoringthroughtopQ}, we have the following homological representation functors for all $i\geq 0$:
\begin{equation}
\label{eq:output_of_general_construction_loop_braids_unlinks_not_preserved}
\begin{split}
L_{i}(\Int{\cF}_{(\underline{k}\bS^{1},\Diff(\underline{\lambda}\bS^{1}),\ell)}(\bD^{3}))\colon\fU\cL\Beta \too \modr[\bZ[\cQ_{(\bS,\lambda,\ell)}(\bD^{3})]]^{\tw};\\
L_{i}(\Int{\cF}^{\unt}_{(\underline{k}\bS^{1},\Diff(\underline{\lambda}\bS^{1}),\ell)}(\bD^{3}))\colon\fU\cL\Beta \too \modr[\bZ[\cQ^{\unt}_{(\bS,\lambda,\ell)}(\bD^{3})]];\\
L_{i}(\Int{\cF}'_{(\underline{k}\bS^{1},\Diff(\underline{\lambda}\bS^{1}),\ell)}(\bD^{3}))\colon\fU\cL\Beta' \too \modr[\bZ[\cQ'_{(\bS,\lambda,\ell)}(\bD^{3})]]^{\tw};\\
L_{i}(\Int{\cF}'^{\unt}_{(\underline{k}\bS^{1},\Diff(\underline{\lambda}\bS^{1}),\ell)}(\bD^{3}))\colon\fU\cL\Beta' \too \modr[\bZ[\cQ'^{\unt}_{(\bS,\lambda,\ell)}(\bD^{3})]].
\end{split}
\end{equation}
The group $\cQ_{(\bS,\lambda,\ell)}(\bD^{3})$ and its variants with ${}^{\unt}$ or $'$ are defined similarly to above, except that the symbol `$\bS$' indicates that we are now in the setting where $Z$ is an \emph{unoriented unlink}.
For $\ell = 2$, we compute that $\cQ_{(\bS,\lambda,2)}(\bD^{3})\cong \bZ^{r}\times (\bZ/2)^{2r'+r(r+1)}$ and $\cQ'_{(\bS,\lambda,2)}(\bD^{3})\cong (\bZ/2)^{2r'+r(r+2)}$; see Lemma~\ref{lem:transformation_groups_ab_quotient_loop_braid_groups}. Moreover, for partitions of the form $\lambda =\boldsymbol{\{}b;\lambda'\boldsymbol{\}}$ with $b\in\{1,2,3\}$ and generically denoting both $\Int{\cF}_{(\underline{k}\bS^{1},\Diff(\underline{\lambda}\bS^{1}),\ell)}$ and $\Int{\cF}'_{(\underline{k}\bS^{1},\Diff(\underline{\lambda}\bS^{1}),\ell)}$ by $\cF_{(\lambda,\ell)}$, we prove in \cite[Prop.~6.2]{PSIN} that $L_{i}(\cF_{(\lambda,\ell)}(\bD^{3}))\neq L_{i}(\cF_{(\lambda,\ell + 1)}(\bD^{3}))$ for each $\ell \geq 2$.

Therefore, because of the computations of their ground rings, all of the representations encoded by the functors \eqref{eq:output_of_general_construction_loop_braids_unlinks} and \eqref{eq:output_of_general_construction_loop_braids_unlinks_not_preserved} appear to be new.

\subsubsection{Mapping class groups of surfaces}\label{ss:mcg_construction}

Although it is best adapted to motion groups, the general construction of \S\ref{ss:homological_representation_functor_motion_groups} may also be used to construct representations of mapping class groups, and recovers several classical constructions. This is detailed in the first paragraph below. We then apply, in the second paragraph, the method of \S\ref{ss:homological_representation_functor_mcg} to define other representations for mapping class groups of surfaces.

\paragraph*{Homological representations from the construction of \S\ref{ss:homological_representation_functor_motion_groups}.}

Let $\fU\cD_{2}^{\emptyset}$ denote the full subcategory of $\fU\cD_{2}$ on the decorated manifolds $(M,A)$ where $A = \emptyset$. We consider the continuous semifunctor $\Int{\cF}_{(\underline{k},\Sym_{\lambda},\ell)}$ of Theorem~\ref{thm:global_functor_motion_groups} where $Z:=\underline{k}$ is a finite set of $k\geq 1$ points, $\sG:=\Sym_{\lambda}$ and $Q:=Q_{\LCS_{\ell}}$ for some $\ell\geq 1$, which we first restrict to $\fU\cD_{2}^{\emptyset}$. We note that, on this subcategory, the short exact sequences of the diagram \eqref{eq:split-short-exact-sequence-quotient} degenerate in the sense that the right-hand side is the trivial group.
Restricting further, we then denote by $\Int{\cF}_{(\underline{k},\Sym_{\lambda},\ell)}(\MCGo)$ and $\Int{\cF}_{(\underline{k},\Sym_{\lambda},\ell)}(\MCGno)$ the restrictions of this semifunctor to the subcategories $\fU\cM_{2}^{+,\mt}$ and $\fU\cM_{2}^{-,\mt}$ of $\fU\cD_{2}^{\emptyset}$ (see \S\ref{sss:category_mcg}).
Note that, by definition of the groupoids $\cM_{2}^{+}$ and $\cM_{2}^{-}$, the posets $\cP(\cM_{2}^{+},\cM_{2}^{+})$ and $\cP(\cM_{2}^{-},\cM_{2}^{-})$ (see Definition~\ref{def:underlying-poset}) are both isomorphic to the poset of natural numbers $\bN$ with its usual order, so in particular they are directed sets.
Hence, using the colimit coefficient systems \eqref{eq:colim_coefficient_system} associated to the semifunctors under consideration, we obtain from Corollary~\ref{coro:def_homological_rep_functor_motion_groups} and Proposition~\ref{prop:factoringthroughtopQtw} the following homological representation functors for all $i\geq 0$:
\begin{equation}
\label{eq:output_of_general_construction_mapping_class_motion}
\begin{split}
L_{i}(\Int{\cF}_{(\underline{k},\Sym_{\lambda},\ell)}(\MCGo)) \colon \fU\cM_{2}^{+} &\too \modr[\bZ[\cQ^{\cF}_{(\lambda,\ell)}(\MCGo)]]^{\tw}; \\
L_{i}(\Int{\cF}_{(\underline{k},\Sym_{\lambda},\ell)}(\MCGno)) \colon \fU\cM_{2}^{-} &\too \modr[\bZ[\cQ^{\cF}_{(\lambda,\ell)}(\MCGno)]]^{\tw}.
\end{split}
\end{equation}

\begin{rmk}
\label{rmk:no-untwisted-versions}
Since the groupoids $\cM_{2}^{+}$ and $\cM_{2}^{-}$ are not motion groupoids, Proposition~\ref{prop:factoringthroughtopQ} does not apply in this setting, so we cannot deduce that the $\Int{\cF}^{\unt}$ variants of Theorem~\ref{thm:global_functor_motion_groups} lead to untwisted versions of \eqref{eq:output_of_general_construction_mapping_class_motion}; thus we do not consider these variants here. This issue is precisely what is solved in the next paragraph below, by applying the construction of \S\ref{ss:homological_representation_functor_mcg} instead of that of \S\ref{ss:homological_representation_functor_motion_groups}.
\end{rmk}

The transformation groups $\cQ^{\cF}_{(\lambda,\ell)}(\MCGo)$ and $\cQ^{\cF}_{(\lambda,\ell)}(\MCGno)$ in \eqref{eq:output_of_general_construction_mapping_class_motion} are defined (via taking a colimit) from the $\ell$-th lower central quotients $\B_{\lambda}(\Sigma_{g,1})/\LCS_{\ell}$ or $\B_{\lambda}(\N_{h,1})/\LCS_{\ell}$ for all $g,h\geq 1$ respectively. These lower central quotients are computed in \cite[\S 6.5--\S 6.6]{DPS}, from which one may deduce complete calculations of all of these transformation groups.

A simple modification of the above construction of \eqref{eq:output_of_general_construction_mapping_class_motion} consists in replacing the lower central quotient $Q = Q_{\LCS_{\ell}}$ with the identity quotient $Q = \id_{\groups}$. The corresponding functors encode the representations given by the natural actions on the homology of the universal covers of $\lambda$-partitioned configuration spaces on the surfaces. In particular, it is easy to check from the definitions that the functor corresponding to $\lambda = \boldsymbol{\{} 1 \boldsymbol{\}}$ and $Q = \id_{\groups}$ encodes the (reduced) \emph{Magnus representations} of the mapping class groups, whose topological interpretation was introduced in \cite{Suzuki}.

Another interesting modification of the construction of \eqref{eq:output_of_general_construction_mapping_class_motion} consists in removing the basepoint $p_{0} \in \partial S$ from the surface and allowing the configuration points to lie in the boundary of the surface. Namely, we replace each surface $S$ by the surface $S\smallsetminus\{ p_{0}\}$ (formally, this means slightly changing the domain category of the homological representation functors, but its automorphism groups are the same up to isomorphism) and we use the closed variant $\cF_{(\underline{k},\Sym_{\lambda},\ell)}$ of \eqref{eq:global-functor2} instead of $\Int{\cF}_{(\underline{k},\Sym_{\lambda},\ell)}$, so that Borel-Moore homology may be applied to this variant by Lemma~\ref{lem:F-proper}. In this way we define homological representation functors $L^{\BM}_{i}(\cF_{(\underline{k},\Sym_{\lambda},\ell)}(\MCGo))$ and $L^{\BM}_{i}(\cF_{(\underline{k},\Sym_{\lambda},\ell)}(\MCGno))$ analogous to those of \eqref{eq:output_of_general_construction_mapping_class_motion}. These alternatives have the advantage of encoding representations endowed with natural free generating sets; see \cite[\S 2.2]{PSIIp}.
In particular, we make here the connection between the representations arising from the functor $L^{\BM}_{k}(\cF_{(\underline{k},0,1)}(\MCGo))$ (i.e.~$i=k$, $\lambda = \boldsymbol{\{} 1^{k}\boldsymbol{\}}$ and $\ell = 1$) with those introduced by Moriyama in \cite{Moriyama}, whose kernels are the $k$-th terms of the \emph{Johnson filtration}. More precisely, for each $g$, the $\MCGo_{g,1}$-representation introduced by \cite{Moriyama} is given by the $\MCGo_{g,1}$-action on the relative homology group $H_{k}(\Sigma_{g,1}^{\times k} , \Delta \cup A_{g};\bZ)$, where $\Delta$ denotes the ``fat diagonal'' of $\Sigma_{g,1}^{\times k}$ where at least two points coincide and $A_{g}$ denotes the subspace of $\Sigma_{g,1}^{\times k}$ where at least one point is equal to $p_{0}$, a chosen basepoint on $\partial \Sigma_{g,1}$.

\begin{prop}\label{prop:Moriyama_recovering}
The restriction of the functor $L^{\BM}_{k}(\cF_{(\underline{k},0,1)}(\MCGo))$ to the mapping class group $\MCGo_{g,1}$ of genus $g$ is isomorphic to the $k$-th Moriyama representation $\MCGo_{g,1}\to\Aut_{\bZ}(H_{k}(\Sigma_{g,1}^{\times k} , \Delta \cup A_{g};\bZ))$.
\end{prop}
\begin{proof}
We write $\Sigma_{g,1}' = \Sigma_{g,1} \smallsetminus \{p_{0}\}$ and $F_{k}(\Sigma_{g,1}')$ for $C_{\boldsymbol{\{} 1^{k}\boldsymbol{\}}}(\Sigma_{g,1}')$. Since $\Sigma_{g,1}^{\times k}$ is a compactification of $F_{k}(\Sigma_{g,1}') = \Sigma_{g,1}^{\times k} \smallsetminus (\Delta \cup A_{g})$, the Borel-Moore homology of $F_{k}(\Sigma_{g,1}')$ is isomorphic to the relative homology $H_*(\Sigma_{g,1}^{\times k} , \Delta \cup A_{g};\bZ)$. That this isomorphism is $\MCGo_{g,1}$-equivariant follows from the fact that the $\MCGo_{g,1}$-actions on $H^{\BM}_{*}(F_{k}(\Sigma_{g,1}');\bZ)$ and $H_*(\Sigma_{g,1}^{\times k} , \Delta \cup A_{g};\bZ)$ are both induced from the diagonal action of $\MCGo_{g,1}$ on $\Sigma_{g,1}^{\times k}$.
\end{proof}

On the other hand, the other representations encoded by the homological representation functors \eqref{eq:output_of_general_construction_mapping_class_motion}, as well as their variants $L^{\BM}_{i}(\cF_{(\underline{k},\Sym_{\lambda},\ell)}(\MCGo))$ and $L^{\BM}_{i}(\cF_{(\underline{k},\Sym_{\lambda},\ell)}(\MCGno))$, do not appear in the literature and thus seem to be new.

\paragraph*{Homological representations from the construction of \S\ref{ss:homological_representation_functor_mcg}.}

As pointed out in Remark~\ref{rmk:no-untwisted-versions}, the construction of \S\ref{ss:homological_representation_functor_motion_groups} typically does not lead to untwisted representations when applied to (full) mapping class groups, since Proposition~\ref{prop:factoringthroughtopQ} requires the groupoid $\cM$ to be a motion groupoid (its automorphism groups must be braided mapping class groups). The purpose of the alternative construction in \S\ref{ss:homological_representation_functor_mcg} is that it solves this problem: the analogue of Proposition~\ref{prop:factoringthroughtopQ} in this setting is Proposition~\ref{prop:factoringthroughcovQ-mcg}, which does not assume that $\cM$ is a motion groupoid. In this section we illustrate this alternative construction with some examples in the context of mapping class groups of surfaces.

Let us consider the continuous semifunctors $\Int{\fF}_{(\underline{k},\Sym_{\lambda},\ell)}$ and $\Int{\fF}^{\unt}_{(\underline{k},\Sym_{\lambda},\ell)}$ of Theorem~\ref{thm:global_functor_mcg}, where $Z:=\underline{k}$ is a finite set of $k\geq 1$ points, $\sG:=\Sym_{\lambda}$ and $Q:=Q_{\LCS_{\ell}}$ for some $\ell\geq 1$.
We denote the restrictions of these semifunctors to the subcategories $\fU\cM_{2}^{+,\mt}$ and $\fU\cM_{2}^{-,\mt}$ of $\fU\cD_{2}$ (see \S\ref{sss:category_mcg}) by $\Int{\fF}_{(\underline{k},\Sym_{\lambda},\ell)}(\MCGo)$, $\Int{\fF}^{\unt}_{(\underline{k},\Sym_{\lambda},\ell)}(\MCGo)$, $\Int{\fF}_{(\underline{k},\Sym_{\lambda},\ell)}(\MCGno)$ and $\Int{\fF}^{\unt}_{(\underline{k},\Sym_{\lambda},\ell)}(\MCGno)$ respectively.
Recall from \S\ref{ss:mcg_construction} that the posets $\cP(\cM_{2}^{+},\cM_{2}^{+})$ and $\cP(\cM_{2}^{-},\cM_{2}^{-})$ are both directed sets.
Using the appropriate colimit coefficient system $V_{\col}(\sF)$ of \S\ref{sss:col_coeff_unwtisted_rep_mcg} in each case, we therefore obtain from Corollary~\ref{coro:def_homological_rep_functor_mcg}, Proposition~\ref{prop:factoringthroughtopQtw} and Proposition~\ref{prop:factoringthroughcovQ-mcg} the following homological representation functors for all $i\geq 0$:
\begin{equation}
\label{eq:output_of_general_construction_mcg_mcg}
\begin{split}
L_{i}(\Int{\fF}_{(\underline{k},\Sym_{\lambda},\ell)}(\MCGo)) \colon \fU\cM_{2}^{+} &\too \modr[\bZ[\cQ^{\fF}_{(\lambda,\ell)}(\MCGo)]]^{\tw}; \\
L_{i}(\Int{\fF}^{\unt}_{(\underline{k},\Sym_{\lambda},\ell)}(\MCGo)) \colon \fU\cM_{2}^{+} &\too \modr[\bZ[\cQ^{\unt,\fF}_{(\lambda,\ell)}(\MCGo)]]; \\
L_{i}(\Int{\fF}_{(\underline{k},\Sym_{\lambda},\ell)}(\MCGno)) \colon \fU\cM_{2}^{-} &\too \modr[\bZ[\cQ^{\fF}_{(\lambda,\ell)}(\MCGno)]]^{\tw}; \\
L_{i}(\Int{\fF}^{\unt}_{(\underline{k},\Sym_{\lambda},\ell)}(\MCGno)) \colon \fU\cM_{2}^{-} &\too \modr[\bZ[\cQ^{\unt,\fF}_{(\lambda,\ell)}(\MCGno)]].
\end{split}
\end{equation}

In contrast with the previous paragraph (constructing homological representations from the construction of \S\ref{ss:homological_representation_functor_motion_groups}), the representations encoded by the functors \eqref{eq:output_of_general_construction_mcg_mcg} for $\ell\geq 2$ are more novel. For instance, for $\ell = 2$, we have $\cQ^{\fF}_{(\lambda,2)}(\MCGo)\cong (\bZ/2)^{r'}$ and $\cQ^{\fF}_{(\lambda,2)}(\MCGno)\cong (\bZ/2)^{r'}\times(\bZ/2)^{r}$; see Corollary~\ref{coro:transformation_group_MCG}. As far as the authors know, there are no representations of the mapping class groups $\MCGo_{g,1}$ and $\MCGno_{h,1}$ in the literature whose ground rings are of the above form. Therefore, many new representations of the mapping class groups arise from the functors \eqref{eq:output_of_general_construction_mcg_mcg}.

\section{Appendix: computations of transformation groups}\label{s:Appendix}

This appendix aims to study some of the transformation groups (i.e.~the groups whose group rings are the ground rings) of the homological representation functors that we construct in \S\ref{s:applications}, in order to offer a more concrete understanding of these representations.
We first deal with some necessary recollections of presentations of surface braid groups in \S\ref{ss:presentation_surface_braids}. Then we compute and present some properties of the transformation groups of the homological representation functors of \S\ref{s:applications} in \S\ref{ss:appendix_computations}.

\subsection{Presentations of surface braid groups}\label{ss:presentation_surface_braids}

Presentations of braid groups on surfaces with one boundary component may be found in \cite[\S 4]{LambropoulouOldenburg} and in \cite[Th.~1.1 and A.2]{Bellingeripresentations}; see also \cite[\S 6.3]{DPS}.
We fix three integers $k\geq 0$, $g\geq 0$ and $h\geq 1$. We shall use the same notation $\Sigma_{g,1}^{k}$ and $\N_{h,1}^{k}$ for surfaces as in \S\ref{sss:category_mcg}.
In the following presentations, we write $x\ \rightleftarrows\ y$ to denote the relation saying that $x$ and $y$ commute.

\begin{prop}\label{prop:presentation_braid_orientable_surfaces}
The braid group on $n$ strands on the orientable surface $\Sigma_{g,1}^{k}$,
denoted by $\mathbf{B}_{n}(\Sigma_{g,1}^{k})$, admits the presentation with generators $\mathscr{S}=\{\sigma_{i}\} _{1\leq i\leq n-1 }$, $A=\{a_{i}\} _{1\leq i\leq g }$, $B=\{b_{i}\} _{1\leq i\leq g }$ and $X=\{\xi_{i}\} _{1\leq i\leq k }$ and relations given by the braid relations for the elements of $\mathscr{S}$, to which are added the following families of relations (where $x$ and $y$ denote either $a$ or $b$, and $1 \leq r, s \leq g$):
\begin{equation}\label{eq:relations_braids_on_orientable_surface}
\begin{cases}
(BS1)\ \textrm{ }\sigma_i\ \rightleftarrows\ x_r & \textrm{for all $r$ and all $1\leq i\leq n-2 $}\ ;\\
(BS2)\ \textrm{ } x_r\ \rightleftarrows\ \sigma_{n-1} y_s \sigma_{n-1}^{-1} & \textrm{for $s < r$}\ ;\\
(BS3)\ \textrm{ }(\sigma_{n-1} x_r)^2 = (x_r \sigma_{n-1})^2 & \textrm{for all $r$}\ ;\\
(BS4)\ \textrm{ }[\sigma_{n-1} b_r \sigma_{n-1}^{-1}, a_r^{-1}] = \sigma_{n-1}^2 & \textrm{for all $r$}\ ;\\
(BS5)\ \textrm{ }\xi_j\ \rightleftarrows\ \sigma_i & \textrm{for all $1\leq j \leq k$ and all $1\leq i\leq n-2 $}\ ;\\
(BS6)\ \textrm{ }x_r\ \rightleftarrows\ \sigma_{n-1} \xi_j \sigma_{n-1}^{-1} & \textrm{for all $1\leq j \leq k$ and all $1\leq r \leq g$}\ ;\\
(BS7)\ \textrm{ }\xi_i\ \rightleftarrows\ \sigma_{n-1} \xi_j \sigma_{n-1}^{-1} & \textrm{for $i < j$}.
\end{cases}
\end{equation}
The braid group on $n$ strands on the non-orientable surface $\N_{h,1}^{k}$,
denoted by $\mathbf{B}_{n}(\N_{h,1}^{k})$, admits
the presentation with generators $\mathscr{S}=\{\sigma_{i}\} _{1\leq i\leq n-1 }$, $C=\{c_{i}\} _{1\leq i\leq h }$ and $X=\{\xi_{i}\} _{1\leq i\leq k }$ and relations given by the braid relations for the elements of $\mathscr{S}$, to which are added the following families of relations (where $1 \leq r, s \leq h$):
\begin{equation}\label{eq:relations_braids_on_non_orientable_surface}
\begin{cases}
(BN1)\ \textrm{ }\sigma_i\ \rightleftarrows\ c_r & \textrm{for all $r$ and all $1\leq i\leq n-2 $}\ ;\\
(BN2)\ \textrm{ } c_r\ \rightleftarrows\ \sigma_{n-1} c_s \sigma_{n-1}^{-1} & \textrm{for $s < r$}\ ;\\
(BN3)\ \textrm{ }[\sigma_{n-1} c_r \sigma_{n-1}^{-1}, c_r^{-1}] = \sigma_{n-1}^2 & \textrm{for all $r$}\ ;\\
(BN4)\ \textrm{ }\xi_j\ \rightleftarrows\ \sigma_i & \textrm{for all $1\leq j \leq k$ and all $1\leq i\leq n-2 $}\ ;\\
(BN5)\ \textrm{ }c_r\ \rightleftarrows\ \sigma_{n-1} \xi_j \sigma_{n-1}^{-1} & \textrm{for all $1\leq j \leq k$ and all $1\leq r\leq h$}\ ;\\
(BN6)\ \textrm{ }\xi_i\ \rightleftarrows\ \sigma_{n-1} \xi_j \sigma_{n-1}^{-1} & \textrm{for $i < j$}\ ;\\
(BN7)\ \textrm{ }(\sigma_{n-1} \xi_j)^2 = (\xi_j \sigma_{n-1})^2 & \textrm{for all $1\leq j \leq k$}.
\end{cases}
\end{equation}
\end{prop}

Now we consider a partition $\lambda = \boldsymbol{\{} \lambda_{1};\ldots;\lambda_{r} \boldsymbol{\}} \vdash k$. There is an isomorphism $\B_{\lambda}(S) \cong \B_{\boldsymbol{\{} \lambda_{1};\ldots;\lambda_{r-1} \boldsymbol{\}}}\bigl(\bD^{2}_{\lambda_{r}}\natural S\bigr) \rtimes \B_{\lambda_{r}}(S)$, which may be deduced from the split short exact sequence \eqref{eq:split-ses-1}.
There is a classical method of constructing a presentation of a group extension from a presentation of the quotient and a presentation of the kernel; see \cite[\S 2.4.3]{HBE} and \cite[Appendix~B]{DPS}.
For instance, a presentation of the group $\B_{k,n}(\Sigma_{g,1})$ is detailed in \cite[Prop.~3.2]{BellingeriGodelleGuaschi} following this method, while a presentation of $\B_{k,n}(\N_{h,1})$ is given in \cite[Prop.~6.58]{DPS}.
It is routine to generalise this to give full presentations for any partition $\lambda$. We thus obtain from Proposition~\ref{prop:presentation_braid_orientable_surfaces} the following result for the partitioned surface braid groups.

\begin{prop}\label{prop:partitioned_braid_groups_generating_sets}
Let $\lambda = \boldsymbol{\{} \lambda_{1};\ldots;\lambda_{r} \boldsymbol{\}}$ be a partition of $k\geq 1$. 
The surface braid group $\mathbf{B}_{\lambda}(S)$ admits a presentation whose generating sets are:
\begin{itemizeb}
    \item $X^{(\rho)}=\{ \xi^{(\rho)}_{i}\mid 1\leq i\leq \varSigma_{\rho}\}$ with $\varSigma_{\rho}:=\sum_{\rho+1\leq l\leq r}\lambda_{l}$, for each block $1\leq \rho\leq r-1$;
    \item $\mathscr{S}^{(\rho')}=\{ \sigma_{i}^{(\rho')}\} _{1\leq i\leq \lambda_{\rho'}-1}$ for each block $1\leq \rho'\leq r$ such that $r_{\rho'}\geq 2$;
    \item if $S=\Sigma_{g,1}$\textup{:} $A^{(\rho)}=\{ a^{(\rho)}_{i}\} _{1\leq i\leq g }$ and $B^{(\rho)}=\{ b^{(\rho)}_{i}\} _{1\leq i\leq g }$ for each block $1\leq \rho\leq r$;
    \item if $S=\N_{h,1}$\textup{:} $C^{(\rho)}=\{ c^{(\rho)}_{i}\} _{1\leq i\leq h }$ for each block $1\leq \rho\leq r$.
\end{itemizeb}
The relations between generators of the same blocks are those of \eqref{eq:relations_braids_on_orientable_surface} and \eqref{eq:relations_braids_on_non_orientable_surface}, while the relations between generators of different blocks are analogous to those of \textup{(c.1)--(c.8)} in \textup{\cite[Prop.~3.2]{BellingeriGodelleGuaschi}}.
\end{prop}

\subsection{Properties of the transformation groups}\label{ss:appendix_computations}

This section deals with some properties of the homological representation functors of \S\ref{s:applications}.
Throughout \S\ref{ss:appendix_computations}, we consider an integer $\ell\geq 1$ corresponding to a lower central series index, and an integer $k\geq 1$ and a partition $\lambda = \boldsymbol{\{} \lambda_{1};\ldots;\lambda_{r} \boldsymbol{\}} \vdash k$.
We also denote by $r'$ the number of indices $1\leq i \leq r$ such that $\lambda_i \geq 2$. For simplicity, we denote the decorated surface $(\bD^{2},\underline{n})$ for each non-negative integer $n$ by $\bD_{n}$ and call it the \emph{$n$-th marked $2$-disc}, and we denote the decorated surface $(S,\underline{n})$ by $S^{(n)}$. If $n=0$, we abbreviate $(S,\underline{0}) = S^{(0)}$ to $S$.

\subsubsection{Surface braid groups}\label{ss:surface_braid_transformation_groups}

We continue to follow the notation of \S\ref{ss:presentation_surface_braids} and consider a compact, connected, smooth surface $S$ with one boundary component.
We deal here with the computation of the transformation groups of the homological representation functors of the form $L_{i}(\Int{\cF}_{(\underline{k},\Sym_{\lambda},\ell)}(S))$ for surface braid groups defined in \S\ref{ss:applications_motion_groups}. We denote by $\cQ_{(\lambda,\ell)}(S^{(n)})$ the group of the form $\cQ(S,\underline{n})$ induced by diagram \eqref{eq:split-short-exact-sequence-quotient} for each integer $n \geq 0$ and by $\cQ_{(\lambda,\ell)}(S)$ the colimit of these groups as $n\to\infty$, which is the colimit transformation group $\cQ_{\col}(T)$ of Notation~\ref{notation-Qcol} in this setting.

\paragraph*{Abelian quotients.}

We start with the homological representation functor defined using the $\LCS_{2}$ term of the lower central series. The abelianisations of the groups $\B_{\lambda,n}(S)$ are explicitly computed in \cite[Props.~3.5 and 6.47]{DPS}, and may be described via the corresponding generating set of Proposition~\ref{prop:partitioned_braid_groups_generating_sets}. In particular, a generating set for $\B_{\lambda,n}(S)^{\ab}$ is given by:
\begin{itemizeb}
    \item the common image $t_{i'}$ of the generators of $\mathscr{S}^{(i')}$ for each $1\leq i'\leq r$ such that $r_{i'}\geq 2$;
    \item if $S=\bD$: the common image $q_{i}$ of all $\xi_{j}^{(i)}\in X^{(i)}$ with $j\geq 1+\sum_{i+1\leq l\leq r}\lambda_{l}$ for each $1\leq i\leq r$, and the common image $s_{i_{1},i_{2}}$ of all $\xi_{j}^{(i_{1})}$ with $j\in \{j'+\sum_{i_{1}+1\leq l\leq i_{2}-1}\lambda_{l}\mid 1\leq j'\leq i_{2} \}$ for each pair $1\leq i_{1}< i_{2}\leq r$;
    \item if $S=\Sigma_{g,1}$: the images $\{ A_{j}^{(i)},B_{j}^{(i)}\} _{1\leq j\leq g }$ of the sets $A^{(\rho)}$ and $B^{(\rho)}$ for each $1\leq i\leq r$;
    \item if $S=\N_{h,1}$: the images $\{ C_{j}^{(i)}\} _{1\leq j\leq h }$ of the set $C^{(\rho)}$ for each $1\leq i\leq r$.
\end{itemizeb}

\begin{lemm}\label{lem:transformation_groups_ab_quotient_surface_braid_groups}
We have $\cQ_{(\lambda,2)}(\bD) \cong \bZ^{r'} \times \bZ^{r(r-1)/2} \times \bZ^{r}$ and $\cQ_{(\lambda,2)}(S) \cong (\bZ/2)^{r'} \times H_{1}(S;\bZ)^{r}$ if $S \not\cong \bD$.
\end{lemm}
\begin{proof}
We compute the group $\cQ_{(\lambda,2)}(S^{(n)})$ for $n\geq 3$ as the kernel of the map $\B_{\lambda,n}(S)^{\ab}\twoheadrightarrow\B_{n}(S)^{\ab}$ given by diagram \eqref{eq:split-short-exact-sequence-quotient}, the abelianisations $\B_{\lambda,n}(S)^{\ab}$ and $\B_{n}(S)^{\ab}$ being calculated in \cite[Props.~3.5 and 6.47]{DPS}. We deduce from these computations that $\cQ_{(\lambda,2)}(S^{(n)}) \cong \cQ_{(\lambda,2)}(S^{(n+1)})$ for all $n\geq 3$ and a fortiori that the colimit $\cQ_{(\lambda,2)}(S)$ is isomorphic to $\cQ_{(\lambda,2)}(S^{(3)})$, thus giving the result.
\end{proof}

\paragraph*{Further $\LCS_{\ell}$-quotients.}

We now consider more generally the surface braid group homological representation $L_{i}(\Int{\cF}_{(\underline{k},\Sym_{\lambda},\ell)}(S))$ defined using any parameter $\ell\geq 2$.

\begin{prop}\label{prop:Q-stab-general}
For any integers $\ell\geq 2$ and $n\geq 4$, we have $\cQ_{(\lambda,\ell)}(S^{(n)})\cong\cQ_{(\lambda,\ell)}(S^{(n+1)})$.
\end{prop}
\begin{proof}
We recall that $\cQ_{(\lambda,\ell)}(S^{(n)})$ is defined as the kernel of the surjection of $\B_{\lambda,n}(S)/\LCS_{\ell}$ onto $\B_{n}(S)/\LCS_{\ell}$. The data which depends on $n$ in the presentation of $\B_{\lambda,n}(S)$ are the set of braid generators $\mathscr{S}^{(n)}$ of the $n$-th block, and, for each block $1\leq \rho\leq r$, the subset of $X^{(\rho)}$ of the pure braid generators $\{ \chi_{i}^{(\rho)}:=\xi^{(\rho)}_{\varSigma_{\rho}+i}\mid 1\leq i\leq n\}$ where $\varSigma_{\rho}$ denotes the sum $\sum_{\rho+1\leq l\leq r}\lambda_{l}$.
For a group $G$, we generically denote by $\gamma_{\ell}$ the projection onto the $\ell$-nilpotent quotient $G/\LCS_{\ell}$.

Since the assignment $S \mapsto \B_{n}(S)$ is functorial with respect to embeddings of surfaces, we have a canonical injection $\B_{n}\hookrightarrow\B_{n}(S)\hookrightarrow \B_{\lambda,n}(S)$. In particular, this morphism sends $\LCS_\infty(\B_{n})$ to $\LCS_\infty(\B_{\lambda,n}(S))$. Since $\LCS_{\infty}(\B_{n}) = \LCS_{2}(\B_{n})$ (see for instance \cite[Ex.~2.3]{DPS}), we know that $\sigma_{i} \sigma_{j}^{-1} \in \LCS_\infty(\B_{n})$ for all $1\leq i,j \leq n-1$. A fortiori, we deduce that $\sigma^{(n)}_i \equiv \sigma^{(n)}_j \pmod{\LCS_\infty(\B_{n}(S))}$ and we denote by $\sigma^{(n)} \in \B_{n}(S)/\LCS_{\ell}$ the common image of all the $\sigma^{(n)}_i$ under $\gamma_{\ell}$.

Furthermore, for each $1\leq \rho\leq r$, we have the relations $\sigma^{(n)}_{i}\chi^{(\rho)}_{i}(\sigma^{(n)}_{i})^{-1}=(\chi^{(\rho)}_{i})^{-1}\chi^{(\rho)}_{i+1}\chi^{(\rho)}_{i}$, $\sigma^{(n)}_{i}\chi^{(\rho)}_{i+1}(\sigma^{(n)}_{i})^{-1}=\chi^{(\rho)}_{i}$ and $\sigma^{(n)}_{i}\chi^{(\rho)}_{j}(\sigma^{(n)}_{i})^{-1}=\chi^{(\rho)}_{j}$ if $j\notin \{i,i+1\}$ by Proposition~\ref{prop:partitioned_braid_groups_generating_sets}. These are the typical relations between pure braids and Artin generators of a given block induced by the injection $\B_{\lambda,n}\hookrightarrow\B_{\lambda,n}(S)$; see \cite[Prop.~3.2]{BellingeriGodelleGuaschi} for the case of $S=\Sigma_{g,1}$ and $\lambda = \boldsymbol{\{} k \boldsymbol{\}}$. Then, we deduce from these relations that:
\begin{itemizeb}
    \item $\gamma_{\ell}(\chi^{(\rho)}_{i+1})
 = \gamma_{\ell}((\sigma^{(n)}_{i})^{-1}\chi^{(\rho)}_{i}\sigma^{(n)}_{i})
 = \gamma_{\ell}((\sigma^{(n)}_{i-1})^{-1}\chi^{(\rho)}_{i}\sigma^{(n)}_{i-1})
 = \gamma_{\ell}(\chi^{(\rho)}_{i})$ for all $2\leq i\leq n-1$;
    \item $\gamma_{\ell}(\chi^{(\rho)}_{2})
 = \gamma_{\ell}((\sigma^{(n)}_{1})^{-1}\chi^{(\rho)}_{1}\sigma^{(n)}_{1})
 = \gamma_{\ell}((\sigma^{(n)}_{3})^{-1}\chi^{(\rho)}_{1}\sigma^{(n)}_{3})
 = \gamma_{\ell}(\chi^{(\rho)}_{1})$.
\end{itemizeb}
We denote by $\chi^{(\rho)} \in \B_{\lambda,n}(S)/\LCS_{\ell}$ the common image of all the $\chi^{(\rho)}_i$ under $\gamma_{\ell}$.

Therefore, the presentation of $\B_{\lambda,n}(S)/\LCS_{\ell}$ is independent of $n$. In particular, it is routine to check that there is a well-defined map $\gamma'_{\ell}\colon\B_{\lambda,n+1}(S)\to \B_{\lambda,n}(S)/\LCS_{\ell}$ defined by $\sigma_{i}^{(n+1)}\mapsto \sigma^{(n)}$, $\chi^{(\rho)}_j\mapsto \chi^{(\rho)}$ (with $1\leq j\leq n+1$) and the assignment of $\gamma_{\ell}$ for the other generators. Then $\gamma'_{\ell}$ induces an inverse to the canonical map $\B_{\lambda,n}(S)/\LCS_{\ell}\to \B_{\lambda,n+1}(S)/\LCS_{\ell}$, which is thus an isomorphism.
We also know the analogous result for $\B_{n}(S)/\LCS_{\ell}$ by \cite[Prop.~3.13]{BellingeriGodelleGuaschi} (see also \cite[Prop.~6.43]{DPS}), whence the result.
\end{proof}

It follows from Proposition~\ref{prop:Q-stab-general} that $\cQ_{(\lambda,\ell)}(S) \cong \cQ_{(\lambda,\ell)}(S^{(n)})$ for each $\ell\geq 2$ and $n\geq 4$. Thus, for each $\ell \geq 2$, we may compute the colimit transformation group from a computation of $\cQ_{(\lambda,\ell)}(S^{(n)})$ for any $n\geq 4$. For $\ell = 2$, this was done in Lemma~\ref{lem:transformation_groups_ab_quotient_surface_braid_groups} above; the following result gathers further results for $\ell \geq 2$.

\begin{prop}\label{prop:transformation_groups_{2}_nilpotent_quotient_surface_braid_groups}
We assume that the partition is such that $\lambda_{l}\geq 3$ for all $1\leq l\leq r$.

For the classical braid groups, we have for all $\ell\geq 2$:
\[
\cQ_{(\boldsymbol{\{}2;\lambda\boldsymbol{\}},\ell)}(\bD) \cong \bZ^{\binom{r+2}{2} - 1} \times ((\bZ^2/2^{\ell-2}\bar{\Delta})^{r+1} \rtimes \bZ),
\]
where $\bar{\Delta} = (1,-1) \in \bZ^2$ and $1 \in \bZ$ acts on each copy of $\bZ^2/2^{\ell-2}\bar{\Delta}$ by swapping coordinates.

For the surfaces different from the disc, we have the following computations for $\ell = 3$. For orientable surfaces, for all $g\geq 1$:
\begin{equation}\label{eq:transformation_group_orientable_partitioned}
\cQ_{(\lambda,3)}(\Sigma_{g,1})\cong ((\bZ^{r(r-1)/2}\times\bZ^{r}\times\bZ^{rg})\rtimes\bZ^{rg})\times \bZ^{r}.
\end{equation}
In more detail, the right-hand side of \eqref{eq:transformation_group_orientable_partitioned} may be written as
\begin{multline*}
\left(\left(\langle s_{r_{1},r_{2}}\rangle_{1\leq r_{1}< r_{2}\leq r} \times \langle t_{1},\ldots,t_{r}\rangle \times\prod_{1\leq \rho\leq r}\langle A_{1}^{(\rho)},\ldots,A_{g}^{(\rho)}\rangle \right)\right. \\
\left. \rtimes \prod_{1\leq \rho\leq r}\langle B_{1}^{(\rho)},\ldots,B_{g}^{(\rho)}\rangle \right) \times \langle q_{1},\ldots,q_{r}\rangle
\end{multline*}
where the action defining the semi-direct product structure is determined by
\begin{itemizeb}
\item $[A^{(\rho)}_{i},B^{(\rho)}_{i}] = t_{\rho}^{2}$ for all $1\leq \rho\leq r$;
\item $[A_{i}^{(r_{1})},B_{i}^{(r_{2})}]= [A_{i}^{(r_{2})},B_{i}^{(r_{1})}]= s_{r_{1},r_{2}}$ for all $1\leq r_{1}< r_{2}\leq r$;
\item all other pairs of generators commute.
\end{itemizeb}
We deduce that $\cQ_{(\lambda,3)}^{\unt}(\Sigma_{g,1}) \cong \cQ_{(\lambda,3)}(\Sigma_{g,1})/\langle q_{1},\ldots,q_{r}\rangle$.
For non-orientable surfaces, for all $h\geq 1$:
\begin{multline}
\label{eq:transformation_group_non-orientable_partitioned}
\cQ_{(\lambda,3)}(\N_{h,1}) \cong \\ \left( (\bZ^{r-1}\times(\bZ^{r-2}\times\cdots\times (\bZ^{2}\times(\bZ\times\bZ^{h})\rtimes \bZ^{h}) \rtimes\cdots\rtimes\bZ^{h})\rtimes\bZ^{h})\rtimes\bZ^{h} \right) \times(\bZ/2)^{r}\times\bZ^{r}.
\end{multline}
In more detail, the right-hand side of \eqref{eq:transformation_group_non-orientable_partitioned} may be written as
\begin{multline*}
\left(\left(\langle s_{1,r_{2}}\rangle_{r_{2}\geq 2} \times\cdots\times \left( \langle s_{r-1,r}\rangle\times\langle C_{1}^{(r)},\ldots,C_{h}^{(r)}\rangle \right) \rtimes \cdots \rtimes \langle C_{1}^{(2)},\ldots,C_{h}^{(2)}\rangle \right) \right. \\
\left. \rtimes \langle C_{1}^{(1)},\ldots,C_{h}^{(1)}\rangle \right) \times \langle t_{1},\ldots,t_{r}\rangle\times \langle q_{1},\ldots,q_{r}\rangle
\end{multline*}
where the action defining the semi-direct product structure is determined by
\begin{itemizeb}
\item $[C_{i}^{(r_{1})},C_{i}^{(r_{2})}]= [C_{i}^{(r_{2})},C_{i}^{(r_{1})}]= s_{r_{1},r_{2}}$ for all $1\leq r_{1}< r_{2}\leq r$;
\item all other pairs of generators commute.
\end{itemizeb}
We deduce that $\cQ_{(\lambda,3)}^{\unt}(\N_{h,1})\cong\cQ_{(\lambda,3)}(\N_{h,1})/\langle q_{1},\ldots,q_{r}\rangle$.
\end{prop}

\begin{proof}
That $\cQ_{(\boldsymbol{\{}2;\lambda\boldsymbol{\}},\ell)}(\bD)=\cQ^{\unt}_{(\boldsymbol{\{}2;\lambda\boldsymbol{\}},\ell)}(\bD)$ and its explicit computation for each $\ell\geq 3$ is done in \cite[\S 4, Cor.~5.4]{PSIN}.

For the surfaces different from the disc, we first recall that the quotient $\B_{n}(\Sigma_{g,1})/\LCS_{3}$ is computed by \cite[Prop.~3.13]{BellingeriGodelleGuaschi}, while $\B_{n}(\N_{h,1})/\LCS_{3}=\B_{n}(\N_{h,1})^{\ab}$ by \cite[Th.~6.42]{DPS}. Now the computation of $\B_{\lambda,n}(S)/\LCS_{3}$ follows the same steps as the proof of \cite[Prop.~6.58]{DPS}, which computes $\B_{k,n}(\N_{h,1})/\LCS_{3}$ and generalises mutatis mutandis as follows.
\begin{itemize}
    \item Using the presentation from Proposition~\ref{prop:partitioned_braid_groups_generating_sets}, let $N$ be the normal closure of the $\sigma^{(\rho)}_{i}(\sigma^{(\rho)}_{i+1})^{-1}$ for $i<\lambda_{\rho}$ together with the $\xi^{(\rho)}_j (\xi^{(\rho)}_{j+1})^{-1}$ for $j < \varSigma_{\rho}$, for each $1\leq \rho\leq n$. Therefore, $N \subseteq \LCS_{3}(\B_{\lambda,n}(S))$ because its generators are in $\LCS_3(\B_{\lambda,n}(S))$; see for instance \cite[Lem.~6.49]{DPS} and its proof. We claim that, in fact, $N = \LCS_{3}(\B_{\lambda,n}(S))$.

    \item Consider the partition $\boldsymbol{\{} \lambda,n\boldsymbol{\}}$ as a partition $\lambda'$ of $k+n$ with length $r+1$ and we make the identification $s_{i,n}:=q_{i}$ for each $1\leq i\leq r$. It is routine, although lengthy (and an inductive generalisation of the analogous point in the proof of \cite[Prop.~6.58]{DPS}), to check from the presentation of Proposition~\ref{prop:partitioned_braid_groups_generating_sets} that the quotient $\B_{\lambda,n}(S)/N$ is
        \begin{multline}\label{eq:B_k,n(S)/LCS_3_orientable}
            \left(\langle s_{r_{1},r_{2}}\rangle_{1\leq r_{1}< r_{2}\leq r+1} \times \langle t_{1},\ldots,t_{r+1}\rangle \times\prod_{1\leq \rho\leq r+1}\langle A_{1}^{(\rho)},\ldots,A_{g}^{(\rho)}\rangle\right) \\
            \rtimes \prod_{1\leq \rho\leq r+1}\langle B_{1}^{(\rho)},\ldots,B_{g}^{(\rho)}\rangle
        \end{multline}
if $S=\Sigma_{g,1}$, and
        \begin{multline}\label{eq:B_k,n(S)/LCS_3_non_orientable}
            \left(\langle s_{1,r_{2}}\rangle_{2\leq r_{2}\leq r+1}\times\cdots\times(\langle s_{r,r+1}\rangle\times\langle C_{1}^{(r+1)},\ldots,C_{h}^{(r+1)}\rangle)\rtimes \cdots \right. \\
            \left. \cdots \rtimes \langle C_{1}^{(1)},\ldots,C_{h}^{(1)}\rangle \right) \times \langle t_{1},\ldots,t_{r+1}\rangle
        \end{multline}
if $S=\N_{h,1}$.

    \item The proof of the claim in the first point thus follows from the observation that, using the second point, $\B_{\lambda,n}(S)/N$ is a $2$-nilpotent group. Indeed, in both cases, the commutator subgroup is generated by the elements $s_{r_{1},r_{2}}$ and $t_{i}^{2}$ (if $S=\Sigma_{g,1}$ for the latter): all these generators are also clearly central in $\B_{\lambda,n}(S)/N$, which proves our claim.
\end{itemize}
The computations of $\cQ_{(\lambda,3)}(S)$ then directly follow from the above descriptions of $\B_{\lambda,n}(S)/\LCS_{3}$ and $\B_{n}(S)/\LCS_{3}$.

Finally, we compute the untwisted quotients $\cQ^{\unt}_{(\lambda,3)}(\Sigma_{g,1})$ and $\cQ^{\unt}_{(\lambda,3)}(\N_{h,1})$ as follows. First we know from the presentation of $\B_{\lambda,n}(S)/\LCS_{3}$ (see \eqref{eq:B_k,n(S)/LCS_3_orientable} and \eqref{eq:B_k,n(S)/LCS_3_non_orientable}) that the only generators of $\cQ_{(\lambda,3)}(S)$ on which the action of $\B_{n}(S)$ (given by conjugation) is not trivial are $A_{i}^{(\rho)}$ and $B_{i}^{(\rho)}$ for all $i$ and each $1\leq\rho\leq r$ if $S=\Sigma_{g,1}$, or the $C_{j}^{(\rho)}$ for all $j$ and each $1\leq\rho\leq r$ if $S=\N_{h,1}$. Then, we deduce from the presentation of Proposition~\ref{prop:partitioned_braid_groups_generating_sets} that for all $1\leq \rho\leq r$:
\[
q_{\rho}B_{i}^{(\rho)}= A_{i}^{(r+1)} B_{i}^{(\rho)} (A_{i}^{(r+1)})^{-1},
\]
and the analogous relation swapping $A$ and $B$, for all $1\leq i\leq g$, and
\[
q_{\rho}C_{j}^{(\rho)}= C_{j}^{(r+1)} C_{j}^{(\rho)} (C_{j}^{(r+1)})^{-1}
\]
for all $1\leq j\leq h$. This proves that the quotienting submodule defining the coinvariants is $\langle q_{1},\ldots,q_{r}\rangle$ in each case.
\end{proof}

\begin{rmk}\label{rmk:limitation_untwisted_transformation_group_surface_braid}
If $\lambda_{i} \leq 2$ for some $1\leq i\leq r$ and $\ell\geq 3$, we do not know whether the quotient $\cQ^{\unt}_{(\lambda,\ell)}(S)$ of $\cQ_{(\lambda,\ell)}(S)$ is a \emph{proper} quotient.
\end{rmk}

\subsubsection{Loop braid groups}\label{ss:loop_braid_transformation_groups}

We now compute the transformation groups of the homological representation functors defined in \S\ref{ss:applications_motion_groups} for loop braid groups with the parameter $\ell = 2$.

\begin{lemm}\label{lem:transformation_groups_ab_quotient_loop_braid_groups}
We have the following descriptions of the colimit transformation groups:
\begin{itemizeb}
    \item for \eqref{eq:output_of_general_construction_loop_braids_points}\textup{:} $\cQ_{(P,\lambda,2)}(\bD^{3})\cong \bZ^{r}\times (\bZ/2)^{r'}$ and $\cQ'_{(P,\lambda,2)}(\bD^{3})\cong (\bZ/2)^{r+r'}$;
    \item for \eqref{eq:output_of_general_construction_loop_braids_unlinks}\textup{:} $\cQ_{(\bS_{+},\lambda,2)}(\bD^{3})\cong \bZ^{r^2+r+r'}\times (\bZ/2)^{r'}$ and $\cQ'_{(\bS_{+},\lambda,2)}(\bD^{3})\cong \bZ^{r^2+r'}\times (\bZ/2)^{r+r'}$;
    \item for \eqref{eq:output_of_general_construction_loop_braids_unlinks_not_preserved}\textup{:} $\cQ_{(\bS,\lambda,2)}(\bD^{3})\cong \bZ^{r}\times (\bZ/2)^{2r'+r(r+1)}$ and $\cQ'_{(\bS,\lambda,2)}(\bD^{3})\cong (\bZ/2)^{2r'+r(r+2)}$.
\end{itemizeb}
\end{lemm}
\begin{proof}
Recall the quotients of embedding spaces $C_{\lambda}(\bD^{3}_{n})$, $U_{\lambda}^{+}(\bD^{3}_{n})$ and $U_{\lambda}(\bD^{3}_{n})$ introduced in \eqref{eq:configurations_loop_braid_groups}; henceforth we denote by $C_{\lambda,n}$ any one of these spaces.
The abelianisations of all of the groups $\pi_{1}(C_{\lambda,n})\rtimes\lB_{n}$, $\lB_{n}$, $\pi_{1}(C_{\lambda,n})\rtimes\lB'_{n}$ and $\lB'_{n}$ are computed in \cite[Props.~4.46 and 5.10]{DPS}.
In particular, via these computations, one sees that $(\pi_{1}(C_{\lambda,n})\rtimes\lB_{n})^{\ab}\cong (\pi_{1}(C_{\lambda,n+1})\rtimes\lB_{n+1})^{\ab}$, $\lB_{n}^{\ab}\cong \lB_{n+1}^{\ab}$, $(\pi_{1}(C_{\lambda,n})\rtimes\lB'_{n})^{\ab}\cong (\pi_{1}(C_{\lambda,n+1})\rtimes\lB'_{n+1})^{\ab}$ and $(\lB'_{n})^{\ab}\cong (\lB'_{n+1})^{\ab}$ for all $n\geq 4$. Therefore, the groups of the statement of the form $\cQ$ given by diagram \eqref{eq:split-short-exact-sequence-quotient} are computed as the kernels of the maps $(\pi_{1}(C_{\lambda,n})\rtimes\lB_{n})^{\ab}\twoheadrightarrow \lB_{n}^{\ab}$ and $(\pi_{1}(C_{\lambda,n})\rtimes\lB'_{n})^{\ab}\twoheadrightarrow(\lB'_{n})^{\ab}$ for some fixed $n\geq 4$.
\end{proof}

\subsubsection{Mapping class groups}

We finally study the transformation groups of the homological representation functors for mapping class group from \S\ref{ss:mcg_construction} with the parameter $\ell = 2$.
For the sake of completeness, we recall that the Lickorish generators together with the Dehn twist along a simple closed curve encircling the boundary component generate the mapping class group $\MCGo_{g,1}$ with $g\geq 1$ (see for instance \cite[\S 4.4]{farbmargalit}), and that generating sets for the mapping class group $\MCGno_{h,1}$ have been worked out by Stukow in \cite[Th.~A.7]{Stukow_Dehn} for $h=2$ and in \cite[Th.~5.2]{Stukow_Generating} for $h\geq 3$, while $\MCGno_{0,1}=\MCGno_{1,1}$ are trivial by \cite{Epstein}.
For brevity, we shall typically use the notation $\MCG(S,\lambda)$ for the mapping class group $\pi_{0}(\diffdec(S,\lambda))$, where $S$ is a compact, connected, smooth surface with one boundary component equipped with $k$ marked points that are fixed setwise, respecting the partition $\lambda\vdash k$.
We first establish the following general decomposition for abelianisations of mapping class groups:

\begin{prop}\label{prop:transformation_group_MCG}
For a compact, connected, smooth, non-planar surface $S$ with one boundary component, we have:
\begin{equation}\label{eq:computation_abel_MCG}
\MCG(\bD_{k}\natural S,\lambda)^{\ab}\cong (\bZ/2)^{r'} \times (H_{1}(S;\bZ)^{r})_{\MCG(S)}\times \MCG(S)^{\ab},
\end{equation}
where each of the first $r'$ $\bZ/2$-summands is generated by the image in the abelianisation $\sigma^{(\rho')}$ (with $1\leq \rho'\leq r$ such that $\lambda_{\rho'}\geq 2$) of a standard braid generator (considered as a mapping class) interchanging two points in the corresponding $\rho'$-th block of the partition.
\end{prop}
\begin{proof}
Considering the split short exact sequence \eqref{eq:split-ses-2} and using the computation of $\B_{\lambda}(S)^{\ab}$ (see for instance \cite[Prop.~6.47]{DPS}), it follows from the general formula for the calculation of the abelianisation of a semi-direct product that $\MCG(\bD_{k}\natural S,\lambda)^{\ab}\cong ((\bZ/2)^{r'} \times H_{1}(S;\bZ)^{r})_{\MCG(S)}\times \MCG(S)^{\ab}$.
We note that the splitting of \eqref{eq:split-ses-2} is induced by the embedding of surfaces $S \hookrightarrow\bD_{k}\natural S$. Therefore, each $\sigma^{(\rho')}$ may be represented as a mapping class supported in the subsurface $\bD_{k}$ on which $\MCG(S)$ acts trivially, whence the result.
\end{proof}

Using this, we deduce the following calculations of the transformation groups of the homological representation functors $L_{i}(\Int{\fF}_{(\underline{k},\Sym_{\lambda},\ell)}(\MCGo))$ and $L_{i}(\Int{\fF}_{(\underline{k},\Sym_{\lambda},\ell)}(\MCGno))$ of \eqref{eq:output_of_general_construction_mcg_mcg} for $\ell = 2$. Let us write $\mathfrak{G} \in \{\MCGo,\MCGno\}$ and correspondingly $S\in\{\Sigma_{1,1},\N_{1,1}\}$. Denote by $\cQ^{\fF}_{(\lambda,2,n)}(\mathfrak{G})$ the group of the form $\cQ(S^{\natural n},\emptyset)$ induced by diagram \eqref{eq:split-short-exact-sequence-quotient-mcg} for each integer $n \geq 0$ and by $\cQ^{\fF}_{(\lambda,2)}(\mathfrak{G})$ the colimit of these groups as $n\to\infty$, which is the colimit transformation group $\cQ_{\col}(T)$ of Notation~\ref{notation-Qcol} in the setting of \eqref{eq:output_of_general_construction_mcg_mcg}.

\begin{coro}\label{coro:transformation_group_MCG}
We have $\cQ^{\fF}_{(\lambda,2)}(\MCGo)\cong(\bZ/2)^{r'}$ and $\cQ^{\fF}_{(\lambda,2)}(\MCGno)\cong(\bZ/2)^{r'}\times(\bZ/2)^{r}$.
\end{coro}
\begin{proof}
Considering the discrete partition $\boldsymbol{\{} 1^{k}\boldsymbol{\}}$, we know from \cite[Th.~5.1]{Korkmaz_low} that the abelianisations of $\MCGo_{g,1}$ and of $\MCGo_{g,1}^{\boldsymbol{\{} 1^{k}\boldsymbol{\}}}$ are trivial for $g\geq 3$, while it follows from \cite[Th.~6.21]{Stukow_Generating} that $\bigl(\MCGno_{h,1}^{\boldsymbol{\{} 1^{k}\boldsymbol{\}}}\bigr)^{\ab}\cong (\bZ/2)^{r}\times (\MCGno_{h,1})^{\ab}$ for $h\geq 7$.
Using Proposition~\ref{prop:transformation_group_MCG} for $\lambda = \boldsymbol{\{} 1^{k}\boldsymbol{\}}$, we deduce that $(H_{1}(\Sigma_{g,1};\bZ)^{r})_{\MCGo_{g,1}}=0$ and that $(H_{1}(\N_{h,1};\bZ)^{r})_{\MCGno_{h,1}}\cong (\bZ/2)^{r}$.
Applying Proposition~\ref{prop:transformation_group_MCG} again, now for arbitrary $\lambda$, we then deduce that $\cQ^{\fF}_{(\lambda,2,g)}(\MCGo)\cong(\bZ/2)^{r'}$ and $\cQ^{\fF}_{(\lambda,2,h)}(\MCGno)\cong(\bZ/2)^{r'}\times(\bZ/2)^{r}$ for all $g \geq 3$ and $h \geq 7$. In particular, these are independent of $g$ and $h$ respectively, so the result follows by taking the colimits as $g\to\infty$ and $h\to\infty$.
\end{proof}

\bibliography{biblio}

\end{document}